\let\chooseClass3   

\ifx\chooseClass1
\IfFileExists{extarticle.cls}{
\documentclass[14pt]{extarticle}
\emergencystretch 7 pt
\setlength{\oddsidemargin}{-17 mm} 
\setlength{\textwidth}{192 mm}     
\setlength{\textheight}{246 mm}    
\setlength{\topmargin}{-31 mm}     
}{
\documentclass[12pt]{article}
\textheight = 8.5in
\textwidth 6.3in
\setlength{\oddsidemargin}{0mm}
\setlength{\topmargin}{-20 mm}
}{}

\makeatletter
\def\@seccntformat#1{\csname the#1\endcsname.\quad}
\renewcommand\section{\@startsection {section}{1}{\z@}%
                                   {-3.5ex \@plus -1ex \@minus -.2ex}%
                                   {2.3ex \@plus.2ex}%
                                   {\normalfont\large\bfseries}}
\renewcommand\subsection{\@startsection{subsection}{2}{\z@}%
                        {3.25ex plus 1ex minus .2ex}{-.5em}%
                        {\normalfont\normalsize\bfseries}}
\renewcommand\subsubsection{\@startsection{subsubsection}{3}{\z@}%
                        {3.25ex plus 1ex minus .2ex}{-.5em}%
                        {\normalfont\normalsize\it}}
\@addtoreset{equation}{section}
\makeatother

\usepackage{amsthm}
\newtheoremstyle{boldhead}
{\topsep}
{\topsep}
{\slshape}
{}
{\bfseries}
{.}
{ }
{\thmname{#1}\thmnumber{ #2}\thmnote{ (#3)}}

\newtheoremstyle{boldremark}
{\topsep}
{\topsep}
{\upshape}
{}
{\bfseries}
{.}
{ }
{\thmname{#1}\thmnumber{ #2}\thmnote{ (#3)}}

\swapnumbers

\theoremstyle{boldhead}
\newtheorem{theorem}[subsection]{Theorem}
\newtheorem{corollary}[subsection]{Corollary}
\newtheorem{lemma}[subsection]{Lemma}
\newtheorem{proposition}[subsection]{Proposition}
\newtheorem{maintheorem}[subsection]{Main Theorem}
\newtheorem{NOproposition}{Proposition}
\newtheorem{NOlemma}[NOproposition]{Lemma}

\theoremstyle{boldremark}
\newtheorem{definition}[subsection]{Definition}
\newtheorem{example}[subsection]{Example}
\newtheorem{examples}[subsection]{Examples}
\newtheorem{remark}[subsection]{Remark}
\newtheorem{problem}[subsection]{Problem}
\newtheorem{question}[subsection]{Question}
\newtheorem{warning}[subsection]{Warning}
\newtheorem{acknowledgement}[subsection]{Acknowledgements}

\fi

\ifx\chooseClass2
\documentclass{tac}
\newtheorem{maintheorem}{Main Theorem}
\newtheoremrm{acknowledgement}{Acknowledgements}
\fi

\ifx\chooseClass3
\documentclass[12pt]{article}
\makeatletter
\def\@seccntformat#1{\csname the#1\endcsname.\quad}
\renewcommand\section{\@startsection {section}{1}{\z@}%
                                   {-3.5ex \@plus -1ex \@minus -.2ex}%
                                   {2.3ex \@plus.2ex}%
                                   {\normalfont\large\bfseries}}
\renewcommand\subsection{\@startsection{subsection}{2}{\z@}%
                        {3.25ex plus 1ex minus .2ex}{-.5em}%
                        {\normalfont\normalsize\bfseries}}
\renewcommand\subsubsection{\@startsection{subsubsection}{3}{\z@}%
                        {3.25ex plus 1ex minus .2ex}{-.5em}%
                        {\normalfont\normalsize\it}}
\@addtoreset{equation}{section}
\makeatother

\usepackage{amsthm}
\newtheoremstyle{boldhead}
{\topsep}
{\topsep}
{\slshape}
{}
{\bfseries}
{.}
{ }
{\thmname{#1}\thmnumber{ #2}\thmnote{ (#3)}}

\newtheoremstyle{boldremark}
{\topsep}
{\topsep}
{\upshape}
{}
{\bfseries}
{.}
{ }
{\thmname{#1}\thmnumber{ #2}\thmnote{ (#3)}}

\swapnumbers

\theoremstyle{boldhead}
\newtheorem{theorem}[subsection]{Theorem}
\newtheorem{corollary}[subsection]{Corollary}
\newtheorem{lemma}[subsection]{Lemma}
\newtheorem{proposition}[subsection]{Proposition}
\newtheorem{maintheorem}[subsection]{Main Theorem}
\newtheorem{NOlemma}{Lemma}

\theoremstyle{boldremark}
\newtheorem{definition}[subsection]{Definition}

\newtheorem{remark}[subsection]{Remark}

\newtheorem{acknowledgement}[subsection]{Acknowledgements}

\fi

\usepackage{    amsmath,
                amsfonts,
                amssymb,
                verbatim,
}

\numberwithin{equation}{subsection}

\IfFileExists{euscript.sty}{\usepackage[mathcal]{euscript}}{}
\IfFileExists{mathrsfs.sty}{\usepackage{mathrsfs}}{\let\mathscr\mathfrak}

\IfFileExists{url.sty}{\usepackage{url}}{}
\providecommand{\url}[1]{{\tt #1}}

\usepackage{ifpdf}
\ifpdf
\usepackage[pdftex]{graphics}
\IfFileExists{hyperref.sty}{\usepackage[pdftex]{hyperref}}{}
\else
\usepackage[dvips]{graphics}
\IfFileExists{hyperref.sty}{\usepackage[hypertex]{hyperref}}{}
\fi

\usepackage{diagrams}
\diagramstyle[height=2em,balance,righteqno,PostScript=dvips,nohug]

\usepackage{t-angles} 
\def\rhaha{\raise.24ex\hbox{$\rightharpoonup$}\kern-1em\lower.24ex\hbox{$\rightharpoondown$}}%
\def\lhaha{\raise.24ex\hbox{$\leftharpoonup$}\kern-1em\lower.24ex\hbox{$\leftharpoondown$}}%
\def\dhaha{\downharpoonleft\kern-.22em\downharpoonright\kern.02em}%
\def\uhaha{\upharpoonleft\kern-.22em\upharpoonright\kern.02em}
\newarrowhead{twoharpoons}\rhaha\lhaha\dhaha\uhaha

\newarrow{Id}===={twoharpoons}
\newarrow{Epi}----{triangle}
\newarrow{MapsTo}{mapsto}---{->}
\newarrow{Mono}{boldhook}---{->}
\newarrow{TTo}----{->}
\newarrow{Twoar}===={=>}

\newcommand\NN{{\mathbb N}}
\newcommand\ZZ{{\mathbb Z}}

\newcommand{\ca}{{\mathcal A}}
\newcommand{\cb}{{\mathcal B}}
\newcommand{\cc}{{\mathcal C}}
\newcommand{\cd}{{\mathcal D}}
\newcommand{\ce}{{\mathcal E}}
\newcommand{\cf}{{\mathcal F}}
\newcommand{\cfd}{{\mathcal F}}

\newcommand{\cm}{{\mathcal M}}
\newcommand{\cn}{{\mathcal N}}

\newcommand{\cq}{{\mathcal Q}}
\newcommand{\ct}{{\mathcal T}}

\newcommand{\cx}{{\mathcal X}}
\newcommand{\cy}{{\mathcal Y}}
\newcommand{\cz}{{\mathcal Z}}

\newcommand{\fA}{{\mathfrak A}}
\newcommand{\opA}{{\mathscr A}}
\newcommand{\fc}{{\mathfrak C}}
\newcommand{\mcM}{{\mathscr M}}
\newcommand{\quo}{\mathsf q}
\newcommand{\Quo}{\mathsf Q}
\newcommand{\SSS}{{\mathfrak S}}
\newcommand{\fu}{{\mathscr U}}
\newcommand{\Yo}{{\mathscr Y}}

\newcommand{\bull}{{\scriptscriptstyle\bullet}}
\newcommand{\Com}{{{\mathsf C}_\kk}}
\newcommand{\uCom}{{\underline{\mathsf C}_\kk}}
\newcommand{\ihat}{\hat{\imath}}
\newcommand\ju{{\overline{\jmath}}}
\newcommand{\tdt}{\tens\dots\tens}
\newcommand{\tr}{{\mathsf{tr}}}
\newcommand{\tree}{{\mathfrak t}}
\newcommand{\uni}{{\mathbf i}}
\newcommand{\unix}{{\sS{_X}\uni^\ca_0}}

\newcommand{\sS}[2]{\vphantom{#2}#1 #2}
\newcommand{\n}[1]{\nobreakdash-\hspace{0pt}}
\newcommand{\ainf}[1]{$A_\infty$\nobreakdash-\hspace{0pt}}
\newcommand{\ainfu}[1]{$A_\infty^u$\nobreakdash-\hspace{0pt}}

\newcommand{\tcolimit}{{\displaystyle\lim_{\longrightarrow}}\raisebox{-1.8mm}{\vphantom{l}}}
\newcommand{\Cat}{{\mathcal C}at}
\newcommand{\Dr}{\mathsf D}
\newcommand{\Rep}{\mathrm{Rep}}
\newcommand{\qverdier}{Q_{\mathrm{Verdier}}}

\let\kk\Bbbk

\let\eps\varepsilon
\let\epsilon\varepsilon
\let\ge\geqslant
\let\le\leqslant
\let\sss\scriptstyle
\let\tens\otimes

\let\wh\widehat
\let\wt\widetilde

\DeclareMathOperator\ad{ad}
\DeclareMathOperator\alve{\overline{Vert}}

\DeclareMathOperator\codom{codom}
\DeclareMathOperator\coev{coev}
\DeclareMathOperator\Coker{Coker}
\DeclareMathOperator\Cone{Cone}

\DeclareMathOperator\dom{dom}

\DeclareMathOperator\ev{ev}
\DeclareMathOperator\gr{gr}
\DeclareMathOperator\Ho{Ho}
\DeclareMathOperator\Hom{Hom}
\DeclareMathOperator\id{id}

\DeclareMathOperator\im{Im}

\DeclareMathOperator\inj{in}
\DeclareMathOperator\inve{Vert}
\DeclareMathOperator\Ker{Ker}
\DeclareMathOperator\Leaf{Leaf}
\DeclareMathOperator\minUV{minUV}
\DeclareMathOperator\modul{-mod}
\DeclareMathOperator\modulo{mod}
\DeclareMathOperator\Ob{Ob}
\newcommand{\op}{{\operatorname{op}}}
\DeclareMathOperator\Out{Out}
\DeclareMathOperator\pr{pr}
\newcommand{\pretr}{{\operatorname{pre-tr}}}
\DeclareMathOperator\restr{restr}
\DeclareMathOperator\UV{UV}

\newcommand{\appref}[1]{Appendix~\ref{#1}}
\newcommand{\corref}[1]{Corollary~\ref{#1}}
\newcommand{\defref}[1]{Definition~\ref{#1}}

\newcommand{\lemref}[1]{Lemma~\ref{#1}}
\newcommand{\propref}[1]{Proposition~\ref{#1}}
\newcommand{\remref}[1]{Remark~\ref{#1}}
\newcommand{\secref}[1]{Section~\ref{#1}}
\newcommand{\thmref}[1]{Theorem~\ref{#1}}

\begin{document}
\bibliographystyle{amsalpha}
\title{Quotients of unital $A_\infty$-categories}
\ifx\chooseClass1
\author{Volodymyr Lyubashenko%
\thanks{Institute of Mathematics,
National Academy of Sciences of Ukraine,
3 Tereshchenkivska st.,
Kyiv-4, 01601 MSP,
Ukraine;
lub@imath.kiev.ua}
\ and Oleksandr Manzyuk%
\thanks{Fachbereich Mathematik,
Postfach 3049,
67653 Kaiserslautern,
Germany;
manzyuk@mathematik.uni-kl.de}
}
\fi
\ifx\chooseClass2
\author{Volodymyr Lyubashenko and Oleksandr Manzyuk}
\address{Institute of Mathematics,
National Academy of Sciences of Ukraine, \\
3 Tereshchenkivska st.,
Kyiv-4, 01601 MSP, Ukraine \\[8pt]
Fachbereich Mathematik,
Postfach 3049,
67653 Kaiserslautern,
Germany \\
}
\copyrightyear{2004}
\eaddress{lub@imath.kiev.ua\CR manzyuk@mathematik.uni-kl.de}
\fi
\ifx\chooseClass3
\author{Volodymyr Lyubashenko%
\thanks{Institute of Mathematics,
National Academy of Sciences of Ukraine,
3 Tereshchenkivska st.,
Kyiv-4, 01601 MSP,
Ukraine;
lub@imath.kiev.ua}
\ and Oleksandr Manzyuk%
\thanks{Fachbereich Mathematik,
Postfach 3049,
67653 Kaiserslautern,
Germany;
manzyuk@mathematik.uni-kl.de}
}
\fi
\maketitle

\begin{abstract}
Assuming that $\cb$ is a full \ainf-subcategory of a unital
\ainf-category $\cc$ we construct the quotient unital \ainf-category
$\cd=$`$\cc/\cb$'. It represents the \ainfu-2-functor
$\ca\mapsto A_\infty^u(\cc,\ca)_{\modulo\cb}$, which associates with a
given unital \ainf-category $\ca$ the \ainf-category of unital
\ainf-functors $\cc\to\ca$, whose restriction to $\cb$ is contractible.
Namely, there is a unital \ainf-functor \(e:\cc\to\cd\) such that the
composition \(\cb\hookrightarrow\cc\overset{e}\longrightarrow\cd\) is
contractible, and for an arbitrary unital \ainf-category $\ca$ the
restriction \ainf-functor
 \((e\boxtimes1)M:A_\infty^u(\cd,\ca)\to
 A_\infty^u(\cc,\ca)_{\modulo\cb}\)
is an equivalence.

Let $\uCom$ be the differential graded category of differential graded
$\kk$\n-modules. We prove that the Yoneda \ainf-functor
$Y:\ca\to A_\infty^u(\ca^\op,\uCom)$ is a full embedding for an
arbitrary unital \ainf-category $\ca$. In particular, such $\ca$ is
\ainf-equivalent to a differential graded category with the same set of
objects.
\end{abstract}

\allowdisplaybreaks[1]

Let $\ca$ be an Abelian category. The question: what is the quotient
\begin{center}
\{category of complexes in $\ca$\}/\{category of acyclic complexes\}?
\end{center}
admits several answers. The first answer -- the derived category of
$\ca$ -- was given by Grothendieck and Verdier~\cite{Verdier2}.

The second answer -- a differential graded category $\cd$ -- is given
by Drinfeld~\cite{Drinf:DGquot}. His article is based on the work of
Bondal and Kapranov~\cite{BondalKapranov:FramedTR} and of
Keller~\cite{MR99m:18012}. The derived category $D(\ca)$ can be
obtained as \(H^0(\cd^\pretr)\).

The third answer -- an \ainf-category of bar-construction type -- is
given by Lyubashenko and Ovsienko~\cite{LyuOvs-iResAiFn}. This
\ainf-category is especially useful when the basic ring $\kk$ is a
field. It is an \ainf-version of one of the constructions of
Drinfeld~\cite{Drinf:DGquot}.

The fourth answer -- an \ainf-category freely generated over the
category of complexes in $\ca$ -- is given in this article. It is
\ainf-equivalent to the third answer and enjoys certain universal
property of the quotient. Thus, it passes this universal property also
to the third answer.

\section{Introduction}
Since \ainf-algebras were introduced by
Stasheff~\cite[II]{Stasheff:HomAssoc} there existed a possibility to
consider \ainf-generalizations of categories. It did not happen until
\ainf-categories were encountered in studies of mirror symmetry by
Fukaya~\cite{Fukaya:A-infty} and
Kontsevich~\cite{Kontsevich:alg-geom/9411018}. \ainf-categories may be
viewed as generalizations of differential graded categories for which
the binary composition is associative only up to a homotopy. The
possibility to define \ainf-functors was mentioned by
Smirnov~\cite{MR90h:55033}, who reformulated one of his results in the
language of \ainf-functors between differential graded categories. The
definition of \ainf-functors between \ainf-categories was published by
Keller~\cite{math.RA/9910179}, who studied their applications to
homological algebra. Homomorphisms of \ainf-algebras (e.g.
\cite{Kadeishvili82}) are particular cases of \ainf-functors.

\ainf-transformations between \ainf-functors are certain coderivations.
Given two \ainf-categories $\ca$ and $\cb$, one can construct a third
\ainf-category $A_\infty(\ca,\cb)$, whose objects are \ainf-functors
$f:\ca\to\cb$, and morphisms are \ainf-transformations
(Fukaya~\cite{Fukaya:FloerMirror-II}, Kon\-tse\-vich and Soibelman
\cite{math.RA/0606241,KonSoi-book},
Le\-f\`evre-Ha\-se\-ga\-wa~\cite{Lefevre-Ainfty-these}, as well as
\cite{Lyu-AinfCat}). For an \ainf-category $\cc$ there is a homotopy
invariant notion of unit elements (identity morphisms)
\cite{Lyu-AinfCat}. They are cycles \(\sS{_X}\uni^\cc_0\in s\cc(X,X)\)
of degree $-1$ such that the maps
\((1\tens\uni^\cc_0)b_2,-(\uni^\cc_0\tens1)b_2:s\cc(X,Y)\to s\cc(X,Y)\)
are homotopic to the identity map. This allows to define the
2\n-category \(\overline{A_\infty^u}\), whose objects are unital
\ainf-categories (those which have units), 1\n-morphisms are unital
\ainf-functors (their first components preserve the units up to a
boundary) and 2\n-morphisms are equivalence classes of natural
\ainf-transformations \cite{Lyu-AinfCat}. We continue to study this
2\n-category. Notations and terminology follow \cite{Lyu-AinfCat},
complemented by \cite{LyuOvs-iResAiFn} and \cite{LyuMan-freeAinf}.

Unital \ainf-categories and unital \ainf-functors can be considered as
strong homotopy generalizations of differential graded categories and
functors. Let us illustrate the notion of \ainf-transformations in a
familiar context.

\subsection{\texorpdfstring{Differential for $A_\infty$-transformations
 compared with the Hochschild differential}{Differential for
 A8-transformations compared with the Hochschild differential}}
Let $\ca$, $\cb$ be ordinary $\kk$\n-linear categories. We consider
$\ca(\_,\_)$ and $\cb(\_,\_)$ as complexes of $\kk$\n-modules
concentrated in degree 0. This turns $\ca$ and $\cb$ into differential
graded categories and, thereby, into unital \ainf-categories. An
\ainf-functor between $\ca$ and $\cb$ is necessarily strict, for
$(s\ca)^{\tens k}=\ca^{\tens k}[k]$ and $s\cb=\cb[1]$ are concentrated
in different degrees if $k>1$. Thus, a unital \ainf-functor
$f:\ca\to\cb$ is the same as an ordinary $\kk$\n-linear functor $f$.
Let $f,g:\ca\to\cb$ be $\kk$\n-linear functors. All complexes
\(\uCom((s\ca)^{\tens k}(X,Y),s\cb(Xf,Yg))\) are concentrated in degree
\(k-1\). Their direct product
\[ \Psi_k =\prod_{X,Y\in\Ob\ca} \uCom((s\ca)^{\tens k}(X,Y),s\cb(Xf,Yg))
\]
is the same, whether taken in the category of $\kk$\n-modules or graded
\(\kk\)\n-modules or complexes of \(\kk\)\n-modules. It is the module
of $k$\n-th components of \ainf-transformations. The graded
$\kk$\n-module of \ainf-transformations \(sA_\infty(\ca,\cb)(f,g)\) is
isomorphic to the direct product $\prod_{k=0}^\infty\Psi_k$ taken in
the category of graded \(\kk\)\n-modules
\cite[Section~2.7]{Lyu-AinfCat}. That is,
\([sA_\infty(\ca,\cb)(f,g)]^n=\prod_{k=0}^\infty\Psi_k^n\), where
\(\Psi_k^n\) is the degree $n$ part of \(\Psi_k\). Therefore, in our
case it simply coincides with the graded \(\kk\)\n-module
\(\Psi[1]:\ZZ\ni n\mapsto\Psi_{n+1}\in\kk\modul\). The graded
$\kk$\n-module \(sA_\infty(\ca,\cb)(f,g)\) is equipped with the
differential $B_1$, \(rB_1=rb^\cb-(-)^rb^\ca r\)
\cite[Proposition~5.1]{Lyu-AinfCat}. Since the only non-vanishing
component of $b$ (resp. $f$) is \(b_2\) (resp. \(f_1\)), the explicit
formula for components of \(rB_1\) is the following:
\[ (rB_1)_{k+1} = (f_1\tens r_k)b_2 + (r_k\tens g_1)b_2
-(-)^{r_k} \sum_{a+c=k-1} (1^{\tens a}\tens b_2\tens1^{\tens c})r_k.
\]
Recalling that $\deg r_k=k-1$, we get the differential in \(\Psi[1]\)
also denoted $B_1$:
\[ r_kB_1 = (f_1\tens r_k)b_2 + (r_k\tens g_1)b_2
+(-)^k \sum_{a+c=k-1} (1^{\tens a}\tens b_2\tens1^{\tens c})r_k,
\]
where \(r_k\in\Psi_k\), \(r_kB_1\in\Psi_{k+1}\). We consider an
isomorphism  of graded \(\kk\)\n-modules
\(\Psi\to\Psi':\ZZ\ni k\mapsto\Psi'_k\) given by
\[ \Psi_k \ni r_k \mapsto (s\tdt s)r_ks^{-1}
= s^{\tens k}r_ks^{-1} \in \Psi'_k \overset{\text{def}}=
\prod_{X,Y\in\Ob\ca} \uCom(\ca^{\tens k}(X,Y),\cb(Xf,Yg)).
\]
Its inverse is
\(\Psi'_k\ni t_k\mapsto(s^{\tens k})^{-1}t_ks\in\Psi_k\). This
isomorphism induces the differential
\[ d: \Psi'_k \to \Psi'_{k+1}, \qquad
t_kd = s^{\tens k+1}\cdot[(s^{\tens k})^{-1}t_ks]B_1\cdot s^{-1}.
\]
The explicit formula for $d$ is
\[ t_kd = (f\tens t_k)m_2
+ \sum_{a+c=k-1} (-1)^{a+1}(1^{\tens a}\tens m_2\tens1^{\tens c})t_k
+ (-1)^{k+1}(t_k\tens g)m_2.
\]
Up to an overall sign this coincides with the differential in the
Hochschild cochain complex $C^\bull(\ca,\sS{_f}\cb_g)$ (cf.
\cite[Section~X.3]{MacLane-Homology}). The \(\ca\)\n-bimodule
\(\sS{_f}\cb_g\) acquires its left \(\ca\)\n-module structure via $f$
and its right \(\ca\)\n-module structure via $g$. Therefore, in our
situation \ainf-transformations are nothing else but Hochschild
cochains. Natural \ainf-transformations \(r:f\to g:\ca\to\cb\) (such
that $\deg r=-1$ and $rB_1=0$) are identified with the Hochschild
cocycles of degree 0, that is, with natural transformations
\(t=rs^{-1}:f\to g:\ca\to\cb\) in the ordinary sense.

When $\ca$, $\cb$ are differential graded categories and
$f,g:\ca\to\cb$ are differential graded functors, we may still
interpret the complex \((sA_\infty(\ca,\cb)(f,g),B_1)\) as the complex
of Hochschild cochains $C^\bull(\ca,\sS{_f}\cb_g)$ for the differential
graded category $\ca$ and the differential graded bimodule
\(\sS{_f}\cb_g\). Indeed, for a homogeneous element
\(r\in sA_\infty(\ca,\cb)(f,g)\) the components of \(rB_1\) are
\begin{multline*}
(rB_1)_k = r_kb_1 +(f_1\tens r_{k-1})b_2 + (r_{k-1}\tens g_1)b_2 \\
-(-)^r \sum_{a+1+c=k} (1^{\tens a}\tens b_1\tens1^{\tens c})r_k
-(-)^r \sum_{a+2+c=k} (1^{\tens a}\tens b_2\tens1^{\tens c})r_{k-1}.
\end{multline*}
For \ainf-functors between differential graded categories or
\ainf-categories the differential \(B_1\) is not interpreted as
Hochschild differential any more. But we may view the complex of
\ainf-transformations as a generalization of the Hochschild cochain
complex.

\subsection{Main result}
By Definition~6.4 of \cite{LyuOvs-iResAiFn} an \ainf-functor
$g:\cb\to\ca$ from a unital \ainf-category $\cb$ is \emph{contractible}
if for all objects $X$, $Y$ of $\cb$ the chain map
\(g_1:s\cb(X,Y)\to s\ca(Xg,Yg)\) is null-homotopic. If $g:\cb\to\ca$ is
a unital \ainf-functor, then it is contractible if and only if for any
$X\in\Ob\cb$ and any $V\in\Ob\ca$ the complexes $s\ca(Xg,V)$ and
$s\ca(V,Xg)$ are contractible. Equivalently,
\(g\uni^\ca\equiv0:g\to g:\cb\to\ca\)
\cite[Proposition~6.1(C5)]{LyuOvs-iResAiFn}. Other equivalent
conditions are listed in Propositions 6.1--6.3 of
\cite{LyuOvs-iResAiFn}.

Let $\cb$ be a full \ainf-subcategory of a unital \ainf-category $\cc$.
Let $\ca$ be an arbitrary unital \ainf-category. Denote by
\(A_\infty^u(\cc,\ca)_{\modulo\cb}\) the full \ainf-subcategory of
\(A_\infty^u(\cc,\ca)\), whose objects are unital \ainf-functors
$\cc\to\ca$, whose restriction to $\cb$ is contractible. We allow
consideration of \ainf-categories with the empty set of objects.

\begin{maintheorem}\label{thm-maintheorem}
In the above assumptions there exists a unital \ainf-category
$\cd=\quo(\cc|\cb)$ and a unital \ainf-functor \(e:\cc\to\cd\) such
that
\begin{enumerate}
\item[1)] the composition \(\cb \rMono \cc \rTTo^e \cd\) is
contractible;

\item[2)] the strict \ainf-functor given by composition with $e$
\[ (e\boxtimes1)M: A_\infty^u(\cd,\ca) \to
A_\infty^u(\cc,\ca)_{\modulo\cb}, \qquad f \mapsto ef,
\]
is an \ainf-equivalence for an arbitrary unital \ainf-category $\ca$.
\end{enumerate}
\end{maintheorem}

\begin{proof}
Let us prove the statement first in a particular case, for a full
subcategory \(\wt{\cb}\) of a strictly unital \ainf-category
\(\wt{\cc}\). Then the representing \ainf-category
\(\cd=\Quo(\wt{\cc}|\wt{\cb})\) is constructed in
\secref{sec-constr-D} as an \ainf-category, freely generated over
\(\wt{\cc}\) by application of contracting homotopies $H$ to morphisms,
whose source or target is in \(\wt{\cb}\). The strict \ainf-functor
\(\wt{e}:\wt{\cc}\to\cd\) is identity on objects and \(\wt{e}_1\) is an
embedding. \thmref{thm-C-strict-D-unital} asserts unitality of
\(\cd=\Quo(\wt{\cc}|\wt{\cb})\). By construction, the \ainf-functor
\(\wt{e}:\wt{\cc}\to\cd\) is unital and
 \(\wt\cb\rMono \wt\cc\rTTo^{\wt{e}} \cd\) is contractible. By
\thmref{thm-restriction-DA-CAmodB-equi} the restriction strict
\ainf-functor
 \(\restr:A_\infty^u(\cd,\ca)\to A_\infty^u(\cc,\ca)_{\modulo\cb}\) is
an \ainf-equivalence. Thus, $\cd$ and \(\wt{e}:\wt{\cc}\to\cd\)
represent the \ainfu-2-functor
 $\ca\mapsto A_\infty^u(\wt{\cc},\ca)_{\modulo\wt{\cb}}$
in the sense of 1), 2), as claimed.

Let now $\cb$ be a full \ainf-subcategory of a unital \ainf-category $\cc$.
There exists a differential graded category \(\wt{\cc}\) with
\(\Ob\wt{\cc}=\Ob\cc\), and quasi-inverse to each other \ainf-functors
\(\wt{Y}:\cc\to\wt{\cc}\), \(\Psi:\wt{\cc}\to\cc\) such that
\(\Ob\wt{Y}=\Ob\Psi=\id_{\Ob\cc}\) (by \remref{rem-Yoneda-tilde} this
follows from the \ainf-version of Yoneda Lemma --
\thmref{thm-Yoneda-equivalence}). Let \(\wt{\cb}\subset\wt{\cc}\) be
the full differential graded subcategory with \(\Ob\wt{\cb}=\Ob\cb\).
By the previous case there is a unital \ainf-category $\cd$ and a
unital \ainf-functor \(\wt{e}:\wt{\cc}\to\cd\) representing the
\ainfu-2-functor
 $\ca\mapsto A_\infty^u(\wt{\cc},\ca)_{\modulo\wt{\cb}}$
in the sense of 1), 2). By considerations in \appref{sec-strict-A2fun}
(Corollaries \ref{cor-Psi1M-Y1M-A-equivalences},
\ref{cor-equi-pairs-equi-ainf2funs}) the pair
 \(\bigl(\cd,e=(\cc \rTTo^{\wt{Y}} \wt{\cc} \rTTo^{\wt{e}} \cd)\bigr)\)
represents \(A_\infty^u(\cc,\ca)_{\modulo\cb}\). Indeed,
\[ (e\boxtimes1)M = \bigl( A_\infty^u(\cd,\ca)
\rTTo^{(\wt{e}\boxtimes1)M} A_\infty^u(\wt{\cc},\ca)_{\modulo\wt{\cb}}
\rTTo^{(\wt{Y}\boxtimes1)M} A_\infty^u(\cc,\ca)_{\modulo\cb} \bigr)
\]
is a composition of two \ainf-equivalences.
\end{proof}

Notation for our quotient constructions is the following. The
construction of \secref{sec-constr-D} is denoted \(\Quo(\_|\_)\). When
it is combined with the Yoneda \ainf-equivalence of
\remref{rem-Yoneda-tilde} we denote it \(\quo(\_|\_)\).

The 2\n-category \(\overline{A_\infty^u}=H^0A_\infty^u\) has unital
\ainf-categories as objects, unital \ainf-functors as 1\n-morphisms and
equivalence classes of natural \ainf-transformations as 2\n-morphisms
\cite{Lyu-AinfCat}. Thus,
 \(\overline{A_\infty^u}(\cc,\ca)(f,g)=[H^0A_\infty^u(\cc,\ca)](f,g)
 =H^0(A_\infty^u(\cc,\ca)(f,g),m_1)\).

A zero object of a category $\ce$ is an object $Z$, which is
simultaneously initial and terminal. For a linear ($Ab$\n-enriched, not
necessarily additive) category $\ce$ this can be formulated as follows:
\(\ce(Z,X)=0\) and \(\ce(X,Z)=0\) for any object $X$ of $\ce$. This
condition is equivalent to the equation \(1_Z=0\in\ce(Z,Z)\).

\begin{corollary}\label{cor-H0Aiu-exact-seq-cats}
The unital \ainf-functor \(e:\cc\to\cd\) from the main theorem has the
following property: composition with $e$ in the sequence of functors
\[ \overline{A_\infty^u}(\cd,\ca) \rTTo^{e\bull}
\overline{A_\infty^u}(\cc,\ca) \rTTo \overline{A_\infty^u}(\cb,\ca)
\]
is an equivalence of the category \(\overline{A_\infty^u}(\cd,\ca)\)
with the full subcategory
\[
\Ker(\overline{A_\infty^u}(\cc,\ca) \to \overline{A_\infty^u}(\cb,\ca))
= H^0(A_\infty^u(\cc,\ca)_{\modulo\cb},m_1)
\subset \overline{A_\infty^u}(\cc,\ca),
\]
consisting of unital \ainf-functors \(f:\cc\to\ca\) such that the
restriction \(f\big|_\cb:\cb\to\ca\) is a zero object of
\(\overline{A_\infty^u}(\cb,\ca)\).
\end{corollary}

\begin{proof}
A unital \ainf-functor \(g:\cb\to\ca\) is contractible if and only if
\(g\uni^\ca\equiv0:g\to g:\cb\to\ca\), that is,
\(1_g=0\in \overline{A_\infty^u}(\cb,\ca)(g,g)\). Thus, for unital $g$
contractibility is equivalent to $g$ being a zero object of
\(\overline{A_\infty^u}(\cb,\ca)\).
\end{proof}

The main theorem asserts that the chain map
 \(e\bull=s(e\boxtimes1)M_{01}s^{-1}:
 A_\infty^u(\cd,\ca)(f,g)\to A_\infty^u(\cc,\ca)(ef,eg)\)
is a homotopy isomorphism, while \corref{cor-H0Aiu-exact-seq-cats}
claims only that it induces isomorphism in 0\n-th homology.

\subsection{\texorpdfstring{Uniqueness of the representing $A_\infty$-category}
 {Uniqueness of the representing A8-category}}
With each strict \ainfu-2-functor $F:A_\infty^u\to A_\infty^u$ is
associated an ordinary strict 2\n-functor
$\overline{F}:\overline{A_\infty^u}\to\overline{A_\infty^u}$,
$\overline{F}\ca=F\ca$ \cite[Section~3.2]{LyuMan-freeAinf}. With a
strict \ainfu-2-transformation
$\lambda=(\lambda_\ca):F\to G:A_\infty^u\to A_\infty^u$ is associated
an ordinary strict 2\n-transformation
 $\overline{\lambda}=(\overline{\lambda_\ca}):
 \overline{F}\to\overline{G}:
 \overline{A_\infty^u}\to\overline{A_\infty^u}$
in cohomology [\textit{ibid}]. Assume that $\lambda$ is a natural
\ainfu-2-equivalence. Since \(\lambda_\ca:F\ca\to G\ca\) are
\ainf-equivalences, the 1\n-morphisms
\(\overline{\lambda_\ca}:\overline{F}\ca\to\overline{G}\ca\) are
equivalences in the 2\n-category \(\overline{A_\infty^u}\). Composing
$\overline{\lambda}$ with the 0\n-th cohomology 2\n-functor
$H^0:\overline{A_\infty^u}\to\Cat$, we get a 2\n-natural equivalence
 \(H^0\overline{\lambda}:H^0\overline{F}\to H^0\overline{G}:
 \overline{A_\infty^u}\to\Cat\),
which consists of equivalences of ordinary categories
\(H^0(\lambda_\ca):H^0(F\ca)\to H^0(G\ca)\). In particular, if
\(F=A_\infty^u(\cd,\_)\) for some unital \ainf-category $\cd$, then
\[ H^0\overline{F} = H^0\overline{A_\infty^u(\cd,\_)}
= \overline{A_\infty^u}(\cd,\_).
\]
Indeed, both the categories \(H^0\overline{A_\infty^u(\cd,\ca)}\) and
\(\overline{A_\infty^u}(\cd,\ca)\) have unital \ainf-functors
$\cd\to\ca$ as objects and equivalence classes of natural
 \ainf-transformations as morphisms. If
$\lambda:A_\infty^u(\cd,\_)\to G:A_\infty^u\to A_\infty^u$ is a natural
\ainfu-2-equivalence ($G$ is \emph{unitally representable} by $\cd$),
then
 \(H^0(\lambda_\ca):\overline{A_\infty^u}(\cd,\ca)\to H^0(G\ca):
 \overline{A_\infty^u}\to\Cat\)
is a 2\n-natural equivalence. Thus, \(H^0\overline{G}\) is represented
by $\cd$ in the 2\n-category sense and $\cd$ is unique up to an
equivalence by \secref{sec-representing-pair-unique}.

In particular, if \(G=A_\infty^u(\cc,\_)_{\modulo\cb}\), then with each
object $e$ of \(G\cd=A_\infty^u(\cc,\cd)_{\modulo\cb}\) is associated a
strict \ainfu-2-transformation
\[ \lambda = (e\boxtimes1)M: A_\infty^u(\cd,\ca) \to
A_\infty^u(\cc,\ca)_{\modulo\cb}: A_\infty^u \to A_\infty^u.
\]
We have identified \(H^0\overline{G}(\ca)\) with
 \(\Ker(\overline{A_\infty^u}(\cc,\ca)
 \to\overline{A_\infty^u}(\cb,\ca))\)
in \corref{cor-H0Aiu-exact-seq-cats}. The strict 2\n-natural
equivalence
 \(H^0\overline{\lambda}:H^0\overline{F}\to H^0\overline{G}:
 \overline{A_\infty^u}\to\Cat\)
identifies with the strict 2\n-transformation
\[ \sS{^{H^0\overline{G}}}\lambda^{\cd,e}:\overline{A_\infty^u}(\cd,\_)
\to H^0\overline{G}: \overline{A_\infty^u}\to\Cat
\]
from \propref{pro-functor-GA-(ACat)}, since
\begin{align*}
(f:\cd\to\ca) & \rMapsTo^{(e\boxtimes1)M} ef = (e)(H^0\overline{G}(f))
\overset{\text{def}}=: (f)\sS{^{H^0\overline{G}}}\lambda^{\cd,e}, \\
(r:f\to g:\cd\to\ca) & \rMapsTo^{(e\boxtimes1)M} er
= (e)(H^0\overline{G}(r))
\overset{\text{def}}=: (r)\sS{^{H^0\overline{G}}}\lambda^{\cd,e}.
\end{align*}
Therefore, the pair $(\cd,e)$ represents the strict 2\n-functor
\(H^0\overline{G}:\overline{A_\infty^u}\to\Cat\) in the sense of
\defref{def-representable-by-pair}.

\begin{corollary}
The pair \((\cd,e:\cc\to\cd)\) is unique up to an equivalence, that is,
for any other quotient \((\cd',e':\cc\to\cd')\) there exists an
\ainf-equivalence \(\phi:\cd\to\cd'\) such that $e\phi$ is isomorphic
to $e'$.
\end{corollary}

The proof immediately follows from results of
\secref{sec-representing-pair-unique}.

The unital \ainf-category $\cd$ obtained in the main theorem can be
replaced with a differential graded category by the \ainf-version of
Yoneda Lemma (\thmref{thm-Yoneda-equivalence}). We may restrict
\corref{cor-H0Aiu-exact-seq-cats} to differential graded categories
$\cb$, $\cc$, $\ca$ for the same reason. Then it becomes parallel to
the second half of main Theorem~1.6.2 of Drinfeld's
work~\cite{Drinf:DGquot}, which asserts exactness of the sequence of
categories
\[ T(\cd,\ca) \to T(\cc,\ca) \to T(\cb,\ca)
\]
in the same sense as in \corref{cor-H0Aiu-exact-seq-cats}. Here $T$ is
a certain 2\n-category whose objects are differential graded
categories. It is not known in general whether categories
\(T(\cc,\ca)\) and \(\overline{A_\infty^u}(\cc,\ca)\) are equivalent.
If $\kk$ is a field, then, as B.~Keller explained to us, equivalence of
\(\overline{A_\infty^u}(\cc,\ca)\) and \(T(\cc,\ca)\) can be deduced
from results of Lef\`evre-Hasegawa
\cite[Section~8.2]{Lefevre-Ainfty-these}.

\subsection{Basic properties of the main construction}
The proof of universality of \(\cd=\Quo(\cc|\cb)\) is based on the fact
that $\cd$ is \emph{relatively free} over $\cc$, that is, it admits a
filtration
\begin{equation*}
\cc = \cd_0 \subset \cq_1 \subset \cd_1 \subset \cq_2
\subset \cd_2 \subset \cq_3 \subset \dots \subset \cd
\end{equation*}
by \ainf-subcategories \(\cd_j\) and differential graded subquivers
\(\cq_j\), such that the graded subquiver \(\cd_j\subset\cq_{j+1}\) has
a direct complement $\cn_{j+1}$ (a graded subquiver of $\cq_{j+1}$),
and such that $\cd_{j+1}$ is generated by $\cn_{j+1}$ over $\cd_j$. The
precise conditions are given in \defref{def-relatively-free-Ainf} (see
also \propref{pro-Rb1-subset(R)}). This filtration allows to write down
a sequence of restriction \ainf-functors
\begin{multline}
A_\infty^u(\cc,\ca)_{\modulo\cb}
\lTTo A_{\infty1}^{\psi u}(\cd_0,\cq_1;\ca)
\lTTo A_\infty^{\psi u}(\cd_1,\ca) \\
\lTTo A_{\infty1}^{\psi u}(\cd_1,\cq_2;\ca)
\lTTo A_\infty^{\psi u}(\cd_2,\ca)
\lTTo A_{\infty1}^{\psi u}(\cd_2,\cq_3;\ca)
\lTTo \dots
\tag{\ref{eq-A(CA)-A(D0Q1A)-A(D1A)-A(D1Q2A)}}
\end{multline}
and to prove that each of these \ainf-functors is an equivalence, surjective
on objects (\thmref{thm-restriction-equivalence-quiver},
Propositions \ref{pro-A(D0Q1A)-A(CA)modB-equivalence},
\ref{pro-restriction-A(EQA)-A(EA)}). The category
\(A_{\infty1}^{\psi u}(\cd_j,\cq_{j+1};\ca)\) is defined via pull-back square
\eqref{dia-A(EQA)}
\begin{diagram}
A_{\infty1}^{\psi u}(\cd_j,\cq_{j+1};\ca) \SEpbk
& \rTTo & A_1(\cq_{j+1},\ca) \\
\dTTo && \dTTo \\
A_\infty^{\psi u}(\cd_j,\ca) & \rTTo & A_1(\cd_j,\ca)
\end{diagram}
The \ainf-categories $\cd_j$ are not unital, but only pseudounital --
there are distinguished cycles \(\uni^\cc_0\in(s\cd_j)^{-1}\), which
are not unit elements of $\cd_j$ if $j>0$. The index $\psi u$ in
\(A_\infty^{\psi u}\) indicates that we consider pseudounital
\ainf-functors -- a generalization of unital ones
(\defref{def-pseudounital}). Their first components preserve the
distinguished cycles up to a boundary. The \ainf-equivalence
\(A_\infty^{\psi u}(\cd,\ca)\to A_\infty^u(\cc,\ca)_{\modulo\cb}\) is
the limit case of \eqref{eq-A(CA)-A(D0Q1A)-A(D1A)-A(D1Q2A)}
(\thmref{thm-restriction-DA-CAmodB-equi}).

The proof of unitality of \(\cd=\Quo(\cc|\cb)\) for strictly unital
$\cc$ is based on the study of the multicategory of \ainf-operations
and contracting homotopies operating in $\cd$
(\thmref{thm-C-strict-D-unital}).

\subsection{Description of various results}
The proof of \thmref{thm-restriction-equivalence-quiver} is based on
description of chain maps \(P\to sA_\infty(\cf\cq,\ca)(\phi,\psi)\) to
the complex of $(\phi,\psi)$-coderivations
(\propref{pro-chain-to-Ainf-FQ-A}), where $\cf\cq$ is the free
\ainf-category, generated by a differential graded quiver $\cq$. A
similar result for the quotient \(\cf\cq/s^{-1}I\) over an \ainf-ideal
$I$ is given in \propref{pro-chain-map-in-subcomplex}. Of course, any
\ainf-category is a quotient of a free one (\propref{pro-FD/(R)-D}). We
also describe homotopies between chain maps
\(P\to sA_\infty(\cf\cq,\ca)(\phi,\psi)\)
(\corref{cor-chain-to-Ainf-FQ-A-null-homotopic}), and generalize the
result to quotients \(\cf\cq/s^{-1}I\)
(\corref{cor-null-homotopic-in-subcomplex}).

In \secref{sec-example-complexes} we consider the example of
differential graded category $\cc=\underline{\mathsf C}(\ca)$ of
complexes in a $\kk$\n-linear Abelian category $\ca$, and the full
subcategory $\cb\subset\cc$ of acyclic complexes. The functor $H^0e$
factors through a functor \(g:D(\ca)\to H^0(\Quo(\cc|\cb))\) by
\corref{cor-DA-H0D}. It is an equivalence, when $\kk$ is a field.

In \appref{ap-Yoneda-Lemma} we define, following
Fukaya~\cite[Lemma~9.8]{Fukaya:FloerMirror-II}, the Yoneda
\ainf-functor $Y:\ca^{\op}\to A_\infty(\ca,\uCom)$, where $\uCom$ is the
differential graded category of complexes of $\kk$\n-modules. We prove
for an arbitrary unital \ainf-category $\ca$ that the Yoneda
\ainf-functor $Y$ is an equivalence of $\ca^\op$ with its image -- full
differential graded subcategory of $A_\infty(\ca,\uCom)$
(\thmref{thm-Yoneda-equivalence}). This is already proven by Fukaya in
the case of strictly unital \ainf-category $\ca$
\cite[Theorem~9.1]{Fukaya:FloerMirror-II}. As a corollary we deduce
that any $\fu$\n-small unital \ainf-category $\ca$ is \ainf-equivalent
to a $\fu$\n-small differential graded category
(\corref{cor-unital-equivalent-differential}).

In \appref{ap-2-Yoneda} we lift the classical Yoneda Lemma one
dimension up -- to strict 2\n-categories, weak 2\n-functors and weak
2\n-transformations. In the completely strict set-up such lifting can
be obtained via enriched category theory, namely, that of
$\Cat$\n-categories, see Street and Walters~\cite{Street:12},
Kelly~\cite{KellyGM:bascec}. The present weak generalization admits a
direct proof.

An important result from another paper is recalled in simplified form,
in which it is used in the present paper:

\begin{corollary}[to Theorem~8.8 \cite{Lyu-AinfCat}]
 \label{cor-Theorem8:8}
Let $\cc$ be an \ainf-category and let $\cb$ be a unital
\ainf-category. Let $\phi:\cc\to\cb$ be an \ainf-functor such that for
all objects $X$, $Y$ of $\cc$ the chain map
$\phi_1:(s\cc(X,Y),b_1)\to(s\cb(X\phi,Y\phi),b_1)$ is homotopy
invertible. If \(\Ob\phi:\Ob\cc\to\Ob\cb\) is surjective, then $\cc$ is
unital and $\phi$ is an \ainf-equivalence.
\end{corollary}

\begin{proof}
Let $h:\Ob\cb\to\Ob\cc$ be an arbitrary mapping such that
\(h\cdot\Ob\phi=\id_{\Ob\cb}\). The remaining data required in
Theorem~8.8 of \cite{Lyu-AinfCat} can be chosen as
\(\sS{_U}r_0=\sS{_U}p_0=\sS{_U}\uni^\cb_0:\kk\to(s\cb)^{-1}(U,U)\) for
all objects $U$ of $\cb$. We conclude by this theorem that there exists
an \ainf-functor $\psi:\cb\to\cc$ with $\Ob\psi=h$, quasi-inverse to
$\phi$.
\end{proof}

Logical dependence of sections is the following. Appendices
\ref{ap-Yoneda-Lemma} and \ref{ap-2-Yoneda} do not depend on other
sections. \appref{sec-strict-A2fun} depends on
\appref{ap-Yoneda-Lemma}. Sections
\ref{sec-Quot-free}--\ref{sec-example-complexes} depend on appendices
and on sections with smaller number. Dependence on the first section
means dependence on \corref{cor-Theorem8:8} and on overall notations
and conventions. Being a summary, the first section depends on all the
rest of the article.

\subsection{Conventions and preliminaries}\label{sec-convent-nota}
We keep the notations and conventions of
\cite{Lyu-AinfCat,LyuOvs-iResAiFn,LyuMan-freeAinf}, sometimes without
explicit mentioning. Some of the conventions are recalled here.

We assume that most quivers, \ainf-categories, etc. are small with
respect to some universe $\fu$. It means that the set of objects and
the set of morphisms are $\fu$\n-small, that is, isomorphic as sets to
an element of $\fu$ \cite[Section~1.0]{GroVer-Prefai-SGA4a}. The
universe $\fu$ is supposed to be an element of a universe $\fu'$, which
in its turn is an element of a universe $\fu''$, and so on. All sets
are supposed to be in bijection with some elements of some of the
universes. Some differential graded categories in this paper will be
$\fu'$\n-small $\fu$\n-categories. A category $\cc$ is a
$\fu$\n-category if all its sets of morphisms $\cc(X,Y)$ are
$\fu$\n-small \cite[Definition~1.1]{GroVer-Prefai-SGA4a}.

The $\fu$\n-small ground ring $\kk$ is a unital associative commutative
ring. A $\kk$\n-module means a $\fu$\n-small $\kk$\n-module.

We use the right operators: the composition of two maps (or morphisms)
$f:X\to Y$ and $g:Y\to Z$ is denoted $fg:X\to Z$; a map is written on
elements as $f:x\mapsto xf=(x)f$. However, these conventions are not
used systematically, and $f(x)$ might be used instead.

The set of non-negative integers is denoted $\NN=\ZZ_{\ge0}$.

We consider only such \ainf-categories $\cc$ that the differential
\(b:Ts\cc\to Ts\cc\) vanishes on $T^0s\cc$, that is, $b_0=0$. We
consider only those \ainf-functors $f:\ca\to\cb$, whose 0\n-th
component $f_0$ vanishes.

\begin{acknowledgement}
We are grateful to all the participants of the \ainf-category seminar
at the Institute of Mathematics, Kyiv, for attention and fruitful
discussions, especially to Yu.~Bespalov and S.~Ovsienko. We thank all
the staff of Max-Planck-Institut f\"ur Mathematik in Bonn for warm
hospitality and support of this research. The main results of this
article were obtained during the stage of the first author in MPIM, and
a short term visit to MPIM of the second author. We are indebted to
M.~Jibladze for a valuable advice to look for an operadic approach to
the quotient category. We thank B.~Keller for the explanation of some
results obtained by K.~Lef\`evre-Hasegawa in his Ph.D. thesis.
\end{acknowledgement}

\section{\texorpdfstring{Quotients of free $A_\infty$-categories}
 {Quotients of free A8-categories}}\label{sec-Quot-free}
\begin{definition}
Let $\cc$ be an \ainf-category, and let \(I\subset s\cc\) be a graded
subquiver with \(\Ob I=\Ob\ca\). The subquiver $I$ is called an
\emph{\ainf-ideal} of $\cc$ if
\[ \im\bigl( b_{\alpha+1+\beta}:
(s\cc)^{\tens\alpha}\tens I\tens(s\cc)^{\tens\beta} \to s\cc\bigr)
\subset I
\]
for all $\alpha$, $\beta\ge0$.
\end{definition}

If \(I\subset s\cc\) is an \ainf-ideal of an \ainf-category $\cc$, then
the quotient graded quiver \(\ce=\cc/s^{-1}I\) with \(\Ob\ce=\Ob\cc\),
\(\ce(X,Y)=\cc(X,Y)/s^{-1}I(X,Y)\) has a unique \ainf-category
structure such that the natural projection
 \(\pi_1:s\cc\to s\ce=s\cc/I\) determines a strict \ainf-functor
\(\pi:\cc\to\ce\). Multiplications \(b^\ce_k\) in $\ce$ are
well-defined for $k\ge1$ by the equation
\begin{multline*}
\bigl[(s\cc)^{\tens k} \rTTo^{b^\cc_k} s\cc \rTTo^{\pi_1} s\cc/I \bigr]
\\
=\Bigl[ (s\cc)^{\tens k} \rEpi
 (s\cc)^{\tens k}\Big/\sum_{\alpha+1+\beta=k}
 (s\cc)^{\tens\alpha}\tens I\tens(s\cc)^{\tens\beta}
\simeq (s\cc/I)^{\tens k} \rTTo^{b^\ce_k} s\cc/I \Bigr].
\end{multline*}

Let $\cq$ be a differential graded $\kk$\n-quiver. There is a free
\ainf-category $\cf\cq$, generated by $\cq$
\cite[Section~2.1]{LyuMan-freeAinf}. Let $R\subset s\cf\cq$ be a graded
subquiver. Denote by $I=(R)\subset s\cf\cq$ the graded subquiver
spanned by multiplying elements of $R$ with some elements of $s\cf\cq$
via several operations $b^{\cf\cq}_k$, $k>1$. It can be described as
follows. Let $t\in\ct^n_{\ge2}$ be a plane rooted tree with $i(t)=n$
input leaves. Each decomposition
\begin{equation}
(t,\le)=(1^{\sqcup\alpha_1}\sqcup\tree_{k_1}\sqcup1^{\sqcup\beta_1})
\cdot(1^{\sqcup\alpha_2}\sqcup\tree_{k_2}\sqcup1^{\sqcup\beta_2})\cdot
\ldots\cdot \tree_{k_N},
\label{eq-ord-tree-decomp-forest}
\end{equation}
of $t$ into the product of elementary forests gives a linear ordering
$\le$ of $t$. Here $\alpha_1+k_1+\beta_1=n$ and $N=|t|$ is the number
of internal vertices of $t$. An operation
\begin{equation}
b^{\cf\cq}_{(t,\le)} = \bigl( s\cf\cq^{\tens n}
\rTTo^{1^{\tens\alpha_1}\tens b^{\cf\cq}_{k_1}\tens1^{\tens\beta_1}}
s\cf\cq^{\tens\alpha_1+1+\beta_1}
\rTTo^{1^{\tens\alpha_2}\tens b^{\cf\cq}_{k_2}\tens1^{\tens\beta_2}}
\dots \rTTo^{b^{\cf\cq}_{k_N}} s\cf\cq \bigr)
\label{eq-bFQ-1bFQ1-bFQ}
\end{equation}
is associated with the linearly ordered tree $(t,\le)$. Different
choices of the ordering of $t$ change only the sign of the above map.
In particular, one may consider the canonical linear ordering $t_<$ of
$t$ \cite[Section~1.7]{LyuMan-freeAinf} and the corresponding map
$b^{\cf\cq}_{t_<}$. So the subquiver $I\subset s\cf\cq$ is defined as
\[ I = (R) = \sum_{t\in\ct_{\ge2}^{\alpha+1+\beta}}
\im(b^{\cf\cq}_{t_<}: s\cq^{\tens\alpha}\tens R\tens s\cq^{\tens\beta}
\to s\cf\cq),
\]
where the summation goes over all $\alpha,\beta\in\ZZ_{\ge0}$ and all
trees with $\alpha+1+\beta$ input leaves.

\begin{proposition}\label{pro-Rb1-subset(R)}
Let $R\subset s\cf\cq$ be a graded subquiver such that
$Rb^{\cf\cq}_1\subset(R)=I$. Then $Ib^{\cf\cq}_1\subset I$, $I$ is an
\ainf-ideal of $\cf\cq$, and $\ce=\cf\cq/s^{-1}I$ is an \ainf-category.
\end{proposition}

\begin{proof}
Clearly, $I$ is closed under multiplications $b_k$, $k>1$, with
elements of $s\cf\cq$.

Let us prove that for all $t\in\ct_{\ge2}$
\begin{equation}
\im(b^{\cf\cq}_{t_<}: s\cq^{\tens\alpha}\tens R\tens s\cq^{\tens\beta}
\to s\cf\cq)b^{\cf\cq}_1 \subset I
\label{eq-Im-bQRQFQ-b1I}
\end{equation}
using induction on $|t|$. This holds for $|t|=0$, $t=\big|$ by
assumption. Let $|t|=N>0$ and assume that \eqref{eq-Im-bQRQFQ-b1I}
holds for all $t'\in\ct_{\ge2}$ with $|t'|<N$. The tree $t$ can be
presented as $t=(t_1\sqcup\dots\sqcup t_k)\tree_k$ for some $k>1$. We
have
 $b^{\cf\cq}_{t_<}=\pm(b^{\cf\cq}_{t_1}\tdt
 b^{\cf\cq}_{t_k})b^{\cf\cq}_k$,
$|t_i|<N$ and
\[ b^{\cf\cq}_kb^{\cf\cq}_1 = - \sum_{a+q+c=k}^{a+c>0}
(1^{\tens a}\tens b^{\cf\cq}_q\tens1^{\tens c})b^{\cf\cq}_{a+1+c}.
\]
One of the $a+1+c$ factors of
\begin{equation}
(s\cq^{\tens\alpha}\tens R\tens s\cq^{\tens\beta})
(b^{\cf\cq}_{t_1}\tdt b^{\cf\cq}_{t_k})
(1^{\tens a}\tens b^{\cf\cq}_q\tens1^{\tens c})
\label{eq-QRQ-bb-1b1}
\end{equation}
is contained in $I$ (for $q=1$ this is the induction assumption).
Hence,
\[ (s\cq^{\tens\alpha}\tens R\tens s\cq^{\tens\beta})
(b^{\cf\cq}_{t_1}\tdt b^{\cf\cq}_{t_k})
(1^{\tens a}\tens b^{\cf\cq}_q\tens1^{\tens c})b^{\cf\cq}_{a+1+c}
\subset I,
\]
and the inclusion $Ib^{\cf\cq}_1\subset I$ follows by induction.

Therefore, $I$ is stable under all $b^{\cf\cq}_k$, $k\ge1$, so it is an
\ainf-ideal, and $\ce=\cf\cq/s^{-1}I$ is an \ainf-category.
\end{proof}

\subsection{\texorpdfstring{$A_\infty$-functors from a quotient of a
 free $A_\infty$-category}
 {A8-functors from a quotient of a free A8-category}}
The following statement is Proposition~2.3 from \cite{LyuMan-freeAinf}.

\begin{proposition}\label{pro-extension-f1}
Let $\cq$ be a differential graded quiver, and let $\ca$ be an
\ainf-category. \ainf-functors $f:\cf\cq\to\ca$ are in bijection with
sequences $(f'_1,f_k)_{k>1}$, where $f'_1:s\cq\to(s\ca,b_1)$ is a chain
morphism of differential graded quivers with the underlying mapping of
objects $\Ob f:\Ob\cq\to\Ob\ca$ and $f_k:T^ks\cf\cq\to s\ca$ are
$\kk$\n-quiver morphisms of degree 0 with the same underlying map
$\Ob f$ for all $k>1$. The morphisms $f_k$ are components of $f$ for
$k>1$. The component $f_1:s\cf\cq\to s\ca$ is an extension of $f'_1$.
\end{proposition}

In fact, it is shown in \cite[Proposition~2.3]{LyuMan-freeAinf} that
such a sequence $(f'_1,f_k)_{k>1}$ extends to a sequence of components
of an \ainf-functor $(f_1,f_k)_{k>1}$ in a unique way.

We are going to extend this description to quotients of free
\ainf-categories. Let a graded subquiver $R\subset s\cf\cq$ satisfy the
assumption $Rb^{\cf\cq}_1\subset(R)=I$. Denote by
$\pi:\cf\cq \rEpi \ce=\cf\cq/s^{-1}I$ the natural projection strict
\ainf-functor.

\begin{proposition}\label{pro-functor-factorizes}
An \ainf-functor $f:\cf\cq\to\ca$ factorizes as
$f=\bigl(\cf\cq \rTTo^\pi \ce \rTTo^{\tilde{f}} \ca\bigr)$ for some
(unique) \ainf-functor $\tilde{f}:\ce\to\ca$ if and only if the
following two conditions are satisfied:
\begin{enumerate}
\item $Rf_1=0$;
\item
 $(s\cf\cq^{\tens\alpha}\tens I\tens s\cf\cq^{\tens\beta})
 f_{\alpha+1+\beta}=0$
for all $\alpha,\beta\in\ZZ_{\ge0}$ such that $\alpha+\beta>0$.
\end{enumerate}
\end{proposition}

\begin{proof}
If $f=\pi\tilde{f}$, then $f_k=\pi_1^{\tens k}\tilde{f}_k$ and the
above conditions are necessary.

Assume that the both conditions are satisfied. Let us prove that
$If_1=0$. We are going to prove that for all $t\in\ct_{\ge2}$
\begin{equation}
\bigl(s\cq^{\tens\alpha}\tens R\tens s\cq^{\tens\beta}
\rTTo^{b^{\cf\cq}_{t_<}} s\cf\cq \rTTo^{f_1} s\ca\bigr) = 0
\label{eq-QRQbFQf1A0}
\end{equation}
by induction on $|t|$. This holds for $|t|=0$, $t=\big|$ by assumption.
Let $|t|=N>0$ and assume that \eqref{eq-QRQbFQf1A0} holds for all
$t'\in\ct_{\ge2}$ with $|t'|<N$. The tree $t$ can be presented as
$t=(t_1\sqcup\dots\sqcup t_k)\tree_k$ for some $k>1$. We have
 $b^{\cf\cq}_{t_<}=\pm(b^{\cf\cq}_{t_1}\tdt
 b^{\cf\cq}_{t_k})b^{\cf\cq}_k$,
$|t_i|<N$ and
\begin{equation}
b^{\cf\cq}_kf_1 = \sum_{i_1+\dots+i_l=k}(f_{i_1}\tdt f_{i_l})b^\ca_l
- \sum_{a+q+c=n}^{a+c>0}
(1^{\tens a}\tens b^{\cf\cq}_q\tens1^{\tens c})f_{a+1+c}.
\label{eq-bkf1-ffb-1b1f}
\end{equation}
One of the $a+1+c$ factors of \eqref{eq-QRQ-bb-1b1} is contained in $I$
and also one of the $l$ factors of
\[ (s\cq^{\tens\alpha}\tens R\tens s\cq^{\tens\beta})
(b^{\cf\cq}_{t_1}\tdt b^{\cf\cq}_{t_k}) (f_{i_1}\tdt f_{i_l})
\]
is contained in $I$ (for $i_j=1$ this is the induction assumption).
Hence, the right hand side of \eqref{eq-bkf1-ffb-1b1f} is contained in
$I$, and $If_1\subset I$ follows by induction. Therefore,
\[ (s\cf\cq^{\tens\alpha}\tens I\tens s\cf\cq^{\tens\beta})f \subset
\sum_{i_1+\dots+i_l=\alpha+1+\beta}
(s\cf\cq^{\tens\alpha}\tens I\tens s\cf\cq^{\tens\beta})
(f_{i_1}\tdt f_{i_l}) = 0
\]
for all $\alpha,\beta\ge0$. Clearly, $f$ factorizes as
$f=\pi\tilde{f}$.
\end{proof}

\subsection{\texorpdfstring{Transformations from a
 free $A_\infty$-category}{Transformations from a free A8-category}}
The following statement is Proposition~2.8 from \cite{LyuMan-freeAinf}.

\begin{proposition}\label{pro-chain-to-Ainf-FQ-A}
Let $\phi,\psi:\cf\cq\to\ca$ be \ainf-functors. For an arbitrary
complex $P$ of $\kk$\n-modules chain maps
$u:P\to sA_\infty(\cf\cq,\ca)(\phi,\psi)$ are in bijection with the
following data: $(u',u_k)_{k>1}$
\begin{enumerate}
\item a chain map $u':P\to sA_1(\cq,\ca)(\phi,\psi)$,

\item $\kk$\n-linear maps
\[ u_k: P \to \prod_{X,Y\in\Ob\cq}
\uCom\bigl((s\cf\cq)^{\tens k}(X,Y),s\ca(X\phi,Y\psi)\bigr)
\]
of degree 0 for all $k>1$.
\end{enumerate}
The bijection maps $u$ to $u_k=u\cdot\pr_k$,
\begin{equation}
u' = \bigl(P \rTTo^u sA_\infty(\cf\cq,\ca)(\phi,\psi)
\rEpi^{\restr_{\le1}} sA_1(\cf\cq,\ca)(\phi,\psi) \rEpi^\restr
sA_1(\cq,\ca)(\phi,\psi) \bigr).
\label{eq-u'-urestr2}
\end{equation}
The inverse bijection can be recovered from the recurrent formula
\begin{multline}
(-)^pb^{\cf\cq}_k(pu_1) = -(pd)u_k + \sum_{a+q+c=k}^{\alpha,\beta}
(\phi_{a\alpha}\tens pu_q\tens\psi_{c\beta})b^\ca_{\alpha+1+\beta} \\
-(-)^p \sum_{\alpha+q+\beta=k}^{\alpha+\beta>0} (1^{\tens\alpha}\tens
b^{\cf\cq}_q\tens1^{\tens\beta})(pu_{\alpha+1+\beta}):
(s\cf\cq)^{\tens k} \to s\ca,
\label{eq-bFQk-pu1}
\end{multline}
where $k>1$, \(p\in P\), and \(\phi_{a\alpha}\), \(\psi_{c\beta}\) are
matrix elements of $\phi$, $\psi$.
\end{proposition}

The following statement is Corollary~2.10 from \cite{LyuMan-freeAinf}.

\begin{corollary}\label{cor-chain-to-Ainf-FQ-A-null-homotopic}
Let $\phi,\psi:\cf\cq\to\ca$ be \ainf-functors. Let $P$ be a complex of
$\kk$\n-modules. Let $w:P\to sA_\infty(\cf\cq,\ca)(\phi,\psi)$ be a
chain map. The set (possibly empty) of homotopies
$h:P\to sA_\infty(\cf\cq,\ca)(\phi,\psi)$, $\deg h=-1$, such that
$w=dh+hB_1$ is in bijection with the set of data $(h',h_k)_{k>1}$,
consisting of
\begin{enumerate}
\item a homotopy $h':P\to sA_1(\cq,\ca)(\phi,\psi)$, $\deg h'=-1$, such
that $dh'+h'B_1=w'$, where
\begin{equation*}
w' = \bigl(P \rTTo^w sA_\infty(\cf\cq,\ca)(\phi,\psi)
\rEpi^{\restr_{\le1}} sA_1(\cf\cq,\ca)(\phi,\psi) \rEpi^\restr
sA_1(\cq,\ca)(\phi,\psi) \bigr);
\end{equation*}

\item $\kk$\n-linear maps
\[ h_k: P \to \prod_{X,Y\in\Ob\cq}
\uCom\bigl((s\cf\cq)^{\tens k}(X,Y),s\ca(X\phi,Y\psi)\bigr)
\]
of degree $-1$ for all $k>1$.
\end{enumerate}
The bijection maps $h$ to $h_k=h\cdot\pr_k$,
\begin{equation}
h' = \bigl(P \rTTo^h sA_\infty(\cf\cq,\ca)(\phi,\psi)
\rEpi^{\restr_{\le1}} sA_1(\cf\cq,\ca)(\phi,\psi) \rEpi^\restr
sA_1(\cq,\ca)(\phi,\psi) \bigr).
\label{eq-h-Ph-Ai(FQA)-A1(FQA)-A1(QA)}
\end{equation}
The inverse bijection can be recovered from the recurrent formula
\begin{multline}
(-)^pb_k(ph_1) = pw_k -(pd)h_k - \sum_{a+q+c=k}^{\alpha,\beta}
(\phi_{a\alpha}\tens ph_q\tens\psi_{c\beta})b_{\alpha+1+\beta} \\
-(-)^p \sum_{a+q+c=k}^{a+c>0}
(1^{\tens a}\tens b_q\tens1^{\tens c})(ph_{a+1+c}):
(s\cf\cq)^{\tens k} \to s\ca,
\label{eq-bk-ph1}
\end{multline}
where $k>1$, \(p\in P\), and \(\phi_{a\alpha}\), \(\psi_{c\beta}\) are
matrix elements of $\phi$, $\psi$.
\end{corollary}

\subsection{\texorpdfstring{Transformations from a quotient of a free
 $A_\infty$-category}{Transformations from a quotient of a free
 A8-category}}\label{sec-tran-quo-free}
We are going to extend the above description to quotients of free
\ainf-categories. Assume that a graded subquiver $R\subset s\cf\cq$
satisfies $Rb^{\cf\cq}_1\subset(R)=I$. Let $\ca$ be an \ainf-category.
Composition with the projection \ainf-functor
$\pi:\cf\cq \rEpi \ce=\cf\cq/s^{-1}I$ gives a strict \ainf-functor
$L^\pi=(\pi\boxtimes1)M:A_\infty(\ce,\ca)\to A_\infty(\cf\cq,\ca)$. It
is injective on objects and morphisms, that is, both maps
$\Ob L^\pi:\phi\mapsto\pi\phi$ and
 $L^\pi_1:sA_\infty(\ce,\ca)(\phi,\psi)\to
 sA_\infty(\cf\cq,\ca)(\pi\phi,\pi\psi)$,
$r\mapsto\pi r$ are injective. We are going to characterize the
subcomplex
 $sA_\infty(\ce,\ca)(\phi,\psi) \rMono
 sA_\infty(\cf\cq,\ca)(\pi\phi,\pi\psi)$.

\begin{proposition}\label{pro-chain-map-in-subcomplex}
Let $P$ be a complex of\/ $\kk$\n-modules, and let
$u:P\to sA_\infty(\cf\cq,\ca)(\pi\phi,\pi\psi)$ be a chain map. Denote
\[ u_k = u\cdot \pr_k : P \to \prod_{X,Y\in\Ob\cq}
\uCom\bigl((s\cf\cq)^{\tens k}(X,Y),s\ca(X\phi,Y\psi)\bigr),
\qquad k \ge 0.
\]
Then the image of $u$ is contained in the subcomplex
$sA_\infty(\ce,\ca)(\phi,\psi)$ if and only if the following two
conditions are satisfied:
\begin{enumerate}
\item $R(pu_1)=\im(pu_1:R(X,Y)\to s\ca(X\phi,Y\psi))=0$ for all
$p\in P$;
\item
$(s\cf\cq^{\tens\alpha}\tens I\tens s\cf\cq^{\tens\beta})(pu_k)=0$ for
all $p\in P$ and all $\alpha,\beta\ge0$ such that $\alpha+1+\beta=k>1$.
\end{enumerate}
\end{proposition}

\begin{proof}
The conditions are obviously necessary. Let us prove that they are
sufficient. First of all, we are going to show that $I(pu_1)=0$.
Namely, we claim that
\begin{equation}
(s\cq^{\tens\alpha}\tens R\tens s\cq^{\tens\beta})
b^{\cf\cq}_{t_<}(pu_1) = 0
\label{eq-QRQbpu10}
\end{equation}
for all trees $t\in\ct_{\ge2}^{\alpha+1+\beta}$. We prove it by
induction on $|t|$. For $|t|=0$, $t=\big|$ this holds by assumption 1.
Let $t\in\ct_{\ge2}$ be a tree with $|t|=N>0$ internal vertices. Assume
that \eqref{eq-QRQbpu10} holds for all $t'\in\ct_{\ge2}$ with $|t'|<N$.
The tree $t$ can be presented as $t=(t_1\sqcup\dots\sqcup t_k)\tree_k$
for some $k>1$, so
 $b^{\cf\cq}_{t_<}=\pm(b^{\cf\cq}_{t_1}\tdt
 b^{\cf\cq}_{t_k})b^{\cf\cq}_k$,
and $|t_i|<N$. Formula~\eqref{eq-bFQk-pu1} in the form
\begin{multline}
(-)^p b^{\cf\cq}_k(pu_1) = -(pd)u_k + \sum_{a+q+c=k}^{m,n}
(\pi_1^{\tens a}\phi_{am}\tens pu_q\tens\pi_1^{\tens c}\psi_{cn})
b^\ca_{m+1+n} \\
-(-)^p \sum_{a+q+c=k}^{a+c>0}
(1^{\tens a}\tens b^{\cf\cq}_q\tens1^{\tens c}) (pu_{a+1+c}):
s\cf_{t_1}\cq\tdt s\cf_{t_k}\cq \to s\ca
\label{eq-bFQk-pu1-pi}
\end{multline}
implies that
\begin{align*}
(s\cq^{\tens\alpha} &\tens R\tens s\cq^{\tens\beta})
b^{\cf\cq}_{t_<}(pu_1) =
(s\cq^{\tens\alpha}\tens R\tens s\cq^{\tens\beta})
(b^{\cf\cq}_{t_1}\tdt b^{\cf\cq}_{t_k})b^{\cf\cq}_k(pu_1) \\
&\subset
(s\cf\cq^{\tens\gamma}\tens I\tens s\cf\cq^{\tens\delta})[(pd)u_k] \\
&+ \sum (s\cq^{\tens\alpha}\tens R\tens s\cq^{\tens\beta})
(b^{\cf\cq}_{t_1}\tdt b^{\cf\cq}_{t_k})
(\pi_1^{\tens a}\phi_{am}\tens pu_q\tens\pi_1^{\tens c}\psi_{cn})
b^\ca_{m+1+n} \\
&+ \sum_{\lambda+\mu>0}
(s\cf\cq^{\tens\lambda}\tens I\tens s\cf\cq^{\tens\mu})
(pu_{\lambda+1+\mu}) = 0.
\end{align*}
Indeed, summands with $q=0$ vanish due to $I\pi_1=0$, summands with
$q=1$ vanish due to induction assumption~\eqref{eq-QRQbpu10}, and other
summands vanish due to condition~2.

Thus condition 2 holds not only for $k>1$ but for $k=1$ as well. At
last
\[ (s\cf\cq^{\tens\alpha}\tens I\tens s\cf\cq^{\tens\beta})(pu) \subset
\sum_{a+q+c=\alpha+1+\beta}^{m,n}
(s\cf\cq^{\tens\alpha}\tens I\tens s\cf\cq^{\tens\beta})
(\pi_1^{\tens a}\phi_{am}\tens pu_q\tens\pi_1^{\tens c}\psi_{cn}) = 0,
\]
and the proof is finished.
\end{proof}

Let us extend the description of homotopies between chain maps to the
case of quotient of a free \ainf-category. We keep the assumptions of
\secref{sec-tran-quo-free}.

\begin{corollary}\label{cor-null-homotopic-in-subcomplex}
Let $P$ be a complex of\/ $\kk$\n-modules, and let
$v:P\to sA_\infty(\ce,\ca)(\phi,\psi)$ be a chain map. Denote
\[ w = \bigl(P \rTTo^v sA_\infty(\ce,\ca)(\phi,\psi) \rMono
sA_\infty(\cf\cq,\ca)(\pi\phi,\pi\psi)\bigr).
\]
Let $h:P\to sA_\infty(\cf\cq,\ca)(\pi\phi,\pi\psi)$ be a homotopy,
$\deg h=-1$, $w=dh+hB_1$. Then the image of $h$ is contained in the
subcomplex $sA_\infty(\ce,\ca)(\phi,\psi)$ if and only if it factorizes
as
\[ h = \bigl(P \rTTo^\eta sA_\infty(\ce,\ca)(\phi,\psi) \rMono
sA_\infty(\cf\cq,\ca)(\pi\phi,\pi\psi)\bigr),
\]
where $\deg\eta=-1$, $v=d\eta+\eta B_1$, or if and only if the
following two conditions are satisfied:
\begin{enumerate}
\item $R(ph_1)=\im(ph_1:R(X,Y)\to s\ca(X\phi,Y\psi))=0$ for all
$p\in P$;
\item
$(s\cf\cq^{\tens\alpha}\tens I\tens s\cf\cq^{\tens\beta})(ph_k)=0$ for
all $p\in P$ and all $\alpha,\beta\ge0$ such that $\alpha+1+\beta=k>1$.
\end{enumerate}
\end{corollary}

\begin{proof}
Given pair \((w,h)\) defines a chain map
 \(\overline{w}:\Cone(\id_P)\to sA_\infty(\cf\cq,\ca)(\pi\phi,\pi\psi)\),
\((q,ps)\mapsto qw+ph\), such that
\[ w = \bigl(P \rTTo^{\inj_1} P\oplus P[1] = \Cone(\id_P)
\rTTo^{\overline{w}} sA_\infty(\cf\cq,\ca)(\pi\phi,\pi\psi)\bigr).
\]
The image of $h$ is contained in \(sA_\infty(\ce,\ca)(\phi,\psi)\) if
and only if the image of $\overline{w}$ is contained in
\(sA_\infty(\ce,\ca)(\phi,\psi)\). By
\propref{pro-chain-map-in-subcomplex} this is equivalent to conditions:
\begin{alignat*}2
&1'.\quad R(qw_1)=0; &\qquad &1''.\quad R(ph_1)=0; \\
&2'.\quad
(s\cf\cq^{\tens\alpha}\tens I\tens s\cf\cq^{\tens\beta})(qw_k)=0;
&\qquad &2''.\quad
(s\cf\cq^{\tens\alpha}\tens I\tens s\cf\cq^{\tens\beta})(ph_k)=0
\end{alignat*}
for all $q,p\in P$ and all $\alpha,\beta\ge0$ such that
$\alpha+1+\beta=k>1$. However, the conditions $1'$ and $2'$ are
satisfied automatically by \propref{pro-chain-map-in-subcomplex}
applied to $w$ and $v$.
\end{proof}

\section{A simple example}\label{sec-simple-example}
We want to consider an example, which is almost tautological. The
non-trivial part of it is the concrete choice of a system of relations
$R$ generating an ideal. Let $\cd$ be an \ainf-category. We view it as
a differential graded quiver and construct the free \ainf-category
$\cf\cd$ out of it. We choose the following subquiver of relations:
$R_\cd=\sum_{n\ge2}\im(b^{\cf\cd}_n-b^\cd_n:s\cd^{\tens n}\to s\cf\cd)$.
More precisely, a map $\delta_n:(s\cd)^{\tens n}\to s\cfd\cd$, $n\ge2$,
is defined as the difference
\[ \delta_n = \bigl((s\cd)^{\tens n} = (s\cfd_|\cd)^{\tens n}
\rTTo^{b_n^{\cfd\cd}}_{s^{-1}} s\cfd_{\tree_n}\cd \rMono s\cfd\cd\bigr)
- \bigl((s\cd)^{\tens n} \rTTo^{b_n^\cd} s\cd = s\cfd_|\cd \rMono
s\cfd\cd\bigr),
\]
and $R_\cd=\sum_{n\ge2}\im(\delta_n)$.

\begin{lemma}\label{lem-I-ideal}
The subquiver $I=(R_\cd)\subset s\cf\cd$ is an \ainf-ideal.
\end{lemma}

\begin{proof}
According to \propref{pro-Rb1-subset(R)} it suffices to check that
$R_\cd b^{\cf\cd}_1\subset(R_\cd)$. As
\[ (b_n^{\cf\cd}-b^\cd_n)b_1^{\cf\cd} = - b^\cd_nb_1^{\cf\cd}
- \sum_{a+k+c=n}^{a+c>0}
(1^{\tens a}\tens b_k^{\cfd\cd}\tens1^{\tens c}) b_{a+1+c}^{\cfd\cd},
\]
we have
\begin{align*}
\im(b^{\cf\cd}_n-b^\cd_n)b^{\cf\cd}_1 &\subset
- \im\Bigl( b^\cd_nb_1^{\cf\cd} + \sum_{a+k+c=n}^{a+c>0}
(1^{\tens a}\tens b_k^\cd\tens1^{\tens c}) b_{a+1+c}^{\cfd\cd}\Bigr)
+(R_\cd) \\
&\subset - \im\Bigl( b^\cd_nb_1^\cd + \sum_{a+k+c=n}^{a+c>0}
(1^{\tens a}\tens b_k^\cd\tens1^{\tens c})b_{a+1+c}^\cd\Bigr)+(R_\cd) =(R_\cd),
\end{align*}
and the lemma is proven.
\end{proof}

Let us describe the ideal $I=(R_\cd)\subset s\cfd\cd$. Let $t'\in\ct$ be a
tree, and let
$t=(1^{\sqcup\alpha}\sqcup\tree_k\sqcup1^{\sqcup\beta})\cdot t'_<$.
Define a map $\delta_{\alpha,k,t'}$ as the difference
\begin{multline*}
\delta_{\alpha,k,t'} = \bigl( (s\cd)^{\tens n}
\rTTo^{1^{\tens\alpha}\tens b_k^{\cfd\cd}\tens1^{\tens\beta}}
(s\cfd_|\cd)^{\tens\alpha}\tens
s\cfd_{\tree_k}\cd\tens(s\cfd_|\cd)^{\tens\beta}
 \rTTo^{b_{t'_<}^{\cfd\cd}} s\cfd_t\cd \rTTo s\cfd\cd\bigr) \\
- \bigl( (s\cd)^{\tens n}
 \rTTo^{1^{\tens\alpha}\tens b_k^\cd\tens1^{\tens\beta}}
(s\cd)^{\tens\alpha+1+\beta} = (s\cfd_|\cd)^{\tens\alpha+1+\beta}
\rTTo^{b_{t'_<}^{\cfd\cd}} s\cfd_{t'}\cd \rTTo s\cfd\cd\bigr)
\end{multline*}
for $k>1$, $\alpha+k+\beta=n$. Clearly,
$\sum\im(\delta_{\alpha,k,t'})=I$. These compositions can be
simplified. Let $h$ be the height of the distinguished vertex $\tree_k$
in $t_<$. Then the above difference equals
\begin{multline*}
\delta_{\alpha,k,t'} = \bigl( (s\cd)^{\tens\alpha+k+\beta}
= (s\cfd_|\cd)^{\tens n} \rTTo^{(-)^{|t|-h}s^{-|t|}}
s\cfd_t\cd \rTTo s\cfd\cd\bigr) \\
- \bigl( (s\cd)^{\tens\alpha+k+\beta}
 \rTTo^{1^{\tens\alpha}\tens b_k^\cd\tens1^{\tens\beta}}
(s\cd)^{\tens\alpha+1+\beta} = (s\cfd_|\cd)^{\tens\alpha+1+\beta}
\rTTo^{s^{-|t'|}} s\cfd_{t'}\cd \rTTo s\cfd\cd\bigr).
\end{multline*}

With the identity map $\id:\cd\to\cd$ is associated a strict
\ainf-functor $\widehat{\id}:\cf\cd\to\cd$
\cite[Section~2.6]{LyuMan-freeAinf}. Its first component equals
\begin{equation}
\widehat{\id}_1 = \bigl(s\cfd_t\cd
\rTTo^{[b^{\cfd\cd}_{(t,\le)}]^{-1}} T^{i(t)}s\cd
\rTTo^{b^\cd_{(t,\le)}} s\cd \bigr)
\label{eq-id1-sFD-TsD-sD}
\end{equation}
for each linear ordering $(t,\le)$ of a plane rooted tree
$t\in\ct^n_{\ge2}$ with $i(t)=n$ input leaves. Note that in the above
formula $b^{\cfd\cd}_{(t,\le)}=\pm s^{-|t|}$. In particular, for the
canonical linear ordering $t_<$ the formula becomes
\[ \widehat{\id}_1 = \bigl(s\cfd_t\cd \rTTo^{s^{|t|}} T^{i(t)}s\cd
\rTTo^{b^\cd_{t_>}} s\cd \bigr).
\]
Here
\[ b^\cd_{(t,\le)} = \bigl( s\cd^{\tens n}
\rTTo^{1^{\tens\alpha_1}\tens b^\cd_{k_1}\tens1^{\tens\beta_1}}
s\cd^{\tens\alpha_1+1+\beta_1}
\rTTo^{1^{\tens\alpha_2}\tens b^\cd_{k_2}\tens1^{\tens\beta_2}} \dots
\rTTo^{b^\cd_{k_N}} s\cd \bigr)
\]
corresponds to the ordered decomposition
\[ (t,\le)=(1^{\sqcup\alpha_1}\sqcup\tree_{k_1}\sqcup1^{\sqcup\beta_1})
\cdot(1^{\sqcup\alpha_2}\sqcup\tree_{k_2}\sqcup1^{\sqcup\beta_2})\cdot
\ldots\cdot \tree_{k_N},
\]
of $(t,\le)$ into the product of forests, $\alpha_1+k_1+\beta_1=n$, and
$b^{\cfd\cd}_{(t,\le)}$ has a similar meaning.

\begin{proposition}\label{pro-FD/(R)-D}
The \ainf-category $\cf\cd/s^{-1}(R_\cd)$ is isomorphic to $\cd$.
\end{proposition}

\begin{proof}
Let $\ce=\cf\cd/s^{-1}(R_\cd)$ be the quotient category. The projection map
$\pi_1:s\cf\cd\to s\ce$ with the underlying map of objects
$\Ob\pi=\id_{\Ob\cd}$ determines a strict \ainf-functor
$\pi:\cf\cd\to\ce$. The embedding
$\iota_1=\bigl(s\cd \rMono s\cf\cd \rEpi^{\pi_1} s\ce\bigr)$ with the
underlying identity map of objects $\Ob\iota=\id_{\Ob\cd}$ determines a
strict \ainf-functor $\iota:\cd\to\ce$. Indeed,
$\iota_1^{\tens n}b^\ce_n=b^\cd_n\iota_1:s\cd^{\tens n}\to s\ce$, for
$\im(b^{\cf\cd}_n-b^\cd_n)=\im\delta_n\subset I\subset s\cf\cd$.

We claim that the \ainf-functor $\wh{\id}:\cf\cd\to\cd$ factorizes as
$\widehat{\id}=\bigl(\cf\cd \rTTo^\pi \ce \rTTo^{\wt{\id}} \cd\bigr)$
for some \ainf-functor $\wt{\id}:\ce\to\cd$. Indeed, both conditions of
\propref{pro-functor-factorizes} are satisfied. The second is obvious
since $\wh{\id}$ is strict. The first condition
$R_\cd\widehat{\id}_1=0$ follows from the computation
 $(b^{\cf\cd}_n-b^\cd_n)\widehat{\id}_1=
 \id_1^{\tens n}b^\cd_n-b^\cd_n\id_1=0:s\cd^{\tens n}\to s\cd$.
Thus, $\wt{\id}:\ce\to\cd$ is a strict \ainf-functor.

Both $\iota$ and $\wt{\id}$ induce the identity map on objects.
Clearly, $\iota\cdot\wt{\id}=\id_\cd$. Furthermore, $\iota$ is
surjective on morphisms because every element of $\cf_t\cd$ reduces to
an element of $\cd$ modulo $(R_\cd)$. Therefore, $\iota$ is invertible
and \ainf-functors $\iota$ and $\wt{\id}$ are strictly inverse to each
other.
\end{proof}

\begin{corollary}\label{cor-(R)id1-0}
$(R_\cd)\widehat{\id}_1=0$.
\end{corollary}

\section{\texorpdfstring{$A_\infty$-categories and quivers}
 {A8-categories and quivers}}
\label{sec-A-categories-and-quivers}
\begin{definition}[Pseudounital $A_\infty$-categories]
 \label{def-pseudounital}
A \emph{pseudounital} structure of an \ainf-category $\cd$ is a choice
of an element $[\iota_X]\in H^{-1}(s\cd(X,X),b_1)$ for each object $X$
of $\cd$. By $\iota_X$ we mean a representative of the chosen
cohomology class, $\deg\iota_X=-1$. An \ainf-functor $f:\cd\to\ca$
between two pseudounital \ainf-categories is called \emph{pseudounital}
if it preserves the distinguished cohomology classes, that is,
$\iota^\ca_{Xf}-\iota^\cd_Xf_1\in\im b_1$ for all objects $X$ of $\cd$.
\end{definition}

The composition of pseudounital \ainf-functors is pseudounital as
well. The full subcategory of pseudounital \ainf-functors is denoted
$A_\infty^{\psi u}(\cd,\ca)\subset A_\infty(\cd,\ca)$.

A unital \ainf-category $\ca$ has a canonical pseudounital structure:
$\iota^\ca_X=\sS{_X}\uni^\ca_0$. An \ainf-functor $f:\cd\to\ca$ between
unital \ainf-categories is unital if and only if it is pseudounital
for the canonical pseudounital structures of $\cd$ and $\ca$
\cite[Definition~8.1]{Lyu-AinfCat}.

Let $\cd$ be a pseudounital \ainf-category, let $\cq$ be a differential
graded quiver with \(\Ob\cq=\Ob\cd\), and let \(\inj^\cd:\cd\to\cq\) be
an $A_1$\n-functor, such that \(\Ob\inj^\cd=\id_{\Ob\cd}\) and
\(\inj^\cd_1:s\cd\hookrightarrow s\cq\) is an embedding. Let $\ca$ be a
pseudounital \ainf-category.

\subsection{\texorpdfstring{$A_{\infty1}$-functors and transformations}
 {A81-functors and transformations}}
We define \ainf-category \(A_{\infty1}^{\psi u}(\cd,\cq;\ca)\) via
pull-back square
\begin{diagram}[LaTeXeqno]
A_{\infty1}^{\psi u}(\cd,\cq;\ca) \SEpbk & \rTTo & A_1(\cq,\ca) \\
\dTTo && \dTTo>{A_1(\inj^\cd,\ca)} \\
A_\infty^{\psi u}(\cd,\ca) & \rTTo^{\restr_{\le1}} & A_1(\cd,\ca)
\label{dia-A(EQA)}
\end{diagram}

In details, the objects of \(A_{\infty1}^{\psi u}(\cd,\cq;\ca)\) are
pairs \((f,f')\), where \(f:\cd\to\ca\) is a pseudounital
\ainf-functor, \(f':\cq\to\ca\) is an $A_1$\n-functor such that
\(\Ob f=\Ob f'\) and \(f_1=f'_1\big|_{s\cd}\). Morphisms of
\(A_{\infty1}^{\psi u}(\cd,\cq;\ca)\) from \((f,f')\) to \((g,g')\) are
pairs \((r,r')\), where \(r\in A_\infty(\cd,\ca)(f,g)\),
\(r'\in A_1(\cq,\ca)(f',g')\) are such that \(\deg r=\deg r'\),
\(r_0=r'_0\) and \(r_1=r'_1\big|_{s\cd}\). For any $n$\n-tuple of
composable morphisms
 \((r^j,p^j)\in
 sA_{\infty1}(\cd,\cq;\ca)((f^{j-1},g^{j-1}),(f^j,g^j))\),
\(1\le j\le n\), their $n$\n-th product is defined as
\[ [(r^1,p^1)\tdt(r^n,p^n)]B_n \overset{\text{def}}=
\bigl((r^1\tdt r^n)B_n,(p^1\tdt p^n)B_n\bigr).
\]
It is well-defined, because formulas for $B_k$ agree for all
$A_N$\n-categories, \(1\le N\le\infty\). The identity $B^2=0$ for
\(A_{\infty1}^{\psi u}(\cd,\cq;\ca)\) follows from that for
\(A_\infty^{\psi u}(\cd,\ca)\) and \(A_1(\cq,\ca)\). Thus,
\(A_{\infty1}^{\psi u}(\cd,\cq;\ca)\) is an \ainf-category.

\begin{proposition}
If $\ca$ is unital, then the \ainf-category
\(A_{\infty1}^{\psi u}(\cd,\cq;\ca)\) is unital as well.
\end{proposition}

\begin{proof}
Let us denote by \((1\boxtimes\uni^\ca)M\) the \ainf-transformation
 \(\id\to\id:A_{\infty1}^{\psi u}(\cd,\cq;\ca)\to
 A_{\infty1}^{\psi u}(\cd,\cq;\ca)\),
whose $n$\n-th component is
\[ (r^1,p^1)\tdt(r^n,p^n) \mapsto
\bigl((r^1\tdt r^n\boxtimes\uni^\ca)M_{n1},
(p^1\tdt p^n\boxtimes\uni^\ca)M_{n1}\bigr).
\]
It is well-defined, since formulas for multiplication $M$ agree for all
$A_N$\n-categories, \(1\le N\le\infty\). In a sense,
\((1\boxtimes\uni^\ca)M\) for \(A_{\infty1}^{\psi u}(\cd,\cq;\ca)\) is
determined by a pair of \ainf-transformations:
\((1\boxtimes\uni^\ca)M\) for \(A_\infty^{\psi u}(\cd,\ca)\) and for
\(A_1(\cq,\ca)\). Since the latter two \ainf-transformations are
natural and satisfy
\[ [(1\boxtimes\uni^\ca)M\tens(1\boxtimes\uni^\ca)M]B_2
\equiv (1\boxtimes\uni^\ca)M,
\]
the \ainf-transformation \((1\boxtimes\uni^\ca)M\) for
\(A_{\infty1}^{\psi u}(\cd,\cq;\ca)\) has the same properties.

We claim that \((1\boxtimes\uni^\ca)M\) is the unit transformation of
\(A_{\infty1}^{\psi u}(\cd,\cq;\ca)\) as defined in
\cite[Definition~7.6]{Lyu-AinfCat}. Indeed, it remains to prove that
chain endomorphisms
\begin{multline*}
a \overset{\text{def}}=
\bigl(1\tens\sS{_{(g,g')}}[(1\boxtimes\uni^\ca)M]_0\bigr)B_2, \quad
c \overset{\text{def}}=
\bigl(\sS{_{(f,f')}}[(1\boxtimes\uni^\ca)M]_0\tens1\bigr)B_2: \\
\Phi \overset{\text{def}}=
\bigl(sA_{\infty1}(\cd,\cq;\ca)((f,f'),(g,g')),B_1\bigr) \to \Phi
\end{multline*}
are homotopy invertible for all pairs \((f,f')\), \((g,g')\) of objects
of \(A_{\infty1}^{\psi u}(\cd,\cq;\ca)\). We have
\begin{align*}
(r,r')a &=\bigl((r\tens g\uni^\ca)B_2,(r'\tens g'\uni^\ca)B_2\bigr), \\
(r,r')c &=\bigl(r(f\uni^\ca\tens1)B_2,r'(f'\uni^\ca\tens1)B_2\bigr).
\end{align*}
As a graded $\kk$\n-module \(\Phi=\prod_{n=0}^\infty V_n\), where
\begin{align}
V_0 &= \prod_{X\in\Ob\cd} s\ca(Xf,Xg), \notag \\
V_1 &= \prod_{X,Y\in\Ob\cd} \uCom(s\cq(X,Y),s\ca(Xf,Yg)), \notag \\
V_n &= \prod_{X,Y\in\Ob\cd} \uCom((s\cd)^{\tens n}(X,Y),s\ca(Xf,Yg))
\qquad \text{for} \qquad n\ge2.
\label{eq-V0-V1-Vn}
\end{align}
Consider the decreasing filtration
 \(\Phi=\Phi_0\supset\Phi_1\supset\dots\supset
 \Phi_n\supset\Phi_{n+1}\supset\dots\)
of the complex $\Phi$, defined by
$\Phi_n=0\times\dots\times0\times\prod_{m=n}^\infty V_m$. We may write
\begin{align*}
\Phi_1 &= \{ (r,r')\in\Phi \mid r_0 = 0 \}, \\
\Phi_n &= \{ (r,0)\in\Phi \mid \forall l<n \quad r_l = 0 \}
\qquad \text{for} \qquad n\ge2.
\end{align*}
The differential \(B_1\) preserves the submodules \(\Phi_n\). It
induces the differential
\[ d = \uCom(1,b_1) - \sum_{\alpha+1+\beta=n}
\uCom(1^{\tens\alpha}\tens b_1\tens1^{\tens\beta},1)
\]
in the quotient \(V_n=\Phi_n/\Phi_{n+1}\). Due to
\[ [(r\tens g\uni^\ca)B_2]_k = \sum_l (r\tens g\uni^\ca)\theta_{kl}b_l
= \sum_{a+p+c+q+e=k}^{\alpha,\gamma,\eps}
(f_{a\alpha}\tens r_p\tens g_{c\gamma}\tens(g\uni^\ca)_q\tens g_{e\eps})
b_{\alpha+\gamma+\eps+2}
\]
and similar formulas, the endomorphisms \(a,c:\Phi\to\Phi\) preserve
the subcomplexes \(\Phi_n\). They induce the endomorphisms
\(\gr a,\gr c:V_n\to V_n\) in the quotient complex $V_n$:
\begin{align*}
(r'_1)\gr a &= \prod_{X,Y\in\Ob\cd}
(\sS{_{X,Y}}r'_1\tens\sS{_{Yg}}\uni^\ca_0)b^\ca_2,\qquad r'_1\in V_1,\\
(r'_1)\gr c &= \prod_{X,Y\in\Ob\cd}
\sS{_{X,Y}}r'_1(\sS{_{Xf}}\uni^\ca_0\tens1)b^\ca_2, \\
(r_n)\gr a &= \prod_{X,Y\in\Ob\cd}
(\sS{_{X,Y}}r_n\tens\sS{_{Yg}}\uni^\ca_0)b^\ca_2,
\qquad r_n\in V_n, \quad n = 0 \text{ or } n \ge 2. \\
(r_n)\gr c &= \prod_{X,Y\in\Ob\cd}
\sS{_{X,Y}}r_n(\sS{_{Xf}}\uni^\ca_0\tens1)b^\ca_2.
\end{align*}
Due to unitality of $\ca$, for each pair $X$, $Y$ of objects of $\cd$
there exist $\kk$\n-linear maps
\(\sS{_{X,Y}}h,\sS{_{X,Y}}h':s\ca(Xf,Yg)\to s\ca(Xf,Yg)\) of degree
$-1$ such that
\begin{align}
(1\tens\sS{_{Yg}}\uni^\ca_0)b^\ca_2 &=
1 + \sS{_{X,Y}}h\cdot d + d\cdot\sS{_{X,Y}}h, \label{eq-1i0b2-hd-dh} \\
(\sS{_{Xf}}\uni^\ca_0\tens1)b^\ca_2 &=
-1 + \sS{_{X,Y}}h'\cdot d + d\cdot\sS{_{X,Y}}h'. \label{eq-i01b2-hd-dh}
\end{align}
We equip the $\kk$\n-modules $\Phi^d=\prod_{n=0}^{\infty}V_n^d$ with
the topology of the product of discrete Abelian groups $V_n^d$. Thus
the \(\kk\)\n-submodules
$\Phi_m^d=0^{m-1}\times\prod_{n=m}^{\infty}V_n^d$ form a basis of
neighborhoods of 0 in $\Phi^d$. Continuous maps \(A:V^d\to V^{d+p}\)
are identified with $\NN\times\NN$-matrices of linear
maps $A_{nm}:V_n^d\to V_m^{d+p}$
with finite number of non-vanishing elements in each column.

In particular, the maps \(B_1:\Phi^d\to\Phi^{d+1}\) are continuous for
all \(d\in\ZZ\). Let us introduce continuous $\kk$\n-linear maps
$H,H':\prod_{n=0}^\infty V_n^d\to\prod_{n=0}^\infty V_n^{d-1}$ by
diagonal matrices $\sS{_{X,Y}}r_n\mapsto\sS{_{X,Y}}r_n\sS{_{X,Y}}h$,
$\sS{_{X,Y}}r_n\mapsto\sS{_{X,Y}}r_n\sS{_{X,Y}}h'$. We may view
$\gr a$, $\gr c$ as diagonal matrices and as the corresponding
continuous endomorphisms of \(\prod_{n=0}^{\infty}V_n^d\). Equations
\eqref{eq-1i0b2-hd-dh}, \eqref{eq-i01b2-hd-dh} can be written as
\[ \gr a = 1 + Hd + dH, \qquad \gr c = -1 + H'd + dH'.
\]
The continuous chain maps
\[ N = a - HB_1 - B_1H - 1, \quad N' = c - H'B_1 - B_1H' + 1:
\prod_{n=0}^{\infty}V_n^d \to \prod_{n=0}^{\infty}V_n^d
\]
have strictly upper triangular matrices. Therefore, the endomorphisms
\(1+N,-1+N':\prod_{n=0}^{\infty}V_n^d\to\prod_{n=0}^{\infty}V_n^d\) are
invertible. Their inverse maps correspond to well-defined
$\NN\times\NN$-matrices $\sum_{i=0}^\infty(-N)^i$ and
$-\sum_{i=0}^\infty(N')^i$. Since $a$, $c$ are homotopic to invertible
maps, they are homotopy invertible and the proposition is proven.
\end{proof}

\subsection{\texorpdfstring{A strict $A_\infty^u$-2-functor}
 {A strict A8u-2-functor}}
Let us describe a strict \ainfu-2-functor
\(F:A_\infty^u\to A_\infty^u\). It maps a unital \ainf-category $\ca$
to the unital \ainf-category \(A_{\infty1}^{\psi u}(\cd,\cq;\ca)\). The
strict unital \ainf-functor
\[ F_{\ca,\cb} = A_{\infty1}^{\psi u}(\cd,\cq;\_): A_\infty^u(\ca,\cb)
\to A_\infty^u\bigl(A_{\infty1}^{\psi u}(\cd,\cq;\ca),
A_{\infty1}^{\psi u}(\cd,\cq;\cb)\bigr)
\]
is specified by the following data. It maps an object \(f:\ca\to\cb\)
(a unital \ainf-functor) to the object
\[ A_{\infty1}^{\psi u}(\cd,\cq;f)
= (A_\infty^{\psi u}(\cd,f),A_1(\cq,f))
= ((1\boxtimes f)M,(1\boxtimes f)M),
\]
a unital \ainf-functor
 \(A_{\infty1}^{\psi u}(\cd,\cq;\ca)\to
 A_{\infty1}^{\psi u}(\cd,\cq;\cb)\),
which sends \((g,g')\in\Ob A_{\infty1}^{\psi u}(\cd,\cq;\ca)\) to
\((gf,g'f)\in\Ob A_{\infty1}^{\psi u}(\cd,\cq;\cb)\). Its $n$\n-th
component is
\[ [A_{\infty1}^{\psi u}(\cd,\cq;f)]_n: (r^1,p^1)\tdt(r^n,p^n) \mapsto
\bigl((r^1\tdt r^n\boxtimes f)M_{n0},
(p^1\tdt p^n\boxtimes f)M_{n0}\bigr).
\]
An \ainf-transformation \(q:f\to g:\ca\to\cb\) is mapped by
\([A_{\infty1}^{\psi u}(\cd,\cq;\_)]_1\) to the \ainf-transformation
\[ A_{\infty1}^{\psi u}(\cd,\cq;q)
= (A_\infty^{\psi u}(\cd,q),A_1(\cq,q))
= ((1\boxtimes q)M,(1\boxtimes q)M),
\]
whose $n$\n-th component is
\[ [A_{\infty1}^{\psi u}(\cd,\cq;q)]_n: (r^1,p^1)\tdt(r^n,p^n) \mapsto
\bigl((r^1\tdt r^n\boxtimes q)M_{n1},
(p^1\tdt p^n\boxtimes q)M_{n1}\bigr).
\]
Thus, a strict \ainf-functor
\(F_{\ca,\cb}=A_{\infty1}^{\psi u}(\cd,\cq;\_)\) is constructed. It is
unital, because the unit element \(f\uni^\cb\) of
\(f\in\Ob A_\infty^u(\ca,\cb)\) is mapped to the unit element
 \(A_{\infty1}^{\psi u}(\cd,\cq;f\uni^\cb)
 =((1\boxtimes f\uni^\cb)M,(1\boxtimes f\uni^\cb)M)\)
of \(A_{\infty1}^{\psi u}(\cd,\cq;f)\in\Ob A_\infty^u(F\ca,F\cb)\).

Necessary equation, given by diagram (3.1.1) of \cite{LyuMan-freeAinf}
follows from the same equation written for
\(A_\infty^{\psi u}(\cd,\_)\) and for \(A_1(\cq,\_)\). Therefore, the
strict \ainfu-2-functor $F$ is constructed.

\subsection{\texorpdfstring{An $A_\infty$-category freely generated
 over an $A_\infty$-category}
 {An A8-category freely generated over an A8-category}}
 \label{sec-Cones-freely-generated}
Let $\cd$ be a pseudounital \ainf-category, let $\cq$ be a differential
graded quiver with \(\Ob\cq=\Ob\cd\), equipped with a chain embedding
\(s\cd \rMono s\cq\), identity on objects. Assume that there exists a
graded $\kk$\n-subquiver $\cn\subset\cq$, \(\Ob\cn=\Ob\cd\), which is a
direct complement of $\cd$. Thus, $\cd\oplus\cn=\cq$ is an isomorphism
of graded $\kk$\n-linear quivers. Then there exists a differential in
the graded quiver \(s\cm=\cn\) and a chain map \(\alpha:s\cm\to s\cd\)
such that \(\Ob\alpha=\id_{\Ob\cd}\) and $\cq=\Cone\alpha$. Indeed, the
embedding \(\inj^\cd\) in the exact sequence
\[ 0 \to s\cd \rTTo^{\inj^\cd} s\cq \rTTo^{\pr^\cn} \cm[2] \to 0
\]
is a chain map. Thus the graded $\kk$\n-quiver
\(\cm[2]=s\cn=\Coker(\inj^\cd)\) acquires a differential
\(d^{\cm[2]}\), such that \(\pr^\cn\) is a chain map. The differential
in \(s\cq\) has the form
\begin{equation}
(t,ms)b^\cq_1 = (tb^\cd_1+m\alpha,msd^{\cm[2]})
= (tb^\cd_1+m\alpha,-md^{\cm[1]}s)
\label{eq-(tm)bQ1-cone}
\end{equation}
for \(t\in s\cd(X,Y)\), \(m\in\cm[1](X,Y)\), where
\(\alpha:s\cm(X,Y)\to s\cd(X,Y)\) are $\kk$\n-linear maps of degree 0.
The condition \((b^\cq_1)^2=0\) is equivalent to $\alpha$ being a chain
map. Therefore, \(s\cq=\Cone(\alpha:s\cm\to s\cd)\).

Define a pseudounital \ainf-category \(\ce=\cf\cq/s^{-1}(R_\cd)\),
where $R_\cd=\sum_{n\ge2}\im(\delta_n)$ for
\[ \delta_n = \bigl((s\cd)^{\tens n} \rMono (s\cf\cq)^{\tens n}
\rTTo^{b_n^{\cf\cq}} s\cf\cq\bigr)
- \bigl((s\cd)^{\tens n} \rTTo^{b_n^\cd} s\cd \rMono s\cf\cq\bigr).
\]
Repeating word by word the proof of \lemref{lem-I-ideal} we deduce that
$J=(R_\cd)\subset s\cfd\cq$ is an \ainf-ideal. The distinguished
elements \(\iota_X\in(s\ce)^{-1}(X,X)\) are those of \(\cd\subset\ce\).
Let $\ca$ be a pseudounital \ainf-category. There is the restriction
strict \ainf-functor
\[ \restr: A_\infty^{\psi u}(\ce,\ca)
\to A_{\infty1}^{\psi u}(\cd,\cq;\ca),
\qquad x \mapsto (x\big|_\cd,x\big|_\cq).
\]
If $\ca$ is unital, then the above \ainf-functor is unital, because the
unit element \(f\uni^\ca\) of \(f\in\Ob A_\infty^{\psi u}(\ce,\ca)\) is
mapped to the unit element
 \(\bigl((f\uni^\ca)\big|_\cd,(f\uni^\ca)\big|_\cq\bigr)
 =\bigl(f\big|_\cd\uni^\ca,f\big|_\cq\uni^\ca\bigr)\)
of \((f\big|_\cd,f\big|_\cq)\in\Ob A_{\infty1}^{\psi u}(\cd,\cq;\ca)\).

\begin{proposition}\label{pro-tilde-f-E-A}
The map
 \(\Ob\restr:\Ob A_\infty^{\psi u}(\ce,\ca)
 \to\Ob A_{\infty1}^{\psi u}(\cd,\cq;\ca)\)
is surjective. An object $(f,f')$ of
 \(A_{\infty1}^{\psi u}(\cd,\cq;\ca)\) is the restriction of a unique
pseudounital \ainf-functor $\wt{f}:\ce\to\ca$ such that
 \(\Ob\wt{f}=\Ob f\), \(\wt{f}_1\big|_{s\cq}=f'_1\),
\(\wt{f}_k\big|_{s\cd^{\tens k}}=f_k\), and \(\wt{f}_k\) vanishes on
all summands of \(T^ks\ce\) containing the factor \(s\cn\) for $k>1$.
\end{proposition}

\begin{proof}
Let us define an \ainf-functor \(\hat{f}:\cf\cq\to\ca\) via
\propref{pro-extension-f1} by the following data. On objects it is
\(\Ob\hat{f}=\Ob f=\Ob f'\), the restriction to $s\cq$ of the first
component is \(\hat{f}_1\big|_{s\cq}=f'_1\). On each direct summand of
$T^ks\cfd\cq$, $k>1$, containing the factor $s\cn$ we set
$\hat{f}_k=0$. On the direct summand $T^ks\cfd\cd$ of $T^ks\cfd\cq$ we
define
\[ \hat{f}_k = \bigl(T^ks\cfd\cd \rTTo^{\widehat{\id}_1^{\tens k}}
T^ks\cd \rTTo^{f_k} s\ca \bigr)
\]
for $k>1$, where \(\widehat{\id}_1\) is defined by
\eqref{eq-id1-sFD-TsD-sD}. These requirements specify $\hat{f}$
completely. It is pseudounital, since \(\hat{f}_1\big|_{s\cd}=f_1\) and
$f$ is pseudounital.

Let us prove that the \ainf-functor $\hat{f}$ factors as
$\hat{f}=\bigl(\cfd\cq \rEpi^\pi \ce \rTTo^{\wt{f}} \ca\bigr)$ for some
unique \ainf-functor $\wt{f}$. Denote \(J=(R_\cd)\subset s\cf\cq\). We
have to check conditions of \propref{pro-functor-factorizes}. The
second condition,
 $(s\cf\cq^{\tens\alpha}\tens J\tens s\cf\cq^{\tens\beta})
 \hat{f}_{\alpha+1+\beta}=0$
if $\alpha+\beta>0$, holds on direct summands of
$s\cf\cq^{\tens\alpha}\tens J\tens s\cf\cq^{\tens\beta}$, which contain
a factor $s\cn$ in some of $\cf\cq$. It holds also on summands of
$J$ of the form
$\im(b^{\cf\cq}_t:\dots\tens s\cn\tdt R_\cd\tens\dots\to s\cf\cq)$ or
$\im(b^{\cf\cq}_t:\dots\tens R_\cd\tdt s\cn\tens\dots\to s\cf\cq)$. We
have to verify that
 $(s\cf\cd^{\tens\alpha}\tens I\tens s\cf\cd^{\tens\beta})
 \hat{f}_{\alpha+1+\beta}=0$
if $\alpha+\beta>0$, where $I\subset s\cf\cd$ is the ideal described in
\lemref{lem-I-ideal}. This equation holds true because
\[ (s\cf\cd^{\tens\alpha}\tens I\tens s\cf\cd^{\tens\beta})
\hat{f}_{\alpha+1+\beta} =
(s\cf\cd^{\tens\alpha}\tens I\tens s\cf\cd^{\tens\beta})
\widehat{\id}_1^{\alpha+1+\beta} f_{\alpha+1+\beta} = 0
\]
due to equation $I\widehat{\id}_1=0$, obtained in
\corref{cor-(R)id1-0}. Therefore, the second condition of
\propref{pro-functor-factorizes} is verified.

Let us observe that the restriction of $\hat{f}$ to $\cfd\cd$ coincides
with $\cfd\cd \rTTo^{\widehat{\id}} \cd \rTTo^f \ca$. Indeed, their
components are equal,
$\hat{f}_k\big|_{\cfd\cd}=\widehat{\id}_1^{\tens k}\cdot f_k$ for
$k>1$, and
 $\hat{f}_1\big|_{s\cfd_|\cd}=f_1=\widehat{\id}_1\cdot f_1:s\cd\to
 s\ca$.
Hence, $\hat{f}\big|_{\cfd\cd}=\widehat{\id}\cdot f$ by
\propref{pro-extension-f1}.

In particular,
 $\hat{f}_1\big|_{s\cfd\cd}=\bigl(s\cfd\cd \rTTo^{\widehat{\id}_1} s\cd
 \rTTo^{f_1} s\ca\bigr)$.
Hence, $R_\cd\hat{f}_1=R_\cd\widehat{\id}_1f_1=0$ due to
\corref{cor-(R)id1-0}. Therefore, the first condition of
\propref{pro-functor-factorizes} is also verified and \(\wt{f}\)
exists. Uniqueness of such \(\wt{f}\) is obvious.
\end{proof}

The projection map $\pi_1:s\cf\cq\to s\ce$ with the underlying map of
objects $\Ob\pi=\id_{\Ob\cd}$ determines a strict \ainf-functor
$\pi:\cf\cq\to\ce$. The embedding
$\iota_1=\bigl(s\cd \rMono^i s\cf\cq \rEpi^{\pi_1} s\ce\bigr)$ with the
underlying identity map of objects $\Ob\iota=\id_{\Ob\cd}$ determines a
strict \ainf-functor $\iota:\cd\to\ce$. Indeed,
$\iota_1^{\tens n}b^\ce_n=b^\cd_n\iota_1:s\cd^{\tens n}\to s\ce$, for
$\im(b^{\cf\cd}_n-b^\cd_n)=\im\delta_n\subset J\subset s\cf\cq$. These
\ainf-functors produce other strict \ainf-functors for an arbitrary
\ainf-category $\ca$. For instance, the functor
\[ (\pi\boxtimes1)M: A_\infty(\ce,\ca) \to A_\infty(\cf\cq,\ca), \qquad
x \mapsto \pi x,
\]
injective on objects and morphisms, and the restriction \ainf-functor
\[ \restr = (\iota\boxtimes1)M: A_\infty(\ce,\ca)\to A_\infty(\cd,\ca),
\qquad y \mapsto \iota y = \overline{y}.
\]

\begin{theorem}\label{thm-restriction-equivalence-quiver}
Let \(\ce=\cf\cq/s^{-1}(R_\cd)\), where $\cd$, $\cq$ satisfy
assumptions of \secref{sec-Cones-freely-generated}. Then the
restriction \ainf-functor
\begin{equation}
\restr: A_\infty^{\psi u}(\ce,\ca)\to A_{\infty1}^{\psi u}(\cd,\cq;\ca)
\label{eq-restr-A(EA)-A(DQA)}
\end{equation}
is an \ainf-equivalence, surjective on objects. The chain surjections
$\restr_1$ admit a chain splitting.
\end{theorem}

\begin{proof}
Let us prove that the chain map
\begin{equation}
\restr_1: sA_\infty^{\psi u}(\ce,\ca)(f,g) \to
sA_{\infty1}^{\psi u}(\cd,\cq;\ca)
\bigl((f|_\cd,f|_\cq),(g|_\cd,g|_\cq)\bigr)
\label{eq-restr1-Apsiu(EA)-Apsiu(DQA)}
\end{equation}
is homotopy invertible for all pairs of pseudounital \ainf-functors
\(f,g:\ce\to\ca\). This will be achieved in a sequence of Lemmata.

The graded $\kk$\n-quiver decomposition \(\cq=\cd\oplus\cn\)
implies that the graded $\kk$\n-quiver \(\cf\cd\) is a direct summand
of \(\cf\cq\). The projection \(\pr^{\cf\cd}:s\cf\cq\to s\cf\cd\)
annihilates all summands with factors \(s\cn\). Define a degree 0 map
\[ \varpi = \bigl(s\cf\cq \rTTo^{\pr^{\cf\cd}} s\cf\cd
\rTTo^{\widehat{\id}_1} s\cd\bigr).
\]

\begin{lemma}\label{lem-unique-u}
There exists a unique chain map
\[ u:
sA_{\infty1}(\cd,\cq;\ca)\bigl((f|_\cd,f|_\cq),(g|_\cd,g|_\cq)\bigr)
\to sA_\infty(\cf\cq,\ca)(\pi f,\pi g)
\]
such that $u'$ obtained from $u$ via formula~\eqref{eq-u'-urestr2}
equals
\begin{equation*}
u' = \restr_\cq:
sA_{\infty1}(\cd,\cq;\ca)\bigl((f|_\cd,f|_\cq),(g|_\cd,g|_\cq)\bigr)
\to sA_1(\cq,\ca)(f|_\cq,g|_\cq), \quad (p,p') \mapsto p',
\end{equation*}
and for $k>1$ the maps $u_k$ are
\begin{multline*}
u_k = \Bigl[
sA_{\infty1}(\cd,\cq;\ca)\bigl((f|_\cd,f|_\cq),(g|_\cd,g|_\cq)\bigr)
\rTTo^{\pr_k} \prod_{X,Y\in\Ob\cd}\uCom(s\cd^{\tens k}(X,Y),s\ca(Xf,Yg))
\\
\rTTo^{\prod\uCom(\varpi^{\tens k},1)}
\prod_{X,Y\in\Ob\cd}\uCom(s\cf\cq^{\tens k}(X,Y),s\ca(Xf,Yg)) \Bigr],
\qquad (p,p') \mapsto p_k \mapsto \varpi^{\tens k}p_k.
\end{multline*}
\end{lemma}

\begin{proof}
Apply \propref{pro-chain-to-Ainf-FQ-A} to
 \(P=
sA_{\infty1}(\cd,\cq;\ca)\bigl((f|_\cd,f|_\cq),(g|_\cd,g|_\cq)\bigr)\).
\end{proof}

\begin{lemma}
The map $u$ from \lemref{lem-unique-u} takes values in
\[ sA_\infty(\ce,\ca)(f,g)\subset sA_\infty(\cf\cq,\ca)(\pi f,\pi g).
\]
\end{lemma}

\begin{proof}
Let us verify conditions of \propref{pro-chain-map-in-subcomplex}. We
have $J\varpi=0$. Indeed, $\varpi$ vanishes on summands of $J=(R_\cd)$ of
the form
$\im(b^{\cf\cq}_t:\dots\tens s\cn\tdt R_\cd\tens\dots\to s\cf\cq)$ or
$\im(b^{\cf\cq}_t:\dots\tens R_\cd\tdt s\cn\tens\dots\to s\cf\cq)$.
Looking at $I=J\cap\cf\cd$ we find that
$J\varpi=I\varpi=I\widehat{\id}_1=0$ by \corref{cor-(R)id1-0}.
Therefore,
\begin{equation}
(s\cf\cq^{\tens\alpha}\tens J\tens s\cf\cq^{\tens\beta})((p,p')u_k) =
(s\cf\cq^{\tens\alpha}\tens J\tens s\cf\cq^{\tens\beta})
\varpi^{\tens k}p_k = 0,
\label{eq-FQJFQ-varpi-p0}
\end{equation}
and the second condition of \propref{pro-chain-map-in-subcomplex} is
verified.

Let us check now that $R_\cd((p,p')u_1)=0$ for any element
\[ (p,p')\in
sA_{\infty1}(\cd,\cq;\ca)\bigl((f|_\cd,f|_\cq),(g|_\cd,g|_\cq)\bigr).
\]
That is,
\begin{equation}
(-)^p(i^{\tens k}b^{\cf\cq}_k - b^\cd_ki) ((p,p')u_1) = 0:
s\cd^{\tens k}(X,Y) \to s\ca(Xf,Yg)
\label{eq-ikbk-bki-pu0}
\end{equation}
for $k>1$, where $i:s\cd \rMono s\cf\cq$ is the embedding of
differential graded $\kk$\n-quivers. By definition \eqref{eq-bFQk-pu1}
\begin{multline*}
(-)^p i^{\tens k}b^{\cf\cq}_k((p,p')u_1)
= -i^{\tens k}[(pB_1,p'B_1)u_k] \\
+ \sum_{a+q+c=k}^{m,n} (i^{\tens a}\pi_1^{\tens a}f_{am}\tens
i^{\tens q}((p,p')u_q)\tens i^{\tens c}\pi_1^{\tens c}g_{cn})
b^\ca_{m+1+n} \\
-(-)^p i^{\tens k} \sum_{a+q+c=k}^{a+c>0}
(1^{\tens a}\tens b^{\cf\cq}_q\tens1^{\tens c}) ((p,p')u_{a+1+c}):
s\cd^{\tens k}(X,Y) \to s\ca(Xf,Yg).
\end{multline*}
For an arbitrary
 $(t,t')\in
 sA_{\infty1}(\cd,\cq;\ca)\bigl((f|_\cd,f|_\cq),(g|_\cd,g|_\cq)\bigr)$,
in particular for $(p,p')$ or $(pB_1,p'B_1)$, we have for $k>1$
\begin{equation}
i^{\tens k}((t,t')u_k) = (t,t')u_k\uCom(i^{\tens k},1)
= t_k\uCom(\varpi^{\tens k},1)\uCom(i^{\tens k},1)
= t_k\uCom(i^{\tens k}\varpi^{\tens k},1) = t_k
\label{eq-ik(tuk)-C2}
\end{equation}
due to relation $i\varpi=\id_{s\cd}$. For $k=0$ or 1 we also have
$i^{\tens k}((t,t')u_k)=t_k$. Indeed,
\begin{gather*}
(t,t')u_0 = (t,t')u'\pr_0 = t'_0 = t_0, \\
i((t,t')u_1) = \inj^\cd((t,t')u')\pr_1 = \inj^\cd t'_1 = t_1.
\label{eq-ik(tuk)-beta0}
\end{gather*}
Notice also that
$\bigl(s\cd \rMono^i s\cf\cq \rEpi^{\pi_1} s\ce\bigr)=\iota_1$, hence,
 $i^{\tens a}\pi_1^{\tens a}f_{am}=\iota_1^{\tens a}f_{am}
 =\overline{f}_{am}\overset{\text{def}}=(f|_\cd)_{am}$.
Due to already proven property \eqref{eq-FQJFQ-varpi-p0} we may replace
$i^{\tens q}b^{\cf\cq}_q$ in the last sum with $b^\cd_qi$. Therefore,
\begin{align*}
(-)^p i^{\tens k}b^{\cf\cq}_k((p&,p')u_1) = -(pB_1)_k
+\sum_{a+q+c=k}^{m,n}(\overline{f}_{am}\tens p_q\tens\overline{g}_{cn})
b^\ca_{m+1+n} \\
&\hspace*{6em}
-(-)^p \sum_{a+q+c=k}^{a+c>0}(1^{\tens a}\tens b^\cd_q\tens1^{\tens c})
i^{\tens a+1+c}((p,p')u_{a+1+c})\\
&= -(pB_1)_k + (pb^\ca)_k -(-)^p \sum_{a+q+c=k}
(1^{\tens a}\tens b^\cd_q\tens1^{\tens c})p_{a+1+c} +(-)^p b^\cd_kp_1\\
&= -\bigl(pb^\ca-(-)^pb^\cd p\bigr)_k + (pb^\ca)_k -(-)^p (b^\cd p)_k
+(-)^p b^\cd_ki((p,p')u_1) \\
&= (-)^p b^\cd_ki((p,p')u_1): s\cd^{\tens k}(X,Y) \to s\ca(Xf,Yg)
\end{align*}
and \eqref{eq-ikbk-bki-pu0} is proven. We conclude by
\propref{pro-chain-map-in-subcomplex} that there is a chain map $\Phi$
such that
\begin{multline}
u = \bigl( sA_{\infty1}^{\psi u}(\cd,\cq;\ca)
\bigl((f|_\cd,f|_\cq),(g|_\cd,g|_\cq)\bigr)
\rTTo^\Phi sA_\infty^{\psi u}(\ce,\ca)(f,g) \\
\rMono^{(\pi\boxtimes1)M_{01}} sA_\infty(\cf\cq,\ca)(\pi f,\pi g)
\bigr),
\label{eq-u-Phi-(pi1)M}
\end{multline}
so the lemma is proven.
\end{proof}

\begin{lemma}\label{lem-Phi-one-sided-inverse}
The map $\Phi$ from \eqref{eq-u-Phi-(pi1)M} is a one-sided inverse to
\(\restr_1\):
\begin{multline*}
\Bigl( sA_{\infty1}^{\psi u}(\cd,\cq;\ca)
\bigl((f|_\cd,f|_\cq),(g|_\cd,g|_\cq)\bigr)
\rTTo^\Phi sA_\infty^{\psi u}(\ce,\ca)(f,g) \\
\rTTo^{\restr_1} sA_{\infty1}^{\psi u}(\cd,\cq;\ca)
\bigl((f|_\cd,f|_\cq),(g|_\cd,g|_\cq)\bigr)
\Bigr) = \id.
\end{multline*}
\end{lemma}

\begin{proof}
Recall that $\restr_1$ is the componentwise map
\begin{align*}
\bigl((\iota\boxtimes1)M_{01},(j\boxtimes1)M_{01}\bigr):
sA_\infty^{\psi u}(\ce,\ca)(f,g) &\to
sA_{\infty1}(\cd,\cq;\ca)\bigl((f|_\cd,f|_\cq),(g|_\cd,g|_\cq)\bigr) \\
r=(r_k)_k &\mapsto (\iota r,jr)
= \bigl((\iota_1^{\tens k}r_k)_{k\ge0},(j_1^{\tens k}r_k)_{k=0,1}\bigr)
\end{align*}
Also
 $(\pi\boxtimes1)M_{01}:sA_\infty^{\psi u}(\ce,\ca)(f,g) \rMono
 sA_\infty(\cf\cq,\ca)(\pi f,\pi g)$,
$r=(r_k)_k\mapsto\pi r=(\pi_1^{\tens k}r_k)_k$ is componentwise.
Introduce another componentwise map of degree 0
\begin{align*}
L^i: sA_\infty(\cf\cq,\ca)(\pi f,\pi g) &\to
sA_{\infty1}(\cd,\cq;\ca)\bigl((f|_\cd,f|_\cq),(g|_\cd,g|_\cq)\bigr),\\
q = (q_k)_k &\mapsto
\bigl((i^{\tens k}q_k)_{k\ge0},(\inj^{\cq\,\tens k}q_k)_{k=0,1}\bigr).
\end{align*}
As $\iota_1=\bigl(s\cd \rMono^i s\cf\cq \rEpi^{\pi_1} s\ce\bigr)$ and
$j_1=\bigl(s\cq \rMono^{\inj^\cq} s\cf\cq \rEpi^{\pi_1} s\ce\bigr)$,
the lower triangle in the following diagram commutes:
\begin{diagram}[nobalance,LaTeXeqno]
sA_{\infty1}(\cd,\cq;\ca)\bigl((f|_\cd,f|_\cq),(g|_\cd,g|_\cq)\bigr)
&& \rLine && \HmeetV \\
& \rdTTo^\Phi && = & \dTTo>u \\
\dTTo<{\Phi\cdot\restr_1}>{\qquad=} && sA_\infty^{\psi u}(\ce,\ca)(f,g)
&& \\
& \ldTTo<{\restr_1} & = & \rdMono^{(\pi\boxtimes1)M_{01}} & \\
sA_{\infty1}(\cd,\cq;\ca)\bigl((f|_\cd,f|_\cq),(g|_\cd,g|_\cq)\bigr)
&& \lTTo^{L^i} && sA_\infty(\cf\cq,\ca)(\pi f,\pi g)
\label{dia-3triangles-Phi}
\end{diagram}
Thus the whole diagram is commutative and $\Phi\cdot\restr_1=uL^i$. We
have proved in \eqref{eq-ik(tuk)-C2} and \eqref{eq-ik(tuk)-beta0} that
for all
 $(p,p')\in
 sA_{\infty1}(\cd,\cq;\ca)\bigl((f|_\cd,f|_\cq),(g|_\cd,g|_\cq)\bigr)$
and all $k\in\ZZ_{\ge0}$ we have $i^{\tens k}((p,p')u_k)=p_k$.
Similarly,
\[ \inj^\cq[(p,p')u_1] = (p,p')u'\pr_1 = p'_1.
\]
Therefore,
\[ (p,p')\Phi\restr_1 =
\bigl( (i^{\tens k}(p,p')u_k)_{k\ge0},
(\inj^{\cq\,\tens k}(p,p')u_k)_{k=0,1} \bigr)
= \bigl((p_k)_{k\ge0},(p'_k)_{k=0,1}\bigr) = (p,p'),
\]
and the equation
 $\Phi\cdot\restr_1
 =\id_{sA_{\infty1}(\cd,\cq;\ca)((f|_\cd,f|_\cq),(g|_\cd,g|_\cq))}$
is proven.
\end{proof}

\begin{lemma}\label{lem-homotopy-h-for-w}
Denote by $v$ the chain map
\[ v = \id - \restr_1\cdot\Phi: sA_\infty(\ce,\ca)(f,g)
\to sA_\infty(\ce,\ca)(f,g).
\]
Denote by $w$ the chain map
\[ w = \bigl[ sA_\infty(\ce,\ca)(f,g) \rTTo^v
sA_\infty(\ce,\ca)(f,g) \rMono^{(\pi\boxtimes1)M_{01}}
sA_\infty(\cf\cq,\ca)(\pi f,\pi g) \bigr].
\]
There exists a unique homotopy
$h:sA_\infty(\ce,\ca)(f,g)\to sA_\infty(\cf\cq,\ca)(\pi f,\pi g)$
of degree $-1$ such that \(w=B_1h+hB_1\),
\begin{multline*}
h' = \bigl(sA_\infty(\ce,\ca)(f,g) \rTTo^h
sA_\infty(\cf\cq,\ca)(\pi f,\pi g) \\
\rEpi^{\restr_{\le1}} sA_1(\cf\cq,\ca)(\pi f,\pi g) \rEpi^\restr
sA_1(\cq,\ca)(jf,jg)\bigr) = 0,
\end{multline*}
\begin{equation*}
h_k = 0: sA_\infty(\ce,\ca)(f,g) \to \prod_{X,Y\in\Ob\cd}
\uCom(s\cf\cq^{\tens k}(X,Y),s\ca(Xf,Yg)), \qquad \text{for } k>1.
\end{equation*}
\end{lemma}

\begin{proof}
We use \corref{cor-chain-to-Ainf-FQ-A-null-homotopic}, setting
$P=sA_\infty(\ce,\ca)(f,g)$. We have
\begin{multline*}
w = \bigl[ sA_\infty(\ce,\ca)(f,g) \rMono^{(\pi\boxtimes1)M_{01}}
sA_\infty(\cf\cq,\ca)(\pi f,\pi g) \bigr] \\
- \bigl[ sA_\infty(\ce,\ca)(f,g) \rTTo^{\restr_1}
sA_{\infty1}(\cd,\cq;\ca)\bigl((f|_\cd,f|_\cq),(g|_\cd,g|_\cq)\bigr) \rTTo^u
sA_\infty(\cf\cq,\ca)(\pi f,\pi g) \bigr]
\end{multline*}
due to \eqref{eq-u-Phi-(pi1)M}. Due to \eqref{eq-u'-urestr2}, $w'$
defined in condition 1 of
\corref{cor-chain-to-Ainf-FQ-A-null-homotopic} is
\begin{multline*}
w' = \bigl[ sA_\infty(\ce,\ca)(f,g) \rMono^{(\pi\boxtimes1)M_{01}}
sA_\infty(\cf\cq,\ca)(\pi f,\pi g) \\
\rEpi^{\restr_{\le1}} sA_1(\cf\cq,\ca)(\pi f,\pi g) \rEpi^\restr
sA_1(\cq,\ca)(jf,jg) \bigr] \\
- \bigl[ sA_\infty(\ce,\ca)(f,g) \rTTo^{\restr_1}
sA_{\infty1}(\cd,\cq;\ca)\bigl((f|_\cd,f|_\cq),(g|_\cd,g|_\cq)\bigr)
\rTTo^{u'}_{\restr_\cq} sA_1(\cq,\ca)(jf,jg) \bigr],
\end{multline*}
where $j:\cq\hookrightarrow\ce$ is the embedding $A_1$\n-functor,
$j_1=\bigl(s\cq\hookrightarrow s\cf\cq \rTTo^{\pi_1} s\ce\bigr)$. We
get
\begin{multline*}
w' = \bigl[ sA_\infty(\ce,\ca)(f,g) \rEpi^{\restr_{\le1}}
sA_1(\ce,\ca)(f,g) \\
\rMono^{(\pi\boxtimes1)M_{01}} sA_1(\cf\cq,\ca)(\pi f,\pi g)
\rEpi^{(\inj^\cq\boxtimes1)M_{01}} sA_1(\cq,\ca)(jf,jg) \bigr] \\
- \bigl[ sA_\infty(\ce,\ca)(f,g)
\rTTo^{((\iota\boxtimes1)M_{01},(j\boxtimes1)M_{01})}
sA_{\infty1}(\cd,\cq;\ca)\bigl((f|_\cd,f|_\cq),(g|_\cd,g|_\cq)\bigr) \\
\hspace{17em} \rTTo^{\restr_\cq} sA_1(\cq,\ca)(jf,jg) \bigr] \\
=\bigl[sA_\infty(\ce,\ca)(f,g) \rEpi^{\restr_{\le1}} sA_1(\ce,\ca)(f,g)
\rEpi^{(j\boxtimes1)M_{01}} sA_1(\cq,\ca)(jf,jg) \bigr] \\
-\bigl[sA_\infty(\ce,\ca)(f,g) \rEpi^{\restr_{\le1}} sA_1(\ce,\ca)(f,g)
\rEpi^{(j\boxtimes1)M_{01}} sA_1(\cq,\ca)(jf,jg) \bigr] = 0.
\end{multline*}
Therefore, $h'=0$ satisfies \(B_1h'+h'B_1=0=w'\). Hence, the unique
homotopy $h$ is constructed by
\corref{cor-chain-to-Ainf-FQ-A-null-homotopic}.
\end{proof}

\begin{remark}\label{rem-expicit-rh1}
In the case of \lemref{lem-homotopy-h-for-w} \(h\cdot\pr_k=h_k=0\) if
$k>1$ or if $k=0$. Indeed, \eqref{eq-h-Ph-Ai(FQA)-A1(FQA)-A1(QA)}
together with $h'=0$ implies that \(h_0=h\pr_0=h'\pr_0=0\), moreover,
\begin{equation}
rh_1|_{s\cq} = rh'\pr_1 =0: s\cq(X,Y) = s\cf_|\cq(X,Y) \to s\ca(Xf,Yg),
\label{eq-rh1Q-rhpr1-0}
\end{equation}
where \(r\in sA_\infty(\ce,\ca)(f,g)\). Therefore, recurrent
formula~\eqref{eq-bk-ph1} simplifies here to
\begin{multline}
(-)^rb^{\cf\cq}_k(rh_1) = rw_k - \sum_{a+1+c=k}^{m,n}
(\pi_1^{\tens a}f_{am}\tens rh_1\tens\pi_1^{\tens c}g_{cn})
b^\ca_{m+1+n}: \\
(s\cf\cq)^{\tens k}(X,Y) \to s\ca(Xf,Yg).
\label{eq-bk-rh1-rwk-pfrh1pg}
\end{multline}
\end{remark}

\begin{lemma}\label{lem-homotopy-h-factorizes}
The homotopy $h$ constructed in \lemref{lem-homotopy-h-for-w}
factorizes as
\[ h = \bigl(sA_\infty(\ce,\ca)(f,g) \rTTo^\eta sA_\infty(\ce,\ca)(f,g)
\rMono^{(\pi\boxtimes1)M_{01}} sA_\infty(\cf\cq,\ca)(\pi f,\pi g)\bigr)
\]
for a unique homotopy $\eta$ of degree $-1$ such that
$v=B_1\eta+\eta B_1$.
\end{lemma}

\begin{proof}
Let us show that $h$ satisfies conditions of
\corref{cor-null-homotopic-in-subcomplex}. The second is obvious. The
first is $R_\cd(rh_1)=0$ for any $r\in sA_\infty(\ce,\ca)(f,g)$, that is,
\begin{equation}
(-)^r(i^{\tens k}b^{\cf\cq}_k - b^\cd_ki) (rh_1) = 0:
s\cd^{\tens k}(X,Y) \to s\ca(Xf,Yg).
\label{eq-ikbk-bki-ph1}
\end{equation}
 From \eqref{eq-bk-rh1-rwk-pfrh1pg} we find the formula for $k>1$
\begin{multline*}
(-)^ri^{\tens k}b^{\cf\cq}_k(rh_1) = i^{\tens k}(rw_k)
- \sum_{a+1+c=k}^{m,n} (\iota_1^{\tens a}f_{am}
\tens i(rh_1)\tens\iota_1^{\tens c}g_{cn}) b^\ca_{m+1+n} \\
= \iota_1^{\tens k}r_k - i^{\tens k}[(\iota r,jr)u_k]
- \sum_{a+1+c=k}^{m,n}
(\overline{f}_{am}\tens i(rh_1)\tens\overline{g}_{cn}) b^\ca_{m+1+n}
:s\cd^{\tens k}(X,Y) \to s\ca(Xf,Yg).
\end{multline*}
A particular case of \eqref{eq-rh1Q-rhpr1-0} is
\begin{equation}
i(rh_1) = [\inj^\cd(rh')]\pr_1 = 0,
\label{eq-irh1-inDrhpr1-0}
\end{equation}
due to $h'=0$, where the $A_1$\n-functor
$\inj^\cd:\cd\hookrightarrow\cq$ is the natural embedding. For our
concrete choice of $u_k$ we get
\begin{equation*}
(-)^ri^{\tens k}b^{\cf\cq}_k(rh_1) = \iota_1^{\tens k}r_k
- i^{\tens k}\varpi^{\tens k}\iota_1^{\tens k}r_k = 0:
s\cd^{\tens k}(X,Y) \to s\ca(Xf,Yg),
\end{equation*}
since $\bigl(s\cd \rTTo^i s\cf\cq \rTTo^\varpi s\cd\bigr)=\id$.
Therefore, $i^{\tens k}b^{\cf\cq}_k(rh_1)=0$ and $b^\cd_ki(rh_1)=0$ due
to \eqref{eq-irh1-inDrhpr1-0}. We conclude that \eqref{eq-ikbk-bki-ph1}
is satisfied, and by \corref{cor-null-homotopic-in-subcomplex} there
exists a homotopy
\[ \eta: sA_\infty(\ce,\ca)(f,g) \to sA_\infty(\ce,\ca)(f,g),
\]
such that $\deg\eta=-1$, $h=\eta\cdot[(\pi\boxtimes1)M_{01}]$ and
$v=B_1\eta+\eta B_1$.
\end{proof}

Lemmata \ref{lem-Phi-one-sided-inverse} and
\ref{lem-homotopy-h-factorizes} show that the maps $\restr_1$ and
$\Phi$ given by \eqref{eq-restr1-Apsiu(EA)-Apsiu(DQA)} and
\eqref{eq-u-Phi-(pi1)M} are homotopy inverse to each other.

The \ainf-functor $\restr$ is surjective on objects by
\propref{pro-tilde-f-E-A}, and its first component is a homotopy
isomorphism. Therefore, it is an \ainf-equivalence by
\corref{cor-Theorem8:8}, and
\thmref{thm-restriction-equivalence-quiver} is proven.
\end{proof}

\begin{corollary}
The collection of \ainf-functors \eqref{eq-restr-A(EA)-A(DQA)} is
natural \ainfu-2-equivalence.
\end{corollary}

\begin{proof}
The restriction \ainf-functors
\(A_\infty^{\psi u}(\ce,\ca)\to A_\infty^{\psi u}(\cd,\ca)\) and
\(A_\infty^{\psi u}(\ce,\ca)\to A_1(\cq,\ca)\) are strict
\ainfu-2-transformations. By \eqref{dia-A(EQA)} the restriction
\ainf-functor
\(A_\infty^{\psi u}(\ce,\ca)\to A_{\infty1}^{\psi u}(\cd,\cq;\ca)\) is
also a strict \ainfu-2-transformation. It is an \ainf-equivalence by
\thmref{thm-restriction-equivalence-quiver}.
\end{proof}

\begin{remark}
The maps $\Phi$, $\eta$ constructed in the proof of
\thmref{thm-restriction-equivalence-quiver} satisfy
\[ \Phi\cdot\eta = 0:
sA_{\infty1}^{\psi u}(\cd,\cq;\ca)
\bigl((f|_\cd,f|_\cq),(g|_\cd,g|_\cq)\bigr)
\to sA_\infty^{\psi u}(\ce,\ca)(f,g).
\]
Indeed, \(\Phi\cdot\eta\) composed with an embedding,
\[ \Phi\cdot\eta\cdot(\pi\boxtimes1)M_{01} = \Phi\cdot h:
sA_{\infty1}^{\psi u}(\cd,\cq;\ca)
\bigl((f|_\cd,f|_\cq),(g|_\cd,g|_\cq)\bigr)
\to sA_\infty(\cf\cq,\ca)(\pi f,\pi g),
\]
is a degree $-1$ homotopy such that
\begin{multline*}
B_1(\Phi h) + (\Phi h)B_1 = \Phi(B_1h+hB_1) = \Phi w \\
= \Phi(\id-\restr_1\Phi)(\pi\boxtimes1)M_{01}
= (\id-\Phi\restr_1)\Phi(\pi\boxtimes1)M_{01} = 0.
\end{multline*}
We have \(\Phi h\pr_k=\Phi\cdot h_k=0\) for $k>1$ and
\[ (\Phi h)' = \Phi h\restr_{\le1}\restr = \Phi\cdot h' = 0.
\]
The 0 homotopy for 0 chain map also has these properties, and by
\corref{cor-chain-to-Ainf-FQ-A-null-homotopic} we conclude that
$\Phi h=0$.
\end{remark}

\begin{remark}
The equation
\[ \eta\cdot\restr_1 = 0: sA_\infty^{\psi u}(\ce,\ca)(f,g) \to
sA_{\infty1}^{\psi u}(\cd,\cq;\ca)
\bigl((f|_\cd,f|_\cq),(g|_\cd,g|_\cq)\bigr)
\]
also holds. Indeed, the decomposition
\(\restr_1=(\pi\boxtimes1)M_{01}\cdot L^i\) from
diagram~\eqref{dia-3triangles-Phi} implies that
\[ \eta\cdot\restr_1 = \eta\cdot(\pi\boxtimes1)M_{01}\cdot L^i
= h\cdot L^i.
\]
For any \(r\in sA_\infty^{\psi u}(\ce,\ca)(f,g)\) all components of the
element
\[ rhL^i =
\bigl((i^{\tens k}(rh_k))_{k\ge0},
(\inj^{\cq\,\tens k}(rh_k))_{k=0,1}\bigr)
\in sA_{\infty1}^{\psi u}(\cd,\cq;\ca)
\bigl((f|_\cd,f|_\cq),(g|_\cd,g|_\cq)\bigr)
\]
vanish except, possibly, those indexed by $k=1$ by
\remref{rem-expicit-rh1}. However, \(i(rh_1)=0\) by
\eqref{eq-irh1-inDrhpr1-0}, and, moreover, \(\inj^\cq(rh_1)=0\) by
\eqref{eq-rh1Q-rhpr1-0}, thus, all the components of \(rhL^i\) vanish.
\end{remark}

\begin{remark}
We have not used in the proof of
\thmref{thm-restriction-equivalence-quiver} the assumption of
pseudounitality of $\cd$ and $\ce$. Its assertion holds without this
property. If $\ca$ is unital, then the restriction \ainf-functor
\[ \restr: A_\infty(\ce,\ca) \to A_{\infty1}(\cd,\cq;\ca)
\]
is an \ainf-equivalence, surjective on objects. Its first component
maps admit a chain splitting. In the particular case \(\cd(X,Y)=0\) for
all $X,Y\in\Ob\cq$  we get Theorem~2.12 of \cite{LyuMan-freeAinf}: the
\ainf-functor
\[ \restr: A_\infty(\cf\cq,\ca) \to A_1(\cq,\ca)
\]
is an \ainf-equivalence.
\end{remark}

\section{\texorpdfstring{Relatively free $A_\infty$-categories}
 {Relatively free A8-categories}}\label{sec-Relatively-free-categories}
Hinich~\cite{Hinich:q-alg/9702015} defines standard cofibrations of
differential graded algebras. This notion is generalized by Drinfeld to
semi-free differential graded categories \cite{Drinf:DGquot}. We give a
definition in the spirit of these two definitions in the framework of
\ainf-categories.

\begin{definition}\label{def-relatively-free-Ainf}
Let \(e:\cc\to\cd\) be a strict \ainf-functor such that \(\Ob e\) is an
isomorphism, and \(e_1:s\cc\to s\cd\) is an embedding. The
\ainf-category $\cd$ is \emph{relatively free} over $\cc$, if it can be
represented as the union of an increasing sequence of its
\ainf-subcategories \(\cd_j\) and differential graded subquivers
\(\cq_j\)
\begin{equation}
\cd_0 \subset \cq_1 \subset \cd_1 \subset \cq_2
\subset \cd_2 \subset \cq_3 \subset \dots \subset \cd
\label{eq-D0-Q1-D1-Q2-D2-Q3}
\end{equation}
with the same set of objects \(\Ob\cd\), such that
\begin{enumerate}
\item \(s\cd_0=(s\cc)e_1\);

\item for each $j\ge0$ the embedding of graded quivers
\(\cd_j \rMono \cq_{j+1}\) admits a splitting map
\(\cq_{j+1} \rEpi \cd_j\) of degree 0;

\item for each $j>0$ the unique strict \ainf-functor
\(\cf\cq_j\to\cd_j\) extending the embedding \(\cq_j \rMono \cd_j\)
factors into the natural projection and an isomorphism
\[ \cf\cq_j \rEpi \cf\cq_j/s^{-1}(R_j) \rTTo^\sim \cd_j,
\]
where the system of relations \(R_j=R_{\cd_{j-1}}\subset s\cf\cq_j\) is
$R_j=\sum_{n\ge2}\im(\delta_n)$ for
\[ \delta_n = \bigl((s\cd_{j-1})^{\tens n} \rMono (s\cf\cq_j)^{\tens n}
\rTTo^{b_n^{\cf\cq_j}} s\cf\cq_j\bigr)
- \bigl((s\cd_{j-1})^{\tens n} \rTTo^{b_n^{\cd_{j-1}}} s\cd_{j-1}
\rMono s\cf\cq_j\bigr).
\]
\end{enumerate}
\end{definition}

When all differential graded quivers
\(\cn_j=\Coker(\cd_{j-1} \rMono \cq_j)\) have zero differential and the
$\kk$\n-modules \(\cn_j^k(X,Y)\) are free for all $j\ge1$, $k\in\ZZ$,
we say that $\cd$ is \emph{semi-free} over $\cc$ in accordance with
terminology of Drinfeld. In fact, if in
\defref{def-relatively-free-Ainf} one replaces \ainf-categories with
differential graded categories and adds the above assumption on
\(\cn_j\), then one recovers Definition~13.4 from \cite{Drinf:DGquot}
of semi-free differential graded categories.

The system of relations $R_j$ is the minimal one that ensures that the
natural embedding \(s\cd_{j-1}\hookrightarrow s\cf\cq_j/(R_j)=s\cd_j\)
is the first component of a strict \ainf-functor. In semi-free case we
may say that \(\cd_j\) is freely generated by \(\cn_j\) over
\(\cd_{j-1}\).

\subsection{The main construction}\label{sec-constr-D}
Let $\cc$ be a unital \ainf-category, and let \(\cb\subset\cc\) be its
full subcategory. The unit \(\uni^\cc_0\) is abbreviated to \(\uni_0\).

A vertex of a tree is $k$\n-ary if it is adjacent to $k+1$ edges. A
unary vertex is a 1\n-ary one.

Define a \emph{labeled tree} \(t=(t;X_0,X_1,\dots,X_n)\) as a non-empty
(non-reduced) plane rooted tree $t$ with $n$ leaves, such that unary
vertices are not joined by an edge, equipped with a sequence
\((X_0,X_1,\dots,X_n)\) of objects of $\cc$.

Let $e$ be an edge of $t$. If $i$ is the smallest number such that
$i$\n-th leaf is above $e$, the \emph{domain} of $e$ is defined as
\(\dom(e)=X_{i-1}\). If $k$ is the biggest number such that $k$\n-th
leaf is above $e$, the \emph{codomain} of $e$ is defined as
\(\codom(e)=X_k\). An \emph{admissible tree} is a labeled tree
\((t;X_0,X_1,\dots,X_n)\) such that for each edge $e$ adjacent to a
unary vertex \(\dom(e)\in\Ob\cb\) or \(\codom(e)\in\Ob\cb\) (or both).

The set of vertices $V(t)$ of a rooted tree $t$ has a canonical
ordering: $x\preccurlyeq y$ iff the minimal path connecting the root
with $y$ contains $x$. A \emph{$\cc$\n-admissible tree} is an admissible
tree \((t;X_0,X_1,\dots,X_n)\) such that top (maximal with respect to
$\preccurlyeq$) internal vertices are unary.

Define a graded quiver $\ce$ with the set of objects \(\Ob\ce=\Ob\cc\).
The $\ZZ$\n-graded $\kk$\n-module of morphisms between $X,Y\in\Ob\ce$
is defined as
\begin{align}
s\ce(X,Y) &= \bigoplus_{n\ge1}
\bigoplus_{\text{admissible }(t;X_0,X_1,\dots,X_n)}^{X_0=X,\;X_n=Y}
s\ce(t)(X,Y), \label{eq-E(XY)-E(t)(XY)} \\
s\ce(t)(X_0,X_n) &= s\ce(t) =
s\cc(X_0,X_1)\tens\dots\tens s\cc(X_{n-1},X_n)\bigl[|t|_1-|t|_>\bigr],
\notag
\end{align}
where $|t|_1$ is the number of unary internal vertices of $t$, and
$|t|_>$ is the number of internal vertices of arity $>1$.

The vertices of arity \(k>1\) are interpreted as $k$\n-ary
multiplications of degree 1. Unary vertices are interpreted as
contracting homotopies $H$ of degree $-1$. Define an \ainf-structure on
$\ce$ in which operations $b_k$, $k>1$, are given by grafting. So for
$k>1$ the operation $b_k$ is a direct sum of maps
\begin{equation*}
b_k = s^{|t_1|}\tens\dots\tens s^{|t_{k-1}|}\tens s^{|t_k|-|t|}:
s\ce(t_1)(Y_0,Y_1)\tens\dots\tens s\ce(t_k)(Y_{k-1},Y_k)
\to s\ce(t)(Y_0,Y_k),
\end{equation*}
where \(|t|=|t|_>-|t|_1\) and
$t=(t_1\sqcup\dots\sqcup t_k)\cdot\tree_k$. In particular,
$|t|=|t_1|+\dots+|t_k|+1$.

Let \(t=(t;X_0,X_1,\dots,X_n)\) be an admissible tree, whose lowest
internal vertex is not unary. In particular, $t$ might be the trivial
tree \(t=(|;X_0,X_1)\). Assume that \(X_0\in\Ob\cb\) or
\(X_n\in\Ob\cb\) (or both). Denote by $H$ the $\kk$\n-linear map
\begin{multline*}
H = s: s\ce(t)(X_0,X_n) =
s\cc(X_0,X_1)\tdt s\cc(X_{n-1},X_n)\bigl[|t|_1-|t|_>\bigr] \\
\to s\cc(X_0,X_1)\tdt s\cc(X_{n-1},X_n)\bigl[1+|t|_1-|t|_>\bigr]
= s\ce(t\cdot\tree_1)(X_0,X_n)
\end{multline*}
of degree $-1$. Here
 \(\tree_1=(\;
 \vstretch40 \begin{tanglec}\n \\ \id\end{tanglec}
 \;;X_0,X_n)\)
is the unary corolla.

The operation $b_1$ is determined by the given differential
$b_1:s\cc\to s\cc$ and by the recursive substitutions
\begin{align}
b_kb_1 &:= - \sum_{\alpha+p+\beta=k}^{\alpha+\beta>0}
(1^{\tens\alpha}\tens b_p\tens1^{\tens\beta})b_{\alpha+1+\beta},
\qquad k>1, \label{eq-bkb1-sum-(1b1)b} \\
Hb_1 &:= 1 - b_1H, \label{eq-Hb1-1-b1H}
\end{align}
where $H$ stands for a unary vertex.
Identity~\eqref{eq-bkb1-sum-(1b1)b} is satisfied for $k>1$ by
definition of $b_1$. We have to prove that \(b_1^2=0\). Assume that
$k>1$, then
\begin{align*}
b_kb_1^2 &= - \sum_{\alpha+p+\beta=k}^{\alpha+\beta>0}
(1^{\tens\alpha}\tens b_p\tens1^{\tens\beta})b_{\alpha+1+\beta}b_1 \\
&= \sum_{\substack{\alpha+p+\beta=k\\ \gamma+q+\delta=\alpha+1+\beta}}
^{\gamma+\delta>0}
(1^{\tens\alpha}\tens b_p\tens1^{\tens\beta})
(1^{\tens\gamma}\tens b_q\tens1^{\tens\delta})b_{\gamma+1+\delta} \\
&= \sum_{\alpha+p+\eps+q+\delta=k}
(1^{\tens\alpha}\tens b_p\tens1^{\tens\eps}\tens b_q\tens1^{\tens\delta})
b_{\alpha+\eps+\delta+2} \\
&\quad- \sum_{\gamma+q+\eta+p+\beta=k}
(1^{\tens\gamma}\tens b_q\tens1^{\tens\eta}\tens b_p\tens1^{\tens\beta})
b_{\gamma+\eta+\beta+2} \\
&\quad+ \sum_{\gamma+r+\delta=k}^{\gamma+\delta>0}
\biggl\{1^{\tens\gamma}\tens \Bigl[\sum_{\kappa+p+\lambda=r}
(1^{\tens\kappa}\tens b_p\tens1^{\tens\lambda})b_{\kappa+1+\lambda}
\Bigr] \tens1^{\tens\delta} \biggr\} b_{\gamma+1+\delta} \\
&= \sum_{\gamma+1+\delta=k}
(1^{\tens\gamma}\tens b_1^2\tens1^{\tens\delta})b_k,
\end{align*}
because the sum in square brackets vanishes for $r>1$ by
\eqref{eq-bkb1-sum-(1b1)b}. We also have
\[ Hb_1^2 = (1-b_1H)b_1 = b_1 - b_1(1-b_1H) = b_1^2H.
\]
By induction the equation \(b_1^2\big|_{s\cc}=0\) implies that
\(b_1^2=0\) on $s\ce$. Therefore, $\ce$ is an \ainf-category.

It has an ideal $(R_\cc)_+$, generated by the $\kk$\n-subquiver
$R_\cc=\sum_{n\ge2}\im(\delta_n)$ for
\[ \delta_n = \bigl((s\cc)^{\tens n} \rMono (s\ce)^{\tens n}
\rTTo^{b_n^{\ce}} s\ce\bigr)
- \bigl((s\cc)^{\tens n} \rTTo^{b_n^\cc} s\cc \rMono s\ce\bigr)
\]
by application of operations
\(1^{\tens\alpha}\tens H\tens1^{\tens\beta}\),
\(1^{\tens\alpha}\tens b_p\tens1^{\tens\beta}\) for $p\ge2$. By
\lemref{lem-I-ideal} $R_\cc b^\ce_1\subset(R_\cc)\subset(R_\cc)_+$,
where $(R_\cc)$ denotes the ideal generated by application of
\(1^{\tens\alpha}\tens b_p\tens1^{\tens\beta}\) ($p\ge2$) only.
Similarly to \propref{pro-Rb1-subset(R)} this implies that
\((R_\cc)_+b^\ce_1\subset(R_\cc)_+\). Indeed, let
\(t=(t;X_0,X_1,\dots,X_n)\) be an admissible tree, whose lowest
internal vertex is not unary. Assume that \(X_0\in\Ob\cb\) or
\(X_n\in\Ob\cb\), so that \(t\cdot\tree_1\) is admissible. For an
arbitrary \(z\in(R_\cc)_+(t)(X_0,X_n)\) there exists
\(zH\in(R_\cc)_+(t\cdot\tree_1)(X_0,X_n)\), and
\(zb_1\in(R_\cc)_+(X_0,X_n)\) by induction. Due to \eqref{eq-Hb1-1-b1H}
\(zHb_1=z-zb_1H\in(R_\cc)_+\), which proves the claim. Therefore, the
ideal \((R_\cc)_+\) is stable with respect to all \ainf-operations,
including $b_1$.

Denote by $\cd=\ce/s^{-1}(R_\cc)_+=\Quo(\cc|\cb)$ the quotient
\ainf-category. It has a direct sum decomposition similar to that of
$\ce$
\begin{align*}
s\cd(X,Y) &= \bigoplus_{n\ge1}
\bigoplus_{\cc\text{-admissible }(t;X_0,X_1,\dots,X_n)}
^{X_0=X,\quad X_n=Y}
s\cd(t)(X,Y), \\
s\cd(t)(X_0,X_n) &= s\ce(t)(X_0,X_n) =
s\cc(X_0,X_1)\tens\dots\tens s\cc(X_{n-1},X_n)\bigl[|t|_1-|t|_>\bigr],
\end{align*}
with the only difference that the sum is taken over $\cc$\n-admissible
trees $t$. We can view $\cd$ as a graded $\kk$\n-subquiver of $\ce$.

The category $\cc$ is embedded in $\cd$ (via a strict \ainf-functor) as
\[ s\cd_0(X,Y) = s\cd(|)(X,Y) = s\cc(X,Y)
\]
for the trivial tree $t=(|;X,Y)$. Let us show that $\cd$ is relatively
free over $\cc$.

Let us define for $j\ge0$ the \ainf-subcategories $\cd_j$ and
differential graded subquivers \(\cq_{j+1}\) of $\cd$ so that
embeddings \eqref{eq-D0-Q1-D1-Q2-D2-Q3} hold. Each leaf $\ell$ and the
root of a tree can be connected by the unique minimal path. We say that
internal vertices occurring at this path are between the root and the
leaf $\ell$. Define for $j\ge0$ the \ainf-subcategory
\(\cd_j=\oplus_t\cd(t)\) of $\cd$, where the summation goes over all
\begin{equation}
\parbox{20em}{$\cc$-admissible trees $t$ with no more than $j$ unary
internal vertices between the root and any leaf.}
\tag{C1}\label{eq-admissible-Dj-C1}
\end{equation}
Define for $j\ge1$ the graded subquiver \(\cn_j=\oplus_t\cd(t)\) of
$\cd$, where the summation goes over all trees $t$ satisfying
\eqref{eq-admissible-Dj-C1} and such that
\begin{equation}
\parbox{30em}{there exists a leaf $\ell$ of $t$ with $j$ unary internal
vertices between the root and $\ell$; the lowest internal vertex
(adjacent to the root) is unary.}
\tag{C2}
\end{equation}
One can easily see that for $j\ge1$
\[ s\cq_j = s\cd_{j-1} \oplus s\cn_j
\]
is a differential graded subquiver of \(\cd_j\subset\cd\). For example,
\(\cd_0=\cd(|)=\cc\),
 \(\cn_1=\cd\bigl(\,
 \vstretch40 \begin{tanglec}\n \\ \id\end{tanglec} \,\bigr)\),
 \(\cq_1=\cd(|)\oplus\cd\bigl(\,
 \vstretch40 \begin{tanglec}\n \\ \id\end{tanglec} \,\bigr)\),
and \(\cd_1=\cd(|)\oplus\oplus_t\cd(t)\), where $t$ runs over
admissible trees with the only unary internal vertex $v$, such that all
other internal vertices lie on the minimal path between the root and
$v$.

The inclusion map of differential graded quivers
\(i:s\cq_j \rMono s\cd_j\) induces a unique strict \ainf-functor
\(\ihat:\cf\cq_j\to\cd_j\) \cite[Corollary~2.4]{LyuMan-freeAinf}.

\begin{proposition}
The map \(\ihat_1:s\cf\cq_j\to s\cd_j\) is surjective and its kernel is
\((R_{\cd_{j-1}})\). Thus it induces an isomorphism
\(\iota_1:s\cf\cq_j/(R_{\cd_{j-1}})\to s\cd_j\).
\end{proposition}

\begin{proof}
The strict \ainf-functor \(\ihat:\cf\cq_j\to\cd_j\) is described in
\cite[Section~2.6]{LyuMan-freeAinf} as follows. Let $t$ be a reduced
labeled tree, with $n$ input leaves, and let $\le$ be a linear order on
\(\inve(t)\), such that $x\preccurlyeq y$ implies $x\le y$ for all
\(x,y\in\inve(t)\). The choice of $\le$ is equivalent to the choice of
decomposition into product of elementary
forests~\eqref{eq-ord-tree-decomp-forest}. The linearly ordered tree
\((t,\le)\) determines the map
\(b^{\cf\cq_j}_{(t,\le)}:(s\cf\cq_j)^{\tens n}\to s\cf\cq_j\) given by
\eqref{eq-bFQ-1bFQ1-bFQ} and a similar map
\(b^{\cd_j}_{(t,\le)}:(s\cd_j)^{\tens n}\to s\cd_j\). In the
commutative diagram from \cite[Section~2.6]{LyuMan-freeAinf}
\begin{diagram}[LaTeXeqno]
(s\cq_j)^{\tens n} & \rTTo_{\pm s^{-|t|}}^{b^{\cf\cq_j}_{(t,\le)}}
& s\cf_t\cq_j \\
\dTTo<{i_1^{\tens n}} && \dTTo>{\ihat_1} \\
(s\cd_j)^{\tens n} & \rTTo^{b^{\cd_j}_{(t,\le)}} & s\cd_j
\label{dia-ihat1-sti1bDj}
\end{diagram}
the top map is invertible, so \(\ihat_1\) is uniquely determined by
this diagram.

Being the first component of a strict \ainf-functor \(\ihat_1\)
satisfies, in particular, the equation
\begin{multline*}
\bigl((s\cd_{j-1})^{\tens n} \rMono (s\cf_|\cq_j)^{\tens n}
\rTTo^{b^{\cf\cq_j}_n} s\cf\cq_j \rTTo^{\ihat_1} s\cd_j\bigr)
= \bigl((s\cd_{j-1})^{\tens n} \rMono^{i_1^{\tens n}}
(s\cd_j)^{\tens n} \rTTo^{b^{\cd_j}_n} s\cd_j\bigr) \\
= \bigl((s\cd_{j-1})^{\tens n} \rTTo^{b^{\cd_{j-1}}_n} s\cd_{j-1}
\rMono^{i_1} s\cd_j\bigr)
= \bigl((s\cd_{j-1})^{\tens n} \rTTo^{b^{\cd_{j-1}}_n} s\cd_{j-1}
\rMono s\cf_|\cq_j \rTTo^{\ihat_1} s\cd_j\bigr).
\end{multline*}
It implies \(R_{\cd_{j-1}}\ihat_1=0\). Since $\ihat$ is strict we have
also \((R_{\cd_{j-1}})\ihat_1=0\). Thus there is a strict \ainf-functor
$\iota$ with the first component
\(\iota_1:s\cf\cq_j/(R_{\cd_{j-1}})\to s\cd_j\), identity on objects.

Let us construct a degree 0 map \(\phi:s\cd_j\to s\cf\cq_j\) for
$j\ge1$. Let $t$ be a tree that satisfies \eqref{eq-admissible-Dj-C1}.
Denote by \(\UV(t)\subset\inve(t)\) the subset of unary internal
vertices. Let \(\minUV(t)\) be the subset of partially ordered set
\((\UV(t),\preccurlyeq)\) consisting of minimal elements. Let
\(L\subset\Leaf(t)\) be the subset of leaves $\ell$ such that between
$\ell$ and the root there are no unary vertices. Let \(\overline{L}\)
be the set of leaf vertices above leaves from $L$. Using the canonical
linear ordering \(t_<=(t,\le)\) of $\alve(t)$
\cite[Section~1.7]{LyuMan-freeAinf} we can write the set
\(\overline{L}\sqcup\minUV(t)\) as \(\{u_1<\dots<u_k\}\). For any
\(1\le p\le k\) denote by $t_p$ the $\cc$\n-admissible subtree of $t$
with
\[ \alve(t_p) = \{y\in\alve(t) \mid y\succcurlyeq u_p\} \sqcup
\{\text{new root vertex }r_p\}.
\]
Edges of $t_p$ are all edges of $t$ above $u_p$ plus a new root edge
between $u_p$ and $r_p$. In particular, if \(u_p\in\overline{L}\), then
\(t_p=|\) is the trivial tree. Denote by $t'$ the reduced labeled tree,
which is $t$ with all vertices and edges above \(\minUV(t)\) removed.
It has precisely $k$ leaves. Thus $t$ is the concatenation of a forest
and $t'$:
\begin{equation}
t = (t_1\sqcup t_2\sqcup\dots\sqcup t_k)\cdot t'.
\label{eq-t-t1-tk-tprim}
\end{equation}
We have
\begin{equation*}
\inve(t) = \inve(t')\sqcup \bigsqcup_{p=1}^k \inve(t_p), \qquad
\Leaf(t) = \bigsqcup_{p=1}^k \Leaf(t_p), \qquad
\Leaf(t') = \bigsqcup_{p=1}^k \Out(t_p).
\end{equation*}
Correspondingly, the labels of the $p$\n-th leaf of $t'$ are those of
\(\Out(t_p)\).

Being simply a shift, the $\kk$\n-linear map
\(b^\cd_{t'_<}:s\cd(t_1)\tdt s\cd(t_k)\to s\cd(t)\) is invertible.
Therefore, for any element \(x\in s\cd(t)\) there exists a unique
tensor \(\sum_iz^i_1\tdt z^i_k\in s\cd(t_1)\tdt s\cd(t_k)\) such that
\(x=\sum_i(z^i_1\tdt z^i_k)b^\cd_{t'_<}\). We have
\(z^i_p\in s\cd(t_p)\subset s\cq_j\), in particular,
\(z^i_p\in s\cd(|)=s\cc\), if \(u_p\in\overline{L}\). Define
\[ x\phi = \sum_i(z^i_1\tdt z^i_k)b^{\cf\cq_j}_{t'_<}
\in s\cf_{t'}\cq_j.
\]
Commutative diagram~\eqref{dia-ihat1-sti1bDj} implies that
\(x\phi\ihat_1=\sum_i(z^i_1\tdt z^i_k)b^{\cd_j}_{t'_<}=x.\) Therefore,
\begin{equation}
\bigl[ s\cd_j \rTTo^\phi s\cf\cq_j \rEpi^{\pi_1}
s\cf\cq_j/(R_{\cd_{j-1}}) \rTTo^{\iota_1} s\cd_j \bigr] = \id.
\label{eq-phi-pi-iota-id}
\end{equation}

Let us prove that \(\phi\pi_1\) preserves the operations $b_n$ for
$n>1$. Indeed,
\[ b^{\cd_j}_n(\phi\pi_1) - (\phi\pi_1)^{\tens n}b_n =
\bigl[ (s\cd_j)^{\tens n}
\rTTo^{b^{\cd_j}_n\phi-\phi^{\tens n}b^{\cf\cq_j}_n} s\cf\cq_j
\rTTo^{\pi_1} s\cf\cq_j/(R_{\cd_{j-1}}) \bigr].
\]
Consider trees $\tau_1$, $\tau_2$, \dots, $\tau_n$ satisfying
condition~\eqref{eq-admissible-Dj-C1}, labeled so that the operation
\(b^{\cd_j}_n:s\cd(\tau_1)\tdt s\cd(\tau_n)\to s\cd_j\) makes sense.
The quiver \((s\cd_j)^{\tens n}\) is a direct sum of such
\(s\cd(\tau_1)\tdt s\cd(\tau_n)\). If some of trees $\tau_p$ are not
trivial, then \(b^{\cd_j}_n\phi=\phi^{\tens n}b^{\cf\cq_j}_n\), because
constructing $\phi$ for
\(\tau=(\tau_1\sqcup\tau_2\sqcup\dots\sqcup\tau_n)\cdot\tree_n\) is
equivalent to decomposing each $\tau_p$ as in \eqref{eq-t-t1-tk-tprim},
collecting the upper parts, and gluing the lowest parts $\tau'_p$ into
\(\tau'=(\tau'_1\sqcup\tau'_2\sqcup\dots\sqcup\tau'_n)\cdot\tree_n\).
If all trees $\tau_p$ are trivial, then \(\cd(\tau_p)=\cd(|)=\cc\) and
\[ (b^\cc_n\phi-\phi^{\tens n}b^{\cf\cq_j}_n)\pi_1 =
(b^\cc_n-b^{\cf\cq_j}_n)\pi_1 = 0:
(s\cc)^{\tens n} \to s\cf\cq_j/(R_{\cd_{j-1}}),
\]
due to \(R_\cc\pi_1\subset R_{\cd_{j-1}}\pi_1=0\).

We claim that
\begin{equation}
\bigl[ s\cf\cq_j \rTTo^{\ihat_1} s\cd_j \rTTo^\phi s\cf\cq_j
\rTTo^{\pi_1} s\cf\cq_j/(R_{\cd_{j-1}}) \bigr]
= \bigl[ s\cf\cq_j \rTTo^{\pi_1} s\cf\cq_j/(R_{\cd_{j-1}}) \bigr].
\label{eq-FQ-D-FQ-FQ/R-FQ-FQ/R}
\end{equation}
First of all, the restriction of this equation to
\(s\cq_j=s\cf_|\cq_j\) holds true:
\begin{equation}
\bigl[ s\cq_j \rTTo^{i_1} s\cd_j \rTTo^\phi s\cf\cq_j \rTTo^{\pi_1}
s\cf\cq_j/(R_{\cd_{j-1}}) \bigr]
= \bigl[ s\cq_j \rMono s\cf\cq_j \rTTo^{\pi_1}
s\cf\cq_j/(R_{\cd_{j-1}}) \bigr].
\label{eq-Q-D-FQ-FQ/R-Q-FQ-FQ/R}
\end{equation}
Indeed, \(s\cq_j=s\cd_{j-1}\oplus s\cn_j\). On the first summand we get
for \(x\in s\cd_{j-1}(t)\)
\begin{multline*}
x \rMapsTo^{i_1} x = \sum_i(z^i_1\tdt z^i_k)b^{\cd_j}_{t'_<}
\rMapsTo^\phi \sum_i(z^i_1\tdt z^i_k)b^{\cf\cq_j}_{t'_<} \\
\rMapsTo^{\pi_1}
\sum_i(z^i_1\tdt z^i_k)b^{\cf\cq_j}_{t'_<} + (R_{\cd_{j-1}}) =
\sum_i(z^i_1\tdt z^i_k)b^{\cd_{j-1}}_{t'_<} + (R_{\cd_{j-1}}) = x\pi_1
\end{multline*}
by \propref{pro-FD/(R)-D} because \(z_p^i\in s\cd_{j-1}\). On the
second summand we get for \(x\in s\cn_j(t)\)
\[ x \rMapsTo^{i_1} x \rMapsTo^\phi x \rMapsTo^{\pi_1} x\pi_1
\]
since \(t'=|\). Thus, \eqref{eq-Q-D-FQ-FQ/R-Q-FQ-FQ/R} is verified.

Now we prove \eqref{eq-FQ-D-FQ-FQ/R-FQ-FQ/R} on the generic summand
\(s\cf_\tau\cq_j\) of \(s\cf\cq_j\), where $\tau$ is a reduced labeled
tree with $n$ leaves. The first map below is an isomorphism:
\begin{align*}
& \bigl[ (s\cq_j)^{\tens n} \rTTo^{b^{\cf\cq_j}_{\tau_<}}
s\cf_\tau\cq_j \rTTo^{\ihat_1} s\cd_j \rTTo^{\phi\pi_1}
s\cf\cq_j/(R_{\cd_{j-1}}) \bigr] \\
&= \bigl[ (s\cq_j)^{\tens n} \rMono^{i_1^{\tens n}} (s\cd_j)^{\tens n}
\rTTo^{b^{\cd_j}_{\tau_<}} s\cd_j \rTTo^{\phi\pi_1}
s\cf\cq_j/(R_{\cd_{j-1}}) \bigr] \\
&= \bigl[ (s\cq_j)^{\tens n} \rMono^{i_1^{\tens n}} (s\cd_j)^{\tens n}
\rTTo^{\phi\pi_1} (s\cf\cq_j/(R_{\cd_{j-1}}))^{\tens n}
\rTTo^{b_{\tau_<}} s\cf\cq_j/(R_{\cd_{j-1}}) \bigr] \\
&= \bigl[ (s\cq_j)^{\tens n} \rMono (s\cf\cq_j)^{\tens n}
\rTTo^{\pi_1^{\tens n}} (s\cf\cq_j/(R_{\cd_{j-1}}))^{\tens n}
\rTTo^{b_{\tau_<}} s\cf\cq_j/(R_{\cd_{j-1}}) \bigr] \\
&= \bigl[ (s\cq_j)^{\tens n} \rTTo^{b^{\cf\cq_j}_{\tau_<}}
s\cf_\tau\cq_j \rTTo^{\pi_1} s\cf\cq_j/(R_{\cd_{j-1}}) \bigr]
\end{align*}
by \eqref{eq-Q-D-FQ-FQ/R-Q-FQ-FQ/R} and by the fact that the considered
maps $\ihat_1$, $\phi\pi_1$ and $\pi_1$ commute with \(b_{\tau_<}\).

Rewriting \eqref{eq-FQ-D-FQ-FQ/R-FQ-FQ/R} in the form
\begin{multline*}
\bigl[ s\cf\cq_j \rEpi^{\pi_1} s\cf\cq_j/(R_{\cd_{j-1}})
\rTTo^{\iota_1} s\cd_j \rTTo^\phi s\cf\cq_j \rTTo^{\pi_1}
s\cf\cq_j/(R_{\cd_{j-1}}) \bigr] \\
= \bigl[ s\cf\cq_j \rEpi^{\pi_1} s\cf\cq_j/(R_{\cd_{j-1}}) \bigr],
\end{multline*}
we find by surjectivity of $\pi_1$ that
\[ \bigl[ s\cf\cq_j/(R_{\cd_{j-1}}) \rTTo^{\iota_1} s\cd_j \rTTo^\phi
s\cf\cq_j \rTTo^{\pi_1} s\cf\cq_j/(R_{\cd_{j-1}}) \bigr] = \id.
\]
Together with \eqref{eq-phi-pi-iota-id} this proves that $\phi\pi_1$ is
an inverse to $\iota_1$.
\end{proof}

\begin{corollary}
The \ainf-category $\cd=\Quo(\cc|\cb)$ is relatively free over $\cc$.
\end{corollary}

\subsection{The first equivalence}
Let $\ca$ be pseudounital, then the restriction functor
\[ A_{\infty1}^{\psi u}(\cd_0,\cq_1;\ca)
\to A_\infty^{\psi u}(\cc,\ca), \qquad (f,f') \mapsto f
\]
takes values in the full subcategory
\(A_\infty^{\psi u}(\cc,\ca)_{\modulo\cb}\). Indeed, let $(f,f')$ be an
object of \(A_{\infty1}^{\psi u}(\cd_0,\cq_1;\ca)\). For arbitrary
objects $X$, $Y$ of $\cb$ we have
\begin{align*}
f_1 &=
\bigl[s\cc(X,Y) \rMono s\cq_1(X,Y) \rTTo^{f'_1} s\ca(Xf,Yf)\bigr] \\
&= \bigl[s\cd(|)(X,Y) \rTTo^{b_1H+Hb_1} s\cq_1(X,Y) \rTTo^{f'_1}
s\ca(Xf,Yf)\bigr] \\
&= \bigl[s\cd(|)(X,Y) \rTTo^{b_1(Hf'_1)+(Hf'_1)b_1} s\ca(Xf,Yf)\bigr],
\end{align*}
where $H$ is the map \(H:s\cd(|) \rTTo^s s\cd(\tree_1) \rMono s\cq_1\).
Hence, the above $f_1$ is null-homo\-to\-pic. By Definition~6.4 of
\cite{LyuOvs-iResAiFn} the \ainf-functor $f\big|_\cb$ is contractible.

A short exact sequence of chain maps of complexes is \emph{semisplit}
(resp. \emph{semisplittable}) if it is split (resp. splittable) as a
sequence of degree 0 maps of graded $\kk$\n-modules.

\begin{lemma}\label{lem-semisplit-split}
Let $0\to C \rTTo^\alpha A \rTTo^\beta B\to0$ be a semisplittable exact
sequence of complexes of $\kk$\n-modules. If $C$ is contractible, then
this sequence is splittable, and the splitting chain map \(\nu:B\to A\) can
be chosen so that $\nu$ is homotopy inverse to $\beta$.
\end{lemma}

\begin{proof}
Let $\phi:A\to C$ be a map of degree 0, such that $\alpha\phi=1_C$.
Assume that $1_C=Hd+dH$ for a homotopy $H:C\to C$ of degree $-1$. Then
$\psi=(\phi H)d=\phi Hd+d\phi H:A\to C$ is a chain map, and
$\alpha\psi=\alpha\phi Hd+d\alpha\phi H=Hd+dH=1_C$. Denote by
\(\nu:B\to A\) the unique $\kk$\n-linear map such that $\nu\psi=0$,
$\nu\beta=1$. The splitting injection $\nu=\ker\psi$ is a chain map.
The sequence looks as follows
\[ 0 \to C \pile{\rTTo^\alpha\\ \lTTo_\psi} A
\pile{\rTTo^\beta\\ \lTTo_\nu} B \to 0.
\]
Since $A$ is a direct sum $C\oplus B$, we have
\[ \id_A - \beta\nu = \psi\alpha = (\phi Hd+d\phi H)\alpha
= (\phi H\alpha)d + d(\phi H\alpha) = \gamma d + d\gamma,
\]
where $\gamma=\bigl(A \rTTo^\phi C \rTTo^H C \rTTo^\alpha A\bigr)$ is a
homotopy.
\end{proof}

\begin{proposition}\label{pro-A(D0Q1A)-A(CA)modB-equivalence}
Let $\ca$ be unital, then the restriction strict \ainf-functor
\begin{equation}
\restr: A_{\infty1}^{\psi u}(\cd_0,\cq_1;\ca)
\to A_\infty^u(\cc,\ca)_{\modulo\cb}
\label{eq-A(D0Q1A)-A(CA)modB}
\end{equation}
is an \ainf-equivalence, surjective on objects. The chain surjections
$\restr_1$ admit a chain splitting.
\end{proposition}

\begin{proof}
First of all, $\restr$ is surjective on objects. Indeed, assume that
\(f:\cc\to\ca\) is unital and  \(\cb \rMono \cc \rTTo^f \ca\) is
contractible. We have to extend the chain maps
\(f_1:s\cc(X,Y)\to s\ca(Xf,Yf)\) to chain maps
\(f'_1:s\cq_1(X,Y)\to s\ca(Xf,Yf)\). If \(X,Y\notin\Ob\cb\), then
\(s\cq_1(X,Y)=s\cc(X,Y)\) and \(f'_1=f_1\). If \(X\in\Ob\cb\) or
\(Y\in\Ob\cb\), then \(s\cq_1(X,Y)=s\cc(X,Y)\oplus s\cd(\tree_1)(X,Y)\)
as a graded quiver, and the complex \(\ca(Xf,Yf)\) is contractible by
Proposition~6.1(C1), (C2) of \cite{LyuOvs-iResAiFn}. Let
\(\chi_{XY}:s\ca(Xf,Yf)\to s\ca(Xf,Yf)\) be a contracting homotopy for
\(\ca(Xf,Yf)\). Define
\[ f'_1 = \bigl(s\cd(\tree_1)(X,Y) \rTTo^{s^{-1}} s\cc(X,Y) \rTTo^{f_1}
s\ca(Xf,Yf) \rTTo^{\chi_{XY}} s\ca(Xf,Yf)\bigr).
\]
Then \(Hf'_1=f_1\chi_{XY}\). Since \(H=s:s\cc\to s\cd(\tree_1)\) is
invertible, the equation
\begin{multline*}
Hf'_1b_1 - Hb_1f'_1 = Hf'_1b_1 + b_1Hf'_1 - f'_1
= f_1\chi_{XY}b_1 + b_1f_1\chi_{XY} - f_1 = 0: \\
s\cc(X,Y) \to s\ca(Xf,Yf)
\end{multline*}
implies that $f'_1$ is a chain map.

Let us prove that the restriction chain map
\[ \restr_1:
sA_{\infty1}^{\psi u}(\cc,\cq_1;\ca)\bigl((f,f'),(g,g')\bigr)
\to sA_\infty(\cc,\ca)(f,g), \qquad (r,r') \mapsto r
\]
is homotopy invertible. This map is a product over \(n\in\ZZ_{\ge0}\)
of the restriction maps \(\rho_n:V_n\to V'_n\) of the graded
$\kk$\n-modules of $n$\n-th components (compare \eqref{eq-V0-V1-Vn}
with analogous decomposition of \(sA_\infty(\cc,\ca)(f,g)\)). Clearly,
the maps \(\rho_n=\id\) for $n=0$ or for $n>1$. On the other hand, for
$n=1$
\begin{multline*}
\rho_1 =\prod \uCom(\inj^\cc,1):
\\
\prod_{X,Y\in\Ob\cc} \uCom(s\cq_1(X,Y),s\ca(Xf,Yg)) \to
\prod_{X,Y\in\Ob\cc} \uCom(s\cc(X,Y),s\ca(Xf,Yg))
\end{multline*}
is surjective with the kernel
\(\Ker\rho_1=\prod_{X,Y\in\Ob\cc}\uCom(s\cd(\tree_1)(X,Y),s\ca(Xf,Yg))\),
because the sequence \(0\to s\cc\to s\cq_1\to s\cd(\tree_1)\to0\) is
semisplit. Since we may write the kernel as
\[ \Ker\rho_1 =
\prod_{\substack{X\in\Ob\cb,Y\in\Ob\cc\\ \text{or }X\in\Ob\cc,Y\in\Ob\cb}}
\uCom(s\cd(\tree_1)(X,Y),s\ca(Xf,Yg)),
\]
it is contractible, because contractibility of
\(f\big|_\cb,g\big|_\cb:\cb\to\ca\) implies contractibility of
complexes \(\ca(Xf,Yg)\) by Proposition~6.1(C1), (C2) of
\cite{LyuOvs-iResAiFn}.

Summing up, the first term of the semisplit exact sequence
\[ 0 \to \Ker\restr_1 \rTTo
sA_{\infty1}(\cc,\cq_1;\ca)\bigl((f,f'),(g,g')\bigr) \rTTo^{\restr_1}
sA_\infty(\cc,\ca)(f,g) \to 0
\]
is contractible. By \lemref{lem-semisplit-split} this sequence admits a
splitting chain map
\[ \nu: sA_\infty(\cc,\ca)(f,g) \to
sA_{\infty1}(\cc,\cq;\ca)\bigl((f,f'),(g,g')\bigr),
\]
and $\nu$ is homotopy inverse to $\restr_1$. Applying
\corref{cor-Theorem8:8} we conclude that \eqref{eq-A(D0Q1A)-A(CA)modB}
is an \ainf-equivalence.
\end{proof}

\begin{corollary}
\ainfu-2-transformation \eqref{eq-A(D0Q1A)-A(CA)modB} is a natural
\ainfu-2-equivalence.
\end{corollary}

\begin{proposition}\label{pro-f-extends-(ff)}
Let $\ca$ be a unital \ainf-category, let $\cd$ be a pseudounital
\ainf-category with distinguished cycles \(\iota^\cd_X\), let $\cq$ be
a differential graded quiver and let $\cn$ be a graded quiver such that
\(\Ob\cd=\Ob\cq=\Ob\cn\) and \(\cq=\cd\oplus\cn\). Suppose that
$\cn(X,Y)\ne0$ implies that \(\iota^\cd_X\in\im b_1\) or
\(\iota^\cd_Y\in\im b_1\). Then an arbitrary pseudounital \ainf-functor
\(f:\cd\to\ca\) extends to an object
\((f,f')\in\Ob A_{\infty1}^{\psi u}(\cd,\cq;\ca)\) for some $f'$.
\end{proposition}

\begin{proof}
Let \(s\cm=\cn\) be the differential graded quiver and let
\(\alpha:s\cm\to s\cd\) be the chain map defined in
\secref{sec-Cones-freely-generated}. There exists a homotopy
\(\overline{f}:s\cm(X,Y)\to s\ca(Xf,Yf)\) of degree $-1$ such that
\[ \alpha f_1 = \overline{f}b_1 + d^{\cm[1]}\overline{f}:
s\cm(X,Y) \rTTo s\ca(Xf,Yf).
\]
Indeed, the case of \(\cm(X,Y)=0\) being obvious, we may assume that
\(\iota_X\in\im b_1\) or \(\iota_Y\in\im b_1\). Then
\(\sS{_{Xf}}\uni^\ca_0\in\iota_Xf_1+\im b_1\subset\im b_1\) or
\(\sS{_{Yf}}\uni^\ca_0\in\iota_Yf_1+\im b_1\subset\im b_1\). Since
$\ca$ is unital, the complex \(s\ca(Xf,Yf)\) is contractible with some
contracting homotopy $\overline{h}$. We may take
\(\overline{f}=\alpha f_1\overline{h}\).

Define a degree 0 map
\[ f'_1 = \bigl( s\cq(X,Y) = s\cd(X,Y)\oplus s\cn(X,Y)
\rTTo^{(f_1,s^{-1}\overline{f})} s\ca(Xf,Yf) \bigr).
\]
For arbitrary \(p\in s\cd(X,Y)\), \(m\in s\cm(X,Y)=\cn(X,Y)\) we have
\begin{align*}
(p,ms)(f'_1b^\ca_1-b^\cq_1f'_1) &= pf_1b^\ca_1 + m\overline{f}b^\ca_1
-(pb^\cd_1+m\alpha)f_1 +(md^{\cm[1]}s)s^{-1}\overline{f} \\
&= p(f_1b^\ca_1-b^\cd_1f_1)
+m(\overline{f}b^\ca_1+d^{\cm[1]}\overline{f}-\alpha f_1) = 0
\end{align*}
by \eqref{eq-(tm)bQ1-cone}. Therefore, $f'_1$ is a chain map and
\(f'_1\big|_{s\cd}=f_1\).
\end{proof}

\begin{proposition}\label{pro-restriction-A(EQA)-A(EA)}
In assumptions of \propref{pro-f-extends-(ff)} the restriction strict
\ainf-functor
\begin{equation}
\restr: A_{\infty1}^{\psi u}(\cd,\cq;\ca)
\to A_\infty^{\psi u}(\cd,\ca), \qquad (x,x') \mapsto x
\label{eq-A(EQA)-A(EA)}
\end{equation}
is an \ainf-equivalence, surjective on objects. The chain surjections
$\restr_1$ admit a chain splitting.
\end{proposition}

\begin{proof}
Let us prove that for an arbitrary pair of objects $(f,f')$, $(g,g')$
of \(A_{\infty1}^{\psi u}(\cd,\cq;\ca)\) the restriction chain map
\[ sA_{\infty1}(\cd,\cq;\ca)\bigl((f,f'),(g,g')\bigr)
\to sA_\infty(\cd,\ca)(f,g)
\]
is homotopy invertible. This map is a product over \(n\in\ZZ_{\ge0}\)
of the restriction maps \(\rho_n:V_n\to V'_n\) of the graded
$\kk$\n-modules of $n$\n-th components (compare \eqref{eq-V0-V1-Vn}
with analogous decomposition of \(sA_\infty(\cd,\ca)(f,g)\)). Clearly,
the maps \(\rho_n=\id\) for $n=0$ or for $n>1$. On the other hand, for
$n=1$
\begin{multline*}
\rho_1 =\prod \uCom(\inj^\cd,1):
\\
\prod_{X,Y\in\Ob\cd} \uCom(s\cq(X,Y),s\ca(Xf,Yg)) \to
\prod_{X,Y\in\Ob\cd} \uCom(s\cd(X,Y),s\ca(Xf,Yg))
\end{multline*}
is surjective with the kernel
\(\Ker\rho_1=\prod_{X,Y\in\Ob\cd}\uCom(s\cn(X,Y),s\ca(Xf,Yg))\). As in
proof of \propref{pro-f-extends-(ff)} $\cn(X,Y)\ne0$ implies that
\(\sS{_{Xf}}\uni^\ca_0\in\iota_Xf_1+\im b_1\subset\im b_1\) or
\(\sS{_{Yg}}\uni^\ca_0\in\iota_Yg_1+\im b_1\subset\im b_1\), hence,
\(s\ca(Xf,Yg)\) is contractible. Therefore, for all objects $X$, $Y$ of
$\cd$ the complex \(\uCom(s\cn(X,Y),s\ca(Xf,Yg))\) is contractible.
Thus, \(\Ker\restr_1=\Ker\rho_1\) is contractible.

Summing up, the first term of the semisplit exact sequence
\[ 0 \to \Ker\restr_1 \rTTo
sA_{\infty1}(\cd,\cq;\ca)\bigl((f,f'),(g,g')\bigr) \rTTo^{\restr_1}
sA_\infty(\cd,\ca)(f,g) \to 0
\]
is contractible. By \lemref{lem-semisplit-split} this sequence admits a
splitting chain map
\[ \nu: sA_\infty(\cd,\ca)(f,g) \to
sA_{\infty1}(\cd,\cq;\ca)\bigl((f,f'),(g,g')\bigr),
\]
and $\nu$ is homotopy inverse to $\restr_1$.

By \propref{pro-f-extends-(ff)} the surjection
\[ \Ob A_{\infty1}^{\psi u}(\cd,\cq;\ca) \ni (f,f')
\mapsto f \in \Ob A_\infty^{\psi u}(\cd,\ca)
\]
admits a splitting \(f\mapsto\wt{f}=(f,f')\). Applying
\corref{cor-Theorem8:8} we conclude that \eqref{eq-A(EQA)-A(EA)} is an
\ainf-equivalence.
\end{proof}

\begin{corollary}
\ainfu-2-transformation \eqref{eq-A(EQA)-A(EA)} is a natural
\ainfu-2-equivalence.
\end{corollary}

An easy converse to \lemref{lem-semisplit-split} is given by

\begin{lemma}\label{lem-split-contractible}
Let \(\beta:A\to B\), \(\nu:B\to A\) be chain maps of complexes of
$\kk$\n-modules, such that \(\nu\beta=\id_B\). Denote \(C=\Ker\nu\),
then \(A\simeq C\oplus B\). If $\beta$ is a homotopy isomorphism, then
the chain complex $C$ is contractible.
\end{lemma}

\begin{proof}
Being one-sided inverse to homotopy isomorphism $\beta$, the map $\nu$
is homotopy inverse to $\beta$. Therefore, \(\id_A-\beta\nu=hd+dh\) for
some homotopy \(h:A\to A\) of degree $-1$. Since $C$ is the image of
the idempotent \(\id_A-\beta\nu=\pr^C\cdot\inj^C\), we have
\[ \id_C = \inj^C\pr^C\inj^C\pr^C = \inj^C(\id_A-\beta\nu)\pr^C
= (\inj^Ch\pr^C)d + d(\inj^Ch\pr^C).
\]
Thus, \(\id_C=Hd+dH\) for \(H=\inj^Ch\pr^C:C\to C\), and $C$ is
contractible.
\end{proof}

\begin{theorem}\label{thm-restriction-DA-CAmodB-equi}
Let $\ca$ be a unital \ainf-category, and let
\(\cd=\cup_{j\ge0}\cd_j=\tcolimit_j\cd_j=\Quo(\cc|\cb)\) be as in
\secref{sec-constr-D}. Then the restriction strict \ainf-functor
\[ \restr: A_\infty^{\psi u}(\cd,\ca)
\to A_\infty^u(\cc,\ca)_{\modulo\cb}
\]
is an \ainf-equivalence, 2\n-natural in $\ca$, surjective on objects.
The chain surjections $\restr_1$ admit a chain splitting.
\end{theorem}

\begin{proof}
All restriction strict \ainf-functors in the sequence
\begin{multline}
A_\infty^u(\cc,\ca)_{\modulo\cb}
\lTTo A_{\infty1}^{\psi u}(\cd_0,\cq_1;\ca)
\lTTo A_\infty^{\psi u}(\cd_1,\ca) \\
\lTTo A_{\infty1}^{\psi u}(\cd_1,\cq_2;\ca)
\lTTo A_\infty^{\psi u}(\cd_2,\ca)
\lTTo A_{\infty1}^{\psi u}(\cd_2,\cq_3;\ca)
\lTTo \dots
\label{eq-A(CA)-A(D0Q1A)-A(D1A)-A(D1Q2A)}
\end{multline}
are \ainf-equivalences (and natural \ainfu-2-equivalences). They are
surjective on objects. The first components are surjective and admit a
chain splitting. For the first functor it follows from
\propref{pro-A(D0Q1A)-A(CA)modB-equivalence}. For other odd-numbered
functors it follows from \propref{pro-restriction-A(EQA)-A(EA)}.
Indeed, if \(X\in\Ob\cb\), then
 $\sS{_X}\uni^\cc_0=\sS{_X}\uni^\cc_0Hb_1\in\im(b_1:
 s\cd_1(X,X)\to s\cd_1(X,X))$.
For even-numbered functors it follows from
\thmref{thm-restriction-equivalence-quiver}.

Let us show that \(A_\infty^{\psi u}(\cd,\ca)\) is the inverse limit of
\eqref{eq-A(CA)-A(D0Q1A)-A(D1A)-A(D1Q2A)} on objects and on morphisms.
There are restriction strict \ainf-functors
\begin{align*}
\restr:A_\infty^{\psi u}(\cd,\ca) &\to A_\infty^u(\cc,\ca)_{\modulo\cb}
\subset A_\infty^u(\cc,\ca), &\quad f &\mapsto f|_\cc, \\
\restr: A_\infty^{\psi u}(\cd,\ca) &\to
A_{\infty1}^{\psi u}(\cd_j,\cq_{j+1};\ca),
&\quad f &\mapsto (f|_{\cd_j},f|_{\cq_{j+1}}), \quad j\ge0, \\
\restr: A_\infty^{\psi u}(\cd,\ca) &\to A_\infty^{\psi u}(\cd_j,\ca),
&\quad f &\mapsto f|_{\cd_j}, \quad j\ge1.
\end{align*}
They agree with the functors $\restr$ from
\eqref{eq-A(CA)-A(D0Q1A)-A(D1A)-A(D1Q2A)} in the sense
\(\restr\cdot\restr=\restr\).

Since \(\cd=\cup_{j\ge0}\cd_j\), pseudounital \ainf-functors
\(f:\cd\to\ca\) are in bijection with sequences $(f^j)_j$ of pseudounital
\ainf-functors \(f^j:\cd_j\to\ca\) such that
\(f^{j+1}\big|_{\cd_j}=f^j\). In other words,
\(\Ob A_\infty^{\psi u}(\cd,\ca)\) is the inverse limit of the sequence
of surjections
\begin{multline*}
\Ob A_\infty^u(\cc,\ca)_{\modulo\cb}
\lEpi \Ob A_{\infty1}^{\psi u}(\cd_0,\cq_1;\ca)
\lEpi \Ob A_\infty^{\psi u}(\cd_1,\ca) \\
\lEpi \Ob A_{\infty1}^{\psi u}(\cd_1,\cq_2;\ca)
\lEpi \Ob A_\infty^{\psi u}(\cd_2,\ca)
\lEpi \Ob A_{\infty1}^{\psi u}(\cd_2,\cq_3;\ca)
\lEpi \dots
\end{multline*}
In particular, the map
\[ \Ob \restr: \Ob A_\infty^{\psi u}(\cd,\ca)
\to \Ob A_\infty^u(\cc,\ca)_{\modulo\cb}
\]
is surjective.

Let \(f,g:\cd\to\ca\) be pseudounital \ainf-functors. Since
 \(\cd=
 \bigl(\cup_{j\ge0}\cd_j\bigr)\bigcup\bigl(\cup_{j\ge1}\cq_j\bigr)\),
\ainf-transformations \(p:f\to g:\cd\to\ca\) are in bijection with
sequences
\[ (p^0,p^{\prime1},p^1,p^{\prime2},p^2,p^{\prime3},\dots)
\]
of \ainf-transformations \(p^j:f|_{\cd_j}\to g|_{\cd_j}:\cd_j\to\ca\),
$j\ge0$, and $A_1$\n-transformations
\(p^{\prime j}:f|_{\cq_j}\to g|_{\cq_j}:\cq_j\to\ca\), $j\ge1$, such
that \(p^{j+1}\big|_{\cd_j}=p^j\), \(p^j\big|_{\cq_j}=p^{\prime j}\).
In other words, \(A_\infty(\cd,\ca)(f,g)\) is the inverse limit of the
sequence of splittable chain surjections
\begin{multline*}
A_\infty(\cc,\ca)(f|_\cc,g|_\cc) \lEpi
A_{\infty1}(\cd_0,\cq_1;\ca)
\bigl((f|_{\cd_0},f|_{\cq_1}),(g|_{\cd_0},g|_{\cq_1})\bigr)
\lEpi \dots \\
\lEpi A_\infty(\cd_j,\ca)(f|_{\cd_j},g|_{\cd_j}) \lEpi
A_{\infty1}(\cd_j,\cq_{j+1};\ca)
\bigl((f|_{\cd_j},f|_{\cq_{j+1}}),(g|_{\cd_j},g|_{\cq_{j+1}})\bigr)
\lEpi \dots
\end{multline*}
Since these surjections are splittable, the above sequence is
isomorphic to the sequence of natural projections
\[ C_0 \lEpi C_0\times C_1 \lEpi C_0\times C_1\times C_2 \lEpi \dots
\lEpi \prod_{m=0}^nC_m \lEpi \prod_{m=0}^{n+1}C_m \lEpi \dots
\]
for some complexes $C_m$ of $\kk$\n-modules. Its inverse limit is
\(\prod_{m=0}^\infty C_m\simeq A_\infty(\cd,\ca)(f,g)\). By
\lemref{lem-split-contractible} all $C_m$ are contractible for $m>0$.
Therefore, \(\prod_{m=1}^\infty C_m\) is contractible. We obtain a
split exact sequence
\[ 0 \to \prod_{m=1}^\infty C_m \rTTo A_\infty(\cd,\ca)(f,g)
\rTTo^\beta A_\infty(\cc,\ca)(f|_\cc,g|_\cc) \to 0
\]
with contractible first term. By \lemref{lem-semisplit-split}
\(\beta=s\restr_1s^{-1}\) is a homotopy isomorphism.

Using \corref{cor-Theorem8:8} we prove the theorem.
\end{proof}

\section{\texorpdfstring{Unitality of $\cd$}{Unitality of D}}
We are going to prove that if $\cc$ is strictly unital, then the
\ainf-category $\cd$ constructed in \secref{sec-constr-D} is not only
pseudounital, but unital with the unit elements
\(\sS{_X}\uni^\cc_0\in s\cd(|)(X,X)\subset s\cd(X,X)\). Let us describe
$\kk$\n-linear maps \(h:s\cd(X,Y)\to s\cd(X,Y)\) of degree $-1$ such
that
\begin{equation}
1 - (1\tens\uni^\cc_0)b_2 = b_1h + hb_1: s\cd(X,Y)\to s\cd(X,Y).
\label{eq-1-1ib-b1-hb1}
\end{equation}
The homotopy $h$ is called the \emph{right unit homotopy}. Let $t$ be a
$\cc$\n-admissible tree. Let $y$ be its rightmost leaf vertex. Let
\(v_0\prec v_1\prec\dots\prec v_p\prec v_{p+1}=y\) be the directed path
connecting the root $v_0$ with $y$, $p\ge0$. Vertices $v_i$ and
$v_{i+1}$ are connected by an edge for all \(0\le i\le p\). Let
\(t_i^+\), $1\le i\le p$ be the tree $t$ with an extra leaf attached on
the right to the vertex $v_i$ if $v_i$ is $n$\n-ary, $n>1$ as in
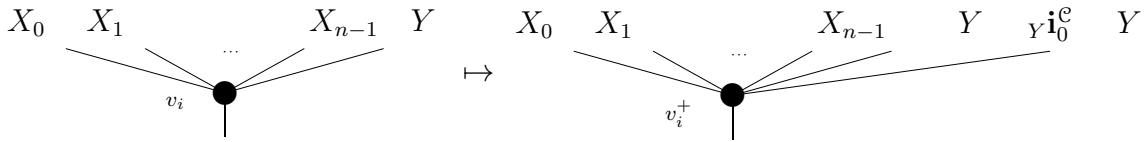
\[
\hstretch360
\begin{tanglec}
\object{X_0}\step\object{X_1}\step[3]\object{X_{n-1}}\step\object{Y}
\\
\nw2\nw1\nodeld{v_i}\node\Put(1,17)[cb]{\scriptscriptstyle\cdots}
\ne1\ne2 \\
\id
\end{tanglec}
\quad\mapsto\quad
\begin{tangle}
\object{X_0}\step\object{X_1}\step[3]\object{X_{n-1}}
\step[1.5]\object{Y}\step\object{\sS{_Y}\uni^\cc_0}\step\object{Y} \\
\hstep\nw2\nw1\node\nodeld{v_i^+}
\Put(1,17)[cb]{\scriptscriptstyle\cdots}\ne1\ne2\step\ne4 \\
\step[2.5]\id
\end{tangle}
\]
If $v_i$ is unary, we attach an extra leaf on the right to $v_i$ and
add two more unary vertices above and below it as in
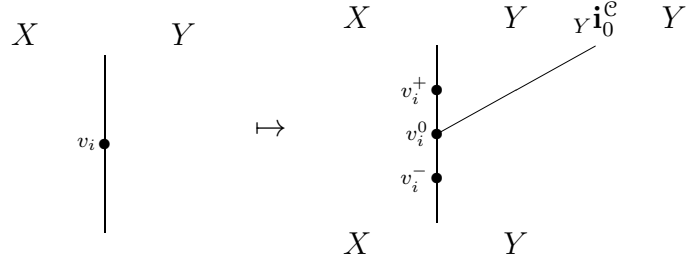
\[
\hstretch180
\vstretch200
\begin{tanglec}
\object{X}\step[4]\object{Y} \\
\nodel{v_i}\n \\
\id
\end{tanglec}
\qquad\mapsto\qquad
\vstretch100
\begin{tangles}{rl}
\object{X}\Step &
\Step\object{Y}\Step\object{\sS{_Y}\uni^\cc_0}\Step\object{Y} \\
& \nodel{v_i^+}\n\step[3]\ne2 \\
& \nodel{v_i^0}\n\step\ne2 \\
& \nodel{v_i^-}\n \\
& \id \\
\object{X}\Step & \Step\object{Y}
\end{tangles}
\]
The obtained trees \(t_i^+\) are $\cc$\n-admissible.

Let $x$ be a homogeneous element of \(s\cd(t)(Z,Y)\). Define
 \(x_i^+=x\tens\sS{_Y}\uni^\cc_0\in s\cd(t_i^+)(Z,Y)=s\cd(t)(Z,Y)\tens
 s\cc(Y,Y)\).
We claim that if $\cc$ is strictly unital, then the map
\begin{equation}
h: x \mapsto \sum_{i=1}^p \pm x_i^+
\label{eq-hx-pmxi}
\end{equation}
with an appropriate choice of signs satisfies \eqref{eq-1-1ib-b1-hb1}.
To describe the signs and to prove the claim we study the set of
operations acting in $\cd$.

\subsection{\texorpdfstring{A multicategory operating in $\cd$}
 {A multicategory operating in D}}
Let $\opA_\infty^{\cc/\cb}$ be the free graded $\kk$\n-linear
(non-symmetric) operad, generated by
\begin{itemize}
\item[---] a 0-ary operation \(\uni_0\in\opA_\infty^{\cc/\cb}(0)\) of
degree $-1$,

\item[---] a unary operation \(H\in\opA_\infty^{\cc/\cb}(1)\) of
degree $-1$,

\item[---] an $n$-ary operation \(b_n\in\opA_\infty^{\cc/\cb}(n)\) of
degree $1$ for each $n>1$.
\end{itemize}
The construction of free non-symmetric operads uses plane trees instead
of abstract trees. Otherwise it is similar to the case of symmetric
operads, see e.g. \cite{MR1898414}. The operad \(\opA_\infty\) of
operations in \ainf-algebras (e.g. \cite{math.AT/9808101}) is a
suboperad of $\opA_\infty^{\cc/\cb}$.

Actually we need a multicategory \cite{LambekJ:dedsc2,LambekJ:mulr} rather
than an operad. A multicategory is a many object version of an operad,
a (non-symmetric) operad is a one-object multicategory. So we define a
graded $\kk$\n-linear multicategory $\opA_\infty^{\cc/\cb}$, whose
objects are pairs of objects of $\cc$, thus
\(\Ob\opA_\infty^{\cc/\cb}=\Ob\cc\times\Ob\cc\). For $n>0$ the graded
$\kk$\n-module of morphisms
\[ \opA_\infty^{\cc/\cb}
\bigl((X'_1,X_1),(X'_2,X_2),\dots,(X'_n,X_n);(Y',Y)\bigr)
\]
is 0 unless \(X'_1=Y'\), \(X_n=Y\) and \(X_i=X'_{i+1}\) for all
\(1\le i<n\). For $n=0$ the graded $\kk$\n-module of morphisms
$\opA_\infty^{\cc/\cb}(;(Y',Y))$ is 0 unless $Y'=Y$.
The morphisms of $\opA_\infty^{\cc/\cb}$ are freely
generated by
\begin{itemize}
\item[---] 0-ary operations
\(\sS{_X}\uni_0\in\opA_\infty^{\cc/\cb}(;(X,X))\), \(X\in\Ob\cc\) of
degree $\deg\sS{_X}\uni_0=-1$,

\item[---] unary operations
\(H=H_{X,Y}\in\opA_\infty^{\cc/\cb}((X,Y);(X,Y))\), \(X,Y\in\Ob\cc\),
where \(X\in\Ob\cb\) or \(Y\in\Ob\cb\), of degree $\deg H=-1$,

\item[---] $n$-ary operations
 \(b_n\in\opA_\infty^{\cc/\cb}\bigl((X_0,X_1),(X_1,X_2),
 \dots,(X_{n-1},X_n);(X_0,X_n)\bigr)\),
of degree $\deg b_n=1$ for $X_0$, \dots, $X_n\in\Ob\cc$, $n>1$.
\end{itemize}
We shall not insist on distinguishing between the operad
\(\opA_\infty^{\cc/\cb}\) and its refinement -- the multicategory
\(\opA_\infty^{\cc/\cb}\), leaving the choice of context to the reader.

Similarly to the operad $\opA_\infty$ \cite{math.AT/9808101},
governing \ainf-algebras, the multicategory \(\opA_\infty^{\cc/\cb}\)
has a differential $d$ -- a derivation of degree 1, such that $d^2=0$.
In general, a derivation $d$ of degree $p$ of a graded $\kk$\n-linear
multicategory $\mcM$ is a collection of $\kk$\n-linear endomorphisms
$d$ of \(\mcM(Z_1,\dots,Z_n;Z)\) of degree $p$, such that all
compositions (which are of degree 0)
\begin{equation}
\mu_i: \mcM(Y_1,\dots,Y_k;Z_i) \tens \mcM(Z_1,\dots,Z_n;Z)
\to \mcM(Z_1,\dots,Z_{i-1},Y_1,\dots,Y_k,Z_{i+1},\dots,Z_n;Z)
\label{eq-mui-MM-M}
\end{equation}
satisfy the equation
\[ \mu_id = (1\tens d + d\tens1)\mu_i
\]
with the sign conventions of this article. If $d^2=0$, we may say that
$\mu_i$ are chain maps.

Since \(\opA_\infty^{\cc/\cb}\) is free, its derivations are uniquely
determined by their values on generators. In particular, the derivation
$d$ of degree 1 is determined by these values:
\begin{align*}
\sS{_X}\uni_0d &= 0, \\
H_{X,Y}d &= 1_{(X,Y)} \qquad (\text{the unit element of }
\opA_\infty^{\cc/\cb}\bigl((X,Y);(X,Y)\bigr)), \\
b_nd &= - \sum_{a+p+c=n}^{p>1,\;a+c>0}
(1^{\tens a}\tens b_p\tens1^{\tens c})b_{a+1+c}, \qquad n>1.
\end{align*}
For \(\opA_\infty^{\cc/\cb}\) we use the notation
\((1^{\tens a}\tens b_p\tens1^{\tens c})b_{a+1+c}\) referring to the
action of \(\opA_\infty^{\cc/\cb}\) in \ainf-algebras as a synonym of
the usual operadic notation
\((b_p\tens b_{a+1+c})\mu_{a+1}=b_p\circ_{a+1}b_{a+1+c}\). Since the
derivation $d$ has odd degree, the $\kk$\n-linear map $d^2$ is also a
derivation. Its value on all generators is 0 (for $b_n$ it follows from
a similar result for \(\opA_\infty\), see e.g.
\cite{math.AT/9808101,MR1898414}). Therefore, $d^2=0$, $d$ is a
differential, and \(\opA_\infty \rMono \opA_\infty^{\cc/\cb}\) is a
chain embedding.

The so defined differential $d$ is distinguished by the following
property. The action maps
\begin{multline}
\alpha: s\ce(X_0,X_1)\tens s\ce(X_1,X_2)\tdt s\ce(X_{n-1},X_n) \\
\tens \opA_\infty^{\cc/\cb}
\bigl((X_0,X_1),(X_1,X_2),\dots,(X_{n-1},X_n);(X_0,X_n)\bigr)
\to s\ce(X_0,X_n)
\label{eq-alpha-EEEA-E}
\end{multline}
are chain maps, where \( s\ce(\_,\_)\) is equipped with the
differential $b_1$. It suffices to check this on generators of
\(\opA_\infty^{\cc/\cb}\). For them the property follows from
\eqref{eq-bkb1-sum-(1b1)b}, \eqref{eq-Hb1-1-b1H} and the equation
\(\sS{_X}\uni^\cc_0b_1=0\).

Let \(\underline{\opA}_\infty^{\cc/\cb}\) be a submulticategory of
\(\opA_\infty^{\cc/\cb}\), generated by $H$ and $b_n$, $n\ge2$, with
the same set of objects
\(\Ob\underline{\opA}_\infty^{\cc/\cb}=\Ob\opA_\infty^{\cc/\cb}\). It
is a differential graded submulticategory without 0\n-ary operations:
\(\underline{\opA}_\infty^{\cc/\cb}(;(X,Y))=0\) for all $X,Y\in\Ob\cc$.
As a $\kk$\n-linear graded multigraph
\(\underline{\opA}_\infty^{\cc/\cb}\) has the following description.
For $n\ge1$
\begin{equation}
\underline{\opA}_\infty^{\cc/\cb}
\bigl((X_0,X_1),(X_1,X_2),\dots,(X_{n-1},X_n);(X_0,X_n)\bigr)
= \bigoplus_{\text{admissible }(t;X_0,X_1,\dots,X_n)} \kk[|t|_1-|t|_>],
\label{eq-opa-k[tt]}
\end{equation}
other $\kk$\n-modules
 \(\underline{\opA}_\infty^{\cc/\cb}
 \bigl((X'_1,X_1),(X'_2,X_2),\dots,(X'_n,X_n);(Y',Y)\bigr)\)
vanish.

The embedding of \(\underline{\opA}_\infty^{\cc/\cb}\) is denoted
\begin{equation}
\iota: \underline{\mcM} \overset{\text{def}}=
\underline{\opA}_\infty^{\cc/\cb} \rMono \mcM \overset{\text{def}}=
\opA_\infty^{\cc/\cb}.
\label{eq-iota-MuA-MA}
\end{equation}
The following general results are valid for an arbitrary embedding
\(\iota:\underline{\mcM} \rMono \mcM\) of differential graded
$\kk$\n-linear multicategories, such that \(\Ob\iota=\id_{\Ob\mcM}\),
and \(\underline{\mcM}(;Z)=0\) for all
\(Z\in\Ob\underline{\mcM}=\Ob\mcM\).

\begin{definition}
A \emph{right derivation} $\delta$ of degree $p$ of the embedding
\(\iota:\underline{\mcM} \rMono \mcM\) is a collection of
$\kk$\n-linear maps
\[ \delta: \underline{\mcM}(Z_1,\dots,Z_n;Z) \to \mcM(Z_1,\dots,Z_n;Z)
\]
of degree $p$ for all \(Z_1,\dots,Z_n,Z\in\Ob\mcM\), such that
compositions $\mu_i$ from \eqref{eq-mui-MM-M} satisfy
\[
\begin{diagram}[inline,nobalance]
\underline{\mcM}\tens\underline{\mcM} & \rTTo^{\mu_n}
& \underline{\mcM} \\
\dTTo<{\iota\tens\delta+\delta\tens\iota} & = & \dTTo>\delta \\
\mcM\tens\mcM & \rTTo^{\mu_n} & \mcM
\end{diagram}
\qquad , \qquad
\begin{diagram}[inline,nobalance]
\underline{\mcM}\tens\underline{\mcM} & \rTTo^{\mu_i}
& \underline{\mcM} \\
\dTTo<{\iota\tens\delta} & = & \dTTo>\delta \\
\mcM\tens\mcM & \rTTo^{\mu_i} & \mcM
\end{diagram}
\]
if \(1\le i<n\).
\end{definition}

Clearly, for $n=0$ the map \(\delta:\underline{\mcM}(;Z)=0\to\mcM(;Z)\)
is 0.

An example of a right derivation is given by an \emph{inner right
derivation}. Let \(\lambda_Z\in\mcM(Z;Z)\) be a family of morphisms of
degree $p$. For $n>0$ define
\begin{align*}
\ad_\lambda: \underline{\mcM}(Z_1,\dots,Z_n;Z) &\to \mcM(Z_1,\dots,Z_n;Z) \\
f &\mapsto (f\tens\lambda_Z)\mu_1 -(-)^{f\lambda}(\lambda_{Z_n}\tens f)\mu_n
= f\circ_1\lambda_Z -(-)^{fp}\lambda_{Z_n}\circ_nf.
\end{align*}
One verifies easily that \(\ad_\lambda\) is a right derivation, which
we call \emph{inner}.

One can show that if $\delta$ is a right derivation of $\iota$, and $d$
is a derivation of $\mcM$ such that
\(\underline{\mcM}d\subset\underline{\mcM}\), then their commutator
\([\delta,d]=\delta d-(-)^{\delta\cdot d}d\delta\) is a right
derivation of $\iota$ as well. In particular, it applies to a
differential $d$ of degree 1. If \(\delta=\ad_\lambda\) is inner, then
\([\delta,d]=\ad_{\lambda d}\) is also inner.

\subsection{Examples of right derivations}
Consider embedding \eqref{eq-iota-MuA-MA} and take the family of
morphisms
\[ \lambda_{(X,Y)} = (1\tens\sS{_Y}\uni_0)b_2 = \sS{_Y}\uni_0\circ_2b_2
\in \opA_\infty^{\cc/\cb}\bigl((X,Y);(X,Y)\bigr).
\]
Since \(\lambda_{(X,Y)}d=0\), we have
\([\ad_\lambda,d]=\ad_{\lambda d}=0\). Thus,
\(\delta=\ad_{(1\tens\uni_0)b_2}\) commutes with $d$.

Let us show that \(\ad_\lambda=[\eta,d]\) for some right derivation
\(\eta:\underline{\opA}_\infty^{\cc/\cb}\to\opA_\infty^{\cc/\cb}\) of
$\iota$. Since \(\underline{\opA}_\infty^{\cc/\cb}\) is a free
multicategory, any such right derivation is uniquely determined by its
value on the generators $H$ and $b_n$, $n>1$. We define a right
derivation $\eta$ of $\iota$ of degree $-1$ by the following
assignment:
\begin{align}
H_{X,Y}\eta &= H_{X,Y}(1\tens\sS{_Y}\uni_0)b_2H_{X,Y}, \notag \\
b_n\eta &= (1^{\tens n}\tens\sS{_{X_n}}\uni_0)b_{n+1}.
\label{eq-bneta-(1i)b}
\end{align}

Any operation \(f\in\underline{\opA}_\infty^{\cc/\cb}\) can be
presented as
\[ f= (g\tens1)(1^{\tens a_1}\tens e_{p_1})(1^{\tens a_2}\tens e_{p_2})
\dots(1^{\tens a_{k-1}}\tens e_{p_{k-1}})e_{p_k}
\]
for some \(g\in\underline{\opA}_\infty^{\cc/\cb}\), where $e_1=H$ and
$e_p=b_p$ if $p>1$. Then
\begin{equation}
f\eta = \sum_{i=1}^k (-)^{k-i} (g\tens1)(1^{\tens a_1}\tens e_{p_1})
\dots(1^{\tens a_i}\tens e_{p_i}\eta)
\dots(1^{\tens a_{k-1}}\tens e_{p_{k-1}})e_{p_k}.
\label{eq-feta-g11e1eeta1ee}
\end{equation}

Let us prove that
\begin{equation}
[\eta,d]  = \eta d + d\eta = \ad_{(1\tens\uni_0)b_2}.
\label{eq-eta-d-ad(1i)b2}
\end{equation}
Since the difference of the both sides is a right derivation of
$\iota$, it suffices to prove \eqref{eq-eta-d-ad(1i)b2} on generators.
First we find that
\begin{multline*}
H_{X,Y}(\eta d+d\eta)
= [H_{X,Y}(1\tens\sS{_Y}\uni_0)b_2H_{X,Y}]d + 1_{(X,Y)}\eta \\
= H_{X,Y}(1\tens\sS{_Y}\uni_0)b_2 - (1\tens\sS{_Y}\uni_0)b_2H_{X,Y}
= H_{X,Y}\ad_{(1\tens\uni_0)b_2}.
\end{multline*}
For $n\ge2$ we have
\begin{align*}
b_n(\eta d+d\eta) &= (1^{\tens n}\tens\uni_0)(b_{n+1}d)
- \sum_{a+p+c=n}^{p>1,\;c>0}
(1^{\tens a}\tens b_p\tens1^{\tens c})(b_{a+1+c}\eta) \\
&\phantom{=\ }
- \sum_{a+p=n}^{p>1,\;a>0} [(1^{\tens a}\tens b_p)b_{a+1}]\eta \\
&= - \sum_{a+p+c=n+1}^{p>1,\;a+c>0} (1^{\tens n}\tens\uni_0)
(1^{\tens a}\tens b_p\tens1^{\tens c})b_{a+1+c} \\
&\phantom{=\ } - \sum_{a+p+c=n}^{p>1,\;c>0}
(1^{\tens a}\tens b_p\tens1^{\tens c})
(1^{\tens a+1+c}\tens\uni_0)b_{a+c+2} \\
&\phantom{=\ } - \sum_{a+p=n}^{p>1,\;a>0}
(1^{\tens a}\tens b_p)(1^{\tens a+1}\tens\uni_0)b_{a+2}
+ \sum_{a+p=n}^{p>1,\;a>0}
[1^{\tens a}\tens(1^{\tens p}\tens\uni_0)b_{p+1}]b_{a+1} \\
&= - \sum_{a+p+c=n+1}^{p>1,\;a+c>0} (1^{\tens n}\tens\uni_0)
(1^{\tens a}\tens b_p\tens1^{\tens c})b_{a+1+c} \\
&\phantom{=\ } + \sum_{a+p+e=n+1}^{p>1,\;e>1} (1^{\tens n}\tens\uni_0)
(1^{\tens a}\tens b_p\tens1^{\tens e})b_{a+1+e} \\
&\phantom{=\ } + \sum_{a+p=n}^{p>1,\;a>0} (1^{\tens n}\tens\uni_0)
(1^{\tens a}\tens b_p\tens1)b_{a+2}
+ \sum_{a+q=n+1}^{q>2,\;a>0} (1^{\tens n}\tens\uni_0)
(1^{\tens a}\tens b_q)b_{a+1} \\
&= - \sum_{0+p+1=n+1} (1^{\tens n}\tens\uni_0)(b_p\tens1)b_2
-\sum_{a+2+0=n+1} (1^{\tens n}\tens\uni_0)(1^{\tens a}\tens b_2)b_{a+1}
\\
&= b_n(1\tens\uni_0)b_2 - [1^{\tens n-1}\tens(1\tens\uni_0)b_2]b_n \\
&= b_n\ad_{(1\tens\uni_0)b_2}.
\end{align*}
Thus equation~\eqref{eq-eta-d-ad(1i)b2} is verified.

Notice that the graded quiver $s\ce$ defined by \eqref{eq-E(XY)-E(t)(XY)}
is a free $\underline{\opA}_\infty^{\cc/\cb}$\n-algebra, generated by
the graded quiver $s\cc$. Indeed, \eqref{eq-E(XY)-E(t)(XY)} can be
written as
\begin{multline*}
s\ce(X,Y) = \bigoplus_{n\ge1}
\bigoplus_{X_1,\dots,X_{n-1}\in\Ob\cc}^{X_0=X,\;X_n=Y}
s\cc(X_0,X_1)\tdt s\cc(X_{n-1},X_n)\tens \\
\tens \underline{\opA}_\infty^{\cc/\cb}
\bigl((X_0,X_1),(X_1,X_2),\dots,(X_{n-1},X_n);(X_0,X_n)\bigr)
\end{multline*}
due to \eqref{eq-opa-k[tt]}. Compare with the usual free algebras over
an operad, e.g. \cite{MR1898414}. The operations $H$, $b_n$ for $n\ge2$
act in $s\ce$ via multicategory compositions in
\(\underline{\opA}_\infty^{\cc/\cb}\).

The multigraph \(\opA_\infty^{\cc/\cb}\) is expressed via
\(\underline{\opA}_\infty^{\cc/\cb}\) as follows:
\begin{multline}
\opA_\infty^{\cc/\cb}
\bigl((X_0,X_1),(X_1,X_2),\dots,(X_{n-1},X_n);(X_0,X_n)\bigr) \\
= \bigoplus_{k_0,\dots,k_n\in\ZZ_{\ge0}}
\underline{\opA}_\infty^{\cc/\cb} \bigl((Y_0,Y_1),(Y_1,Y_2),
\dots,(Y_{p-1},Y_p);(X_0,X_n)\bigr)[k_0+\dots+k_n],
\label{eq-Acb-undAcb}
\end{multline}
\[ (Y_0,Y_1,\dots,Y_p) =
(\underbrace{X_0,\dots,X_0}_{k_0+1},\underbrace{X_1,\dots,X_1}_{k_1+1},
\dots,\underbrace{X_n,\dots,X_n}_{k_n+1}),
\]
where \(p=k_0+\dots+k_n+n\). Indeed, 0\n-ary operations can be
performed first. The summand of \eqref{eq-Acb-undAcb} corresponds to
insertion of $k_0$ symbols \(\sS{_{X_0}}\uni_0\), $k_1$ symbols
\(\sS{_{X_1}}\uni_0\), and so on. In terms of trees such summand is
described by concatenation of the forest
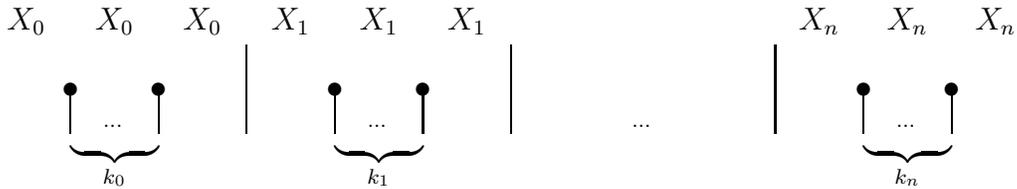
\[
\hstretch200
\begin{tangles}{ccccrclcc}
\object{X_0}\step & \object{X_0} & \step\object{X_0}
\Step\object{X_1}\step & \object{X_1} & \step\object{X_1}\step &&
\step\object{X_n}\step & \object{X_n} & \step\object{X_n} \\
&& \id && \id && \id && \\
& \s\step\nodeu{\dots}\step\s & \id & \s\step\nodeu{\dots}\step\s & \id
& \step[3]\nodeu{\dots}\step[3] & \id & \s\step\nodeu{\dots}\step\s & \\
& \object{\underbrace{\Step}_{k_0}} && \object{\underbrace{\Step}_{k_1}}
&&&& \object{\underbrace{\Step}_{k_n}} &
\end{tangles}
\]
with an admissible tree \((t;Y_0,Y_1,\dots,Y_p)\).

The action \eqref{eq-alpha-EEEA-E} of \(\opA_\infty^{\cc/\cb}\) in
$s\ce$ is described as follows. An element
\[ (\sS{_{X_0}}\uni_0^{\tens k_0}\tens1\tens
\sS{_{X_1}}\uni_0^{\tens k_1}\tens1\tdt1\tens
\sS{_{X_n}}\uni_0^{\tens k_n})\cdot f \in
\opA_\infty^{\cc/\cb}
\bigl((X_0,X_1),\dots,(X_{n-1},X_n);(X_0,X_n)\bigr),
\]
where
 \(f\in\underline{\opA}_\infty^{\cc/\cb}
 \bigl((Y_0,Y_1),(Y_1,Y_2),\dots,(Y_{p-1},Y_p);(X_0,X_n)\bigr)\),
acts by the map
\begin{multline*}
s\ce(X_0,X_1)\tdt s\ce(X_{n-1},X_n)
\rTTo^{\sS{_{X_0}^\cc}\uni_0^{\tens k_0}\tens1\tens
\sS{_{X_1}^\cc}\uni_0^{\tens k_1}\tens1\tdt1\tens
\sS{_{X_n}^\cc}\uni_0^{\tens k_n}} \\
s\cc(X_0,X_0)^{\tens k_0}\tens s\ce(X_0,X_1)\tens
s\cc(X_1,X_1)^{\tens k_1}\tdt
s\ce(X_{n-1},X_n)\tens s\cc(X_n,X_n)^{\tens k_n}
\end{multline*}
followed by the action of $f$ via multiplication in multicategory
\(\underline{\opA}_\infty^{\cc/\cb}\), that is, via grafting of trees.

Define a $\kk$\n-linear map \(\bar{h}:s\ce(X,Y)\to s\ce(X,Y)\) of
degree $-1$ for all objects $X$, $Y$ of $\cc$ as follows:
\begin{diagram}[height=2.9em]
s\ce(X,Y) \\
\dEq \\
\bigoplus_{X_1,\dots,X_{n-1}\in\Ob\cc}^{n,\;X_0=X,\;X_n=Y}
s\cc(X_0,X_1)\tdt s\cc(X_{n-1},X_n)\tens
\underline{\opA}_\infty^{\cc/\cb}
\bigl((X_0,X_1),\dots,(X_{n-1},X_n);(X_0,X_n)\bigr) \\
\dTTo<{-1^{\tens n}\tens\eta} \\
\bigoplus_{X_1,\dots,X_{n-1}\in\Ob\cc}^{n,\;X_0=X,\;X_n=Y}
s\cc(X_0,X_1)\tdt s\cc(X_{n-1},X_n)\tens
\opA_\infty^{\cc/\cb}
\bigl((X_0,X_1),\dots,(X_{n-1},X_n);(X_0,X_n)\bigr) \\
\dTTo<\alpha \\
s\ce(X,Y)
\end{diagram}
The concrete choice \eqref{eq-bneta-(1i)b} of value of $\eta$ on
generators shows that on the summand
\[ s\cc(X_0,X_1)\tdt s\cc(X_{n-1},X_n)\tens
\underline{\opA}_\infty^{\cc/\cb}
\bigl((X_0,X_1),\dots,(X_{n-1},X_n);(X_0,X_n)\bigr)
\]
the map $\bar{h}$ takes values in
\[ s\cc(X_0,X_1)\tdt s\cc(X_{n-1},X_n)\tens s\cc(X_n,X_n)\tens
\underline{\opA}_\infty^{\cc/\cb}
\bigl(\dots,(X_{n-1},X_n),(X_n,X_n);(X_0,X_n)\bigr)
\]
Due to \eqref{eq-feta-g11e1eeta1ee} the explicit formula for $\bar{h}$
has the form \eqref{eq-hx-pmxi} with the concrete choice of signs.

Let us compute the commutator
\begin{align*}
-\bar{h}b^\ce_1 -b^\ce_1\bar{h} &= (1^{\tens n}\tens\eta)\alpha b^\ce_1
+ b^\ce_1(1^{\tens n}\tens\eta)\alpha \\
&= (1^{\tens n}\tens\eta)\Bigl(\sum_{a+1+c=n}
(1^{\tens a}\tens b_1\tens1^{\tens c})\tens1+1^{\tens n}\tens d
\Bigr)\alpha \\
&\phantom{=\ }+ \Bigl(\sum_{a+1+c=n}
(1^{\tens a}\tens b_1\tens1^{\tens c})\tens1+1^{\tens n}\tens d
\Bigr)(1^{\tens n}\tens\eta)\alpha \\
&= [1^{\tens n}\tens(\eta d+d\eta)]\alpha
= [1^{\tens n}\tens\ad_{(1\tens\uni_0)b_2}]\alpha:
\end{align*}
\begin{equation*}
s\cc(X_0,X_1)\tdt s\cc(X_{n-1},X_n)\tens
\underline{\opA}_\infty^{\cc/\cb}
\bigl((X_0,X_1),\dots,(X_{n-1},X_n);(X_0,X_n)\bigr) \to
\end{equation*}
\begin{equation*}
s\cc(X_0,X_1)\tdt s\cc(X_{n-1},X_n)\tens s\cc(X_n,X_n)\tens
\underline{\opA}_\infty^{\cc/\cb}
\bigl(\dots,(X_{n-1},X_n),(X_n,X_n);(X_0,X_n)\bigr).
\end{equation*}

We write an element \(z_1\tdt z_n\tens f\) of the source, which is a
direct summand of \(s\ce(X_0,X_n)\), in the form \((z_1\tdt z_n)f\),
meaning that \(z_i\in s\cc(X_{i-1},X_i)\) and $f$ is a composition of
expressions \(1^{\tens a}\tens b_p\tens1^{\tens c}\) and
\(1^{\tens a}\tens H\tens1^{\tens c}\) ending in $b_p$ or $H$. Then
\begin{multline}
-[(z_1\tdt z_n)f](\bar{h}b^\ce_1+b^\ce_1\bar{h})
= [(z_1\tdt z_n)(f.\ad_{(1\tens\uni_0)b_2})]\alpha \\
= \bigl\{[(z_1\tdt z_n)f]\tens\sS{_{X_n}}\uni^\cc_0\bigr\}b^\ce_2
- (z_1\tdt z_n\tens\sS{_{X_n}}\uni^\cc_0)(1^{\tens n-1}\tens b^\ce_2)f.
\label{eq-zzf-hbbh-zzfib-zzi1bf}
\end{multline}

\begin{proposition}
Assume that $\cc$ is a unital \ainf-category with the unit elements
satisfying equations
\begin{align*}
(1\tens\uni^\cc_0)b_2 &= 1, \\
(1^{\tens n}\tens\uni^\cc_0)b_{n+1} &= 0, \qquad \text{if } n > 1.
\end{align*}
Then there exists a map \(h:s\cd\to s\cd\) such that
\[ \bigl(s\ce \rTTo^{\bar{h}} s\ce \rEpi^\pi s\cd\bigr)
= \bigl(s\ce \rEpi^\pi s\cd \rTTo^h s\cd\bigr).
\]
The map \(h:s\cd\to s\cd\) is a right unit homotopy for $\cd$.
\end{proposition}

\begin{proof}
We have to prove that \((R)_+\bar{h}\subset(R)_+\). The left hand side
is the sum of images of maps
\begin{equation}
[(1^{\tens a}\tens b^\ce_n\tens1^{\tens c})f]\eta
- (1^{\tens a}\tens b^\cc_n\tens1^{\tens c})(f\eta):
s\cc^{\tens k} \to s\ce.
\label{eq-1b1feta-1b1feta}
\end{equation}
If $c>0$, this map equals
\[ (1^{\tens a}\tens b^\ce_n\tens1^{\tens c})(f\eta)
- (1^{\tens a}\tens b^\cc_n\tens1^{\tens c})(f\eta),
\]
and the claim holds. If $c=0$, expression \eqref{eq-1b1feta-1b1feta} is
\begin{align*}
(1^{\tens a}\tens b^\ce_n)(f\eta) &- (1^{\tens a}\tens b^\cc_n)(f\eta)
+ (-)^f(1^{\tens a}\tens(1^{\tens n}\tens\uni^\cc_0)b^\ce_{n+1})f \\
&= (1^{\tens a}\tens b^\ce_n)(f\eta) -(1^{\tens a}\tens b^\cc_n)(f\eta)
\\
&\phantom{=\ }+
(-)^f(1^{\tens a+n}\tens\uni^\cc_0)(1^{\tens a}\tens b^\ce_{n+1})f
- (-)^f(1^{\tens a+n}\tens\uni^\cc_0)(1^{\tens a}\tens b^\cc_{n+1})f \\
&\phantom{=\ }+ (-)^f(1^{\tens a}\tens(1^{\tens n}\tens\uni^\cc_0)
b^\cc_{n+1})f: s\cc^{\tens k} \to s\ce
\end{align*}
and the last summand equals 0. So the claim holds, and $h$ exists.

Property~\eqref{eq-zzf-hbbh-zzfib-zzi1bf} turns into
\[ hb^\cd_1 + b^\cd_1h = 1 - (1\tens\uni^\cc_0)b_2:
s\cd(X,Y)\to s\cd(X,Y).
\]
Therefore, $h$ is a right unit homotopy for $\cd$.
\end{proof}

\begin{theorem}\label{thm-C-strict-D-unital}
Assume that $\cc$ is a unital \ainf-category with the unit elements
satisfying equations
\begin{align*}
(1\tens\uni^\cc_0)b_2 &= 1, &\qquad (\uni^\cc_0\tens1)b_2 &= -1, \\
(1^{\tens n}\tens\uni^\cc_0)b_{n+1} &= 0, &\qquad
(\uni^\cc_0\tens1^{\tens n})b_{n+1} &= 0, \qquad \text{if } n > 1.
\end{align*}
Then the \ainf-category $\cd$ is unital.
\end{theorem}

\begin{proof}
Besides constructing $\cd$ from the pair \((\cc,\cb)\), we may apply
the construction to the pair \((\cc^\op,\cb^\op)\), and we get an
\ainf-category isomorphic to $\cd^\op$. The opposite \ainf-category
$\ca^\op$ to an \ainf-category $\ca$ is the opposite quiver, equipped
with operations \(b^\op_k\), see \defref{def-opposite-ainf-category}.
In particular, \(b^\op_1=b_1\) and
\((x\tens\uni_0)b^\op_2=-x(\uni_0\tens1)b_2\). Thus we may use
\(h_\op=h\) for \(\cd^\op\) in place of $h'$ for \(\cd\). Thus,
\ainf-category $\cd$ is unital.
\end{proof}

\begin{corollary}\label{cor-C-strictly-unital-D-unital}
If $\cc$ is strictly unital, then the \ainf-category $\cd$ is unital.
\end{corollary}

\begin{corollary}
The \ainfu-2-functor \(\ca\mapsto A_\infty^u(\cc,\ca)_{\modulo\cb}\) is
unitally representable for an arbitrary unital \ainf-category $\cc$ by
\[ \bigl(\cc,e:\cc\to\quo(\cc|\cb)\bigr) \overset{\text{def}}=
\bigl(\cc,\cc \rTTo^{\wt{Y}} \wt{\cc} \rTTo^{\wt{e}}
\Quo(\wt{\cc}|\wt{\cb})\bigr),
\]
where \(\wt{Y}:\cc\to\wt{\cc}\) is the Yoneda \ainf-equivalence
identity on objects from \remref{rem-Yoneda-tilde}.
\end{corollary}

\begin{proof}
By \corref{cor-equi-pairs-equi-ainf2funs} it suffices to prove unital
representability for differential graded categories \(\wt{\cc}\) in
place of $\cc$. In this case the representing unital \ainf-category
$\Quo(\wt{\cc}|\wt{\cb})$ exists by
\corref{cor-C-strictly-unital-D-unital}.
\end{proof}

\section{\texorpdfstring{Equivalence of two quotients of
 $A_\infty$-categories}{Equivalence of two quotients of A8-categories}}
 \label{sec-Equivalence-two-quotients-Ainfty-categories}
Let $\cb$ be a full \ainf-subcategory of a unital \ainf-category $\cc$.
By \remref{rem-Yoneda-tilde} there exists a differential graded
category \(\wt{\cc}\) with \(\Ob\wt{\cc}=\Ob\cc\), and quasi-inverse to
each other \ainf-functors \(\wt{Y}:\cc\to\wt{\cc}\),
\(\Psi:\wt{\cc}\to\cc\) such that \(\Ob\wt{Y}=\Ob\Psi=\id_{\Ob\cc}\).
Let \(\wt{\cb}\subset\wt{\cc}\) be the full differential graded
subcategory with \(\Ob\wt{\cb}=\Ob\cb\). The quotient \ainf-category
\(\Quo(\wt{\cc}|\wt{\cb})\) and the quotient \ainf-functor
\(\wt{e}:\wt{\cc}\to\cd\) is constructed in
\thmref{thm-restriction-DA-CAmodB-equi}. The same
\(\Quo(\wt{\cc}|\wt{\cb})\) denoted also \(\quo(\cc|\cb)\) with the
quotient \ainf-functor
 \(e=\bigl(\cc \rTTo^{\wt{Y}} \wt{\cc} \rTTo^{\wt{e}} \cd\bigr)\)
represents the \ainfu-2-functor
\(\ca\mapsto A_\infty^u(\cc,\ca)_{\modulo\cb}\).

There is also a construction of \cite{LyuOvs-iResAiFn} which gives a
unital \ainf-category \(\Dr(\cc|\cb)\) and, in particular, a
differential graded category \(\Dr(\wt{\cc}|\wt{\cb})\). These are
smaller than \(\Quo(\wt{\cc}|\wt{\cb})\), however, we are going to
prove that all these three \ainf-categories are equivalent. Thus a
simpler construction \(\Dr(\cc|\cb)\) enjoys the same universal
properties as \(\quo(\cc|\cb)\) does. As a graded $\kk$\n-quiver
$\ca=\Dr(\cc|\cb)$ has the set of objects $\Ob\ca=\Ob\cc$, the
morphisms for $X,Y\in\Ob\ca$ are
\begin{equation*}
s\ca(X,Y) = \oplus_{C_1,\dots,C_{n-1}\in\cb} s\cc(X,C_1)\tens
s\cc(C_1,C_2)\tdt s\cc(C_{n-2},C_{n-1})\tens s\cc(C_{n-1},Y),
\end{equation*}
where the summation extends over all sequences of objects
$(C_1,\dots,C_{n-1})$ of $\cb$. To the empty sequence ($n=1$)
corresponds the summand $s\cc(X,Y)$. The operations
$\bar{b}_n:s\ca^{\tens n}\to s\ca$ are restrictions of maps
$\bar{b}_0=0$, $\bar{b}_1=b$ and for $n\ge2$
\begin{equation*}
\bar{b}_n = \mu^{(n)}
\sum_{m;q<k;t<l} 1^{\tens q}\tens b_m\tens1^{\tens t}:
T^ks\cc\tens(T^{\ge1}s\cc)^{\tens n-2}\tens T^ls\cc \to T^{\ge1}s\cc.
\end{equation*}
via the natural embedding $s\ca\subset T^{\ge1}s\cc$ of graded
$\kk$\n-quivers \cite[Proposition~2.2]{LyuOvs-iResAiFn}. Here
\(\mu^{(k)}:T^kT^{\ge1}s\cc\to T^{\ge1}s\cc\), $k\ge1$, is the
multiplication in the tensor algebra.

Denote by $\Quo(\cc|\cb)=\cd$ the quotient \ainf-category, constructed
in \thmref{thm-restriction-DA-CAmodB-equi}.

\begin{lemma}\label{lem-quiver-map-psi}
There is a chain quiver map \(\psi:s\Dr(\cc|\cb)\to s\Quo(\cc|\cb)\),
whose summands
\[ \psi_n: s\cc(X,C_1)\tens s\cc(C_1,C_2)\tdt s\cc(C_{n-2},C_{n-1})
\tens s\cc(C_{n-1},Y) \to s\Quo(\cc|\cb)(X,Y)
\]
for \(C_i\in\Ob\cb\) are defined by recurrent relation:
\(\psi_1=e_1:s\cc(X,Y) \rMono s\Quo(\cc|\cb)(X,Y)\) is the natural
embedding, and for $n>1$
\begin{equation}
\psi_n = - \sum_{k=1}^{n-1} (e_1^{\tens k}\tens\psi_{n-k}H)b_{k+1}.
\label{eq-psin-e1kpsiHb}
\end{equation}
\end{lemma}

For example, \(\psi_2=-(e_1\tens e_1H)b^\cd_2\),
\[ \psi_3 = (e_1\tens e_1\tens e_1H)(1\tens b^\cd_2H)b^\cd_2
- (e_1\tens e_1\tens e_1H)b^\cd_3.
\]
In general, expansion of \eqref{eq-psin-e1kpsiHb} contains $2^{n-2}$
summands.

\begin{proof}
We have to prove equation
\[ \psi_nb_1 =b\psi: s\cc(X,C_1)\tens s\cc(C_1,C_2)\tdt s\cc(C_{n-1},Y)
\to s\Quo(\cc|\cb)(X,Y)
\]
for all $n\ge1$. It is obvious for $n=1$. Let us prove it by induction.
Assume that it holds for number of factors smaller than $n$. Then
\begin{multline*}
\psi_nb_1 = - \sum_{k=1}^{n-1}(e_1^{\tens k}\tens\psi_{n-k}H)b_{k+1}b_1
=\sum_{\substack{0<k<n\\k+1=a+p+c}}^{a+c>0}(e_1^{\tens k}\tens\psi_{n-k}H)
(1^{\tens a}\tens b_p\tens1^{\tens c})b_{a+1+c}
\\
\quad = - \sum_{\substack{0<k<n\\k+1=a+p+c}}^{a\ge0,\,c>0}
(1^{\tens a}\tens b_p\tens1^{\tens n-a-p})
(e_1^{\tens a+c}\tens\psi_{n-k}H)b_{a+c+1} \hfill \\
+\sum_{\substack{0<k<n\\k+1=a+p}}^{a>0,\,p>1}(e_1^{\tens k}\tens\psi_{n-k}H)
(1^{\tens a}\tens b_p)b_{a+1}
+ \sum_{0<a<n}(e_1^{\tens a}\tens\psi_{n-a})b_{a+1}
- \sum_{0<k<n}(e_1^{\tens k}\tens\psi_{n-k}b_1H)b_{k+1}.
\end{multline*}
The second sum in the right hand side can be presented as
\[ \sum_{0<a<n-1} \Bigl[e_1^{\tens a}\tens \sum_{0<k-a<n-a}
(e_1^{\tens k-a}\tens\psi_{n-a-(k-a)}H)b_{k-a+1}\Bigr] b_{a+1}.
\]
Thus, it nearly cancels with the third sum except for one summand,
corresponding to $a=n-1$. In the fourth sum in the right hand side we
replace \(\psi_{n-k}b_1\) with $b\psi$ by induction assumption, and we
get
\begin{multline*}
\psi_nb_1 = - \sum_{\substack{0<k<n\\k+1=a+p+c}}^{a\ge0,\,c>0}
(1^{\tens a}\tens b_p\tens1^{\tens n-a-p})
(e_1^{\tens a+c}\tens\psi_{n-k}H)b_{a+c+1}
+ (e_1^{\tens n-1}\tens\psi_1)b_n \\
- \sum_{\substack{0<k<n\\k+\alpha+\beta+\gamma=n}}
(1^{\tens k+\alpha}\tens b_\beta\tens1^{\tens\gamma})
(e_1^{\tens k}\tens\psi_{\alpha+1+\gamma}H)b_{k+1} = b\psi,
\end{multline*}
since
\((e_1^{\tens n-1}\tens\psi_1)b_n=e_1^{\tens n}b_n=b_ne_1=b_n\psi_1\).
\end{proof}

According to \thmref{thm-restriction-DA-CAmodB-equi} for any unital
\ainf-category $\ca$ the map
\begin{multline*}
\restr: A_\infty^u(\cd,\ca) \to A_\infty^u(\cc,\ca)_{\modulo\cb} \\
= \bigl\{ f \in A_\infty^u(\cc,\ca) \mid f|_\cb =
\bigl(\cb \rMono \cc \rTTo^f \ca\bigr) \text{ is contractible } \bigr\}
\end{multline*}
is surjective. A splitting of this surjection is defined recurrently in
Propositions \ref{pro-A(D0Q1A)-A(CA)modB-equivalence},
\ref{pro-tilde-f-E-A}, \ref{pro-f-extends-(ff)}. Let us describe
another splitting map which differs from the mentioned one and is more
suitable for our purposes.

\begin{proposition}\label{pro-fCA-extends-to-tildefDA}
Let \(f:\cc\to\ca\) be a unital \ainf-functor such that $f|_\cb$ is
contractible. For each pair $X$, $Y$ of objects of $\cc$ such that
\(X\in\Ob\cb\) or \(Y\in\Ob\cb\) choose a contracting homotopy
\(\chi_{XY}:s\ca(Xf,Yf)\to s\ca(Xf,Yf)\), thus, \(\chi b_1+b_1\chi=1\).
Let \(\tilde{f}_k:T^ks\cd\to s\ca\), $k>1$, be $\kk$\n-quiver morphisms
of degree 0 which extend the components $f_k:T^ks\cc\to s\ca$. Then
there exists a unique extension of $f_1:s\cc\to s\ca$ to a quiver
morphism $\tilde{f}_1:s\cd\to s\ca$ such that
$(\tilde{f}_1,\tilde{f}_2,\dots)$ are components of a unital
\ainf-functor $\tilde{f}:\cd\to\ca$ and
\begin{equation}
H\tilde{f}_1 = \tilde{f}_1\chi: s\cd(X,Y) \to s\ca(Xf,Yf),
\label{eq-Hf1-f1chi}
\end{equation}
whenever \(X\in\Ob\cb\) or \(Y\in\Ob\cb\).
\end{proposition}

\noindent
    {\sl Warning:}
Extensions $\tilde{f}:\cd\to\ca$ of \(f:\cc\to\ca\) constructed in
Sections \ref{sec-A-categories-and-quivers},
\ref{sec-Relatively-free-categories} do not, in general, satisfy
condition \eqref{eq-Hf1-f1chi}.

\begin{proof}
Let us extend $f$ to an \ainf-functor \(\hat{f}:\ce\to\ca\) such that
 \(\hat{f}_k=\bigl(s\ce^{\tens k} \rEpi s\cd^{\tens k}
 \xrightarrow{\tilde{f}_k} \ca\bigr)\)
and \(H\hat{f}_1=\hat{f}_1\chi\), whenever \(X\in\Ob\cb\) or
\(Y\in\Ob\cb\). Suppose that $t_1$, \dots, $t_n$ are trees, $n\ge1$,
and $\hat{f}_1:s\ce(t_i)\to s\ca$ is already defined for all
 $1\le i\le n$. Then there is only one way to extend $\hat{f}$ on
\(\ce(t)\) for $t=(t_1\sqcup\dots\sqcup t_n)\cdot\tree_n$, where
$\tree_n$ is the corolla with $n$ leaves. Since
\begin{alignat*}2
b_n = s^{|t_1|}\tdt s^{|t_{n-1}|}\tens s^{|t_n|-|t|} &:
s\ce(t_1)\tdt s\ce(t_n) \to s\ce(t), &\qquad &\text{ for } n>1, \\
H = s &: s\ce(t_1) \to s\ce(t), &\qquad &\text{ for } n=1
\end{alignat*}
is invertible, we find that, respectively,
\begin{align*}
\hat{f}_1 &= \bigl( s\ce(t) \rTTo^{b_n^{-1}} s\ce(t_1)\tdt s\ce(t_n)
\rTTo^{\sum(\hat{f}_{i_1}\tdt \hat{f}_{i_l})b_l -
\sum_{\alpha+k+\beta=n}^{\alpha+\beta>0}
(1^{\tens\alpha}\tens b_k\tens1^{\tens\beta})\hat{f}_{\alpha+1+\beta}}
s\ca \bigr),
\\
\hat{f}_1 &= \bigl( s\ce(t) \rTTo^{H^{-1}} s\ce(t_1) \rTTo^{\hat{f}_1}
s\ca(Xf,Yf) \rTTo^{\chi_{XY}} s\ca(Xf,Yf)\bigr) \qquad \text{ for }n=1.
\end{align*}

Let us prove that coalgebra homomorphism $\hat{f}:\ce\to\ca$ with
recursively defined components $(\hat{f}_1,\hat{f}_2,\dots)$ is an
\ainf-functor. Equation
\begin{equation}
b_n\hat{f}_1 =
\sum_{i_1+\dots+i_l=n}(\hat{f}_{i_1}\tdt \hat{f}_{i_l})b_l
- \sum_{\alpha+k+\beta=n}^{\alpha+\beta>0}
(1^{\tens\alpha}\tens b_k\tens1^{\tens\beta})\hat{f}_{\alpha+1+\beta}:
T^ns\ce\to s\ca
\label{eq-bnf1-def-f1}
\end{equation}
is satisfied for all $n>1$ by construction of $\hat{f}_1$. So it
remains to prove that $\hat{f}_1$ is a chain map. To prove by induction
on the number of vertices of $t$ that
$\hat{f}_1b_1=b_1\hat{f}_1:s\ce(t)\to s\ca$, it suffices to show that
$b_n\hat{f}_1b_1=b_nb_1\hat{f}_1$ for all $n>1$ and that
\(H\hat{f}_1b_1=Hb_1\hat{f}_1\) due to invertibility of $b_n$ and $H$.
The first assertion is proven in
\cite[Proposition~2.3]{LyuOvs-iResAiFn}. The second follows from the
computation
\begin{multline*}
H(\hat{f}_1b_1-b_1\hat{f}_1) = H\hat{f}_1b_1 - Hb_1\hat{f}_1
= \hat{f}_1\chi b_1 - \hat{f}_1 + b_1H\hat{f}_1 \\
= \hat{f}_1 - \hat{f}_1b_1\chi - \hat{f}_1 + b_1\hat{f}_1\chi
= -(\hat{f}_1b_1-b_1\hat{f}_1)\chi,
\end{multline*}
which vanishes by induction assumption.

Using \eqref{eq-Hf1-f1chi} and \eqref{eq-bnf1-def-f1} one can prove
that the ideal \((R_\cc)_+\) is mapped by \(\hat{f}_1\) to 0.
Therefore, $\hat{f}$ factorizes as
\(\ce \rEpi \cd \rTTo^{\tilde{f}} \ca\) for a unique \ainf-functor
$\tilde{f}$. It is unital since the unit elements of $\cc$ are also the
unit elements of $\cd$.
\end{proof}

\begin{lemma}\label{lem-iso2iso}
Let \(f:\ca\to\cb\) be an \ainf-equivalence. Let objects $Xf$, $Yf$ of
$\cb$ be isomorphic via inverse to each other isomorphisms
 \(r\in s\cb(Xf,Yf)\), \(p\in s\cb(Yf,Xf)\) (that is,
\([rs^{-1}]\in H^0\cb(Xf,Yf)\), \([ps^{-1}]\in H^0\cb(Yf,Xf)\) are
inverse to each other in the ordianry category $H^0\cb$). Then the
objects $X$, $Y$ of $\ca$ are isomorphic via inverse to each other
isomorphisms \(q\in s\ca(X,Y)\), \(t\in s\ca(Y,X)\) such that
\(qf_1-r\in\im b_1\), \(tf_1-p\in\im b_1\).
\end{lemma}

\begin{proof}
Let chain maps \(g_{X,Y}:s\cb(Xf,Yf)\to s\ca(X,Y)\),
\(g_{Y,X}:s\cb(Yf,Xf)\to s\ca(Y,X)\) be homotopy inverse to maps
\(f_1:s\ca(X,Y)\to s\cb(Xf,Yf)\), \(f_1:s\ca(Y,X)\to s\cb(Yf,Xf)\).
Define $q=rg_{X,Y}$, $t=pg_{Y,X}$. Then
\begin{multline*}
[(q\tens t)b_2-\sS{_X}\uni^\ca_0]f_1
= (rg\tens pg)b_2f_1 -\sS{_X}\uni^\ca_0f_1
\equiv (rg\tens pg)(f_1\tens f_1)b_2 -\sS{_{Xf}}\uni^\cb_0 \\
= [(r+vb_1)\tens(p+wb_1)]b_2 -\sS{_{Xf}}\uni^\cb_0 \equiv
(r\tens p)b_2 -\sS{_{Xf}}\uni^\cb_0 \equiv 0 \pmod{\im b_1}.
\end{multline*}
Hence,
\[ (q\tens t)b_2-\sS{_X}\uni^\ca_0 \equiv
[(q\tens t)b_2-\sS{_X}\uni^\ca_0]f_1g_{X,X} \equiv 0 \pmod{\im b_1}.
\]
By symmetry, \((t\tens q)b_2-\sS{_Y}\uni^\ca_0\in\im b_1\). Other
properties are easy to verify.
\end{proof}

We are going to apply \propref{pro-fCA-extends-to-tildefDA} to the
unital quotient \ainf-functor $\ju:\cc\to\Dr(\cc|\cb)$, constructed in
\cite{LyuOvs-iResAiFn}. When restricted to $\cb$ the \ainf-functor
$\ju$ becomes contractible, therefore, there exists a unital
\ainf-functor \(f:\quo(\cc|\cb)\to\Dr(\cc|\cb)\) (unique up to an
isomorphism) such that $\ju$ is isomorphic to the composition
\(\cc\rTTo^e \quo(\cc|\cb)\rTTo^f \Dr(\cc|\cb)\).

\begin{proposition}
The \ainf-functor \(f:\quo(\cc|\cb)\to\Dr(\cc|\cb)\) (defined up to an
isomorphism) is an equivalence.
\end{proposition}

\begin{proof}
In the following diagram the top and the bottom rows compose to
contractible \ainf-functors.
\begin{diagram}
\cb & \rMono & \cc && \rTTo^{\ju^\cc} && \Dr(\cc|\cb) \\
&&& \rdTTo^e & \dTwoar<\alpha & \ruTTo^f \ruTwoar(1,3)>\beta & \\
\dTTo<{\wt{Y}} &=& \dTTo<{\wt{Y}} &=&
\Quo(\wt{\cc}|\wt{\cb}) && \dTTo>{\Dr(\wt{Y})} \\
&&& \ruTTo_{\wt{e}} &=& \rdTTo_g & \\
\wt{\cb} & \rMono & \wt{\cc} &&
\rTTo^{\ju^{\wt\cc}} && \Dr(\wt{\cc}|\wt{\cb})
\end{diagram}
Here the \ainf-functor $f$ and an isomorphism $\alpha$ exist due to $e$
being quotient \ainf-functor. The existence of \ainf-functor $g$ such
that \(\wt{e}g=\ju^{\wt\cc}\) follows from
\thmref{thm-restriction-DA-CAmodB-equi}. The isomorphism
\[ \alpha\Dr(\wt{Y}): eg = \ju^\cc\Dr(\wt{Y})
\to ef\Dr(\wt{Y}): \cc \to \Dr(\wt{\cc}|\wt{\cb})
\]
is equivalent by \lemref{lem-iso2iso} to $e\beta$ for some isomorphism
\[ \beta: g \to f\Dr(\wt{Y}):
\Quo(\wt{\cc}|\wt{\cb}) \to \Dr(\wt{\cc}|\wt{\cb}).
\]
The \ainf-functor \(\Dr(\wt{Y})\) is an equivalence by corollary~4.9
and section~5 of \cite{LyuOvs-iResAiFn}. Therefore, if we prove that
$g$ is an \ainf-equivalence, then $f$ is an \ainf-equivalence as well.
So in the following we consider only the case of strictly unital
\ainf-category $\cc$, and we are proving that the natural \ainf-functor
\(f:\Quo(\cc|\cb)\to\Dr(\cc|\cb)\) (defined up to an isomorphism) is an
equivalence.

Let $\ce$ be an arbitrary unital \ainf-category. Let \(\cf\subset\ce\)
be its full contractible subcategory, that is, complexes
\((s\ce(X,X),b_1)\) are contractible for all objects $X$ of $\cf$. Let
\(e:\cc\to\ce\) be a unital \ainf-functor such that
\(e(\Ob\cb)\subset\Ob\cf\).
Then there is a unital \ainf-functor
\(\Dr(e):\Dr(\cc|\cb)\to\Dr(\ce|\cf)\)
\cite[Corollary~5.6]{LyuOvs-iResAiFn}. There is a unital \ainf-functor
$\pi^\ce:\Dr(\ce|\cf)\to\ce$, quasi-inverse to the canonical strict
embedding $\ju^\ce:\ce\to\Dr(\ce|\cf)$ and such that
$\ju^\ce\pi^\ce=\id_\ce$ \cite[Proposition~7.4]{LyuOvs-iResAiFn}. In
particular, \(\Ob\pi^\ce=\id_\ce\). The diagram
\begin{diagram}
\cb & \rMono & \cc & \rTTo^{\ju^\cc} & \Dr(\cc|\cb) && \\
\dTTo<e && \dTTo>e && \dTTo>{\Dr(e)} && \\
\cf & \rMono & \ce & \rTTo^{\ju^\ce} & \Dr(\ce|\cf) & \rTTo^{\pi^\ce}
& \ce
\end{diagram}
is commutative due to [\textit{ibid}, Corollary~3.2]. Thus the
composition \(h=\Dr(e)\pi^\ce:\Dr(\cc|\cb)\to\ce\) is a unital
\ainf-functor such that \(\Ob h=\Ob e\) and
\(\ju^\cc h=\ju^\cc\Dr(e)\pi^\ce=e\ju^\ce\pi^\ce=e\).

Now we apply these considerations to the quotient functor
\(e:\cc\to\ce=\cd=\Quo(\cc|\cb)\). Define \(\cf\subset\Quo(\cc|\cb)\)
as a full \ainf-subcategory with \(\Ob\cf=\Ob\cb\). It is contractible.
In the commutative diagram
\begin{diagram}
\cb & \rMono & \cc & \rTTo^{\ju^\cc} & \Dr(\cc|\cb) && \\
\dTTo<e && \dTTo<e && \dTTo<{\Dr(e)} & \rdTTo^h & \\
\cf & \rMono & \cd & \rTTo^{\ju^\cd} & \Dr(\cd|\cf) & \rTTo^{\pi^\cd}
& \cd
\end{diagram}
the composition \(h=\Dr(e)\pi^\cd:\Dr(\cc|\cb)\to\cd\) satisfies
equation \(\ju^\cc h=e\).

The restriction of the strict \ainf-functor
\(\ju^\cc:\cc\to\Dr(\cc|\cb)\) to $\cb$ is contractible by
\cite[Example~6.6]{LyuOvs-iResAiFn}. Choose the maps
\[ \chi = - \sS{_X}\uni^\cc_0\tens1 =
[(\uni^\cc_0\tens\uni^\cc_0)\tens1]\bar{b}_2:
s\Dr(\cc|\cb)(X,Y) \to s\Dr(\cc|\cb)(X,Y)
\]
as contracting homotopies if \(X\in\Ob\cb\). Indeed,
\((\uni^\cc_0\tens\uni^\cc_0)\bar{b}_1=\uni^\cc_0\) and
\[ [(\uni^\cc_0\tens\uni^\cc_0)\tens1]\bar{b}_2\bar{b}_1
+ \bar{b}_1[(\uni^\cc_0\tens\uni^\cc_0)\tens1]\bar{b}_2
= - [(\uni^\cc_0\tens\uni^\cc_0)\bar{b}_1\tens1]\bar{b}_2
= - (\uni^\cc_0\tens1)\bar{b}_2 = 1,
\]
because \(\uni^\cc_0\) is the strict unit of \(\ca=\Dr(\cc|\cb)\). If
\(X\notin\Ob\cb\), but \(Y\in\Ob\cb\) we choose the contracting
homotopies
\[ \chi = 1\tens\sS{_Y}\uni^\cc_0:
s\Dr(\cc|\cb)(X,Y) \to s\Dr(\cc|\cb)(X,Y).
\]
Using \propref{pro-fCA-extends-to-tildefDA} we extend $\ju^\cc$ to a
unique unital strict \ainf-functor \(f:\cd\to\ca=\Dr(\cc|\cb)\) which
satisfies the equation
\[ \bigl(s\cd \rTTo^H s\cd \rTTo^{f_1} s\Dr(\cc|\cb) \bigr)
= \bigl(s\cd \rTTo^{f_1} s\Dr(\cc|\cb) \rTTo^\chi s\Dr(\cc|\cb) \bigr)
\]
whenever the left hand side is defined. In particular, \(ef=\ju^\cc\)
and \(\Ob f=\id_{\Ob\cc}\). The composition \(fh:\cd\to\cd\) satisfies
equation \(efh=\ju^\cc h=e=e\id_\cd:\cc\to\cd\). According to
\thmref{thm-maintheorem} the strict \ainf-functor given by composition
with $e$
\[ (e\boxtimes1)M: A_\infty^u(\cd,\cd) \to
A_\infty^u(\cc,\cd)_{\modulo\cb}, \qquad f \mapsto ef,
\]
is an \ainf-equivalence. Therefore, \(fh\simeq\id_\cd\) due to
\lemref{lem-iso2iso}. We conclude that \(f_1h_1\) is homotopy
invertible, and $f_1$ is homotopy invertible on the right.

We claim that
 \(\bigl(s\ca(X,Y) \rTTo^\psi s\ca(X,Y) \rTTo^{f_1} s\ca(X,Y)\bigr)
 =\id\),
where $\psi$ is constructed in \lemref{lem-quiver-map-psi}. Indeed,
consider this equation on the summand
 \(s\cc(X,C_1)\tens s\cc(C_1,C_2)\tdt s\cc(C_{n-2},C_{n-1})\tens
 s\cc(C_{n-1},Y)\)
of \(s\ca(X,Y)\). For $n=1$ the equation \(e_1f_1=1\) is obvious. Let
us prove it by induction on $n$. If $n>1$, then
\begin{multline*}
\psi_nf_1 = -\sum_{k=1}^{n-1} (e_1^{\tens k}\tens\psi_{n-k}H)b_{k+1}f_1
= - \sum_{k=1}^{n-1} (1^{\tens k}\tens\chi)\bar{b}_{k+1}
= - (1\tens\chi)\bar{b}_2 \\
= [1\tens(\uni^\cc_0\tens1^{\tens n-1})]\bar{b}_2
= (1\tens\uni^\cc_0)b_2\tens1^{\tens n-1} = 1.
\end{multline*}
Since $f_1$ is homotopy invertible on the right and on the left, it is
homotopy invertible. Since \(\Ob f=\id_{\Ob\cc}\), the \ainf-functor
\(f:\Quo(\cc|\cb)\to\Dr(\cc|\cb)\) is an equivalence.
\end{proof}

\section{The example of complexes}\label{sec-example-complexes}
Let $\ca$ be a $\kk$\n-linear Abelian category, let
$\cc=\underline{\mathsf C}(\ca)$ be the differential graded category of
complexes in $\ca$, and let $\cb$ be its full subcategory of acyclic
complexes. Let $\cd=\Quo(\cc|\cb)$ be the quotient unital
\ainf-category. The embedding $e:\cc\hookrightarrow\cd$ induces a
$\kk$\n-linear functor of homotopy categories
\[ H^0e: \Ho\cc = H^0(\cc,m_1) \to H^0(\cd,m_1) = \Ho\cd,
\]
where \(m_1=sb_1s^{-1}\). Morphisms of $\Ho\cc$ are homotopy
equivalence classes $[q]$ of chain morphisms \(q:X\to Y\).

\begin{proposition}
For a quasi-isomorphism $q$ its image \([qe]\in\Ho\cd(X,Y)\) is
invertible.
\end{proposition}

\begin{proof}
If \(q:X\to Y\) is a quasi-isomorphism, then \(C=\Cone q\) is acyclic.
The complex $C$ is the graded object \(Y\oplus X[1]\) with the
differential given by the matrix
\[ d^C =
\begin{pmatrix}
d^Y & 0 \\
s^{-1}q & d^{X[1]}
\end{pmatrix}
=
\begin{pmatrix}
d^Y & 0 \\
s^{-1}q & -s^{-1}d^Xs
\end{pmatrix}
    .
\]
There is a semisplit exact sequence of chain morphisms
$0\to Y \rTTo^n C \rTTo^k X[1]\to0$, where $n=\inj^Y=(1,0)$,
 \(k=\pr^{X[1]}=\bigl(
\begin{smallmatrix}
0 \\ 1
\end{smallmatrix}
 \bigr)\).
Thus \(n\in\cc(Y,C)^0\), \(k\in\cc(C,X[1])^0\) are cycles, $nm_1=0$,
$km_1=0$. The morphisms \(s\in\cc(X,X[1])^{-1}\),
\(s^{-1}\in\cc(X[1],X)^1\) are also cycles, because in our conventions
\(sm_1=sd^{X[1]}+d^Xs=0\), similarly \(s^{-1}m_1=0\). Hence,
\(ks^{-1}\in\cc(C,X)^1\) is also a cycle, \((ks^{-1})m_1=0\).

Since \(C\in\Ob\cb\), we have a map
\[ \underline{H} = sHs^{-1}: \cc(C,C) \rTTo \cd(\tree_1)(C,C) \rMono
\cd(C,C), \qquad f \mapsto fsHs^{-1} = f\underline{H}.
\]
It satisfies the equation
\(\underline{H}m_1+m_1\underline{H}=e:\cc(C,C)\to\cd(C,C)\). In
particular, there is an element \(1_C\underline{H}\in\cd(C,C)\). Define
an element
\[ p =(n\tens1_C\underline{H}\tens ks^{-1})(1\tens m_2)m_2 \in\cd(Y,X),
\]
where \(m_2=(s\tens s)b_2s^{-1}\). More generally,
\(m_n=s^{\tens n}b_ns^{-1}\). We have \(\deg p=0\) and
\[ pm_1 = -(n\tens1_C\tens ks^{-1})(1\tens m_2)m_2 = nks^{-1} = 0.
\]

Let us show that \([p]\in H^0\cd(Y,X)\) is inverse to
\([q]\in H^0\cd(X,Y)\). Define
\(h=\bigl(X \rTTo^s X[1] \rTTo^{\inj^{X[1]}} C\bigr)\), then
\(h\in\cc(X,C)^{-1}\), and
\[ hm_1 = hd^C + d^Xh = (0,s)
\begin{pmatrix}
d^Y & 0 \\
s^{-1}q & -s^{-1}d^Xs
\end{pmatrix}
+ d^X(0,s) = (q,0) = qn.
\]
One can check that
\begin{multline*}
(h\tens1_C\underline{H}\tens ks^{-1})(1\tens m_2)m_2m_1
-(q\tens n\tens1_C\underline{H}\tens ks^{-1})(1\tens1\tens m_2)m_3m_1 \\
= -hks^{-1}
+(q\tens n\tens1_C\underline{H}\tens ks^{-1})(1\tens1\tens m_2)
(1\tens m_2)m_2
= -1_X + (q\tens p)m_2,
\end{multline*}
because \(hk=s\).

Denote by $z$ the morphism \(z=\pr^Y:C\to Y\), then \(z\in\cc(C,Y)^0\)
and
\[ zm_1 = zd^Y - d^Cz = \binom10 d^Y -
\begin{pmatrix}
d^Y & 0 \\
s^{-1}q & d^{X[1]}
\end{pmatrix}
\binom10 = - \binom0{s^{-1}q} = -ks^{-1}q.
\]
One can check that
\begin{multline*}
(n\tens1_C\underline{H}\tens z)(1\tens m_2)m_2m_1
-(n\tens1_C\underline{H}\tens ks^{-1}\tens q)(1\tens m_3)m_2m_1 \\
-(n\tens1_C\underline{H}\tens ks^{-1}\tens q)(1\tens m_2\tens1)m_3m_1 \\
= nz
-(n\tens1_C\underline{H}\tens ks^{-1}\tens q)(1\tens m_2\tens1)
(m_2\tens1)m_2
= 1_Y -(p\tens q)m_2,
\end{multline*}
because $nz=1_Y$. Therefore, the cycles \(p\in\cd(Y,X)^0\) and
\(q\in\cd(X,Y)^0\) are inverse to each other modulo boundaries.
\end{proof}

\begin{corollary}\label{cor-DA-H0D}
The functor $H^0e$ factors as
 \(H^0\cc\rTTo^\qverdier H^0\cc/H^0\cb\rTTo^g H^0(\Quo(\cc|\cb))\),
where the Verdier quotient \(H^0\cc/H^0\cb=D(\ca)\) is the derived
category of $\ca$, and the functor \(g:D(\ca)\to H^0\cd\) is identity
on objects.
\end{corollary}

\subsection{Consequences of further research}
We shall use the results of the forthcoming book \cite{BesLyuMan-book}
to draw more conclusions for the example of complexes. The above
differential graded categories of complexes $\cc$, $\cb$ are
pretriangulated, see \cite{BondalKapranov:FramedTR}. Therefore,
\(\Quo(\cc|\cb)\) is a pretriangulated \ainf-category by results of
\cite[Chapters 15, 16]{BesLyuMan-book}.
The \ainf-equivalent to it (see
\secref{sec-Equivalence-two-quotients-Ainfty-categories}) differential
graded category \(\Dr(\cc|\cb)\) is pretriangulated as well by
[loc.~cit.].
The differential graded category \(\Dr(\cc|\cb)\) is precisely the
Drinfeld's quotient \(\cc/\cb\) introduced in
\cite[Section~3.1]{Drinf:DGquot}.

The isomorphism of \ainf-functors from
\secref{sec-Equivalence-two-quotients-Ainfty-categories} yields an
isomorphism of triangulated functors
\begin{equation}
H^0\ju \simeq \bigl[ H^0\cc \rTTo^{H^0e} H^0(\Quo(\cc|\cb))
\rTTo^{H^0f} H^0(\Dr(\cc|\cb)) \bigr]
\label{eq-H0j-HoeHof}
\end{equation}
by \cite[Chapter~18]{BesLyuMan-book}.
Notice that $H^0e$ and \(H^0\ju\) take objects of $H^0\cb$ to zero
objects. Hence, \eqref{eq-H0j-HoeHof} can be presented as the pasting
\begin{diagram}[LaTeXeqno]
\HmeetV &\rTTo &H^0(\Quo(\cc|\cb)) &\rLine &\HmeetV
\\
&= &\uTTo<g &\simeq &
\\
\uLine<{H^0e} &&H^0\cc/H^0\cb &&\dTTo<{\text{equiv}}>{H^0f}
\\
&\ruTTo^\qverdier &= &\rdTTo^\Psi &
\\
H^0\cc &&\rTTo^{H^0\ju} &&H^0(\Dr(\cc|\cb))
\label{dia-qverdier-pasting}
\end{diagram}
for a unique triangulated functor $\Psi$ by the properties of the
Verdier quotient/localization $\qverdier$
\cite[Section~2.2]{MR1453167}.

Denote by $\ce^\tr$ the pretriangulated envelope of an \ainf-category
$\ce$ \cite[Chapter~18]{BesLyuMan-book}. The \ainf-functor
\(u_\tr:\ce\to\ce^\tr\) is the natural embedding. The commutative
square
\begin{diagram}
\cc &\rTTo^\ju &\Dr(\cc|\cb)
\\
\dTTo<{u_\tr} &&\dTTo>{u_\tr}
\\
\cc^\tr &\rTTo^{\ju^\tr} &\Dr(\cc|\cb)^\tr
\end{diagram}
whose vertical arrows are \ainf-equivalences implies the commutative
diagram
\begin{diagram}[width=5.5em]
H^0\cc &\rTTo^\qverdier &H^0\cc/H^0\cb &\rTTo^\Psi &H^0(\Dr(\cc|\cb))
\\
\dTTo<{H^0(u_\tr)} &= &\dTTo<{H^0(u_\tr)} &= &\dTTo>{H^0(u_\tr)}
\\
H^0(\cc^\tr) &\rTTo^\qverdier &H^0(\cc^\tr)/H^0(\cb^\tr) &\rTTo^\Phi
&H^0(\Dr(\cc|\cb)^\tr)
\end{diagram}
whose rows compose to $H^0\ju$ and \(H^0(\ju^\tr)\), and columns are
equivalences.

When $\kk$ is a field, the functor $\Phi$ is an equivalence by
Theorem~3.4 of Drinfeld~\cite{Drinf:DGquot}. In this case $\Psi$ is an
equivalence, as well as $H^0f$ from
diagram~\eqref{dia-qverdier-pasting}. Hence, the triangulated functor
\(g:D(\ca)=H^0\cc/H^0\cb\to H^0(\Quo(\cc|\cb))\) from
\corref{cor-DA-H0D} is also an equivalence.

\appendix
\section{\texorpdfstring{The Yoneda Lemma for unital $A_\infty$-categories}
 {The Yoneda Lemma for unital A8-categories}}\label{ap-Yoneda-Lemma}

\subsection{Basic identities in symmetric closed monoidal category of
 complexes}
We want to work out in detail a system of notations suitable for
computations in symmetric closed monoidal categories. Actually we need
only the category of $\ZZ$\n-graded $\kk$\n-modules with a differential
of degree 1. The corresponding system of notations was already used in
\cite{Lyu-AinfCat,LyuOvs-iResAiFn}.

There exists a $\fu'$\n-small set $S$ of $\fu$\n-small $\kk$\n-modules
such that any $\fu$\n-small $\kk$\n-module $M$ is isomorphic to some
$\kk$\n-module $N\in S$ (due to presentations
 $\kk^{(P)}\to\kk^{(Q)}\to M\to0$). We turn $S$ into a category of
$\kk$\n-modules $\kk\modul$ with $\Ob\kk\modul=S$. Thus $\kk\modul$ is
an Abelian $\kk$\n-linear symmetric closed monoidal $\fu'$\n-small
$\fu$\n-category. The category $\Com={\mathsf C}(\kk\modul)$ of
complexes in $\kk\modul$ inherits all these properties from
$\kk\modul$, except that the symmetry becomes $c:X\tens Y\to Y\tens X$,
 $x\tens y\mapsto(-)^{xy}y\tens x=(-)^{\deg x\cdot\deg y}y\tens x$.
Therefore, we may consider the category of complexes enriched in $\Com$
(a differential graded category), and it is denoted by $\uCom$ in this
case. The (inner) hom-object between complexes $X$ and $Y$ is the
complex
\begin{equation*}
\uCom(X,Y)^n = \prod_{i\in\ZZ}\Hom_{\kk}(X^i,Y^{i+n}), \qquad
(f^i)_{i\in\ZZ}d = (f^id^{i+\deg f} -(-)^{\deg f}d^if^{i+1})_{i\in\ZZ}.
\end{equation*}

The product $X=\prod_{\iota\in I}X_\iota$ in the category of complexes
of $\kk$\n-modules of the family of objects $(X_\iota)_{\iota\in I}$
coincides with the product in the category of $\ZZ$\n-graded
$\kk$\n-modules (and differs from the product in the category of
$\kk$\n-modules). It is given by $X^m=\prod_{\iota\in I}X_\iota^m$.

Given a complex $Z$ and an element $a\in\uCom(X,Y)$, we assign to it
elements $1\tens a\in\uCom(Z\tens X,Z\tens Y)$,
$(z\tens x)(1\tens a)=z\tens xa$, and
$a\tens1\in\uCom(X\tens Z,Y\tens Z)$,
$(x\tens z)(a\tens1)=(-)^{za}xa\tens z$. Clearly,
$(1\tens a)c=c(a\tens1)\in\uCom(Z\tens X,Y\tens Z)$ and
$(a\tens1)c=c(1\tens a)\in\uCom(X\tens Z,Z\tens Y)$. If $g\in\uCom(Z,W)$,
then we have
$(1\tens a)(g\tens1)=(-)^{ag}(g\tens1)(1\tens a)\in\uCom(Z\tens X,W\tens Y)$
(Koszul's rule).

For any pair of complexes $X$, $Y\in\Ob\Com$ denote by
\begin{equation*}
\ev_{X,Y}:X\tens\uCom(X,Y)\to Y, \qquad
\coev_{X,Y}:Y\to\uCom(X,X\tens Y)
\end{equation*}
the canonical evaluation and coevaluation maps respectively. Then
the adjunction isomorphisms are explicitly given as follows:
\begin{align}
\Com(Y,\uCom(X,Z)) & \longleftrightarrow \Com(X\tens Y,Z), \notag \\
f &\rMapsTo (1_X\tens f)\ev_{X,Z}, \notag \\
\coev_{X,Y}\uCom(X,g) &\lMapsTo g.
\label{eq-g-to-coev-C1g}
\end{align}

Given a complex $Z$ and an element $a\in\uCom(X,Y)$, we assign to it the
element $\uCom(1,a)=\uCom(Z,a)=a\phi$ of $\uCom(\uCom(Z,X),\uCom(Z,Y))$
obtained from the equation
\[ m_2^\uCom = \bigl( \uCom(Z,X)\tens\uCom(X,Y)
\rTTo^{1\tens\phi}_{\exists!} \uCom(Z,X)\tens\uCom(\uCom(Z,X),\uCom(Z,Y))
\rTTo^\ev \uCom(Z,Y) \bigr),
\]
which holds for a unique chain map $\phi$. Despite that the map $a$ is
not a chain map we write this element as $a:X\to Y$, and we write
$a\phi$ as
\[ \uCom(1,a): \uCom(Z,X) \to \uCom(Z,Y), \qquad
(f^i)_{i\in\ZZ} \mapsto (f^ia^{i+\deg f})_{i\in\ZZ}.
\]
Similarly, given a complex $X$ and an element $g\in\uCom(W,Z)$, we
assign to it the element
$\uCom(g,1)=\uCom(g,X)=g\psi\in\uCom(\uCom(Z,X),\uCom(W,X))$ obtained from
the diagram
\begin{diagram}
\uCom(Z,X)\tens\uCom(W,Z) & \rTTo^c & \uCom(W,Z)\tens\uCom(Z,X) \\
\dTTo<{\exists!}>{1\tens\psi} && \dTTo>{m_2^\uCom} \\
\uCom(Z,X)\tens\uCom(\uCom(Z,X),\uCom(W,X)) & \rTTo^\ev & \uCom(W,X)
\end{diagram}
commutative for a unique chain map $\psi$. Although the map $\uCom(g,1)$
is not a chain map we write it as
\[ \uCom(g,1): \uCom(Z,X) \to \uCom(W,X), \qquad
(f^i)_{i\in\ZZ} \mapsto
\bigl((-)^{\deg f\cdot\deg g}g^if^{i+\deg g}\bigr)_{i\in\ZZ}.
\]

For each pair of homogeneous elements $a\in\uCom(X,Y)$, $g\in\uCom(W,Z)$
we have
\begin{multline*}
\bigl( \uCom(Z,X) \rTTo^{\uCom(Z,a)} \uCom(Z,Y) \rTTo^{\uCom(g,Y)}
\uCom(W,Y) \bigr)
\\
= (-)^{ag} \bigl( \uCom(Z,X) \rTTo^{\uCom(g,X)} \uCom(W,X)
\rTTo^{\uCom(W,a)} \uCom(W,Y) \bigr).
\end{multline*}
This equation follows from one of the standard identities in symmetric
closed monoidal categories \cite{EK}, and can be verified directly. We
also have $\uCom(1,a)\uCom(1,h)=\uCom(1,ah)$ and
$\uCom(g,1)\uCom(e,1)=(-)^{ge}\uCom(eg,1)$, whenever these maps are
defined.

One easily sees that $m^{\uCom}_1=d:\uCom(X,Y)\to\uCom(X,Y)$ coincides
with $\uCom(1,d_Y)-\uCom(d_X,1)$.

Let $f:A\tens X\to B$, $g:B\tens Y\to C$ be two homogeneous
$\kk$\n-linear maps of arbitrary degrees. Then the following holds:
\begin{multline}\label{eq-identity-m2}
\bigl( X\tens Y \rTTo^{\coev_{A,X}\tens\coev_{B,Y}}
\uCom(A,A\tens X)\tens\uCom(B,B\tens Y)
\\
\hfill \rTTo^{\uCom(A,f)\tens\uCom(B,g)} \uCom(A,B)\tens\uCom(B,C)
\rTTo^{m_2} \uCom(A,C) \bigr) \quad
\\
= \bigl( X\tens Y \rTTo^{\coev_{A,X\tens Y}} \uCom(A,A\tens X\tens Y)
\rTTo^{\uCom(A,f\tens1)} \uCom(A,B\tens Y) \rTTo^{\uCom(A,g)}
\uCom(A,C) \bigr).
\end{multline}
Indeed,
$(\coev_{A,X}\tens\coev_{B,Y})(\uCom(A,f)\tens\uCom(B,g))
=(\coev_{A,X}\uCom(A,f)\tens\coev_{B,Y}\uCom(B,g))$,
for $\coev$ has degree $0$. Denote $\bar{f}=\coev_{A,X}\uCom(A,f)$,
$\bar{\vphantom{f}g}=\coev_{B,Y}\uCom(B,g)$. The morphisms $\bar f$ and
$\bar{\vphantom{f}g}$ correspond to $f$ and $g$ by adjunction. Further,
the morphism $m_2$ comes by adjunction from the following map:
\[ A\tens\uCom(A,B)\tens\uCom(B,C) \rTTo^{\ev_{A,B}\tens1}
B\tens\uCom(B,C) \rTTo^{\ev_{B,C}} C,
\]
in particular the following diagram commutes:
\begin{diagram}
A\tens\uCom(A,B)\tens\uCom(B,C) & \rTTo^{\ev_{A,B}\tens 1} &
B\tens\uCom(B,C) \\
\dTTo<{1\tens m_2} && \dTTo>{\ev_{B,C}} \\
A\tens\uCom(A,C) & \rTTo^{\ev_{A,C}} & C
\end{diagram}
Thus we have a commutative diagram
\begin{diagram}
A\tens X\tens Y & \rTTo^{1\tens\bar f\tens 1} &
A\tens\uCom(A,B)\tens Y & \rTTo^{\ev_{A,B}\tens 1} & B\tens Y\\
&& \dTTo<{1\tens 1\tens\bar{\vphantom{f}g}} &&
\dTTo>{1\tens\bar{\vphantom{f}g}} \\
&& A\tens\uCom(A,B)\tens\uCom(B,C) & \rTTo^{\ev_{A,B}\tens 1} &
B\tens\uCom(B,C) \\
&& \dTTo<{1\tens m_2} && \dTTo>{\ev_{B,C}} \\
&& A\tens\uCom(A,C) & \rTTo^{\ev_{A,C}} & C
\end{diagram}
The top row composite coincides with $f\tens1$ and the right-hand side
vertical composite coincides with $g$ (by adjunction). Thus
$(1\tens\bar f\tens 1)(1\tens 1\tens\bar{\vphantom{f} g})
(1\tens m_2)\ev_{A,C}$
coincides with $(f\tens1)g$. But the mentioned morphism comes from
$(\bar f\tens\bar{\vphantom{f}g})m_2$ by adjunction, so that
$(\bar f\tens\bar{\vphantom{f}g})m_2=
\coev_{A,X\tens Y}\uCom(A,(f\tens1)g)$
(the latter morphism is the image of $(f\tens1)g$ under the
adjunction), and we are done.

One verifies similarly the following assertion \cite{EK}: given
$f\in\uCom(A,B)$, then the diagram
\begin{diagram}[LaTeXeqno]
X & \rTTo^{\coev_{A,X}} & \uCom(A,A\tens X) \\
\dTTo<{\coev_{B,X}} && \dTTo>{\uCom(A,f\tens 1)} \\
\uCom(B,B\tens X) & \rTTo^{\uCom(f,B\tens X)} & \uCom(A,B\tens X)
\label{dia-coev-coev-ff}
\end{diagram}
commutes.

If $P$ is a $\ZZ$\n-graded $\kk$\n-module, then $sP=P[1]$ denotes the
same $\kk$\n-module with the grading $(sP)^d=P^{d+1}$. The ``identity''
map $P\to sP$ of degree $-1$ is also denoted $s$. The map $s$ commutes
with the components of the differential $b$ in an \ainf-category
(\ainf-algebra) in the following sense: $s^{\tens n}b_n=m_ns$. The main
identity $b^2=0$ written in components takes the form
\begin{equation}
\sum_{r+n+t=k} (1^{\tens r}\tens b_n\tens1^{\tens t})b_{r+1+t} = 0 :
T^ks\ca \to s\ca.
\label{eq-b-b-0}
\end{equation}

\subsection{\texorpdfstring{$A_\infty$-functor $h^X$}{A8-functor hX}}
Let $\ca$ be an \ainf-category. Following Fukaya
\cite[Definition~7.28]{Fukaya:FloerMirror-II} for any object $X$ of
$\ca$ we define a cocategory homomorphism $h^X:Ts\ca\to Ts\uCom$ as
follows. It maps an object $Z$ to the complex $h^XZ=(s\ca(X,Z),-b_1)$.
The minus sign is explained by the fact that
\(H^XZ\overset{\text{def}}=(h^XZ)[-1]=(\ca(X,Z),m_1)\).
Actually, $h^XZ$ is some fixed complex of
$\kk$\n-modules from $\Ob\Com$, with a fixed isomorphism
$h^XZ\rTTo^\sim (s\ca(X,Z),-b_1)$. These isomorphisms enter implicitly
into the structure maps of $h^X$, however, we shall pretend that they
are identity morphisms. The closed monoidal structure of $\Com$ gives
us the right to omit these isomorphisms in all the formulae.

We require $h^X$ to be pointed, that is,
 $(T^0s\ca)h^X\subset T^0s\uCom$. Therefore, $h^X$ is completely
specified by its components $h^X_k$ for $k\ge1$:
\begin{multline}
h^X_k = \left[s\ca(Z_0,Z_1)\tdt s\ca(Z_{k-1},Z_k) \rTTo^\coev
\uCom(h^XZ_0,h^XZ_0\tens h^{Z_0}Z_1\tdt h^{Z_{k-1}}Z_k)\right.
\\
\left. \rTTo^{\uCom(1,b_{1+k})} \uCom(s\ca(X,Z_0),s\ca(X,Z_k)) \rTTo^s
s\uCom(h^XZ_0,h^XZ_k)\right].
 \label{eq-hXk-components}
\end{multline}
The composition \(H^X=h^X\cdot[-1]\) is described in
\cite[eq.~(A.1)]{LyuMan-AmodSerre}. It is proven in this work that
$H^X:\ca\to\uCom$ is an \ainf-functor. Therefore,
$h^X=H^X\cdot[1]:\ca\to\uCom$ is an \ainf-functor as well. A similar
statement is known from Fukaya's work
\cite[Proposition~7.18]{Fukaya:FloerMirror-II} under slightly more
restrictive general assumptions.

\ifx\chooseClass1
\begin{proof}
Nevertheless, for the sake of completeness, we give a full proof in our
own notations. We have to verify the following equation for an
arbitrary $k\ge 1$:
\begin{multline*}
\sum_{p+t+q=k}(1^{\tens p}\tens b_t\tens 1^{\tens q})h^X_{p+1+q}
=h^X_kb^{\uCom}_1 + \sum_{i+j=k}(h^X_i\tens h^X_j)b^{\uCom}_2: \\
s\ca(Z_0,Z_1)\tdt s\ca(Z_{k-1},Z_k) \to s\uCom(h^XZ_0,h^XZ_k).
\end{multline*}
Let us gradually expand each summand. First of all,
\begin{align*}
h^X_kb^{\uCom}_1 &= \coev\uCom(1,b_{1+k})sb^{\uCom}_1
\\
&= \coev\uCom(1,b_{1+k})m^{\uCom}_1s
\\
&= \coev\uCom(1,b_{1+k})\bigl(-\uCom(1,b_1)+\uCom(b_1,1)\bigr)s
\\
&= -\coev\uCom(1,b_{1+k}b_1)s -\coev\uCom(b_1,1)\uCom(1,b_{1+k})s
\\
&= -\coev\uCom(1,b_{1+k}b_1)s-\coev\uCom(1,(b_1\tens1^{\tens k})b_{1+k})s
\end{align*}
by \eqref{dia-coev-coev-ff}. Similarly,
\begin{align}
(h^X_i\tens h^X_j)b^{\uCom}_2 &=
(\coev\uCom(1,b_{1+i})s\tens\coev\uCom(1,b_{1+j})s)b_2^\uCom \notag
\\
&= -(\coev\tens\coev)(\uCom(1,b_{1+i})\tens\uCom(1,b_{1+j}))m_2s \notag
\\
&= -\coev\uCom(1,(b_{1+i}\tens1^{\tens j})b_{1+j})s.
\label{eq-compare-hhb}
\end{align}
The sign is due to the fact $s$ and $b$ have odd degree. The last line
is due to~\eqref{eq-identity-m2}.

Further, given $\alpha:Y\to U$, then the following diagram commutes
\begin{diagram}
X\tens Y\tens Z & \rTTo^{1\tens\alpha\tens 1} & X\tens U\tens Z \\
\dTTo<{\coev_{A,X\tens Y\tens Z}} && \dTTo>{\coev_{A,X\tens U\tens Z}}
\\
\uCom(A,A\tens X\tens Y\tens Z)
& \rTTo^{\uCom(A,1\tens 1\tens\alpha\tens 1)}
& \uCom(A,A\tens X\tens U\tens Z),
\end{diagram}
for $\coev$ is a functor morphism. Applying this to $\alpha=b_t$ we
obtain
\begin{align}
(1^{\tens p}\tens b_t\tens 1^{\tens q})h^X_{p+1+q} &=
 (1^{\tens p}\tens b_t\tens 1^{\tens  q})\coev\uCom(1,b_{p+2+q})s \notag
\\
&= \coev\uCom(1,(1^{\tens1+p}\tens b_t\tens 1^{\tens q})b_{p+2+q})s.
\label{eq-compare2-bh}
\end{align}
Thus,
\[ \sum_{p+t+q=k}(1^{\tens p}\tens b_t\tens 1^{\tens q})h^X_{p+1+q}
=\coev\uCom\bigl(1,\sum_{\substack{p+t+q=k+1,\\ p\ge 1}}
 (1^{\tens p}\tens b_t\tens 1^{\tens q})b_{p+1+q}\bigr)s,
\]
so that
\begin{multline*}
\sum_{p+t+q=k}(1^{\tens p}\tens b_t\tens 1^{\tens q})h^X_{p+1+q}
-h^X_kb^{\uCom}_1
-\sum_{\substack{i+j=k,\\ i,j\ge1}}(h^X_i\tens h^X_j)b^{\uCom}_2
\\
=\coev\uCom\bigl(1,\sum_{\substack{p+t+q=k+1,\\ p\ge 1}}
(1^{\tens p}\tens b_t\tens 1^{\tens q})b_{p+1+q} +b_{1+k}b_1
+(b_1\tens1^{\tens k})b_{1+k}
+\sum_{\substack{i+j=k,\\ i,j\ge 1}}(b_{1+i}\tens 1^{\tens j})
b_{1+j}\bigr)s.
\end{multline*}
The sum in the right-hand side of the above formula vanishes by
equation~\eqref{eq-b-b-0}.
\end{proof}
\fi

Let $\ca$ be a unital \ainf-category. Then for each object $X$ of $\ca$
the \ainf-functor $H^X:\ca\to\uCom$ is unital by
\cite[Remark~5.19]{LyuMan-AmodSerre}. Hence,
the \ainf-functor $h^X:\ca\to\uCom$ is unital as well.

\ifx\chooseClass1
\begin{proof}
For arbitrary objects $Y$, $Z$ of $\ca$ the map $h^X_1$ is given by
\[ h^X_1 = \bigl[s\ca(Y,Z) \rTTo^{\coev} \uCom(h^XY,h^XY\tens h^YZ)
\rTTo^{\uCom(1,b_2)} \uCom(h^XY,h^XZ) \rTTo^s s\uCom(h^XY,h^XZ)\bigr].
\]
We need this formula for $Z=Y$. We have to find such
$v\in(s\uCom)^{-2}(h^XY,h^XY)$ that
\[ \sS{_Y}\uni^\ca_0h^X_1 - \sS{_{h^XY}}\uni^\uCom_0 = vb^\uCom_1:
\kk \to (s\uCom)^{-1}(h^XY,h^XY).
\]
We have $\sS{_{h^XY}}\uni^\uCom_0=1_{h^XY}s$. We seek for an element $v$
in the form $v=h's$, where $h'\in\uCom^{-1}(h^XY,h^XY)$. The required
equation is:
\begin{multline*}
\bigl[\kk \rTTo^{\sS{_Y}\uni^\ca_0} (s\ca)^{-1}(Y,Y) \rTTo^{\coev}
\uCom^{-1}(h^XY,h^XY\tens h^YY) \rTTo^{\uCom(1,b_2)}
\uCom^0(h^XY,h^XY)\bigr] \\
= \bigl[\kk \rTTo^{1_{h^XY}+h'm^\uCom_1} \uCom^0(h^XY,h^XY)\bigr].
\end{multline*}
By functoriality of $\coev$ we have
\(\sS{_Y}\uni^\ca_0\coev_{h^XY,h^YY}
=\coev_{h^XY,\kk}\uCom(h^XY,1\tens\sS{_Y}\uni^\ca_0)\),
and the above equation reduces to
\[ \bigl(h^XY = h^XY\tens\kk \rTTo^{1\tens\sS{_Y}\uni^\ca_0}
h^XY\tens h^YY \rTTo^{b_2} h^XY\bigr) =
\bigl(h^XY \rTTo^{1-h'b_1-b_1h'} h^XY\bigr).
\]
The existence of such a homotopy $h'$ immediately follows from the
definition of a unital \ainf-category, see
\cite[Lemma~7.4]{Lyu-AinfCat}.
\end{proof}
\fi

\subsection{\texorpdfstring{The opposite $A_\infty$-category}
 {The opposite A8-category}}
Let $\ca$ be a graded $\kk$\n-quiver. Then its \emph{opposite quiver}
$\ca^\op$ is defined as the quiver with the same class of objects
$\Ob\ca^\op=\Ob\ca$, and with graded $\kk$\n-modules of morphisms
$\ca^\op(X,Y)=\ca(Y,X)$.

Let $\gamma:Ts\ca^\op\to Ts\ca$ denote the following cocategory
anti-isomorphism:
\begin{equation}
\gamma= (-1)^k\omega^0_c: s\ca^\op(X_0,X_1)\tdt s\ca^\op(X_{k-1},X_k)
\to s\ca(X_k,X_{k-1})\tdt s\ca(X_1,X_0),
\label{eq-gamma-anti}
\end{equation}
where the permutation
$\omega^0=\bigl(
\begin{smallmatrix}
1 & 2 & \dots & k-1 & k \\
k &k-1& \dots &  2  & 1
\end{smallmatrix}
\bigr)$
is the longest element of $\SSS_k$, and $\omega^0_c$ is the
corresponding signed permutation, the action of $\omega^0$ in tensor
products via standard symmetry. Clearly,
$\gamma\Delta=\Delta(\gamma\tens\gamma)c=\Delta c(\gamma\tens\gamma)$,
which is the anti-isomorphism property. Notice also that
$(\ca^\op)^\op=\ca$ and $\gamma^2=\id$.

When $\ca$ is an \ainf-category with the codifferential
$b:Ts\ca\to Ts\ca$, then $\gamma b\gamma:Ts\ca^\op\to Ts\ca^\op$ is
also a codifferential. Indeed,
\begin{align*}
\gamma b\gamma\Delta &= \gamma b\Delta c(\gamma\tens\gamma)
= \gamma\Delta(1\tens b+b\tens1)c(\gamma\tens\gamma)
= \Delta(\gamma\tens\gamma)c(1\tens b+b\tens1)c(\gamma\tens\gamma) \\
&= \Delta(\gamma\tens\gamma)(b\tens1+1\tens b)(\gamma\tens\gamma)
= \Delta(\gamma b\gamma\tens1+1\tens\gamma b\gamma).
\end{align*}

\begin{definition}[cf. Fukaya \cite{Fukaya:FloerMirror-II}
Definition~7.8]\label{def-opposite-ainf-category}
The \emph{opposite \ainf-category} $\ca^\op$ to an \ainf-category $\ca$
is the opposite quiver, equipped with the codifferential
$b^\op=\gamma b\gamma:Ts\ca^\op\to Ts\ca^\op$.
\end{definition}

The components of $b^\op$ are computed as follows:
\begin{multline*}
b^\op_k = (-)^{k+1}\bigl[s\ca^\op(X_0,X_1)\tdt s\ca^\op(X_{k-1},X_k) \\
\rTTo^{\omega^0_c} s\ca(X_k,X_{k-1})\tdt s\ca(X_1,X_0) \rTTo^{b_k}
s\ca(X_k,X_0) = s\ca^\op(X_0,X_k) \bigr].
\end{multline*}
The sign $(-1)^k$ in \eqref{eq-gamma-anti} ensures that the above
definition agrees with the definition of the opposite usual category.

\subsection{\texorpdfstring{The Yoneda $A_\infty$-functor}
 {The Yoneda A8-functor}}
Since the considered \ainf-category $\ca$ is $\fu$\n-small, and $\Com$
is a $\fu'$\n-small $\fu$\n-category, $A_\infty(\ca,\uCom)$ is a
$\fu'$\n-small differential graded $\fu$\n-category. Indeed, every its
set of morphisms $sA_\infty(\ca,\uCom)(f,g)$ is isomorphic to the
product of graded $\kk$\n-modules
\[ \prod_{k=0}^\infty \prod_{X,Y\in\Ob\ca}
\uCom\bigl(T^ks\ca(X,Y),s\uCom(Xf,Yg)\bigr),
\]
that is $\fu$\n-small.

The Yoneda \ainf-functor exists in two versions: $Y$ and
\(\Yo:\ca^\op\to A_{\infty}(\ca,\uCom)\) which differ by a shift:
\(Y=\Yo\cdot A_\infty(1,[1])\). The pointed cocategory homomorphism
$Y:Ts\ca^{\op}\to TsA_{\infty}(\ca,\uCom)$ is given as follows: on
objects $X\mapsto h^X$, the components
\begin{equation}
Y_n: s\ca^{\op}(X_0,X_1)\tdt s\ca^{\op}(X_{n-1},X_n) \to
sA_{\infty}(\ca,\uCom)(h^{X_0},h^{X_n})
\label{eq-components-Yn}
\end{equation}
are determined by the following formulas. The composition of $Y_n$ with
\[ \pr_k: sA_{\infty}(\ca,\uCom)(h^{X_0},h^{X_n}) \to
\uCom(h^{Z_0}Z_1\tdt h^{Z_{k-1}}Z_k,s\uCom(h^{X_0}Z_0,h^{X_n}Z_k))
\]
(that is, the $k$\n-th component of the coderivation
 $(x_1\tdt x_n)Y_n$) is given by the formula
\begin{multline*}
Y_{nk} = (-)^n \bigl[s\ca(X_1,X_0)\tdt s\ca(X_n,X_{n-1}) \rTTo^\coev \\
\uCom(h^{X_0}Z_0\tens h^{Z_0}Z_1\tdt h^{Z_{k-1}}Z_k,h^{X_0}Z_0\tens
h^{Z_0}Z_1\tdt h^{Z_{k-1}}Z_k\tens h^{X_1}X_0\tdt h^{X_n}X_{n-1}) \\
\rTTo^{\uCom(1,\tau_cb_{k+n+1})}
\uCom(h^{X_0}Z_0\tens h^{Z_0}Z_1\tdt h^{Z_{k-1}}Z_k,h^{X_n}Z_k) \\
\rTTo^\sim
\uCom(h^{Z_0}Z_1\tdt h^{Z_{k-1}}Z_k,\uCom(h^{X_0}Z_0,h^{X_n}Z_k)) \\
\rTTo^{\uCom(1,s)}
\uCom(h^{Z_0}Z_1\tdt h^{Z_{k-1}}Z_k,s\uCom(h^{X_0}Z_0,h^{X_n}Z_k))
\bigr],
\end{multline*}
where the permutation $\tau\in\SSS_{k+n+1}$ is
 $\tau^{k,n}=\bigl(
\begin{smallmatrix}
0 & 1 & \dots & k &k+1& \dots &k+n \\
n &1+n& \dots &k+n&n-1& \dots & 0
\end{smallmatrix}
 \bigr)$.

In other words, the coderivation $(x_1\tdt x_n)Y_n$ has components
$s\ca(Z_0,Z_1)\tdt s\ca(Z_{k-1},Z_k)\to s\uCom(h^{X_0}Z_0,h^{X_n}Z_k)$,
$(z_1\tdt z_k)\mapsto (z_1\tdt z_k\tens x_1\tdt x_n)Y'_{nk}$, where
\begin{multline}
Y'_{nk} = (-)^n\bigl[ s\ca(Z_0,Z_1)\tdt s\ca(Z_{k-1},Z_k)\tens
s\ca(X_1,X_0)\tdt s\ca(X_n,X_{n-1}) \\
\rTTo^{\coev} \uCom(h^{X_0}Z_0,h^{X_0}Z_0\tens
h^{Z_0}Z_1\tdt h^{Z_{k-1}}Z_k\tens h^{X_1}X_0\tdt h^{X_n}X_{n-1}) \\
\rTTo^{\uCom(1,\tau_cb_{k+n+1})} \uCom(h^{X_0}Z_0,h^{X_n}Z_k)
\rTTo^s s\uCom(h^{X_0}Z_0,h^{X_n}Z_k) \bigr].
\label{eq-Y-nk-prime}
\end{multline}

The pointed cocategory homomorphism
 $Y:\ca^{\op}\to A_{\infty}(\ca,\uCom)$ is an \ainf-functor. An
equivalent statement is already proved by Fukaya
\cite[Lemma~9.8]{Fukaya:FloerMirror-II} and by the authors in
\cite[Section~5.5]{LyuMan-AmodSerre}.

\ifx\chooseClass1
\begin{proof}
The equation to prove is the following:
\begin{multline*}
\sum_{l>0;i_1+\dots+i_l=n} (Y_{i_1}\tens Y_{i_2} \tdt Y_{i_l}) B_l =
\sum_{\alpha+p+\beta=n}
(1^{\tens\alpha}\tens b^\op_p\tens1^{\tens\beta})Y_{\alpha+1+\beta}: \\
s\ca^{\op}(X_0,X_1)\tdt s\ca^{\op}(X_{n-1},X_n)
\to sA_{\infty}(\ca,\uCom)(h^{X_0},h^{X_n}).
\end{multline*}
We compare the $k$\n-th components of the both sides:
\begin{multline*}
[Y_nB_1]_k + \sum_{p+q=n} [(Y_p\tens Y_q)B_2]_k =
\sum_{\alpha+p+\beta=n}
[(1^{\tens\alpha}\tens b^\op_p\tens1^{\tens\beta})
Y_{\alpha+1+\beta}]_k: \\
s\ca(X_1,X_0)\tdt s\ca(X_n,X_{n-1})
\to sA_{\infty}(\ca,\uCom)(h^{X_0},h^{X_n}) \\
\rTTo^{\pr_k}
\uCom(h^{Z_0}Z_1\tdt h^{Z_{k-1}}Z_k,s\uCom(h^{X_0}Z_0,h^{X_n}Z_k)).
\end{multline*}

Let us denote by $\bar r$ an element of
$s\ca(X_1,X_0)\tdt s\ca(X_n,X_{n-1})$. The terms in the left hand side
are expanded to
\begin{multline*}
[(\bar rY_n)B_1]_k = (\bar rY_n)_kb_1^\uCom \\
+ \sum_{i+j=k}^{i>0}(h^{X_0}_i\tens(\bar rY_n)_j)b_2^\uCom
+ \sum_{i+j=k}^{j>0}((\bar rY_n)_i\tens h^{X_n}_j)b_2^\uCom
-(-)^r \sum_{\alpha+p+\beta=k}
(1^{\tens\alpha}\tens b_p\tens1^{\tens\beta})
(\bar rY_n)_{\alpha+1+\beta},
\end{multline*}
\[ [(Y_p\tens Y_q)B_2]_k =
\sum_{i+j=k}^{p,q>0} (Y'_{pi}\tens Y'_{qj})b_2^\uCom.
\]
To reduce further the number of terms we introduce the notation
$Y'_{0k}=h^X_k$. Indeed, when $n=0$, formulas \eqref{eq-Y-nk-prime} and
\eqref{eq-hXk-components} coincide. So the left hand side expands to
the following expression:
\begin{multline*}
[(\bar rY_n)B_1]_k + \sum_{p+q=n} [(Y_p\tens Y_q)B_2]_k \\
= (\bar rY_n)_kb_1^\uCom
+ \sum_{i+j=k}^{p,q\ge0} (Y'_{pi}\tens Y'_{qj})b_2^\uCom
-(-)^r \sum_{\alpha+p+\beta=k}
(1^{\tens\alpha}\tens b_p\tens1^{\tens\beta})
(\bar rY_n)_{\alpha+1+\beta},
\end{multline*}
Thus the equation to prove is this one:
\begin{multline*}
\bigl[ h^{Z_0}Z_1\tdt h^{Z_{k-1}}Z_k\tens h^{X_1}X_0\tdt h^{X_n}X_{n-1}
\rTTo^{Y'_{nk}} s\uCom(h^{X_0}Z_0,h^{X_n}Z_k) \\
 \rTTo^{b_1^\uCom} s\uCom(h^{X_0}Z_0,h^{X_n}Z_k) \bigr]
+ \sum_{\substack{p+q=n\\i+j=k}}^{\substack{p,q,i,j\ge0\\p+i>0,q+j>0}}
\bigl[ h^{Z_0}Z_1\tdt h^{Z_{k-1}}Z_k\tens h^{X_1}X_0\tdt h^{X_n}X_{n-1}
\rTTo^{1\tens c\tens1} \\
 h^{Z_0}Z_1\tdt h^{Z_{i-1}}Z_i\tens h^{X_1}X_0\tdt h^{X_p}X_{p-1}\tens
h^{Z_i}Z_{i+1}\tdt h^{Z_{k-1}}Z_k\tens h^{X_{p+1}}X_p\tdt h^{X_n}X_{n-1} \\
 \rTTo^{Y'_{pi}\tens Y'_{qj}}
s\uCom(h^{X_0}Z_0,h^{X_p}Z_i)\tens s\uCom(h^{X_p}Z_i,h^{X_n}Z_k)
\rTTo^{b_2^\uCom} s\uCom(h^{X_0}Z_0,h^{X_n}Z_k) \bigr] \\
 - \sum_{\alpha+p+\beta=k}
\bigl[ h^{Z_0}Z_1\tdt h^{Z_{k-1}}Z_k\tens h^{X_1}X_0\tdt h^{X_n}X_{n-1}
 \rTTo^{1^{\tens\alpha}\tens b_p\tens1^{\tens\beta+n}} \\
h^{Z_0}Z_1\tdt h^{Z_{\alpha-1}}Z_\alpha\tens
h^{Z_\alpha}Z_{\alpha+p}\tens h^{Z_{\alpha+p}}Z_{\alpha+p+1}\tdt
h^{Z_{k-1}}Z_k\tens h^{X_1}X_0\tdt h^{X_n}X_{n-1}
\\
\hfill \rTTo^{Y'_{n,\alpha+1+\beta}} s\uCom(h^{X_0}Z_0,h^{X_n}Z_k)
\bigr] \quad
\\
\quad = \sum_{\alpha+p+\beta=n} \bigl[
h^{Z_0}Z_1\tdt h^{Z_{k-1}}Z_k\tens h^{X_1}X_0\tdt h^{X_n}X_{n-1}
\rTTo^{1^{\tens k+\alpha}\tens b^\op_p\tens1^{\tens\beta}} \hfill
\\
 h^{Z_0}Z_1\tdt h^{Z_{k-1}}Z_k\tens h^{X_1}X_0\tdt
h^{X_\alpha}X_{\alpha-1}\tens h^{X_{\alpha+p}}X_\alpha\tens
h^{X_{\alpha+p+1}}X_{\alpha+p}\tdt h^{X_n}X_{n-1} \\
 \rTTo^{Y'_{\alpha+1+\beta,k}} s\uCom(h^{X_0}Z_0,h^{X_n}Z_k) \bigr].
\end{multline*}
Substituting explicit expression for $Y'_{nk}$ we get the following
equation to verify:
\begin{multline*}
 (-)^n \bigl[
h^{Z_0}Z_1\tdt h^{Z_{k-1}}Z_k\tens h^{X_1}X_0\tdt h^{X_n}X_{n-1} \\
\rTTo^\coev \uCom(h^{X_0}Z_0,h^{X_0}Z_0\tens
h^{Z_0}Z_1\tdt h^{Z_{k-1}}Z_k\tens h^{X_1}X_0\tdt h^{X_n}X_{n-1})
\\
\rTTo^{\uCom(1,\tau^{n,k}_cb_{k+n+1})} \uCom(h^{X_0}Z_0,h^{X_n}Z_k)
\rTTo^{-\uCom(1,b_1)+\uCom(b_1,1)} \uCom(h^{X_0}Z_0,h^{X_n}Z_k)
\xrightarrow{s} s\uCom(h^{X_0}Z_0,h^{X_n}Z_k) \bigr]
\\
 - \sum_{\substack{p+q=n\\i+j=k}}^{\substack{p,q,i,j\ge0\\p+i>0,q+j>0}}
(-)^{p+q} \bigl[
h^{Z_0}Z_1\tdt h^{Z_{k-1}}Z_k\tens h^{X_1}X_0\tdt h^{X_n}X_{n-1}
\rTTo^{1^{\tens i}\tens c\tens1^{\tens n-p}} \\
 h^{Z_0}Z_1\tdt h^{Z_{i-1}}Z_i\tens h^{X_1}X_0\tdt h^{X_p}X_{p-1}\tens
h^{Z_i}Z_{i+1}\tdt h^{Z_{k-1}}Z_k\tens h^{X_{p+1}}X_p\tdt h^{X_n}X_{n-1}\\
\rTTo^{\coev\tens\coev}
\uCom(h^{X_0}Z_0,h^{X_0}Z_0\tens h^{Z_0}Z_1\tdt
h^{Z_{i-1}}Z_i\tens h^{X_1}X_0\tdt h^{X_p}X_{p-1})\\
\tens \uCom(h^{X_p}Z_i,h^{X_p}Z_i\tens h^{Z_i}Z_{i+1}\tdt
 h^{Z_{k-1}}Z_k\tens h^{X_{p+1}}X_p\tdt h^{X_n}X_{n-1}) \\
\rTTo^{\uCom(1,\tau^{i,p}_cb_{1+i+p})\tens\uCom(1,\tau^{j,q}_cb_{1+j+q})}
\uCom(h^{X_0}Z_0,h^{X_p}Z_i)\tens\uCom(h^{X_p}Z_i,h^{X_n}Z_k)
\rTTo^{m^\uCom_2} \uCom(h^{X_0}Z_0,h^{X_n}Z_k)
\\
\hfill \rTTo^s s\uCom(h^{X_0}Z_0,h^{X_n}Z_k) \bigr] \quad
\\
\quad -(-)^n \sum_{\alpha+p+\beta=k}
\bigl[ h^{Z_0}Z_1\tdt h^{Z_{k-1}}Z_k\tens h^{X_1}X_0\tdt h^{X_n}X_{n-1}
 \rTTo^{1^{\tens\alpha}\tens b_p\tens1^{\tens\beta+n}} \hfill
\\
h^{Z_0}Z_1\tdt h^{Z_{\alpha-1}}Z_\alpha\tens
h^{Z_\alpha}Z_{\alpha+p}\tens h^{Z_{\alpha+p}}Z_{\alpha+p+1}\tdt
h^{Z_{k-1}}Z_k\tens h^{X_1}X_0\tdt h^{X_n}X_{n-1} \\
\rTTo^\coev \uCom(h^{X_0}Z_0,h^{X_0}Z_0\tens h^{Z_0}Z_1\tdt
h^{Z_{\alpha-1}}Z_\alpha\tens h^{Z_\alpha}Z_{\alpha+p}
 \tens h^{Z_{\alpha+p}}Z_{\alpha+p+1}\tens\dots\tens h^{Z_{k-1}}Z_k\\
\tens h^{X_1}X_0\tdt h^{X_n}X_{n-1})
\rTTo^{\uCom(1,\tau^{\alpha+1+\beta,n}_cb_{n+\alpha+\beta+2})}
\uCom(h^{X_0}Z_0,h^{X_n}Z_k) \rTTo^s s\uCom(h^{X_0}Z_0,h^{X_n}Z_k)
\bigr] \\
 = \sum_{\alpha+p+\beta=n} (-1)^{\alpha+\beta+p}\bigl[
h^{Z_0}Z_1\tdt h^{Z_{k-1}}Z_k\tens h^{X_1}X_0\tdt h^{X_n}X_{n-1}
\rTTo^{1^{\tens k+\alpha}\tens(\omega^0_cb_p)\tens1^{\tens\beta}} \\
 h^{Z_0}Z_1\tdt h^{Z_{k-1}}Z_k\tens h^{X_1}X_0\tdt
h^{X_\alpha}X_{\alpha-1}\tens h^{X_{\alpha+p}}X_\alpha\tens
h^{X_{\alpha+p+1}}X_{\alpha+p}\tdt h^{X_n}X_{n-1} \\
 \rTTo^\coev \uCom(h^{X_0}Z_0,h^{X_0}Z_0\tens h^{Z_0}Z_1\tdt
h^{Z_{k-1}}Z_k\tens h^{X_1}X_0\tdt h^{X_\alpha}X_{\alpha-1}\tens
h^{X_{\alpha+p}}X_\alpha\tens \\
  h^{X_{\alpha+p+1}}X_{\alpha+p}\tdt h^{X_n}X_{n-1})
\rTTo^{\uCom(1,\tau^{k,\alpha+1+\beta}_cb_{k+\alpha+\beta+2})}
\uCom(h^{X_0}Z_0,h^{X_n}Z_k) \\
\rTTo^s s\uCom(h^{X_0}Z_0,h^{X_n}Z_k) \bigr].
\end{multline*}
Using functoriality of $\coev$, identities \eqref{eq-identity-m2} and
\eqref{dia-coev-coev-ff} we reduce all terms of this equation to the
form
\begin{multline*}
\bigl[ h^{Z_0}Z_1\tdt h^{Z_{k-1}}Z_k\tens h^{X_1}X_0\tdt h^{X_n}X_{n-1}
\\
 \rTTo^\coev \uCom(h^{X_0}Z_0,h^{X_0}Z_0\tens h^{Z_0}Z_1\tdt
h^{Z_{k-1}}Z_k\tens h^{X_1}X_0\tdt h^{X_n}X_{n-1}) \\
\rTTo \dots \rTTo \uCom(h^{X_0}Z_0,h^{X_n}Z_k) \rTTo^s
s\uCom(h^{X_0}Z_0,h^{X_n}Z_k) \bigr].
\end{multline*}
Cancelling the above $\coev$ and $s$, we obtain an equation
between complexes and maps of the form $\uCom(h^{X_0}Z_0,\_)$.
Simplifying it we get an equation between maps
$\_:h^{X_0}Z_0\tens h^{Z_0}Z_1\tdt
h^{Z_{k-1}}Z_k\tens h^{X_1}X_0\tdt h^{X_n}X_{n-1}\to h^{X_n}Z_k$,
whose corollary is the required identity. Multiplying it with the
composite symmetry
\begin{multline*}
\bigl(
\begin{smallmatrix}
  1  & \dots &n-1& n &n+1& \dots &n+k+1 \\
k+n+1& \dots &k+3&k+2& 1 & \dots & k+1
\end{smallmatrix}
\bigr)_c:
h^{X_n}X_{n-1}\tdt h^{X_1}X_0\tens h^{X_0}Z_0\tens h^{Z_0}Z_1\tdt
h^{Z_{k-1}}Z_k \\
\to h^{X_0}Z_0\tens h^{Z_0}Z_1\tdt
h^{Z_{k-1}}Z_k\tens h^{X_1}X_0\tdt h^{X_n}X_{n-1}
\end{multline*}
and with the common sign $(-1)^n$ we get the equation
\begin{multline*}
b_{n+1+k}b_1 + (1^{\tens n}\tens b_1\tens1^{\tens k})b_{n+1+k}
+ \sum_{\substack{p+q=n\\i+j=k}}^{\substack{p,q,i,j\ge0\\p+i>0,q+j>0}}
(1^{\tens q}\tens b_{p+1+i}\tens1^{\tens j})b_{q+1+j} \\
+ \sum_{\alpha+p+\beta=k}
(1^{\tens n+1+\alpha}\tens b_p\tens1^{\tens\beta})b_{n+1+\alpha+1+\beta}
+ \sum_{\alpha+p+\beta=n}
(1^{\tens\alpha}\tens b_p\tens1^{\tens\beta+1+k})b_{\alpha+1+\beta+1+k}
= 0.
\end{multline*}
This is identity~\eqref{eq-b-b-0}, so the proof is finished.
\end{proof}
\fi

\subsection{The Yoneda embedding}
We claim that the Yoneda \ainf-functor $Y$ is an equivalence of
$\ca^\op$ with its image. This is already proven by Fukaya in the case
of strictly unital \ainf-category $\ca$
\cite[Theorem~9.1]{Fukaya:FloerMirror-II}. This result extends to
arbitrary unital \ainf-categories as follows.

\ifx\chooseClass1
\begin{NOproposition}\label{pro-Y1-homotopy-invertible}
Let $\ca$ be a unital \ainf-category. Let $X$, $W$ be objects of $\ca$.
Then the map $Y_1:s\ca^{\op}(X,W)\to sA_{\infty}(\ca,\uCom)(h^X,h^W)$ is
homotopy invertible.
\end{NOproposition}

\begin{proof}
The composition of $Y_1$ with
\[ \pr_k: sA_{\infty}(\ca,\uCom)(h^X,h^W) \to
\uCom(h^{Z_0}Z_1\tdt h^{Z_{k-1}}Z_k,s\uCom(h^XZ_0,h^WZ_k)) \]
is given by the formula
\begin{multline*}
Y_{1k} = - \bigl[s\ca(W,X) \rTTo^\coev \\
\uCom(h^XZ_0\tens h^{Z_0}Z_1\tdt h^{Z_{k-1}}Z_k,h^XZ_0\tens
h^{Z_0}Z_1\tdt h^{Z_{k-1}}Z_k\tens h^WX) \\
\rTTo^{\uCom(1,\tau^{k,1}_cb_{k+2})}
\uCom(h^XZ_0\tens h^{Z_0}Z_1\tdt h^{Z_{k-1}}Z_k,h^WZ_k) \\
\rTTo^\sim
\uCom(h^{Z_0}Z_1\tdt h^{Z_{k-1}}Z_k,\uCom(h^XZ_0,h^WZ_k)) \\
\rTTo^{\uCom(1,s)}
\uCom(h^{Z_0}Z_1\tdt h^{Z_{k-1}}Z_k,s\uCom(h^XZ_0,h^WZ_k))
\bigr],
\end{multline*}
where the permutation $\tau^{k,1}\in\SSS_{k+2}$ is
$\tau^{k,1}=\bigl(
\begin{smallmatrix}
0 & 1 & \dots & k &k+1 \\
1 & 2 & \dots &k+1& 0
\end{smallmatrix}
\bigr)$.

Define a degree 0 map
$\alpha:sA_{\infty}(\ca,\uCom)(h^X,h^W)\to s\ca^{\op}(X,W)$ as follows:
\begin{multline*}
\alpha = \bigl[sA_{\infty}(\ca,\uCom)(h^X,h^W) \rTTo^{\pr_0}
s\uCom(h^XX,h^WX) \\
\rTTo^{s^{-1}} \uCom(h^XX,h^WX) \rTTo^{\uCom(\sS{_X}\uni^\ca_0,1)}
\uCom(\kk,s\ca(W,X)) =s\ca^{\op}(X,W) \bigr].
\end{multline*}
It is a chain map. Indeed, $\pr_0$ is a chain map, and so is
$s^{-1}\uCom(\sS{_X}\uni^\ca_0,1)$, since
\begin{multline*}
s^{-1}\uCom(\sS{_X}\uni^\ca_0,1)b^{\ca^{\op}}_1
= s^{-1}\uCom(\sS{_X}\uni^\ca_0,1)\uCom(1,b_1)
= s^{-1}[-\uCom(1,b_1)+\uCom(b_1,1)]\uCom(\sS{_X}\uni^\ca_0,1)
\\
= s^{-1}m^\uCom_1\uCom(\sS{_X}\uni^\ca_0,1)
= b^\uCom_1s^{-1}\uCom(\sS{_X}\uni^\ca_0,1).
\end{multline*}

Let us prove that $\alpha$ is homotopy inverse to $Y_1$. First of all,
\begin{equation*}
Y_{10} = - \bigl[s\ca(W,X) \rTTo^\coev \uCom(h^XZ,h^XZ\tens h^WX)
\rTTo^{\uCom(1,cb_2)} \uCom(h^XZ,h^WZ) \xrightarrow{s} s\uCom(h^XZ,h^WZ)
\bigr].
\end{equation*}
Using this expression for $Z=X$ we compute
\begin{align*}
Y_1\alpha &= - \bigl[s\ca(W,X) \rTTo^\coev \uCom(h^XX,h^XX\tens h^WX)
\rTTo^{\uCom(1,cb_2)} \uCom(h^XX,h^WX) \\
&\hspace*{19em} \rTTo^{\uCom(\sS{_X}\uni^\ca_0,1)}
\uCom(\kk,h^WX) = s\ca(W,X) \bigr]
\\
&= \bigl[s\ca(W,X) \rTTo^\coev \uCom(h^XX,h^XX\tens h^WX)
\rTTo^{\uCom(\sS{_X}\uni^\ca_0,1)} \uCom(\kk,h^XX\tens h^WX) \\
&\hspace*{19em} \rTTo^{\uCom(1,cb_2)} \uCom(\kk,h^WX) = s\ca(W,X) \bigr]
\\
&= \bigl[s\ca(W,X) \rTTo^\coev \uCom(\kk,\kk\tens h^WX)
\rTTo^{\uCom(1,\sS{_X}\uni^\ca_0\tens1)} \uCom(\kk,h^XX\tens h^WX) \\
&\hspace*{19em} \rTTo^{\uCom(1,cb_2)} \uCom(\kk,h^WX) = s\ca(W,X) \bigr]
\\
&= \bigl[s\ca(W,X) = \kk\tens h^WX \rTTo^{\sS{_X}\uni^\ca_0\tens1}
h^XX\tens h^WX \rTTo^{cb_2} s\ca(W,X) \bigr] \\
&= \bigl[s\ca(W,X) = h^WX\tens\kk \rTTo^{1\tens\sS{_X}\uni^\ca_0}
h^WX\tens h^XX \rTTo^{b_2} s\ca(W,X) \bigr].
\end{align*}
Therefore, the map $Y_1\alpha$ is homotopic to identity map
$\id:s\ca(W,X)\to s\ca(W,X)$ by \cite[Lemma~7.4]{Lyu-AinfCat}. Now we
are going to prove that $\alpha Y_1$ is homotopy invertible.

The graded $\kk$-module $sA_{\infty}(\ca,\uCom)(h^X,h^W)$ is
$V=\prod_{n=0}^{\infty}V_n$, where
\[ V_n=\prod_{Z_0,\,\dots,\,Z_n\in\Ob\ca}
\uCom(s\ca(Z_0,Z_1)\tdt s\ca(Z_{n-1},Z_n),s\uCom(h^XZ_0,h^WZ_n))
\]
and all products are taken in the category of graded $\kk$\n-modules.
In other terms, for $d\in\ZZ$
\[ V^d = \prod_{n=0}^{\infty}V_n^d, \qquad
V_n^d = \prod_{Z_0,\,\dots,\,Z_n\in\Ob\ca}
\uCom(h^{Z_0}Z_1\tdt h^{Z_{n-1}}Z_n,s\uCom(h^XZ_0,h^WZ_n))^d.
\]
We consider $V_n^d$ as Abelian groups with discrete topology. The
Abelian group $V^d$ is equipped with the topology of the product. Thus,
its basis of neighborhoods of 0 is given by \(\kk\)\n-submodules
$\Phi_m^d=0^{m-1}\times\prod_{n=m}^{\infty}V_n^d$. They form a
filtration $V^d=\Phi_0^d\supset\Phi_1^d\supset\Phi_2^d\supset\dots$. We
call a \(\kk\)\n-linear map $a:V\to V$ of degree $p$ continuous if the
induced maps \(a^{d,d+p}=a\big|_{V^d}:V^d\to V^{d+p}\) are continuous
for all $d\in\ZZ$. This holds if and only if for any $d\in\ZZ$ and
$m\in\NN\overset{\text{def}}=\ZZ_{\ge0}$ there exists an integer
$\kappa=\kappa(d,m)\in\NN$ such that
$(\Phi_\kappa^d)a\subset\Phi_m^{d+p}$. We may assume that
\begin{equation}
m'<m'' \quad \text{ implies } \quad \kappa(d,m') \le \kappa(d,m'').
\label{eq-mm-kk-inequa}
\end{equation}
Indeed, a given function $m\mapsto\kappa(d,m)$ can be replaced with the
function $m\mapsto\kappa'(d,m)=\min_{n\ge m}\kappa(d,n)$ and
\(\kappa'\) satisfies condition~\eqref{eq-mm-kk-inequa}. Continuous
linear maps $a:V\to V$ of degree $p$ are in bijection with families of
$\NN\times\NN$-matrices $(A^{d,d+p})_{d\in\ZZ}$ of linear maps
$A_{nm}^{d,d+p}:V_n^d\to V_m^{d+p}$ with finite number of non-vanishing
elements in each column of $A^{d,d+p}$. Indeed, to each continuous map
\(a^{d,d+p}:V^d\to V^{d+p}\) corresponds the inductive limit over $m$
of $\kappa(d,m)\times m$-matrices of maps
$V^d/\Phi_{\kappa(d,m)}^d\to V^{d+p}/\Phi_m^{d+p}$. On the other hand,
to each family $(A^{d,d+p})_{d\in\ZZ}$ of $\NN\times\NN$-matrices with
finite number of non-vanishing elements in each column correspond
obvious maps $a^{d,d+p}:V^d\to V^{d+p}$, and they are continuous. Thus,
\(a=(a^{d,d+p})_{d\in\ZZ}\) is continuous. A continuous map
\(a:V\to V\) can be completely recovered from one $\NN\times\NN$-matrix
\((a_{nm})_{n,m\in\NN}\) of maps
$a_{nm}=(A_{nm}^{d,d+p})_{d\in\ZZ}:V_n\to V_m$ of degree $p$.
Naturally, not any such matrix determines a continuous map, however, if
the number of non-vanishing elements in each column of \((a_{nm})\) is
finite, then this matrix does determine a continuous map.

The differential $D\overset{\text{def}}=B_1:V\to V$,
$r\mapsto(r)B_1=rb-(-)^rbr$ is continuous and the function \(\kappa\)
for it is simply \(\kappa(d,m)=m\). Its matrix is given by
\begin{align*}
D = B_1 =
\begin{bmatrix}
    D_{0,0} & D_{0,1} & D_{0,2} & \dots\\
     0 & D_{1,1} & D_{1,2} & \dots\\
     0 & 0 & D_{2,2} & \dots\\
     \vdots & \vdots & \vdots &\ddots
\end{bmatrix}
,
\end{align*}
where
\begin{gather*}
D_{k,k} = \uCom(1,b^\uCom_1) -
\uCom\bigl(\sum_{\alpha+1+\gamma=k}
1^{\tens\alpha}\tens b_1\tens1^{\tens\gamma},1\bigr): V_k \to V_k, \\
r_kD_{k,k} = r_kb_1^{\uCom} -(-)^r\sum_{\alpha+1+\gamma=k}
(1^{\tens\alpha}\tens b_1\tens 1^{\tens\gamma})r_k,
\end{gather*}
(one easily recognizes the differential in the complex $V_k$),
\[ r_kD_{k,k+1} = (r_k\tens h^W_1)b^{\uCom}_2
+(h_1^X\tens r_k)b^{\uCom}_2 -(-)^r\sum_{\alpha+\gamma=k-1}
(1^{\tens\alpha}\tens b_2\tens1^{\tens\gamma})r_k.
\]
Further we will see that we do not need to compute other components.

Composition of $\alpha Y_1$ with $\pr_k$ equals
\begin{multline*}
\alpha Y_{1k} = - \bigl[sA_{\infty}(\ca,\uCom)(h^X,h^W) \rTTo^{\pr_0}
s\uCom(h^XX,h^WX) \rTTo^{s^{-1}} \uCom(h^XX,h^WX) \\
\rTTo^{\uCom(\sS{_X}\uni^\ca_0,1)} \uCom(\kk,h^WX) = h^WX \rTTo^\coev \\
\uCom(h^XZ_0\tens h^{Z_0}Z_1\tdt h^{Z_{k-1}}Z_k,h^XZ_0\tens
h^{Z_0}Z_1\tdt h^{Z_{k-1}}Z_k\tens h^WX) \\
\rTTo^{\uCom(1,\tau^{k,1}_cb_{k+2})}
\uCom(h^XZ_0\tens h^{Z_0}Z_1\tdt h^{Z_{k-1}}Z_k,h^WZ_k) \\
\rTTo^\sim
\uCom(h^{Z_0}Z_1\tdt h^{Z_{k-1}}Z_k,\uCom(h^XZ_0,h^WZ_k)) \\
\rTTo^{\uCom(1,s)}
\uCom(h^{Z_0}Z_1\tdt h^{Z_{k-1}}Z_k,s\uCom(h^XZ_0,h^WZ_k))
\bigr].
\end{multline*}
Clearly, $\alpha Y_1$ is continuous (take $\kappa(d,m)=1$). Its
$\NN\times\NN$-matrix has the form
\[ \alpha Y_1 =
\begin{bmatrix}
    * & * & * &\dots \\
    0 & 0& 0& \dots\\
    0 & 0 & 0 &\dots\\
    \vdots & \vdots &\vdots & \ddots
\end{bmatrix}
 .
\]

The reader might skip the proof of the following lemma at first reading.

\begin{NOlemma}\label{lem-alphaY1-homotopy-upper-triangular}
The map $\alpha Y_1:V\to V$ is homotopic to a continuous map $V\to V$,
whose $\NN\times\NN$-matrix is upper-triangular with the identity maps
$\id:V_k\to V_k$ on the diagonal.
\end{NOlemma}

\begin{proof}
Define a continuous $\kk$\n-linear map
$H:sA_{\infty}(\ca,\uCom)(h^X,h^W)\to sA_{\infty}(\ca,\uCom)(h^X,h^W)$ of
degree $-1$ by its matrix
\[ H =
\begin{bmatrix}
    0 & 0 & 0 &\dots \\
    H_{1,0} & 0& 0& \dots\\
    0 & H_{2,1} & 0 &\dots\\
    \vdots & \vdots &\vdots & \ddots
\end{bmatrix}
,
\]
so $\kappa(d,m)=m+1$, where $H_{k+1,k}$ maps the factor indexed by
$(X,Z_0,\dots,Z_k)$ to the factor indexed by $(Z_0,\dots,Z_k)$ as
follows:
\begin{align*}
H_{k+1,k} = \bigl[
&\uCom(s\ca(X,Z_0)\tens s\ca(Z_0,Z_1)\tdt s\ca(Z_{k-1},Z_k),s\uCom(h^XX,h^WZ_k))
\\
\rTTo^{\uCom(1,s^{-1})} &
\uCom(s\ca(X,Z_0)\tens s\ca(Z_0,Z_1)\tdt s\ca(Z_{k-1},Z_k),\uCom(h^XX,h^WZ_k))
\\
\rTTo^{\uCom(1,\uCom(\unix,1))} &
\uCom(s\ca(X,Z_0)\tens s\ca(Z_0,Z_1)\tdt s\ca(Z_{k-1},Z_k),\uCom(\kk,h^WZ_k))
\\
\rTTo^{\sim}& \uCom(s\ca(Z_0,Z_1)\tdt s\ca(Z_{k-1},Z_k),\uCom(h^XZ_0,h^WZ_k))
\\
\rTTo^{\uCom(1,s)}&
\uCom(s\ca(Z_0,Z_1)\tdt s\ca(Z_{k-1},Z_k),s\uCom(h^XZ_0,h^WZ_k)) \bigr].
\end{align*}
Other factors are ignored. Here the adjunction isomorphism is that of
\eqref{eq-g-to-coev-C1g}.

The composition of continuous maps \(V\to V\) is continuous as well. In particular,
one finds the matrices of $B_1H$ and $HB_1$:
\begin{gather*}
B_1H =
\begin{bmatrix}
D_{0,1}H_{1,0} & D_{0,2}H_{2,1} & D_{0,3}H_{3,2} &\dots\\
D_{1,1}H_{1,0} & D_{1,2}H_{2,1} & D_{1,3}H_{3,2} &\dots\\
0 & D_{2,2}H_{2,1} & D_{2,3}H_{3,2} &\dots\\
0 & 0 &D_{3,3}H_{3,2}& \dots\\
\vdots & \vdots & \vdots &\ddots
\end{bmatrix}
,\\
HB_1 =
\begin{bmatrix}
0 & 0 & 0 & \dots \\
H_{1,0}D_{0,0} & H_{1,0}D_{0,1} & H_{1,0}D_{0,2} &\dots\\
0 & H_{2,1}D_{1,1} & H_{2,1}D_{1,2} &\dots\\
0 & 0 & H_{3,2}D_{2,2} &\dots\\
\vdots & \vdots & \vdots &\ddots
\end{bmatrix}
 .
\end{gather*}
We have $D_{k+1,k+1}H_{k+1,k}+H_{k+1,k}D_{k,k}=0$ for all $k\ge0$.
Indeed, conjugating the expanded left hand side with $\uCom(1,s)$
we come to the following identity
\begin{multline*}
\bigl[
\uCom(h^{X}Z_0\tens h^{Z_0}Z_1\tdt h^{Z_{k-1}}Z_k,\uCom(h^XX,h^WZ_k))
\rTTo^{\uCom(1,m^\uCom_1)+\uCom(\sum_{\alpha+\gamma=k}
1^{\tens\alpha}\tens b_1\tens1^{\tens\gamma},1)} \\
\uCom(h^{X}Z_0\tens h^{Z_0}Z_1\tdt h^{Z_{k-1}}Z_k,\uCom(h^XX,h^WZ_k))
\rTTo^{\uCom(1,\uCom(\unix,1))} \\
\uCom(h^{X}Z_0\tens h^{Z_0}Z_1\tdt h^{Z_{k-1}}Z_k,\uCom(\kk,h^WZ_k))
\overset\sim\to \uCom(h^{Z_0}Z_1\tdt h^{Z_{k-1}}Z_k,\uCom(h^XZ_0,h^WZ_k))
\bigr]
\\
+ \bigl[
\uCom(h^{X}Z_0\tens h^{Z_0}Z_1\tdt h^{Z_{k-1}}Z_k,\uCom(h^XX,h^WZ_k))
\rTTo^{\uCom(1,\uCom(\unix,1))} \\
\uCom(h^{X}Z_0\tens h^{Z_0}Z_1\tdt h^{Z_{k-1}}Z_k,\uCom(\kk,h^WZ_k))
\overset\sim\to \uCom(h^{Z_0}Z_1\tdt h^{Z_{k-1}}Z_k,\uCom(h^XZ_0,h^WZ_k))
\\
\rTTo^{\uCom(1,m^\uCom_1)+\uCom(\sum_{\alpha+\gamma=k-1}
1^{\tens\alpha}\tens b_1\tens1^{\tens\gamma},1)}
\uCom(h^{Z_0}Z_1\tdt h^{Z_{k-1}}Z_k,\uCom(h^XZ_0,h^WZ_k)) \bigr] = 0.
\end{multline*}
After reducing all terms to the common form, beginning with
$\uCom(1,\uCom(\unix,1))$, all terms cancel each other, so the identity
is proven.

Therefore, the chain map $a=\alpha Y_1+B_1H+HB_1$ is also represented
by an upper-triangular matrix. Its diagonal elements are chain maps
$a_{kk}:V_k\to V_k$. We are going to show that they are homotopic to
identity maps.

First of all, let us compute the matrix element
$a_{00}:V_0\to V_0=\prod_{Z\in\Ob\ca}s\uCom(h^XZ,h^WZ)$. We have
\[ a_{00}\pr_Z = (\alpha Y_{10}+B_1H_{1,0})\pr_Z:
V_0\to V_0 \rTTo^{\pr_Z} s\uCom(h^XZ,h^WZ).
\]
In the expanded form these terms are as follows:
\begin{multline*}
\alpha Y_{10}\pr_Z = - \bigl[V_0 \rTTo^{\pr_X} s\uCom(h^XX,h^WX)
\rTTo^{s^{-1}} \uCom(h^XX,h^WX) \rTTo^{\uCom(\sS{_X}\uni^\ca_0,1)}
\uCom(\kk,h^WX)
\\
= h^WX \rTTo^\coev \uCom(h^XZ,h^XZ\tens h^WX) \rTTo^{\uCom(1,cb_2)}
\uCom(h^XZ,h^WZ) \rTTo^s s\uCom(h^XZ,h^WZ) \bigr],
\end{multline*}
\[ B_1H_{1,0}\pr_Z = [(1\tens h^W_1)b_2+(h^X_1\tens1)b_2]H_{1,0}\pr_Z,
\]
\begin{multline*}
(1\tens h^W_1)b_2 = \bigl[V_0 \rTTo^{\pr_X} s\uCom(h^XX,h^WX)
\rTTo^\coev \uCom(h^XZ,h^XZ\tens s\uCom(h^XX,h^WX)) \\
\rTTo^{\uCom(1,c(1\tens\coev))}
\uCom(h^XZ,s\uCom(h^XX,h^WX)\tens\uCom(h^WX,h^WX\tens h^XZ)) \\
\rTTo^{\uCom(1,1\tens\uCom(1,b_2))}
\uCom(h^XZ,s\uCom(h^XX,h^WX)\tens\uCom(h^WX,h^WZ))
\rTTo^{\uCom(1,1\tens s)} \\
\uCom(h^XZ,s\uCom(h^XX,h^WX)\tens s\uCom(h^WX,h^WZ))
\rTTo^{\uCom(1,b^\uCom_2)}
\uCom(h^XZ,s\uCom(h^XX,h^WZ)) \bigr],
\end{multline*}
\begin{multline*}
(h^X_1\tens1)b_2 = \bigl[V_0 \rTTo^{\pr_Z} s\uCom(h^XZ,h^WZ)
\rTTo^\coev \uCom(h^XZ,h^XZ\tens s\uCom(h^XZ,h^WZ)) \\
\rTTo^{\uCom(1,\coev\tens1)}
\uCom(h^XZ,\uCom(h^XX,h^XX\tens h^XZ)\tens s\uCom(h^XZ,h^WZ)) \\
\rTTo^{\uCom(1,\uCom(1,b_2)\tens1)}
\uCom(h^XZ,\uCom(h^XX,h^XZ)\tens s\uCom(h^XZ,h^WZ))
\rTTo^{\uCom(1,s\tens1)} \\
\uCom(h^XZ,s\uCom(h^XX,h^XZ)\tens s\uCom(h^XZ,h^WZ))
\rTTo^{\uCom(1,b^\uCom_2)}
\uCom(h^XZ,s\uCom(h^XX,h^WZ)) \bigr],
\end{multline*}
\begin{multline*}
H_{1,0}\pr_Z = \bigl[V_1 \rTTo^{\pr_{X,Z}} \uCom(h^XZ,s\uCom(h^XX,h^WZ))
\rTTo^{\uCom(1,s^{-1})} \uCom(h^XZ,\uCom(h^XX,h^WZ)) \\
\rTTo^{\uCom(1,\uCom(\unix,1))} \uCom(h^XZ,\uCom(\kk,h^WZ))
= \uCom(h^XZ,h^WZ) \rTTo^s s\uCom(h^XZ,h^WZ) \bigr].
\end{multline*}
We claim that in the sum
\[ a_{00}\pr_Z = \alpha Y_{10}\pr_Z + (1\tens h^W_1)b_2H_{1,0}\pr_Z
+ (h^X_1\tens1)b_2H_{1,0}\pr_Z
\]
the first two summands cancel each other, while the last,
$(h^X_1\tens1)b_2H_{1,0}$ is homotopic to identity. Indeed,
$\alpha Y_{10}\pr_Z+(1\tens h^W_1)b_2H_{1,0}\pr_Z$ factors through
\begin{multline*}
- \bigl[ s\uCom(h^XX,h^WX) \rTTo^{s^{-1}} \uCom(h^XX,h^WX)
\rTTo^{\uCom(\unix,1)} \uCom(\kk,h^WX) = h^WX \\
\hspace*{12em} \rTTo^\coev \uCom(h^XZ,h^XZ\tens h^WX)
\rTTo^{\uCom(1,cb_2)} \uCom(h^XZ,h^WZ) \bigr] \\
+ \bigl[ s\uCom(h^XX,h^WX) \rTTo^\coev
\uCom(h^XZ,h^XZ\tens s\uCom(h^XX,h^WX)) \hspace*{12em} \\
\rTTo^{\uCom(1,c(1\tens\coev))}
\uCom(h^XZ,s\uCom(h^XX,h^WX)\tens\uCom(h^WX,h^WX\tens h^XZ))
\\
\rTTo^{\uCom(1,1\tens\uCom(1,b_2))}
\uCom(h^XZ,s\uCom(h^XX,h^WX)\tens\uCom(h^WX,h^WZ))
\rTTo^{\uCom(1,1\tens s)}
\\
\uCom(h^XZ,s\uCom(h^XX,h^WX)\tens s\uCom(h^WX,h^WZ))
\rTTo^{\uCom(1,b^\uCom_2)} \uCom(h^XZ,s\uCom(h^XX,h^WZ))
\rTTo^{\uCom(1,s^{-1})}
\\
\uCom(h^XZ,\uCom(h^XX,h^WZ))
\rTTo^{\uCom(1,\uCom(\unix,1))} \uCom(h^XZ,\uCom(\kk,h^WZ))
= \uCom(h^XZ,h^WZ) \bigr].
\end{multline*}
Composing this sum with $s$ we get
\begin{multline*}
- \bigl[ \uCom(h^XX,h^WX) \rTTo^{\uCom(\unix,1)} \uCom(\kk,h^WX) = h^WX
\\
\hfill \rTTo^\coev \uCom(h^XZ,h^XZ\tens h^WX) \rTTo^{\uCom(1,cb_2)}
\uCom(h^XZ,h^WZ) \bigr] \quad
\\
\quad - \bigl[ \uCom(h^XX,h^WX) \rTTo^\coev
\uCom(h^XZ,h^XZ\tens\uCom(h^XX,h^WX)) \rTTo^{\uCom(1,c(1\tens\coev))}
\hfill
\\
\uCom(h^XZ,\uCom(h^XX,h^WX)\tens\uCom(h^WX,h^WX\tens h^XZ))
\rTTo^{\uCom(1,1\tens\uCom(1,b_2))}
\\
\uCom(h^XZ,\uCom(h^XX,h^WX)\tens\uCom(h^WX,h^WZ))
\rTTo^{\uCom(1,m^\uCom_2)} \uCom(h^XZ,\uCom(h^XX,h^WZ)) \\
\rTTo^{\uCom(1,\uCom(\unix,1))} \uCom(h^XZ,\uCom(\kk,h^WZ))
= \uCom(h^XZ,h^WZ) \bigr].
\end{multline*}
This expression vanishes by the properties of closed monoidal category
$\Com$. Therefore, $\alpha Y_{10}\pr_Z+(1\tens h^W_1)b_2H_{1,0}\pr_Z$
also vanishes.

Now we prove that $(h^X_1\tens1)b_2H_{1,0}$ is homotopic to identity.
It maps each factor $s\uCom(h^XZ,h^WZ)$ into itself via the following
map:
\begin{multline*}
\bigl[ s\uCom(h^XZ,h^WZ) \rTTo^\coev
\uCom(h^XZ,h^XZ\tens s\uCom(h^XZ,h^WZ)) \\
\rTTo^{\uCom(1,\coev\tens1)}
\uCom(h^XZ,\uCom(h^XX,h^XX\tens h^XZ)\tens s\uCom(h^XZ,h^WZ)) \\
\rTTo^{\uCom(1,\uCom(1,b_2)\tens1)}
\uCom(h^XZ,\uCom(h^XX,h^XZ)\tens s\uCom(h^XZ,h^WZ)) \\
\rTTo^{\uCom(1,s\tens1)}
\uCom(h^XZ,s\uCom(h^XX,h^XZ)\tens s\uCom(h^XZ,h^WZ)) \\
\rTTo^{\uCom(1,b^\uCom_2)} \uCom(h^XZ,s\uCom(h^XX,h^WZ))
\rTTo^{\uCom(1,s^{-1})} \uCom(h^XZ,\uCom(h^XX,h^WZ)) \\
\rTTo^{\uCom(1,\uCom(\unix,1))} \uCom(h^XZ,\uCom(\kk,h^WZ))
= \uCom(h^XZ,h^WZ) \rTTo^s s\uCom(h^XZ,h^WZ) \bigr]
\end{multline*}
\begin{multline*}
= \bigl[ s\uCom(h^XZ,h^WZ) \rTTo^{s^{-1}} \uCom(h^XZ,h^WZ)
\rTTo^\coev \uCom(h^XZ,h^XZ\tens\uCom(h^XZ,h^WZ))
\\
\rTTo^{\uCom(1,\coev\tens1)}
\uCom(h^XZ,\uCom(h^XX,h^XX\tens h^XZ)\tens\uCom(h^XZ,h^WZ))
\rTTo^{\uCom(1,\uCom(1,b_2)\tens1)}
\\
\uCom(h^XZ,\uCom(h^XX,h^XZ)\tens\uCom(h^XZ,h^WZ))
\rTTo^{\uCom(1,m^\uCom_2)} \uCom(h^XZ,\uCom(h^XX,h^WZ))
\\
\rTTo^{\uCom(1,\uCom(\unix,1))} \uCom(h^XZ,\uCom(\kk,h^WZ))
= \uCom(h^XZ,h^WZ) \rTTo^s s\uCom(h^XZ,h^WZ) \bigr]
\end{multline*}
\begin{multline*}
= - \bigl[ s\uCom(h^XZ,h^WZ) \rTTo^{s^{-1}} \uCom(h^XZ,h^WZ)
\rTTo^\coev \uCom(h^XZ,h^XZ\tens\uCom(h^XZ,h^WZ))
\\
\rTTo^{\uCom(1,(\unix\tens1)b_2\tens1)}
\uCom(h^XZ,h^XZ\tens\uCom(h^XZ,h^WZ))
\\
\rTTo^{\uCom(1,\ev)} \uCom(h^XZ,h^WZ) \rTTo^s s\uCom(h^XZ,h^WZ) \bigr].
\end{multline*}
Due to \cite[Lemma~7.4]{Lyu-AinfCat} there exists a homotopy
$h'':h^XZ\to h^XZ$ (a map of degree $-1$), such that
$(\unix\tens1)b_2=-1+h''b_1+b_1h''$. Therefore, the map considered
above equals to
\begin{multline*}
\id_{s\uCom(h^XZ,h^WZ)} - \bigl[ s\uCom(h^XZ,h^WZ) \rTTo^{s^{-1}}
\uCom(h^XZ,h^WZ) \rTTo^\coev \uCom(h^XZ,h^XZ\tens\uCom(h^XZ,h^WZ)) \\
 \rTTo^{\uCom(1,(h''b_1+b_1h'')\tens1)}
\uCom(h^XZ,h^XZ\tens\uCom(h^XZ,h^WZ)) \rTTo^{\uCom(1,\ev)} \uCom(h^XZ,h^WZ)
\xrightarrow{s} s\uCom(h^XZ,h^WZ) \bigr]
\\
= \bigl( \id_{s\uCom(h^XZ,h^WZ)} +b^\uCom_1H'_{00} +H'_{00}b^\uCom_1
\bigr),
\end{multline*}
where
\[ H'_{00} = \bigl[ s\uCom(h^XZ,h^WZ) \rTTo^{s^{-1}}
\uCom(h^XZ,h^WZ) \rTTo^{\uCom(h'',1)} \uCom(h^XZ,h^WZ) \xrightarrow{s}
s\uCom(h^XZ,h^WZ) \bigr].
\]
Therefore, $(h^X_1\tens1)b_2H_{1,0}$ and $a_{00}$ are homotopic to
identity.

Now we are proving that diagonal elements $a_{kk}:V_k\to V_k$ are
homotopic to identity maps for $k>0$. Using the explicit formula for
adjunction isomorphism we can re-write the formula for $H$ as follows.
An element $r_{k+1}\in V_{k+1}$ is mapped to direct product over
$Z_0,\dots,Z_k\in\Ob\ca$ of
\begin{multline*}
r_{k+1}H_{k+1,k} = \coev_{h^XZ_0,*}
\uCom(h^XZ_0,r_{k+1}^{X,Z_0,\dots,Z_k}s^{-1}\uCom(\unix,1))s:\\
h^{Z_0}Z_1\tdt h^{Z_{k-1}}Z_k \to s\uCom(h^XZ_0,h^WZ_k).
\end{multline*}
Here
$r_{k+1}^{X,Z_0,\dots,Z_k}:
h^XZ_0\tens h^{Z_0}Z_1\tdt h^{Z_{k-1}}Z_k \to s\uCom(h^XX,h^WZ_k)$
is one of the coordinates of $r_{k+1}$.

Thus, $r_{k}D_{k,k+1}H_{k+1,k}$ is the sum of three terms
\eqref{eq-rBHa}--\eqref{eq-rBHc}:
\begin{subequations}
\begin{multline}
\coev_{h^XZ_0,*}\uCom(h^XZ_0,(r_k^{X,Z_0,\dots,Z_{k-1}}\tens h^W_1)
b_2^{\uCom}s^{-1}\uCom(\unix,1))s: \\
h^{Z_0}Z_1\tdt h^{Z_{k-1}}Z_k \to s\uCom(h^XZ_0,h^WZ_k),
\label{eq-rBHa}
\end{multline}
where
$r_k^{X,Z_0,\dots,Z_{k-1}}:
h^XZ_0\tens h^{Z_0}Z_1\tdt h^{Z_{k-2}}Z_{k-1} \to s\uCom(h^XX,h^WZ_{k-1})$.
Since $b_2^{\uCom}=(s\tens s)^{-1}m_2s=-(s^{-1}\tens s^{-1})m_2s$,
we have
\[ (r_k\tens h^W_1)b^{\uCom}_2s^{-1}
=-(r_ks^{-1}\tens h^W_1s^{-1})m_2^\uCom.
\]
Further, for any $a:X\to A$ the following equation holds:
\begin{multline*}
\bigl(\uCom(A,B)\tens\uCom(B,C) \rTTo^{m_2} \uCom(A,C)
\rTTo^{\uCom(a,C)} \uCom(X,C)\bigr) \\
= \bigl(\uCom(A,B)\tens\uCom(B,C) \rTTo^{\uCom(a,B)\tens1}
\uCom(X,B)\tens\uCom(B,C) \rTTo^{m_2} \uCom(X,C)\bigr).
\end{multline*}
Thus,
\begin{multline*}
(r_k\tens h^W_1)b^{\uCom}_2s^{-1}\uCom(\unix,1) =
-(r_ks^{-1}\tens h^W_1s^{-1})m_2^\uCom\uCom(\unix,1) \\
= -(r_ks^{-1}\tens h^W_1s^{-1})(\uCom(\unix,1)\tens 1)m_2^\uCom
= (r_ks^{-1}\uCom(\unix,1)\tens h^{W}_1s^{-1})m_2^\uCom
\end{multline*}
(we have used the fact that $\uCom(\unix,1)$ has degree $-1$ and
$h^{W}_1s^{-1}$ has degree $1$).
\begin{multline}
\coev_{h^XZ_0,*}\uCom(h^XZ_0,(h^{X}_1\tens r_k^{Z_0,\dots,Z_k})
b_2^{\uCom}s^{-1}\uCom(\unix,1))s: \\
h^{Z_0}Z_1\tdt h^{Z_{k-1}}Z_k \to s\uCom(h^XZ_0,h^WZ_k),
\label{eq-rBHb}
\end{multline}
where
$r_k^{Z_0,\dots,Z_k}:
h^{Z_0}Z_1\tdt h^{Z_{k-1}}Z_k \to s\uCom(h^XZ_0,h^WZ_k)$.
Similarly to above
\begin{equation*}
(h^{X}_1\tens r_k)b^{\uCom}_2s^{-1}
= (-)^{r+1}(h^X_1s^{-1}\tens r_ks^{-1})m_2^\uCom,
\end{equation*}
so that
\begin{multline*}
(h^{X}_1\tens r_k)b^{\uCom}_2s^{-1}\uCom(\unix,1)
= (-)^{r+1}(h^X_1s^{-1}\tens r_ks^{-1})
(\uCom(\unix,1)\tens1)m_2^\uCom \\
= (h^X_1s^{-1}\uCom(\unix,1)\tens r_ks^{-1})m_2^\uCom
\end{multline*}
($r_ks^{-1}$ has degree $\deg r+1$ and $\uCom(\unix,1)$ has degree
$-1$).

For each $\alpha$, $\gamma$, such that $\alpha+\gamma=k-1$
\begin{multline}
\coev_{h^XZ_0,*}
\uCom(h^XZ_0,(1^{\tens\alpha}\tens b_2\tens1^{\tens\gamma})r_k
s^{-1}\uCom(\unix,1))s: \\
h^{Z_0}Z_1\tdt h^{Z_{k-1}}Z_k \to s\uCom(h^XZ_0,h^WZ_k),
\label{eq-rBHc}
\end{multline}
where $r_k$ means
$r_k^{X,Z_0,\dots,Z_{\alpha-1},Z_{\alpha+1},\dots,Z_k}:
h^XZ_0\tens h^{Z_0}Z_1\tdt h^{Z_{\alpha-1}}Z_{\alpha+1}\tdt h^{Z_{k-1}}Z_k
\to s\uCom(h^XX,h^WZ_k)$,
and $Z_{-1}=X$.
\end{subequations}

Thus,
$r_kD_{k,k+1}H_{k+1,k}=\coev_{h^XZ_0,*}\uCom(h^XZ_0,\Sigma_1)s$,
where
\begin{multline*}
\Sigma_1 = (r_k^{X,Z_0,\dots,Z_{k-1}}s^{-1}
\uCom(\unix,1)\tens h^{W}_1s^{-1})m_2^\uCom
+(h^X_1s^{-1}\uCom(\unix,1)\tens r_k^{Z_0,\dots,Z_k}s^{-1})
m_2^\uCom \\
-(-)^r\sum_{\alpha+\gamma=k-1}
(1^{\tens\alpha}\tens b_2\tens1^{\tens\gamma})
r_k^{X,Z_0,\dots,Z_{\alpha-1},Z_{\alpha+1},\dots,Z_k}
s^{-1}\uCom(\unix,1): \\
h^XZ_0\tens h^{Z_0}Z_1\tdt h^{Z_{k-1}}Z_k \to h^{W}Z_k.
\end{multline*}

Similarly, $r_kH_{k,k-1}D_{k-1,k}$ is the sum of three terms
\eqref{eq-rHBa}--\eqref{eq-rHBc}.
\begin{subequations}
\begin{multline}
\bigl(\coev_{h^XZ_0,*}\uCom(h^XZ_0,r_ks^{-1}\uCom(\unix,1))s
\tens\coev_{h^WZ_{k-1},h^{Z_{k-1}}Z_k}\uCom(h^WZ_{k-1},b_2)s\bigr)
b^{\uCom}_2 \\
=-\bigl(\coev_{h^XZ_0,*}\uCom(h^XZ_0,r_ks^{-1}\uCom(\unix,1))\tens
\coev_{h^WZ_{k-1},h^{Z_{k-1}}Z_k}\uCom(h^WZ_{k-1},b_2)\bigr)m_2s \\
=-\coev_{h^XZ_0,*}\uCom\bigl(h^XZ_0,(r_ks^{-1}\uCom(\unix,1)\tens1)
b_2\bigr)s: \\
h^{Z_0}Z_1\tdt h^{Z_{k-1}}Z_k \to s\uCom(h^XZ_0,h^WZ_k),
\label{eq-rHBa}
\end{multline}
where $r_k$ means
$r_k^{X,Z_0,\dots,Z_{k-1}}:
h^XZ_0\tens h^{Z_0}Z_1\tdt h^{Z_{k-2}}Z_{k-1} \to s\uCom(h^XX,h^WZ_{k-1})$.
Here we use \eqref{eq-identity-m2} (compare with
\eqref{eq-compare-hhb}).
\begin{multline}
\bigl(\coev_{h^XZ_0,h^{Z_0}Z_1}\uCom(h^XZ_0,b_2)s\tens
\coev_{h^XZ_1,*}\uCom(h^XZ_1,r_ks^{-1}\uCom(\unix,1))
s\bigr)b^{\uCom}_2 \\
=(-)^r\bigl(\coev_{h^XZ_0,h^{Z_0}Z_1}\uCom(h^XZ_0,b_2)\tens
\coev_{h^XZ_1,*}\uCom(h^XZ_1,r_ks^{-1}\uCom(\unix,1))\bigr)
m^\uCom_2s \\
=(-)^r\coev_{h^XZ_0,*}\uCom\bigl(h^XZ_0,
(b_2\tens1^{\tens k-1})r_ks^{-1}\uCom(\unix,1)\bigr)s: \\
h^{Z_0}Z_1\tdt h^{Z_{k-1}}Z_k \to s\uCom(h^XZ_0,h^WZ_k),
\label{eq-rHBb}
\end{multline}
where $r_k$ means the coordinate
$r_k^{X,Z_1,\dots,Z_k}:
h^XZ_1\tens h^{Z_1}Z_2\tdt h^{Z_{k-1}}Z_k \to s\uCom(h^XX,h^WZ_k)$.
We have used that $r_ks^{-1}\uCom(\unix,1)$ has degree $\deg r$ and
formula~\eqref{eq-identity-m2}.
\begin{multline}
(1^{\tens\alpha}\tens b_2\tens1^{\tens\gamma})
\coev_{h^XZ_0,*}\uCom\bigl(h^XZ_0,r_ks^{-1}\uCom(\unix,1)\bigr)s \\
= \coev_{h^XZ_0,*}\uCom(h^XZ_0,1^{\tens\alpha+1}\tens b_2\tens1^{\tens\gamma})
\uCom\bigl(h^XZ_0,r_ks^{-1}\uCom(\unix,1)\bigr)s \\
=\coev_{h^XZ_0,*}\uCom\bigl(h^XZ_0,(1^{\tens\alpha+1}\tens b_2\tens1^{\tens\gamma})
r_ks^{-1}\uCom(\unix,1)\bigr)s: \\
h^{Z_0}Z_1\tdt h^{Z_{k-1}}Z_k \to s\uCom(h^XZ_0,h^WZ_k),
\label{eq-rHBc}
\end{multline}
where $r_k$ means the coordinate
$r_k^{X,Z_0,\dots,Z_{\alpha},Z_{\alpha+2},\dots,Z_k}:
h^XZ_0\tens h^{Z_0}Z_1\tdt h^{Z_{\alpha}}Z_{\alpha+2}\tdt h^{Z_{k-1}}Z_k
\to s\uCom(h^XX,h^WZ_k)$.
Here we use functoriality of $\coev$ (compare with
\eqref{eq-compare2-bh}).
\end{subequations}

Thus,
$r_kH_{k,k-1}D_{k-1,k}=\coev_{h^XZ_0,*}\uCom(h^XZ_0,\Sigma_2)s$,
where
\begin{multline*}
\Sigma_2 =-(r_k^{X,Z_0,\dots,Z_{k-1}}s^{-1}\uCom(\unix,1)\tens1)b_2
+(-)^r(b_2\tens1^{\tens k-1})r_k^{X,Z_1,\dots,Z_k}s^{-1}
\uCom(\unix,1) \\
-(-)^{r-1}\sum_{\alpha+\gamma=k-2}
(1^{\tens\alpha+1}\tens b_2\tens1^{\tens\gamma})
r_k^{X,Z_0,\dots,Z_{\alpha},Z_{\alpha+2},\dots,Z_k}
s^{-1}\uCom(\unix,1): \\
h^XZ_0\tens h^{Z_0}Z_1\tdt h^{Z_{k-1}}Z_k \to h^{W}Z_k.
\end{multline*}
The element $rH$ has degree $\deg r-1$, so the sign $(-)^{r-1}$ arises.
Combining this with the expression for $r_kD_{k,k+1}H_{k+1,k}$ we
obtain
\[ r_kD_{k,k+1}H_{k+1,k}+r_kH_{k,k-1}D_{k-1,k}
= \coev_{h^XZ_0,*}\uCom(h^XZ_0,\Sigma)s,
\]
where $\Sigma=\Sigma_1+\Sigma_2$. We claim that
$\Sigma=(h^X_1s^{-1}\uCom(\unix,1)\tens
r_k^{Z_0,\dots,Z_k}s^{-1})m_2$. Indeed, first of all,
\begin{multline*}
(b_2\tens1^{\tens k-1})r_k^{X,Z_1,\dots,Z_k}s^{-1}\uCom(\unix,1) \\
+\sum_{\alpha+\gamma=k-2}(1^{\tens\alpha+1}\tens b_2\tens1^{\tens\gamma})
r_k^{X,Z_0,\dots,Z_\alpha,Z_{\alpha+2},\dots,Z_k}s^{-1}\uCom(\unix,1)\\
=\sum_{\alpha+\gamma=k-1}(1^{\tens\alpha}\tens b_2\tens1^{\tens\gamma})
r_k^{X,Z_0,\dots,Z_{\alpha-1},Z_{\alpha+1},\dots,Z_k}
s^{-1}\uCom(\unix,1),
\end{multline*}
so that
\begin{multline*}
\Sigma = (r_k^{X,Z_0,\dots,Z_{k-1}}s^{-1}\uCom(\unix,1)\tens
h^W_1s^{-1})m^\uCom_2
+(h^X_1s^{-1}\uCom(\unix,1)\tens r_k^{Z_0,\dots,Z_k}s^{-1})
m^\uCom_2 \\
-(r_k^{X,Z_0,\dots,Z_{k-1}}s^{-1}\uCom(\unix,1)\tens 1)b_2:
h^XZ_0\tens h^{Z_0}Z_1\tdt h^{Z_{k-1}}Z_k \to h^{W}Z_k.
\end{multline*}
Further, for all $y\in h^XZ_0$, $z_i\in h^{Z_{i-1}}Z_i$ we have the
action of the third summand
\begin{multline*}
(y\tens z_1\tdt z_{k-1}\tens z_k)(r_k^{X,Z_0,\dots,Z_{k-1}}
s^{-1}\uCom(\unix,1)\tens1)b_2 \\
= (-)^{rz_k}((y\tens z_1\tdt z_{k-1})r_k^{X,Z_0,\dots,Z_{k-1}}
s^{-1}\uCom(\unix,1)\tens z_k)b_2 \\
=(-)^{rz_k+(r+1)+(y+z_1+\dots+z_{k-1})}((\unix)(y\tens z_1\tdt z_{k-1})
r_k^{X,Z_0,\dots,Z_{k-1}}s^{-1}\tens z_k)b_2.
\end{multline*}
Notice that $r_ks^{-1}\uCom(\unix,1)$ has degree $\deg r$, and
$(y\tens z_1\tdt z_{k-1})r_ks^{-1}$ has degree
$(r+1)+(y+z_1+\dots+z_{k-1})$. Similarly, the action of the first
summand is
\begin{multline*}
(y\tens z_1\tdt z_{k-1}\tens z_k)(r_k^{X,Z_0,\dots,Z_{k-1}}
s^{-1}\uCom(\unix,1)\tens h^W_1s^{-1})m^\uCom_2 \\
=(-)^{rz_k+(r+1)+(y+z_1+\dots+z_k)}((\unix)(y\tens z_1\tdt z_{k-1})
r_k^{X,Z_0,\dots,Z_{k-1}}s^{-1}\tens z_kh^W_1s^{-1})m^\uCom_2,
\end{multline*}
and since
$z_kh^W_1s^{-1}=z_k\coev_{h^WZ_{k-1},h^{Z_{k-1}}Z_k}
\uCom(h^WZ_{k-1},b_2)=(t\mapsto(t\tens z_k)b_2)$,
we obtain
\begin{multline*}
(y\tens z_1\tdt z_{k-1}\tens z_k)(r_k^{X,Z_0,\dots,Z_{k-1}}
s^{-1}\uCom(\unix,1)\tens h^W_1s^{-1})m^\uCom_2 \\
=(-)^{rz_k+(r+1)+(y+z_1+\dots+z_k)}((\unix)(y\tens z_1\tdt z_{k-1})
r_k^{X,Z_0,\dots,Z_{k-1}}s^{-1}\tens z_k)b_2.
\end{multline*}
We see that the first and the third terms in $\Sigma$ cancel each
other. Hence, only the second summand remains in
 $\Sigma=(h^X_1s^{-1}\uCom(\unix,1)\tens
 r_k^{Z_0,\dots,Z_k}s^{-1})m^\uCom_2$,
\begin{multline*}
(y\tens z_1\tdt z_k)(h^X_1s^{-1}\uCom(\unix,1)\tens
r_k^{Z_0,\dots,Z_k}s^{-1})m^\uCom_2 \\
= (yh^X_1s^{-1}\uCom(\unix,1)\tens(z_1\tdt z_k)
r_k^{Z_0,\dots,Z_k}s^{-1})m^\uCom_2.
\end{multline*}
We have
\begin{equation*}
yh^X_1s^{-1} = (y)\coev_{h^XX,h^XZ_0}\uCom(h^XX,b_2)
= (t\mapsto (t\tens y)b_2)
\end{equation*}
so that
$yh^X_1s^{-1}\uCom(\unix,1)=(-)^{y+1}(\unix\tens y)b_2
=-(y)(\unix\tens 1)b_2=y-yh''b_1-yb_1h''$.
Then
\begin{equation*}
(z_1\tdt z_k)\coev_{h^XZ_0,*}\uCom(h^XZ_0,\Sigma) =
\bigl(y\mapsto(y-yh''b_1-yb_1h'')(z_1\tdt z_k)(r_k^{Z_0,\dots,Z_k}s^{-1})\bigr)
\end{equation*}
and
\begin{equation*}
(z_1\tdt z_k)\coev_{h^XZ_0,*}\uCom(h^XZ_0,\Sigma)s =
\bigl(y\mapsto(y-yh''b_1-yb_1h'')(z_1\tdt z_k)r_k^{Z_0,\dots,Z_k}\bigr).
\end{equation*}

We see that for $k>0$
\begin{multline*}
a_{kk} = D_{k,k+1}H_{k+1,k} + H_{k,k-1}D_{k-1,k} = 1+f: \\
\uCom(h^{Z_0}Z_1\tdt h^{Z_{k-1}}Z_k,s\uCom(h^XZ_0,h^WZ_k))
\to \uCom(h^{Z_0}Z_1\tdt h^{Z_{k-1}}Z_k,s\uCom(h^XZ_0,h^WZ_k)),
\end{multline*}
where
 $f=-\uCom(h^{Z_0}Z_1\tdt h^{Z_{k-1}}Z_k,s^{-1}\uCom(h''b_1+b_1h'',1)s)$.
More precisely, $D_{k,k+1}H_{k+1,k}+H_{k,k-1}D_{k-1,k}$ is a diagonal
map, whose components are $1+f$. We claim that
$f=m^\uCom_1H'_{kk}+H'_{kk}m^\uCom_1$, where
 $H'_{kk}=\uCom(h^{Z_0}Z_1\tdt h^{Z_{k-1}}Z_k,s^{-1}\uCom(h'',1)s)$.
Indeed,
 $m^\uCom_1=\uCom(1,b^{\uCom}_1)-\uCom(\sum_{\alpha+\gamma=k-1}
 (1^{\tens\alpha}\tens b_1\tens1^{\tens\gamma}),1)$,
so that
\begin{align*}
m^\uCom_1H'_{kk} &= \bigl(\uCom(1,b^\uCom_1)-\uCom(\sum_{\alpha+\gamma=k-1}
1^{\tens\alpha}\tens b_1\tens1^{\tens\gamma},1)\bigr)
\uCom\bigl(1,s^{-1}\uCom(h'',1)s\bigr) \\
&= \uCom\bigl(1,b^{\uCom}_1s^{-1}\uCom(h'',1)s\bigr)
+\uCom\bigl(1,s^{-1}\uCom(h'',1)s\bigr)
\uCom\bigl(\sum_{\alpha+\gamma=k-1}
1^{\tens\alpha}\tens b_1\tens1^{\tens\gamma},1\bigr), \\
H'_{kk}m^\uCom_1 &= \uCom\bigl(1,s^{-1}\uCom(h'',1)s\bigr)
\bigl(\uCom(1,b^{\uCom}_1)-\uCom(\sum_{\alpha+\gamma=k-1}
1^{\tens\alpha}\tens b_1\tens 1^{\tens\gamma},1)\bigr) \\
&= \uCom\bigl(1,s^{-1}\uCom(h'',1)s b^{\uCom}_1\bigr)
-\uCom\bigl(1,s^{-1}\uCom(h'',1)s\bigr)
\uCom\bigl(\sum_{\alpha+\gamma=k-1}
1^{\tens\alpha}\tens b_1\tens 1^{\tens\gamma},1\bigr).
\end{align*}
Therefore,
\begin{align*}
H'_{kk}m^\uCom_1 + m^\uCom_1H'_{kk}
&= \uCom\bigl(1,b^{\uCom}_1s^{-1}\uCom(h'',1)s
+s^{-1}\uCom(h'',1)s b^{\uCom}_1\bigr)
\\
&= \uCom\bigl(1,s^{-1}(m^\uCom_1\uCom(h'',1)
+\uCom(h'',1)m^\uCom_1)s\bigr).
\end{align*}
Since $m^\uCom_1=-\uCom(1,b_1)+\uCom(b_1,1)$, we have
\begin{multline*}
m^\uCom_1\uCom(h'',1)+\uCom(h'',1)m^\uCom_1
\\
= -\uCom(1,b_1)\uCom(h'',1) +\uCom(b_1,1)\uCom(h'',1)
-\uCom(h'',1)\uCom(1,b_1) +\uCom(h'',1)\uCom(b_1,1)
\\
= \uCom(h''b_1+b_1h'',1),
\end{multline*}
so that
$H'_{kk}m^\uCom_1+m^\uCom_1H'_{kk}
=-\uCom(1,s^{-1}\uCom(h''b_1+b_1h'',1)s)=f$.

Summing up, we have proved that
\[ a = \alpha Y_1 + B_1H + HB_1 = 1 + B_1H' + H'B_1 + N,
\]
where
 $H':sA_{\infty}(\ca,\uCom)(h^X,h^W)\to sA_{\infty}(\ca,\uCom)(h^X,h^W)$
is a continuous $\kk$\n-linear map of degree $-1$ determined by a
diagonal matrix with the matrix elements
\[ H'_{kk}=\uCom(h^{Z_0}Z_1\tdt h^{Z_{k-1}}Z_k,
s^{-1}\uCom(h'',1)s): V_k \to V_k,
\]
$H'_{kl}=0$ for $k\ne l$, and the matrix of the remainder $N$ is
strictly upper-triangular: $N_{kl}=0$ for all $k\ge l$.
\end{proof}

The continuous map of degree 0
\[ \alpha Y_1 + B_1(H-H') + (H-H')B_1 = 1 + N: V \to V,
\]
obtained in \lemref{lem-alphaY1-homotopy-upper-triangular}, is
invertible (its inverse is determined by the upper-triangular matrix
$\sum_{i=0}^\infty(-N)^i$, which is well-defined). Therefore,
$\alpha Y_1$ is homotopy invertible. We have proved earlier that
$Y_1\alpha$ is homotopic to identity. Viewing $\alpha$, $Y_1$ as
morphisms of homotopy category $H^0(\uCom)$, we see that both of them
are homotopy invertible. Hence, they are homotopy inverse to each
other.
\end{proof}
\fi

Let us define a full subcategory $\Rep A_\infty^u(\ca,\uCom)$ of the
$\fu'$\n-small differential graded $\fu$\n-category
$A_\infty^u(\ca,\uCom)$ as follows. Its objects are all \ainf-functors
$h^X:\ca\to\uCom$ for $X\in\Ob\ca$. As we know, they are unital. The
differential graded category $\Rep A_\infty^u(\ca,\uCom)$ is
$\fu$\n-small. Thus, the Yoneda \ainf-functor
 $Y:\ca^{\op}\to A_{\infty}(\ca,\uCom)$ takes values in the
$\fu$\n-small subcategory $\Rep A_\infty^u(\ca,\uCom)$.

\begin{theorem}\label{thm-Yoneda-equivalence}
Let $\ca$ be a unital \ainf-category. Then the restricted Yoneda
\ainf-functor $Y:\ca^{\op}\to\Rep A_\infty^u(\ca,\uCom)$ is an
equivalence.
\end{theorem}

The theorem follows immediately from Corollary~A.9 of
\cite{LyuMan-AmodSerre} which states that
 $\Yo:\ca^{\op}\to\Rep A_\infty^u(\ca,\uCom)$ is an \ainf-equivalence.
This is a corollary to a much stronger result, the \ainf-version of the
Yoneda Lemma \cite[Theorem~A.1]{LyuMan-AmodSerre}.

\ifx\chooseClass1
\begin{proof}
Let us verify the assumptions of \corref{cor-Theorem8:8}.
\propref{pro-Y1-homotopy-invertible} shows that the maps $Y_1$ are
homotopy invertible. The map
 $\Ob Y:\Ob\ca\to\Ob\Rep A_\infty^u(\ca,\uCom)$, $X\mapsto h^X$ is
surjective (in our conventions it may happen that $X\ne Y$, but
$h^X=h^Y$). Thus, assumptions of \corref{cor-Theorem8:8} are satisfied.
We deduce by it that $Y$ is an equivalence.
\end{proof}
\fi

\begin{corollary}\label{cor-unital-equivalent-differential}
Each $\fu$\n-small unital \ainf-category $\ca$ is \ainf-equivalent to a
$\fu$\n-small differential graded category
$\Rep A_\infty^u(\ca,\uCom)^\op$.
\end{corollary}

\begin{remark}\label{rem-Yoneda-tilde}
We may use the surjective map
 $\Ob Y:\Ob\ca\to\Ob\Rep A_\infty^u(\ca,\uCom)$ to transfer the
differential graded category structure of $\Rep A_\infty^u(\ca,\uCom)$
to $\Ob\ca$. This new $\fu$\n-small differential graded category is
denoted $\wt{\Rep}A_\infty^u(\ca,\uCom)$. Thus, its set of objects is
$\Ob\ca$, the sets of morphisms are
\[ \wt{\Rep}A_\infty^u(\ca,\uCom)(X,Y) =A_\infty(\ca,\uCom)(h^X,h^Y),
\]
and the operations are those of $A_\infty(\ca,\uCom)$. It is equivalent
to $\Rep A_\infty^u(\ca,\uCom)$ by surjectivity of $\Ob Y$. The Yoneda
\ainf-functor can be presented as an \ainf-equivalence
\(\wt{Y}:\ca^{\op}\to\wt{\Rep}A_\infty^u(\ca,\uCom)\), identity on
objects, whose components $\wt{Y}_n=Y_n$ are given by
\eqref{eq-components-Yn}. A quasi-inverse to \(\wt{Y}\) equivalence
\(\wt{\Rep}A_\infty^u(\ca,\uCom)\to\ca^{\op}\) can be chosen so that it
induces the identity map on objects as well by \corref{cor-Theorem8:8}.
\end{remark}

\section{\texorpdfstring{Strict $A_\infty^u$-2-functor}
 {Strict A8u-2-functor}}\label{sec-strict-A2fun}
The goal of this section is to show that the problem of representing
the \ainfu-2-functor \(\ca\mapsto A^u_\infty(\cc,\ca)_{\modulo\cb}\) for
a pair $(\cc,\cb)$ of a unital \ainf-category $\cc$ and its full
subcategory $\cb$ reduces to the case of differential graded $\cc$.

\subsection{\texorpdfstring{An $A_\infty$-functor}{An A8-functor}}
For arbitrary \ainf-categories $\cx$, $\cy$, $\cz$ the left hand side
of the equation
\begin{multline}
\bigl[ TsA_\infty(\cy,\cz)\boxtimes TsA_\infty(\cx,\cy) \rTTo^c
TsA_\infty(\cx,\cy)\boxtimes TsA_\infty(\cy,\cz) \rTTo^M
TsA_\infty(\cx,\cz) \bigr] \\
= \bigl[ TsA_\infty(\cy,\cz)\boxtimes TsA_\infty(\cx,\cy)
\rTTo^{1\boxtimes A_\infty(\_,\cz)}
TsA_\infty(\cy,\cz)\boxtimes
TsA_\infty(A_\infty(\cy,\cz),A_\infty(\cx,\cz)) \\
\rTTo^\alpha TsA_\infty(\cx,\cz) \bigr]
\label{eq-Ain(Z)-def}
\end{multline}
is an \ainf-functor. Therefore, by Proposition~5.5 of
\cite{Lyu-AinfCat} there exists a unique \ainf-functor
\[ A_\infty(\_,\cz): A_\infty(\cx,\cy) \to
A_\infty(A_\infty(\cy,\cz),A_\infty(\cx,\cz))
\]
in the right hand side, which makes equation~\eqref{eq-Ain(Z)-def} hold
true. The proof of Proposition~3.4 of \cite{Lyu-AinfCat} contains a
recipe for finding the components of \(A_\infty(\_,\cz)\). Namely, the
equation
\begin{equation}
(p\boxtimes1)M = [p.A_\infty(\_,\cz)]\theta
\label{eq-p1M-pA(Z)theta}
\end{equation}
has to hold for all \(p\in TsA_\infty(\cx,\cy)\). In particular,
\begin{align*}
f.A_\infty(\_,\cz) = (f\boxtimes1)M = (f\boxtimes1_\cz)M &:
A_\infty(\cy,\cz) \to A_\infty(\cx,\cz)
\quad \text{for } f\in\Ob A_\infty(\cx,\cy), \\
r.A_\infty(\_,\cz)_1 = (r\boxtimes1)M = (r\boxtimes1_\cz)M &:
(f\boxtimes1)M \to (g\boxtimes1)M
\text{ for } r\in sA_\infty(\cx,\cy)(f,g).
\end{align*}
Other components of \(A_\infty(\_,\cz)\) are obtained from the
recurrent relation, which is equation~\eqref{eq-p1M-pA(Z)theta} written
for \(p=p^1\tens\dots\tens p^n\):
\begin{multline}
(p^1\tens\dots\tens p^n)A_\infty(\_,\cz)_n
= (p^1\tens\dots\tens p^n\boxtimes1)M \\
- \sum_{i_1+\dots+i_l=n}^{l>1}
[(p^1\tens\dots\tens p^n).(A_\infty(\_,\cz)_{i_1}\tens
A_\infty(\_,\cz)_{i_2}\tens\dots\tens A_\infty(\_,\cz)_{i_l})]\theta.
\label{eq-A(Z)-components}
\end{multline}
In particular, for \(r\tens t\in T^2sA_\infty(\cx,\cy)\) we get
\[ (r\tens t)A_\infty(\_,\cz)_2 = (r\tens t\boxtimes1)M
- [(r\boxtimes1)M\tens(t\boxtimes1)M]\theta.
\]
Given $g^0 \rTTo^{p^1} g^1 \rTTo^{p^2} \dots g^{n-1} \rTTo^{p^n} g^n$
we find from \eqref{eq-A(Z)-components} the components of the
\ainf-transformation
\[ (p^1\tens\dots\tens p^n)A_\infty(\_,\cz)_n \in
sA_\infty(A_\infty(\cy,\cz),A_\infty(\cx,\cz))
((g^0\boxtimes1)M,(g^n\boxtimes1)M)
\]
in the form
\[ [(p^1\tens\dots\tens p^n)A_\infty(\_,\cz)_n]_m =
(p^1\tens\dots\tens p^n\boxtimes1)M_{nm}.
\]
So they vanish for $m>1$.

If $\cz$ is unital, then the \ainf-functor \(A_\infty(\_,\cz)\) takes
values in the subcategory
\[ A_\infty^u(A_\infty(\cy,\cz),A_\infty(\cx,\cz)),
\]
because the \ainf-functors \((f\boxtimes1_\cz)M\) commute with
\((1\boxtimes\uni^\cz)M\), so they are unital.

\begin{proposition}
For arbitrary \ainf-categories $\cx$, $\cy$ and unital \ainf-categories
$\cc$, $\cd$ we have
\begin{multline}
\bigl[ TsA_\infty^u(\cc,\cd)\boxtimes TsA_\infty(\cx,\cy) \rTTo^c
TsA_\infty(\cx,\cy)\boxtimes TsA_\infty^u(\cc,\cd)
\rTTo^{A_\infty(\_,\cc)\boxtimes A_\infty(\cx,\_)} \\
TsA_\infty^u(A_\infty(\cy,\cc),A_\infty(\cx,\cc))\boxtimes
TsA_\infty^u(A_\infty(\cx,\cc),A_\infty(\cx,\cd)) \\
\hspace{15em} \rTTo^M
TsA_\infty^u(A_\infty(\cy,\cc),A_\infty(\cx,\cd)) \bigr] \\
= \bigl[ TsA_\infty^u(\cc,\cd)\boxtimes TsA_\infty(\cx,\cy)
\rTTo^{A_\infty(\cy,\_)\boxtimes A_\infty(\_,\cd)} \hspace{15em} \\
TsA_\infty^u(A_\infty(\cy,\cc),A_\infty(\cy,\cd))\boxtimes
TsA_\infty^u(A_\infty(\cy,\cd),A_\infty(\cx,\cd)) \\
\rTTo^M TsA_\infty^u(A_\infty(\cy,\cc),A_\infty(\cx,\cd)) \bigr].
\label{eq-A(-C)-A(X-)-A(Y-)-A(-D)}
\end{multline}
The same statement holds true if one removes the unitality superscript
$u$, and do not assume $\cc$, $\cd$ unital. The same equation holds
true if all four \ainf-categories $\cx$, $\cy$, $\cc$, $\cd$ are unital
and all \ainf-categories $A_\infty(,)$ are replaced with their
subcategories $A_\infty^u(,)$.
\end{proposition}

\begin{proof}
Due to Proposition~3.4 of \cite{Lyu-AinfCat}
equation~\eqref{eq-A(-C)-A(X-)-A(Y-)-A(-D)} is equivalent to the
following one:
\begin{multline*}
\bigl[
TsA_\infty(\cy,\cc)\boxtimes TsA_\infty^u(\cc,\cd)\boxtimes TsA_\infty(\cx,\cy)
\rTTo^{1\boxtimes c}
TsA_\infty(\cy,\cc)\boxtimes TsA_\infty(\cx,\cy)\boxtimes TsA_\infty^u(\cc,\cd)
\\
\rTTo^{1\boxtimes A_\infty(\_,\cc)\boxtimes A_\infty(\cx,\_)}
TsA_\infty(\cy,\cc)\boxtimes TsA_\infty^u(A_\infty(\cy,\cc),A_\infty(\cx,\cc))
\boxtimes TsA_\infty^u(A_\infty(\cx,\cc),A_\infty(\cx,\cd)) \\
\hspace{5em} \rTTo^{1\boxtimes M}
TsA_\infty(\cy,\cc)\boxtimes TsA_\infty^u(A_\infty(\cy,\cc),A_\infty(\cx,\cd))
\rTTo^\alpha TsA_\infty(\cx,\cd) \bigr] \\
= \bigl[
TsA_\infty(\cy,\cc)\boxtimes TsA_\infty^u(\cc,\cd)\boxtimes TsA_\infty(\cx,\cy)
\rTTo^{1\boxtimes A_\infty(\cy,\_)\boxtimes A_\infty(\_,\cd)} \hspace{9em} \\
TsA_\infty(\cy,\cc)\boxtimes TsA_\infty^u(A_\infty(\cy,\cc),A_\infty(\cy,\cd))
\boxtimes TsA_\infty^u(A_\infty(\cy,\cd),A_\infty(\cx,\cd)) \\
\rTTo^{1\boxtimes M}
TsA_\infty(\cy,\cc)\boxtimes TsA_\infty^u(A_\infty(\cy,\cc),A_\infty(\cx,\cd))
\rTTo^\alpha TsA_\infty(\cx,\cd) \bigr].
\end{multline*}
Using the definition of $M$ \cite[diagram~(4.0.1)]{Lyu-AinfCat} we
transform this equation into
\begin{multline*}
\bigl[
TsA_\infty(\cy,\cc)\boxtimes TsA_\infty^u(\cc,\cd)\boxtimes TsA_\infty(\cx,\cy)
\rTTo^{1\boxtimes c}
TsA_\infty(\cy,\cc)\boxtimes TsA_\infty(\cx,\cy)\boxtimes TsA_\infty^u(\cc,\cd)
\\
\rTTo^{1\boxtimes A_\infty(\_,\cc)\boxtimes A_\infty(\cx,\_)}
TsA_\infty(\cy,\cc)\boxtimes TsA_\infty^u(A_\infty(\cy,\cc),A_\infty(\cx,\cc))
\boxtimes TsA_\infty^u(A_\infty(\cx,\cc),A_\infty(\cx,\cd)) \\
\hspace{5em} \rTTo^{\alpha\boxtimes1}
TsA_\infty(\cx,\cc)\boxtimes TsA_\infty^u(A_\infty(\cx,\cc),A_\infty(\cx,\cd))
\rTTo^\alpha TsA_\infty(\cx,\cd) \bigr] \\
= \bigl[
TsA_\infty(\cy,\cc)\boxtimes TsA_\infty^u(\cc,\cd)\boxtimes TsA_\infty(\cx,\cy)
\rTTo^{1\boxtimes A_\infty(\cy,\_)\boxtimes A_\infty(\_,\cd)} \hspace{9em} \\
TsA_\infty(\cy,\cc)\boxtimes TsA_\infty^u(A_\infty(\cy,\cc),A_\infty(\cy,\cd))
\boxtimes TsA_\infty^u(A_\infty(\cy,\cd),A_\infty(\cx,\cd)) \\
\rTTo^{\alpha\boxtimes1}
TsA_\infty(\cy,\cd)\boxtimes TsA_\infty^u(A_\infty(\cy,\cd),A_\infty(\cx,\cd))
\rTTo^\alpha TsA_\infty(\cx,\cd) \bigr].
\end{multline*}
Using the definitions of \(A_\infty(\_,\cc)\) and \(A_\infty(\cy,\_)\)
\cite[(6.1.2)]{Lyu-AinfCat} we rewrite this equation as follows:
\begin{multline*}
\bigl[
TsA_\infty(\cy,\cc)\boxtimes TsA_\infty^u(\cc,\cd)\boxtimes TsA_\infty(\cx,\cy)
\rTTo^{c_{(123)}}
TsA_\infty(\cx,\cy)\boxtimes TsA_\infty(\cy,\cc)\boxtimes TsA_\infty^u(\cc,\cd)
\\
\rTTo^{M\boxtimes A_\infty(\cx,\_)}
TsA_\infty(\cx,\cc)\boxtimes TsA_\infty^u(A_\infty(\cx,\cc),A_\infty(\cx,\cd))
\rTTo^\alpha TsA_\infty(\cx,\cd) \bigr] \\
= \bigl[
TsA_\infty(\cy,\cc)\boxtimes TsA_\infty^u(\cc,\cd)\boxtimes TsA_\infty(\cx,\cy)
\rTTo^{M\boxtimes A_\infty(\_,\cd)} \hspace{9em} \\
TsA_\infty(\cy,\cd)\boxtimes TsA_\infty^u(A_\infty(\cy,\cd),A_\infty(\cx,\cd))
\rTTo^\alpha TsA_\infty(\cx,\cd) \bigr].
\end{multline*}
Now we use definitions of \(A_\infty(\cx,\_)\) and \(A_\infty(\_,\cd)\)
to get an equivalent form of the required equation:
\begin{multline*}
\bigl[
TsA_\infty(\cy,\cc)\boxtimes TsA_\infty^u(\cc,\cd)\boxtimes TsA_\infty(\cx,\cy)
\rTTo^{c_{(123)}}
TsA_\infty(\cx,\cy)\boxtimes TsA_\infty(\cy,\cc)\boxtimes TsA_\infty^u(\cc,\cd)
\\
\rTTo^{M\boxtimes1} TsA_\infty(\cx,\cc)\boxtimes TsA_\infty^u(\cc,\cd)
\rTTo^M TsA_\infty(\cx,\cd) \bigr] \\
= \bigl[
TsA_\infty(\cy,\cc)\boxtimes TsA_\infty^u(\cc,\cd)\boxtimes TsA_\infty(\cx,\cy)
\rTTo^{M\boxtimes1} TsA_\infty(\cy,\cd)\boxtimes TsA_\infty(\cx,\cy) \\
\rTTo^c TsA_\infty(\cx,\cy)\boxtimes TsA_\infty(\cy,\cd)
\rTTo^M TsA_\infty(\cx,\cd) \bigr].
\end{multline*}
This equation holds true due to associativity of $M$, since $M$ has
degree 0.

Other statements are similar or follow from the already proven one.
\end{proof}

\subsection{\texorpdfstring{An $A_\infty^u$-2-functor}
 {An A8u-2-functor}}
Let \ainf-category $\ca$ be unital. The \ainf-category $\ce$ is
pseudounital with distinguished elements equal to the unit elements of
$\cc$.

Strict \ainfu-2-functors are defined in
\cite[Definition~3.1]{LyuMan-freeAinf}. There is a strict
\ainfu-2-functor $F$, given by the following data:
\begin{enumerate}
\item the mapping of objects $F:\Ob A_{\infty}^u\to\Ob A_{\infty}^u$,
\( \ca\mapsto F\ca = A_\infty(\cc,\ca)_{\modulo\cb} \)
($F\ca$ is a full subcategory of the unital \ainf-category
$A_{\infty}(\cc,\ca)$, hence it is unital as well);

\item the strict unital \ainf-functor
 $F=F_{\ca_1,\ca_2}:A_{\infty}^u(\ca_1,\ca_2)\to
 A_{\infty}^u(F\ca_1,F\ca_2)$
for each pair of unital \ainf-categories $\ca_1$, $\ca_2$ given as
follows:
\[ \Ob F: g\mapsto gF = (1\boxtimes g)M|_{F\ca_1},
\]
where $(1\boxtimes g)M:A_\infty(\cc,\ca_1)\to A_\infty(\cc,\ca_2)$.
Indeed, if $\cb \rMono \cc \rTTo^{f} \ca_1$ is a contractible
\ainf-functor, then so is $\cb \rMono \cc \rTTo^{fg} \ca_2$. Actually,
if ${}_X\uni_0f_1=w_Xb_1$ for some $w_X\in(s\ca_1)^{-2}(Xf,Xf)$,
$X\in\Ob\cb$, then ${}_X\uni_0f_1g_1=w_Xb_1g_1=(w_Xg_1)b_1$, where
$w_Xg_1\in(s\ca_2)^{-2}(Xfg,Xfg)$. Furthermore,
\begin{align*}
F_1:A_\infty^u(\ca_1,\ca_2)(g,h) &\to A_\infty^u(F\ca_1,F\ca_2)(gF,hF),
\\
(r:g\to h:\ca_1\to\ca_2) &\mapsto rF_1 = (1\boxtimes r)M|_{F\ca_1},
\end{align*}
or more precisely,
\end{enumerate}
\begin{multline*}
\left[sF\ca_1(f_0,f_1)\tens\dots\tens sF\ca_1(f_{n-1},f_n)
\rTTo^{[rF_1]_n} sF\ca_2(f_0g,f_nh)\right] \\
= \left[sA_{\infty}(\cc,\ca_1)(f_0,f_1)\tens\dots\tens
sA_{\infty}(\cc,\ca_1)(f_{n-1},f_n) \rTTo^{(1\boxtimes r)M_{n1}}
sA_{\infty}(\cc,\ca_2)(f_0g,f_nh)
\right].
\end{multline*}

The necessary equations for $F_1$ are consequences of those for
$A_{\infty}(\cc,\_)$ \cite[Proposition~6.2]{Lyu-AinfCat}. Clearly, $F$
is a unital \ainf-functor. Let us check that the following diagram
commutes:
\begin{diagram}
Ts A_\infty^u(\ca_0,\ca_1)\boxtimes Ts A_\infty^u(\ca_1,\ca_2) &
\rTTo^M & Ts A_\infty^u(\ca_0,\ca_2) \\
\dTTo<{F_{\ca_0,\ca_1}\tens F_{\ca_1,\ca_2}} & = &
\dTTo>{F_{\ca_0,\ca_2}} \\
Ts A_\infty^u(F\ca_0,F\ca_1)\boxtimes Ts A_\infty^u(F\ca_1,F\ca_2) &
\rTTo^M & Ts A_\infty^u(F\ca_0,F\ca_2).
\end{diagram}
It follows from a similar diagram for $A_{\infty}(\cc,\_)$ in place of
$F$ (equation~(3.3.1) of \cite{LyuMan-freeAinf}). The commutativity is
clear on objects; since both sides of the required identity are
cocategory homomorphisms it suffices to show that
\begin{multline*}
(F_{\ca_0,\ca_1}\boxtimes F_{\ca_1,\ca_2})M\pr_1
= MF_{\ca_0,\ca_1}\pr_1: \\
TsA_\infty^u(\ca_0,\ca_1)\boxtimes Ts A_\infty^u(\ca_1,\ca_2)\to
sA_\infty^u(F\ca_0,F\ca_2).
\end{multline*}
Since $F_{-,-}$ are strict \ainf-functors, we must show that for any
non-negative integers $n$, $m$
\begin{multline*}
((F_{\ca_0,\ca_1})_1^{\tens n}\boxtimes
(F_{\ca_1,\ca_2})_1^{\tens m})M_{nm} = M_{nm}(F_{\ca_0,\ca_2})_1: \\
T^nsA_\infty^u(\ca_0,\ca_1)\boxtimes T^ms A_\infty^u(\ca_1,\ca_2)
\to sA_\infty^u(F\ca_0,F\ca_2).
\end{multline*}
Since $M_{nm}$ vanishes whenever $m>1$, we restrict our attention to
$m=0$ and $m=1$. These cases are similar and we will give verification
in the case $m=1$. Given diagrams of \ainf-functors and
\ainf-transformations
\begin{align*}
f^0\rTTo^{r^1}f^1\rTTo\dots\rTTo^{r^n}f^n:\ca_0\to\ca_1, \qquad
g^0\rTTo^{t^1}g^1:\ca_1\to\ca_2,
\end{align*}
we must show that
 $(1\boxtimes(r^1\tens\dots\tens r^n\boxtimes t^1)M)M|_{F\ca_0}
 =((1\boxtimes r^1)M|_{F\ca_0}\tens\dots\tens(1\boxtimes r^n)M|_{F\ca_0}
 \boxtimes(1\boxtimes t^1)M|_{F\ca_1})M$.
This is a particular case of equation~(3.3.1) of \cite{LyuMan-freeAinf}
since it coincides with a similar equation for $A_{\infty}(\cc,\_)$ in
place of $F$.

Let $\cd$ be the \ainf-category defined in \secref{sec-constr-D}. We
claim that $\restr:A_{\infty}(\cd,\ca)\to F\ca$ is an
\ainfu-2-transformation as defined in
\cite[Definition~3.2]{LyuMan-freeAinf}. The strict \ainf-functor
$\restr$ is induced by the inclusion $\iota:\cc \rMono \cd$:
\begin{gather*}
f\mapsto (\cc \rMono \cd \rTTo^{f} \ca) = f|_{\cc}, \\
\restr_1: A_{\infty}(\cd,\ca)(f,g) \to
F\ca(f|_{\cc},g|_{\cc}) = A_{\infty}(\cc,\ca)(f|_{\cc},g|_{\cc}).
\end{gather*}
The restriction $f|_{\cb}$ is a contractible \ainf-functor, for
${}_X\uni_0f_1=\eps_Xb_1f_1=(\eps_Xf_1)b_1$ for $X\in\cb$. Let us check
that the following diagram of \ainf-functors commutes:
\begin{diagram}
A_\infty^u(\ca_1,\ca_2) & \rTTo^{A_\infty(\cd,-)} &
A_\infty^u(A_\infty(\cd,\ca_1),A_\infty(\cd,\ca_2)) \\
\dTTo<{F} & = & \dTTo>{(1\boxtimes\restr_{\ca_2})M} \\
A_\infty^u(F\ca_1,F\ca_2) & \rTTo^{(\restr_{\ca_1}\boxtimes1)M} &
A_\infty^u(A_\infty(\cd,\ca_1),F\ca_2)
\end{diagram}
All functors in the diagram above are strict (the proof is given in
\cite[Section~3.4]{LyuMan-freeAinf}). We must verify the equation
\[ A_{\infty}(\cd,-)(1\boxtimes\restr_{\ca_2})M
= F(\restr_{\ca_1}\boxtimes1)M.
\]
On objects: given a unital \ainf-functor $g:\ca_1\to\ca_2$, we are
going to check that
\[ [(1\boxtimes g)M]_n\cdot\restr_1
= \restr^{\tens n}_1\cdot[(1\boxtimes g)M]_n
\]
for any $n\ge 1$. Indeed, for any $n$\n-tuple of composable
\ainf-transformations
\[ f^0 \rTTo^{r^1} f^1 \rTTo \dots \rTTo^{r^n} f^n:\cd\to\ca_1
\]
we have
\begin{align*}
\{(r^1\tens\dots\tens r^n)[(1\boxtimes g)M]_n\}_k\big|_\cc
&= [(r^1\tdt r^n\boxtimes g)M_{n0}]_k\big|_\cc \\
&= \bigl\{\sum_{l}(r^1\tens\dots\tens r^n)\theta_{kl}g_l\bigr\}
\big|_\cc \\
&= \sum_{l}(r^1|_\cc\tens\dots\tens r^n|_\cc)\theta_{kl}g_l|_\cc \\
&= \{(r^1|_\cc\tens\dots\tens r^n|_\cc)[(1\boxtimes g|_\cc)M]_n\}_k.
\end{align*}

The coincidence of \ainf-transformations
$((1\boxtimes t)M)\cdot\restr_\cd=\restr_\cc\cdot((1\boxtimes t)M)$
follows similarly from the computation:
\begin{align*}
\{(r^1\tens\dots\tens r^n)[(1\boxtimes t)M]_n\}_k\big|_\cc
&= [(r^1\tdt r^n\boxtimes t)M_{n1}]_k\big|_\cc \\
&= \bigl\{\sum_{l}(r^1\tens\dots\tens r^n)\theta_{kl}t_l\bigr\}
\big|_\cc \\
&= \sum_{l}(r^1|_\cc\tens\dots\tens r^n|_\cc)\theta_{kl}t_l|_\cc \\
&= \{(r^1|_\cc\tens\dots\tens r^n|_\cc)[(1\boxtimes t|_\cc)M]_n\}_k.
\end{align*}

We are more interested in the following \ainfu-2-subfunctor of $F$,
denoted $G:\Ob A_{\infty}^u\to\Ob A_{\infty}^u$,
\(\ca\mapsto G\ca=A_\infty^u(\cc,\ca)_{\modulo\cb}\subset F\ca\). All
the structure data of $G$ are restrictions of those of $F$. Hence, the
conclusion that $\restr:A_\infty^{\psi u}(\cd,\ca)\to G\ca$ is an
\ainfu-2-natural transformation remains valid. We prove in
\thmref{thm-restriction-DA-CAmodB-equi} that this \ainf-functor is an
equivalence. Therefore, this restriction \ainf-functor is a 2\n-natural
\ainf-equivalence. If $\cc$ is strictly unital, then
$\cd=\Quo(\cc|\cb)$ is unital by \thmref{thm-C-strict-D-unital}.
Whenever $\cd$ is unital, we say that $\cd$ \emph{unitally} represents
$G$.

We are going to discuss how an \ainfu-2-functor represented by an
\ainf-category $\cx$ depends on $\cx$.

\begin{proposition}\label{pro-(f1)M-Ainf2transformation}
Let \(f:\cx\to\cy\) be an \ainf-functor. Then
\[ \lambda_\ca = (f\boxtimes1)M : A_\infty(\cy,\ca) \to
A_\infty(\cx,\ca), \qquad g \mapsto fg
\]
is a strict \ainfu-2-transformation between two \ainfu-2-functors of
\(\ca\in\Ob A_\infty^u\). If $f:\cx\to\cy$ is a unital \ainf-functor,
then
\[ (f\boxtimes1)M : A_\infty^u(\cy,\ca) \to A_\infty^u(\cx,\ca)
\]
is also a strict \ainfu-2-transformation.
\end{proposition}

\begin{proof}
The \ainf-functor \((f\boxtimes1)M\) strictly commutes with the unit
transformations \((1\boxtimes\uni^\ca)M\) in \(A_\infty(\cy,\ca)\) and
\(A_\infty(\cx,\ca)\):
\[ (f\boxtimes1_\ca)M \cdot (1_\cx\boxtimes\uni^\ca)M
= (1_\cy\boxtimes\uni^\ca)M \cdot (f\boxtimes1_\ca)M
\]
due to associativity of $M$. Therefore, \((f\boxtimes1)M\) is unital.

We have to prove that the diagram of \ainf-functors
\begin{diagram}[LaTeXeqno]
A_\infty^u(\cc,\cd) & \rTTo^{A_\infty(\cy,\_)}
& A_\infty^u(A_\infty(\cy,\cc),A_\infty(\cy,\cd)) \\
\dTTo<{A_\infty(\cx,\_)} && \dTTo>{(1\boxtimes\lambda_\cd)M} \\
A_\infty^u(A_\infty(\cx,\cc),A_\infty(\cx,\cd))
& \rTTo^{(\lambda_\cc\boxtimes1)M}
& A_\infty^u(A_\infty(\cy,\cc),A_\infty(\cx,\cd))
\label{dia-ainf-2-transf-A()}
\end{diagram}
commutes. All four \ainf-functors in this diagram are strict. So it
suffices to check commutativity on objects and for the first
components.

If \(h:\cc\to\cd\) is a unital \ainf-functor, then
\[ h.A_\infty(\cy,\_)(1\boxtimes\lambda_\cd)M
= (1\boxtimes h)M\cdot(f\boxtimes1)M
= (f\boxtimes1\boxtimes h)(1\boxtimes M)M
\]
equals to
\[ h.A_\infty(\cx,\_)(\lambda_\cc\boxtimes1)M
= (f\boxtimes1)M\cdot(1\boxtimes h)M
= (f\boxtimes1\boxtimes h)(M\boxtimes1)M
\]
due to associativity of $M$.

If \(t\in sA_\infty^u(\cc,\cd)(g,h)\), then it is mapped by the first
components of \ainf-functors in diagram~\eqref{dia-ainf-2-transf-A()}
to
\[ t.A_\infty(\cy,\_)_1(1\boxtimes\lambda_\cd)M_{10}
= (1\boxtimes t)M\cdot(f\boxtimes1)M
= (f\boxtimes1\boxtimes t)(1\boxtimes M)M
\]
and
\[ t.A_\infty(\cx,\_)_1(\lambda_\cc\boxtimes1)M_{01}
= (f\boxtimes1)M\cdot(1\boxtimes t)M
= (f\boxtimes1\boxtimes t)(M\boxtimes1)M.
\]
These expressions are equal due to associativity of $M$.

The case of unital $f$ and \ainfu-2-subfunctors \(A_\infty^u(\cy,\ca)\),
\(A_\infty^u(\cx,\ca)\) follows from the general case.
\end{proof}

With two strict \ainfu-2-transformations $\lambda$, $\mu$ in
\(F \rTTo^\lambda G \rTTo^\mu H:A_\infty^u\to A_\infty^u\) is
associated the third \ainfu-2-transformation
\(\lambda\mu:F\to H:A_\infty^u\to A_\infty^u\) -- their composition,
specified by the family of unital \ainf-functors
\[ (\lambda\mu)_\cc = \lambda_\cc\mu_\cc : F\cc \to H\cc,
\qquad \cc \in \Ob A_\infty^u.
\]
In order to verify that \(\lambda\mu\) is indeed a strict
\ainfu-2-transformation, we have to check equation~(3.2.1) of
\cite{LyuMan-freeAinf}:
\[ F\cdot(1\boxtimes(\lambda\mu)_\cd)M
= H\cdot((\lambda\mu)_\cc\boxtimes1)M:
A_\infty^u(\cc,\cd) \to A_\infty^u(F\cc,H\cd).
\]
We do it as follows:
\begin{align*}
\bigl[ A_\infty^u(\cc,\cd) &\rTTo^F A_\infty^u(F\cc,F\cd)
\rTTo^{(1\boxtimes\lambda_\cd\mu_\cd)M}
_{(1\boxtimes\lambda_\cd)M\cdot(1\boxtimes\mu_\cd)M}
A_\infty^u(F\cc,H\cd) \bigr] \\
= \bigl[ A_\infty^u(\cc,\cd) &\rTTo^G A_\infty^u(G\cc,G\cd)
\rTTo^{(\lambda_\cc\boxtimes1)M\cdot(1\boxtimes\mu_\cd)M}
_{(1\boxtimes\mu_\cd)M\cdot(\lambda_\cc\boxtimes1)M}
A_\infty^u(F\cc,H\cd) \bigr] \\
= \bigl[ A_\infty^u(\cc,\cd) &\rTTo^H A_\infty^u(H\cc,H\cd)
\rTTo^{(\mu_\cc\boxtimes1)M\cdot(\lambda_\cc\boxtimes1)M}
_{(\lambda_\cc\mu_\cc\boxtimes1)M}
A_\infty^u(F\cc,H\cd) \bigr].
\end{align*}

\begin{definition}
A strict \emph{modification}
$m:\lambda\to\mu:F\to G:A_\infty^u\to A_\infty^u$ of strict
$A_\infty^u$-2-transformations $\lambda$, $\mu$ is
\begin{enumerate}
\item a family of \ainf-transformations
$m_\cc:\lambda_\cc\to\mu_\cc:F\cc\to G\cc$ for $\cc\in\Ob A_\infty^u$
\newline
such that

\item for any pair of unital \ainf-categories $\cc$, $\cd$ the
\ainf-transformations
\end{enumerate}
\[ (F_{\cc,\cd}\boxtimes\lambda_\cd)\cdot M
= (\lambda_\cc\boxtimes G_{\cc,\cd})\cdot M
\rTTo^{(m_\cc\boxtimes G_{\cc,\cd})\cdot M}
(\mu_\cc\boxtimes G_{\cc,\cd})\cdot M:
A_\infty^u(\cc,\cd) \to A_\infty^u(F\cc,G\cd),
\]
\[ (F_{\cc,\cd}\boxtimes\lambda_\cd)\cdot M
\rTTo^{(F_{\cc,\cd}\boxtimes m_\cd)\cdot M}
(F_{\cc,\cd}\boxtimes\mu_\cd)\cdot M
= (\mu_\cc\boxtimes G_{\cc,\cd})\cdot M:
A_\infty^u(\cc,\cd) \to A_\infty^u(F\cc,G\cd)
\]
\qquad are equal, in short,
\[ G_{\cc,\cd}(m_\cc\boxtimes1)M = F_{\cc,\cd}(1\boxtimes m_\cd)M:
F_{\cc,\cd}(1\boxtimes\lambda_\cd)M\to G_{\cc,\cd}(\mu_\cc\boxtimes1)M.
\]
A modification $m$ is \emph{natural} if all \ainf-transformations
$m_\cc$ are natural.
\end{definition}

\subsection{Example of strict modification}
We claim that for an arbitrary \ainf-transformation
\(r:f\to g:\cx\to\cy\) the family of \ainf-transformations
\[ m_\cc = (r\boxtimes1)M:
\lambda_\cc = (f\boxtimes1)M \to \mu_\cc = (g\boxtimes1)M:
F\cc = A_\infty(\cy,\cc) \to G\cc = A_\infty(\cx,\cc)
\]
is a strict modification. In order to prove it, we notice first of all
that \((r\boxtimes1)M\) is a
\(((f\boxtimes1)M,(g\boxtimes1)M)\)-coderivation for an arbitrary
$\cc$, since $M$ is a cocategory homomorphism. When we want to indicate
the category $\cc$, we write this coderivation as
 \((r\boxtimes1_\cc)M\in
 sA_\infty^u(A_\infty(\cy,\cc),A_\infty(\cx,\cc))
 ((f\boxtimes1_\cc)M,(g\boxtimes1_\cc)M)\).
The equation to verify is
\begin{multline*}
\bigl[ TsA_\infty^u(\cc,\cd)
\rTTo^{(r\boxtimes1_\cc)M\boxtimes A_\infty(\cx,\_)}
T^1sA_\infty^u(A_\infty(\cy,\cc),A_\infty(\cx,\cc))\boxtimes
TsA_\infty^u(A_\infty(\cx,\cc),A_\infty(\cx,\cd)) \\
\rTTo^M TsA_\infty^u(A_\infty(\cy,\cc),A_\infty(\cx,\cd)) \bigr] \\
= \bigl[ TsA_\infty^u(\cc,\cd)
\rTTo^{A_\infty(\cy,\_)\boxtimes(r\boxtimes1_\cd)M}
TsA_\infty^u(A_\infty(\cy,\cc),A_\infty(\cy,\cd))\boxtimes
T^1sA_\infty^u(A_\infty(\cy,\cd),A_\infty(\cx,\cd)) \\
\rTTo^M TsA_\infty^u(A_\infty(\cy,\cc),A_\infty(\cx,\cd)) \bigr].
\end{multline*}
It is an immediate consequence of
equation~\eqref{eq-A(-C)-A(X-)-A(Y-)-A(-D)} restricted to the element
\(r\in T^1sA_\infty(\cx,\cy)\).

\begin{lemma}\label{lem-mB1-modif}
If $m:\lambda\to\mu:F\to G:A_\infty^u\to A_\infty^u$ is a strict
modification, then so is
$mB_1:\lambda\to\mu:F\to G:A_\infty^u\to A_\infty^u$, where
$(mB_1)_\cc=m_\cc B_1:\lambda_\cc\to\mu_\cc:F\cc\to G\cc$ for all
$\cc\in\Ob A_\infty^u$. The $\kk$\n-linear map
\[ sA_\infty(\cx,\cy)(f,g) \ni r \mapsto (r\boxtimes1)M \in
sA_\infty^u(A_\infty(\cy,\cc),A_\infty(\cx,\cc))
((f\boxtimes1)M,(g\boxtimes1)M)
\]
is a chain map.
\end{lemma}

\begin{proof}
Let us prove that the family \(m_\cc B_1\) constitutes a strict
modification. The identity \((1\boxtimes B+B\boxtimes1)M=MB\) implies
that
\begin{multline*}
(F_{\cc,\cd}\boxtimes m_\cd.B_1)M = (F_{\cc,\cd}\boxtimes m_\cd)MB
-(-)^mB(F_{\cc,\cd}\boxtimes m_\cd)M
= [(F_{\cc,\cd}\boxtimes m_\cd)M].B_1 \\
= [(m_\cc\boxtimes G_{\cc,\cd})M].B_1 = (m_\cc\boxtimes G_{\cc,\cd})MB
-(-)^mB(m_\cc\boxtimes G_{\cc,\cd})M = (m_\cc.B_1\boxtimes G_{\cc,\cd})M.
\end{multline*}
Here we use the fact that \(m_\cd.B=m_\cd.B_1\) due to
\(m_\cd\in T^1sA_\infty(F\cd,G\cd)(\lambda_\cd,\mu_\cd)\), and,
similarly, \(m_\cc.B=m_\cc.B_1\). The equation
\begin{equation}
[(r\boxtimes1)M].B_1 = (r\boxtimes1)MB -(-)^rB(r\boxtimes1)M
= (r.B\boxtimes1)M = (r.B_1\boxtimes1)M
\label{eq-r1MB1-rB11M}
\end{equation}
proves that \(r\mapsto(r\boxtimes1)M\) is a chain map.
\end{proof}

\begin{lemma}
Let $m$, $n$ in
\(\lambda \rTTo^m \mu \rTTo^n \nu:F\to G:A_\infty^u\to A_\infty^u\) be
strict modifications. Then
\begin{multline*}
F_{\cc,\cd}[1\boxtimes(m_\cd\tens n_\cd)B_2]M
= G_{\cc,\cd}[(m_\cc\tens n_\cc)B_2\boxtimes1]M \\
+G_{\cc,\cd}[(m_\cc\tens n_\cc)(1\tens B_1+B_1\tens1)A_\infty(\_,G\cd)_2]
- G_{\cc,\cd}[(m_\cc\tens n_\cc)A_\infty(\_,G\cd)_2B_1]
\end{multline*}
\end{lemma}

\begin{proof}
Since \(A_\infty(\_,G\cd)\) is an \ainf-functor we have
\begin{align}
&[ (m_\cc\boxtimes1)M\tens(n_\cc\boxtimes1)M]B_2
+ (m_\cc\tens n_\cc).A_\infty(\_,G\cd)_2B_1 \notag \\
&= [m_\cc.A_\infty(\_,G\cd)_1\tens n_\cc.A_\infty(\_,G\cd)_1]B_2
+ (m_\cc\tens n_\cc).A_\infty(\_,G\cd)_2B_1 \notag \\
&= (m_\cc\tens n_\cc)B_2.A_\infty(\_,G\cd)_1
+ (m_\cc\tens n_\cc)(1\tens B_1+B_1\tens1)A_\infty(\_,G\cd)_2 \notag \\
&= [(m_\cc\tens n_\cc)B_2\boxtimes1]M
+ (m_\cc\tens n_\cc)(1\tens B_1+B_1\tens1)A_\infty(\_,G\cd)_2.
\label{eq-Ainf(-GD)2}
\end{align}
Using this identity we find
\begin{multline*}
F_{\cc,\cd}[1\boxtimes(m_\cd\tens n_\cd)B_2]M
= F_{\cc,\cd}[(1\boxtimes m_\cd)M\tens(1\boxtimes n_\cd)M]B_2 \\
= [F_{\cc,\cd}(1\boxtimes m_\cd)M\tens F_{\cc,\cd}(1\boxtimes n_\cd)M]B_2
= [G_{\cc,\cd}(m_\cc\boxtimes1)M\tens G_{\cc,\cd}(n_\cc\boxtimes1)M]B_2
\\
= G_{\cc,\cd}[(m_\cc\boxtimes1)M\tens(n_\cc\boxtimes1)M]B_2
= G_{\cc,\cd}[(m_\cc\tens n_\cc)B_2\boxtimes1]M \hspace{7.8em} \\
+G_{\cc,\cd}[(m_\cc\tens n_\cc)(1\tens B_1+B_1\tens1)A_\infty(\_,G\cd)_2]
- G_{\cc,\cd}[(m_\cc\tens n_\cc)A_\infty(\_,G\cd)_2B_1],
\end{multline*}
so the lemma is proven.
\end{proof}

\begin{proposition}
If unital \ainf-categories $\cx$, $\cy$ are equivalent, then
\ainfu-2-functors \(\ca\mapsto A_\infty^u(\cx,\ca)\) and
\(\ca\mapsto A_\infty^u(\cy,\ca)\) are naturally \ainfu-2-equivalent.
\end{proposition}

\begin{proof}
Let \(\phi:\cx\to\cy\), \(\psi:\cy\to\cx\) be \ainf-equivalences,
quasi-inverse to each other. Then there are natural
\ainf-transformations
\[ r: \phi\psi \to \id_\cx: \cx \to \cx, \qquad
p: \id_\cx \to \phi\psi: \cx \to \cx,
\]
inverse to each other, that is,
\[ (r\tens p)B_2 \equiv \phi\psi\uni^\cx,
\qquad (p\tens r)B_2 \equiv \uni^\cx.
\]
The \ainf-transformations
\begin{align}
(r\boxtimes1)M &: (\psi\boxtimes1)M(\phi\boxtimes1)M \to \id:
A_\infty^u(\cx,\ca) \to A_\infty^u(\cx,\ca),
\label{eq-(r1)M-psi1Mphi1M-id} \\
(p\boxtimes1)M &: \id \to (\psi\boxtimes1)M(\phi\boxtimes1)M:
A_\infty^u(\cx,\ca) \to A_\infty^u(\cx,\ca)
\label{eq-(p1)M-id-psi1Mphi1M}
\end{align}
are also natural by \eqref{eq-r1MB1-rB11M}. We are going to prove that
these \ainf-transformations are inverse to each other.

The natural \ainf-transformation
\[ (1\boxtimes\uni^\ca)M: \id \to \id:
A_\infty^u(\cx,\ca) \to A_\infty^u(\cx,\ca)
\]
is a unit transformation of \ainf-category \(A_\infty^u(\cx,\ca)\) by
Proposition~7.7 of \cite{Lyu-AinfCat}. Another unit transformation is
given by
\begin{equation}
(\uni^\cx\boxtimes1)M: \id \to \id:
A_\infty^u(\cx,\ca) \to A_\infty^u(\cx,\ca).
\label{eq-iX1M-id-Au(XA)}
\end{equation}
Indeed, values of 0\n-th components of the both natural transformations
on an object $f$ of \(A_\infty^u(\cx,\ca)\) differ by a boundary since
\[ \sS{_f}[(\uni^\cx\boxtimes1)M]_0 = (\uni^\cx\boxtimes f)M_{10}
= \uni^\cx f \equiv f\uni^\ca = (f\boxtimes\uni^\ca)M_{01} =
\sS{_f}[(1\boxtimes\uni^\ca)M]_0.
\]
Therefore, \([1\tens(\uni^\cx\boxtimes1)M_{10}]B_2\) and
\([(\uni^\cx\boxtimes1)M_{10}\tens1]B_2\) are homotopy invertible.

Furthermore, by \eqref{eq-Ainf(-GD)2}
\[ [(\uni^\cx\boxtimes1)M\tens(\uni^\cx\boxtimes1)M]B_2 \equiv
[(\uni^\cx\tens\uni^\cx)B_2\boxtimes1]M \equiv (\uni^\cx\boxtimes1)M
\]
due to \lemref{lem-mB1-modif}. Therefore, homotopy
idempotent~\eqref{eq-iX1M-id-Au(XA)} is a unit transformation of
\(A_\infty^u(\cx,\ca)\) by \cite[Definition~7.6]{Lyu-AinfCat}. Since
unit transformation is unique up to equivalence by
\cite[Corollary~7.10]{Lyu-AinfCat} we have
\[ (\uni^\cx\boxtimes1)M \equiv (1\boxtimes\uni^\ca)M: \id \to \id:
A_\infty^u(\cx,\ca) \to A_\infty^u(\cx,\ca).
\]
Composing natural \ainf-transformations \((r\boxtimes1)M\) and
\((p\boxtimes1)M\) given by \eqref{eq-(r1)M-psi1Mphi1M-id} and
\eqref{eq-(p1)M-id-psi1Mphi1M} we get
\begin{multline*}
[(r\boxtimes1)M\tens(p\boxtimes1)M]B_2 \equiv
[(r\tens p)B_2\boxtimes1]M \equiv (\phi\psi\uni^\cx\boxtimes1)M \\
\equiv (\uni^\cx\phi\psi\boxtimes1)M
= (\psi\boxtimes1)M(\phi\boxtimes1)M(\uni^\cx\boxtimes1)M,
\end{multline*}
\[ [(p\boxtimes1)M\tens(r\boxtimes1)M]B_2 \equiv
[(p\tens r)B_2\boxtimes1]M \equiv (\uni^\cx\boxtimes1)M
\]
by \eqref{eq-Ainf(-GD)2} and \lemref{lem-mB1-modif}. Since
\((\uni^\cx\boxtimes1)M\) is a unit transformation, the
\ainf-transformations \((r\boxtimes1)M\) and \((p\boxtimes1)M\) are
inverse to each other.

The obtained statement together with one more statement in which $\phi$
and $\psi$ exchange their places implies that \ainf-functors
\begin{align*}
A_\infty^u(\phi,\ca) = (\phi\boxtimes1)M &:
A_\infty^u(\cy,\ca) \to A_\infty^u(\cx,\ca), \\
A_\infty^u(\psi,\ca) = (\psi\boxtimes1)M &:
A_\infty^u(\cx,\ca) \to A_\infty^u(\cy,\ca)
\end{align*}
are quasi-inverse to each other.

They form strict \ainfu-2-transformations by
\propref{pro-(f1)M-Ainf2transformation}. Therefore,
\(A_\infty^u(\phi,\ca)\) and \(A_\infty^u(\psi,\ca)\) are natural
\ainfu-2-equivalences.
\end{proof}

Given a pair \((\cc,\cb)\) consisting of a unital \ainf-category $\cc$
and its full subcategory $\cb$, we shall construct another pair
\((\wt{\cc},\wt{\cb})\) consisting of a differential graded category
\(\wt{\cc}\) and its full subcategory \(\wt{\cb}\) as follows. Set
\(\wt{\cc}\) to be the differential graded category
\(\wt{\Rep}A_\infty^u(\cc^\op,\uCom)\) of \ainf-functors, represented
by objects of $\cc$, see \remref{rem-Yoneda-tilde}. Thus,
\(\Ob\wt{\cc}=\Ob\cc\). The category \(\wt{\cc}\) is equivalent to
$\cc$, the Yoneda \ainf-equivalence \(\wt{Y}:\cc\to\wt{\cc}\) and its
quasi-inverse \(\Psi:\wt{\cc}\to\cc\) are identity on objects by
\remref{rem-Yoneda-tilde}. We take \(\wt{\cb}\) to be the full
subcategory of \(\wt{\cc}\) with the set of objects
\(\Ob\wt{\cb}=\Ob\cb\). Therefore, the \ainf-functors \(\wt{Y}\) and
$\Psi$ can be restricted to quasi-inverse to each other
\ainf-equivalences \(Y'=\wt{Y}\big|_\cb:\cb\to\wt{\cb}\) and
\(\Psi'=\Psi\big|_{\wt{\cb}}:\wt{\cb}\to\cb\).

\begin{corollary}\label{cor-Psi1M-Y1M-A-equivalences}
Let $\ca$ be a unital \ainf-category. The \ainf-functors
\begin{align*}
(\Psi\boxtimes1)M: A^u_\infty(\cc,\ca) &\to
A^u_\infty(\wt{\cc},\ca), \qquad f \mapsto \Psi f, \\
(\wt{Y}\boxtimes1)M: A^u_\infty(\wt{\cc},\ca) &\to
A^u_\infty(\cc,\ca), \qquad g \mapsto \wt{Y}g,
\end{align*}
are  quasi-inverse to each other \ainf-equivalences. The first maps
objects of \(A^u_\infty(\cc,\ca)_{\modulo\cb}\) to objects of
\(A^u_\infty(\wt{\cc},\ca)_{\modulo\wt{\cb}}\), the second does vice
versa. Therefore, their restrictions determine quasi-inverse to each
other \ainf-equivalences
\begin{align*}
(\Psi\boxtimes1)M: A^u_\infty(\cc,\ca)_{\modulo\cb} &\to
A^u_\infty(\wt{\cc},\ca)_{\modulo\wt{\cb}},
\\
(\wt{Y}\boxtimes1)M: A^u_\infty(\wt{\cc},\ca)_{\modulo\wt{\cb}} &\to
A^u_\infty(\cc,\ca)_{\modulo\cb}.
\end{align*}
\end{corollary}

\begin{proof}
Let \(f:\cc\to\ca\) be a unital \ainf-functor such that \(f\big|_\cb\)
is contractible. Then for each object $X$ of $\cb$ the complex
\((s\ca(Xf,Xf),b_1)\) is contractible
\cite[Proposition~6.1]{LyuOvs-iResAiFn}. Equivalently we may say that
for each object $Z$ of $\wt{\cb}$ the complex
 \((s\ca(Z\Psi f,Z\Psi f),b_1)\) is contractible, as the following
commutative diagram shows:
\begin{diagram}
\wt{\cb} & \rMono & \wt{\cc} & \rTTo^{\Psi f} & \ca \\
\dTTo<{\Psi'} && \dTTo<\Psi && \dEq \\
\cb & \rMono & \cc & \rTTo^f & \ca
\end{diagram}
This implies contractibility of the \ainf-functor
\(\Psi f\big|_{\wt{\cb}}\). Similarly, if
\(g\in A^u_\infty(\wt{\cc},\ca)_{\modulo\wt{\cb}}\), then
\(\wt{Y}g\in A^u_\infty(\cc,\ca)_{\modulo\cb}\).
\end{proof}

\begin{corollary}\label{cor-equi-pairs-equi-ainf2funs}
The \ainfu-2-functors \(\ca\mapsto A^u_\infty(\cc,\ca)_{\modulo\cb}\)
and \(\ca\mapsto A^u_\infty(\wt{\cc},\ca)_{\modulo\wt{\cb}}\) are
naturally \ainfu-2-equivalent. Therefore, if one of them is
representable, then so is the other.
\end{corollary}

\section{The Yoneda Lemma for 2-categories and bicategories}
    \label{ap-2-Yoneda}
In this article we deal with bicategories of a particular kind --
strict 2\n-categories. However, 2\n-functors and their transformations
need to be weak for our purposes.

\subsection{2-categories}
We recall the definitions of strict 2\n-categories and associated weak
notions originating in \cite{Benabou}. We use the form and the notation
of \cite{Lyu-sqHopf}.

\begin{definition}
A (strict) \emph{2-category} $\fA$ consists of
\begin{enumerate}
\item a set of objects $\Ob\fA$;
\item for any pair of objects $X,Y\in\Ob\fA$ a category $\fA(X,Y)$;
\item
\begin{enumerate}
\renewcommand{\labelenumii}{(\alph{enumii})}
\item for any object $X\in\Ob\fA$ an object $1_X$ of $\fA(X,X)$;
\item for any triple of objects $X,Y,Z\in\Ob\fA$ a functor
\[ \bull:\fA(X,Y)\times\fA(Y,Z) \to \fA(X,Z),\qquad (F,G)\mapsto
FG = F\bull G = G\circ F; \]
\end{enumerate}
such that the following functors are equal
\item $F\bull1=F=1\bull F$, $F(G H)=(F G) H$.
\end{enumerate}
\end{definition}

The 2\n-category of ($\fu$\n-small) categories is denoted $\Cat$.

\begin{definition}\label{def-weak-2-func}
A \emph{weak 2-functor} (a homomorphism in \cite{Benabou}) between
2\n-categories $\fA$ and $\fc$ consists of
\begin{enumerate}
\item a function $F:\Ob\fA\to\Ob\fc$;
\item a functor $F=F_{X,Y}:\fA(X,Y)\to\fc(FX,FY)$ for each pair of
objects $X,Y\in\Ob\fA$;
\item
\begin{enumerate}
\renewcommand{\labelenumii}{(\alph{enumii})}
\item an isomorphism $\phi_0:1_{FX}\to F1_X$;
\item an invertible (natural) transformation
\begin{diagram}
\fA(X,Y)\times\fA(Y,Z) & \rTTo^{\bull} & \fA(X,Z) \\
\dTTo<{F_{X,Y}\times F_{Y,Z}} & \ruTwoar^{\phi_2} & \dTTo>{F_{X,Z}} \\
\fc(FX,FY)\times\fc(FY,FZ) & \rTTo^{\bull} & \fc(FX,FZ)
\end{diagram}
for each triple $X,Y,Z\in\Ob\fA$;
\end{enumerate}
such that
\item
\begin{enumerate}
\renewcommand{\labelenumii}{(\alph{enumii})}
\item for any object $M\in\fA(X,Y)$ the composites
\begin{gather}
FM = FM\bull1_{FY} \rTTo^{FM\bull\phi_0} FM\bull F1_Y \rTTo^{\phi_2}
F(M\bull1_Y) = FM
\label{eq-w-2-fun-a1} \\
FM = 1_{FX}\bull FM \rTTo^{\phi_0\bull FM} F1_X\bull FM \rTTo^{\phi_2}
F(1_X\bull M) = FM
\label{eq-w-2-fun-a2}
\end{gather}
are identity morphisms in $\fc(FM,FM)$;
\item For any objects $W,X,Y,Z\in\Ob\fA$ and any object
\[ (K,L,M)\in\fA(W,X)\times\fA(X,Y)\times\fA(Y,Z) \]
there is an equation
\begin{multline*}
\bigl(FK\bull(FL\bull FM) \rTTo^{FK\bull\phi_2} FK\bull F(L\bull M)
\rTTo^{\phi_2} F(K\bull(L\bull M))\bigr) \\
= \bigl((FK\bull FL)\bull FM \rTTo^{\phi_2\bull FM} F(K\bull L)\bull FM
\rTTo^{\phi_2} F((K\bull L)\bull M)\bigr).
\end{multline*}
\end{enumerate}
\end{enumerate}
If $\phi_2$ and $\phi_0$ are identity isomorphisms, $F$ is called a
\emph{strict 2\n-functor}.
\end{definition}

\begin{definition}\label{def-weak-2-trans}
A \emph{weak 2-transformation} (pseudo-natural transformation
\cite{Gray})
$\lambda:(F,\phi_2,\phi_0)\to(G,\psi_2,\psi_0):\fA\to\fc$ is
\begin{enumerate}
\item a family of 1-morphisms $\lambda_X:FX\to GX$, $X\in\Ob\fA$;
\item for any 1-morphism $f:X\to Y$ in $\fA$ a 2\n-isomorphism in $\fc$
\[ \lambda_f: Ff\bull \lambda_Y \rTTo^\sim \lambda_X\bull Gf: FX\to GY,
\]
which is an isomorphism of functors
\[ \lambda_-:F-\bull\lambda_Y \to \lambda_X\bull G-: \fA(X,Y)\to
\fc(FX,GY), \]
that is, for any 2-morphism $\xi:f\to g:X\to Y$
\[
\begin{diagram}
FX & \rTTo^{Ff} & FY \\
& \ldTwoar_{\lambda_f} & \\
\dTTo<{\lambda_X} && \dTTo>{\lambda_Y} \\
GX &\pile{\rTTo^{Gf}\\ \Downarrow \scriptstyle G\xi\\ \rTTo_{Gg}} & GY
\end{diagram}
=
\begin{diagram}
FX & \pile{\rTTo^{Ff}\\ \Downarrow \scriptstyle F\xi\\ \rTTo_{Fg} }
& FY \\
&& \\
\dTTo<{\lambda_X} & \ldTwoar^{\lambda_g} & \dTTo>{\lambda_Y} \\
GX &\rTTo_{Gg} & GY
\end{diagram},
\]
such that
\item
\begin{enumerate}
\renewcommand{\labelenumii}{(\alph{enumii})}
\item for any object $X\in\Ob\fA$
\[
\begin{diagram}[inline]
FX & \pile{\rTTo^{1_{FX}}\\ \Downarrow\phi_0\\ \rTTo_{F1_X}} & FX \\
\dTTo<{\lambda_X} && \dTTo>{\lambda_X} \\
& \ldTwoar^{\lambda_{1_X}} & \\
GX & \rTTo_{G1_X} & GX
\end{diagram}
\quad = \quad
\begin{diagram}[inline]
FX & \rTTo^{1_{FX}} & FX \\
\dTTo<{\lambda_X} & = & \dTTo>{\lambda_X} \\
GX & \pile{\rTTo^{1_{GX}}\\ \Downarrow\psi_0\\ \rTTo_{G1_X}} & GX
\end{diagram}
\quad;
\]
\item for any pair of composable 1\n-morphisms $f,g\in\fA_1$
\begin{equation*}
\begin{diagram}[inline]
&& FY && \\
& \ruTTo^{Ff} & \dTwoar>{\phi_2} & \rdTTo^{Fg} & \\
FX && \rTTo_{F(f g)} && FZ \\
\dTTo<{\lambda_X} &&& \ldTwoar(4,2)_{\lambda_{f g}} &
\dTTo>{\lambda_Z} \\
GX && \rTTo_{G(f g)} && GZ
\end{diagram}
\quad = \quad
\begin{diagram}[inline]
&& FY && \\
& \ruTTo^{Ff} \ldTwoar(2,4)[nohug]^{\lambda_f} & \dTTo>{\lambda_Y} &
\rdTTo^{Fg} & \\
FX && GY & \lTwoar^{\lambda_g} & FZ \\
\dTTo<{\lambda_X} & \ruTTo_{Gf} & \dTwoar>{\psi_2} & \rdTTo_{Gg} &
\dTTo>{\lambda_Z} \\
GX && \rTTo_{G(f g)} && GZ
\end{diagram}
\ .
\end{equation*}
\end{enumerate}
\end{enumerate}
If $\lambda_f$ are identity isomorphisms, $\lambda$ is called a
\emph{strict 2\n-transformation}. A weak 2\n-transformation
$\lambda=(\lambda_X)$ for which $\lambda_X$ are equivalences is called
a \emph{2\n-natural equivalence}.
\end{definition}

\begin{definition}\label{def-modif}
A \emph{modification} $m:\lambda\to\mu:F\to G:\fA\to\fc$ is
\begin{enumerate}
\item a family of 2-morphisms $m_X:\lambda_X\to\mu_X$, $X\in\Ob\fA$
\newline
such that
\item for any 1-morphism $f:X\to Y$ in $\fA$
\[
\begin{diagram}
FX & \rTTo^{Ff} & FY \\
\dTTo<{\mu_X} \overset{\sss m_X}{\Leftarrow} \dTTo>{\lambda_X} &
\ldTwoar^{\lambda_f} & \dTTo>{\lambda_Y} \\
GX & \rTTo^{Gf} & GY
\end{diagram}
\quad = \quad
\begin{diagram}
FX & \rTTo^{Ff} & FY \\
\dTTo<{\mu_X} & \ldTwoar^{\mu_f} & \dTTo<{\mu_Y}
\overset{\sss m_Y}{\Leftarrow} \dTTo>{\lambda_Y}\\
GX & \rTTo^{Gf} & GY
\end{diagram}
\]
\end{enumerate}
\end{definition}

\begin{proposition}[Invertibility of 2-natural equivalences]
 \label{pro-2-natural-quasi-inverse}
Let $\lambda:F\to G:\fA\to\fc$ be a 2-natural equivalence. Then there
exist a weak 2-transformation $\mu:G\to F:\fA\to\fc$ and invertible
modifications $\eps:\lambda\mu\to 1_F:F\to F:\fA\to\fc$ and
$\eta:1_G\to\mu\lambda:G\to G:\fA\to\fc$. Thus, $\mu$ is quasi-inverse
to $\lambda$.
\end{proposition}

\ifx\chooseClass2
The proof is left to the reader as an exercise.
    \else
\begin{proof}
Since $\lambda_X:FX\to GX$ is an equivalence for every $X\in\Ob\fA$, we
obtain: for every $X\in\Ob\fA$ there exist a 1-morphism
$\mu_X:GX\to FX$ and invertible 2-morphisms
$\eps_X:\lambda_X\mu_X\to 1_{FX}:FX\to FX$,
$\beta_X:1_{GX}\to\mu_X\lambda_X:GX\to GX$.

\begin{NOlemma}
There exist such invertible 2\n-morphisms
$\eta_X:1_{GX}\to\mu_X\lambda_X:GX\to GX$ that the following equations
hold true:
\begin{gather}
\bigl(\mu_X \rTTo^{\eta_X\bull\mu_X} \mu_X\lambda_X\mu_X
\rTTo^{\mu_X\bull\eps_X} \mu_X\bigr) = 1_{\mu_X}, \label{adj1} \\
\bigl(\lambda_X \rTTo^{\lambda_X\bull\eta_X} \lambda_X\mu_X\lambda_X
\rTTo^{\eps_X\bull\lambda_X} \lambda_X\bigr) = 1_{\lambda_X}.
\label{adj2}
\end{gather}
\end{NOlemma}

\begin{proof}
Consider the following functors:
\begin{align}
\fc(GX,GX) &\longrightarrow \fc(GX,FX),
& \fc(GX,GX) &\longrightarrow \fc(FX,GX), \label{eq-C(GXGX)-fun} \\
f &\longmapsto f\mu_X, & f &\longmapsto \lambda_Xf, \notag \\
\phi:f\to g &\longmapsto \phi\bull\mu_X:f\mu_X\to g\mu_X, &
\phi:f\to g &\longmapsto \lambda_X\bull\phi:\lambda_Xf\to\lambda_Xg,
\label{eq-phi-fg2} \\
\fc(GX,FX) &\longrightarrow \fc(GX,GX),
& \fc(FX,GX) &\longrightarrow \fc(GX,GX), \notag \\
h &\longmapsto h\lambda_X, & h &\longmapsto \mu_Xh, \notag \\
\chi:h\to k &\longmapsto \chi\bull\lambda_X:h\lambda_X\to k\lambda_X, &
\chi:h\to k &\longmapsto \mu_X\bull\chi:\mu_Xh\to \mu_Xk. \notag
\end{align}
These functors are faithful. Let us prove it for the first one. Indeed,
if $\phi\bull\mu_X=\psi\bull\mu_X:f\mu_X\to g\mu_X$, then
 $\phi\bull\mu_X\bull\lambda_X=\psi\bull\mu_X\bull\lambda_X:
 f\mu_X\lambda_X\to g\mu_X\lambda_X$
and
 $(f\bull\beta_X)(\phi\bull\mu_X\bull\lambda_X)
 =(f\bull\beta_X)(\psi\bull\mu_X\bull\lambda_X):f\to g\mu_X\lambda_X$,
i.e., $\phi\bull(g\beta_X)=\psi\bull(g\beta_X):f\to g\mu_X\lambda_X$,
hence $\phi=\psi$, because $\beta_X$ is invertible. Similarly or by
symmetry the other 3 functors are also faithful.

Functors~\eqref{eq-C(GXGX)-fun} are full. Let us prove it for the first
one. Given $\psi:f\mu_X\to g\mu_X$, we set
$\phi=(f\bull\beta)(\psi\bull\lambda_X)(g\bull\beta^{-1}_X):f\to g$.
Then
 $(f\bull\beta_X)(\phi\bull\mu_X\bull\lambda_X)(g\bull\beta^{-1}_X)
 =(f\bull\beta_X)(f\bull\beta^{-1}_X)\phi=\phi
 =(f\bull\beta_X)(\psi\bull\lambda_X)(g\bull\beta^{-1}_X)$.
This yields $\phi\bull\mu_X\bull\lambda_X=\psi\bull\lambda_X$. Since
multiplication with $\lambda_X$ is a faithful functor, we obtain
$\phi\bull\mu_X=\psi$. Notice that if $\psi$ is a 2\n-isomorphism, then
the composition of 2\n-isomorphisms $\phi$ is a 2\n-isomorphism as
well. Similarly we prove that the second functor of
\eqref{eq-C(GXGX)-fun} is full.

Now let us take $f=1_{GX}$, $g=\mu_X\lambda_X$. By the first bijection
of \eqref{eq-phi-fg2} we find a 2\n-isomorphism
$\eta_X:1_{GX}\to\mu_X\lambda_X$ corresponding to the 2\n-isomorphism
$\psi=\mu_X\bull\eps^{-1}_X:\mu_X\to\mu_X\lambda_X\mu_X$. Then
$\eta_X\bull\mu_X=\mu_X\bull\eps^{-1}_X$, that is,
$(\eta_X\bull\mu_X)(\mu_X\bull\eps_X)=1_{\mu_X}$. By the second
bijection of \eqref{eq-phi-fg2} we find a 2\n-isomorphism
$\gamma_X:1_{GX}\to\mu_X\lambda_X$ corresponding to the 2\n-isomorphism
$\xi=\eps^{-1}_X\bull\lambda_X:\lambda_X\to\lambda_X\mu_X\lambda_X$.
Then \(\lambda_X\bull\gamma_X=\eps^{-1}_X\bull\lambda_X\), that is,
\((\lambda_X\bull\gamma_X)(\eps_X\bull\lambda_X)=1_{\lambda_X}\). We
have
\begin{multline*}
\gamma = \bigl(1_{GX} \rTTo^\gamma \mu\lambda \rTTo^{1_\mu\lambda}
\mu\lambda\bigr)
= \bigl(1_{GX} \rTTo^\gamma \mu\lambda \rTTo^{\eta\mu\lambda}
\mu\lambda\mu\lambda \rTTo^{\mu\eps\lambda} \mu\lambda\bigr) \\
= \bigl(1_{GX} \rTTo^\eta \mu\lambda \rTTo^{\mu\lambda\gamma}
\mu\lambda\mu\lambda \rTTo^{\mu\eps\lambda} \mu\lambda\bigr)
= \bigl(1_{GX} \rTTo^\eta \mu\lambda \rTTo^{\mu1_\lambda}
\mu\lambda\bigr) = \eta.
\end{multline*}
Therefore, the 2\n-isomorphism $\eta$ fulfills both equations
\eqref{adj1} and \eqref{adj2}.
\end{proof}

For any 1-morphism $f:X\to Y$ in $\fA$ we define a 2-isomorphism
$\mu_f$ as the pasting of
\[
\begin{diagram}[tight,height=1em,width=1.5em]
GX & & \rTTo^{\mu_X} & & FX & & \rTTo^{Ff} & &  FY & & & & \\
   &\rdTTo(4,4)_{1_{GX}} &       &\ruTwoar_{\eta_X} & & &
&\ruTwoar(4,4)_{\lambda^{-1}_f} & &\rdTTo(4,4)^{1_{FY}}  & & & \\
   & & & & \dTTo>{\lambda_X} & & & & \dTTo<{\lambda_Y} & &  & & \\
   & & & & & & & & & \ruTwoar^{\eps_Y} & & & \\
   & & & & GX & & \rTTo_{Gf} & & GY & &\rTTo_{\mu_Y} & & FY
\end{diagram}.
\]
Since $\lambda_f$ determine an isomorphism of functors
\[
\lambda_{-}:F-\bull\lambda_Y\to\lambda_X\bull G-:
\fA(X,Y)\to\fc(FX,GY),
\]
for any 2-morphism $\xi:f\to g:X\to Y$ equation 2 of
\defref{def-weak-2-trans} holds. It implies
\[
\begin{diagram}[tight,height=1em,width=1.4em]
GX & & \rTTo^{\mu_X} & & FX & &
\pile{\rTTo^{Fg}\\ \Uparrow\scriptstyle F\xi \\ \rTTo_{Ff}}
& &  FY & & & & \\
   &\rdTTo(4,4)_{1_{GX}} &       &\ruTwoar_{\eta_X} & & &
&\ruTwoar(4,4)_{\lambda^{-1}_f} & &\rdTTo(4,4)^{1_{FX}}  & & & \\
   & & & & \dTTo>{\lambda_X} & & & & \dTTo<{\lambda_Y} & &  & & \\
   & & & & & & & & & \ruTwoar^{\eps_Y} & & & \\
   & & & & GX & & \rTTo_{Gf} & & GY & &\rTTo_{\mu_Y} & & FY
\end{diagram}
=
\begin{diagram}[tight,height=1em,width=1.4em]
GX & & \rTTo^{\mu_X} & & FX & & \rTTo^{Fg} & &  FY & & & & \\
   &\rdTTo(4,4)_{1_{GX}} &       &\ruTwoar_{\eta_X} & & &
&\ruTwoar(4,4)^{\lambda^{-1}_g} & &\rdTTo(4,4)^{1_{FX}}  & & & \\
   & & & & \dTTo>{\lambda_X} & & & & \dTTo<{\lambda_Y} & &  & & \\
   & & & & & & & & & \ruTwoar^{\eps_Y} & & & \\
   & & & & GX & &
\pile{\rTTo^{Gg}\\ \Uparrow \scriptstyle G\xi\\ \rTTo_{Gf}}
& & GY & &\rTTo_{\mu_Y} & & FY
\end{diagram}.
\]
This shows that the collection of 2-isomorphisms $\mu_f$ determines an
isomorphism of functors
\[
\mu_{-}:G-\bull\mu_Y\to\mu_X\bull F-: \fA(X,Y)\to\fc(GX, FY).
\]

Let us check conditions 3(a),(b) of \defref{def-weak-2-trans} for
$\mu$. For 3(a) we have
\[
\begin{diagram}[tight,height=1em,width=1.4em]
GX & & \rTTo^{\mu_X} & & FX & & \rTTo^{F1_X} & &  FX & & & & \\
   &\rdTTo(4,4)_{1_{GX}} &       &\ruTwoar_{\eta_X} & & &
&\ruTwoar(4,4)^{\lambda^{-1}_{1_X}} & &\rdTTo(4,4)^{1_{FX}}  & & & \\
   & & & & \dTTo>{\lambda_X} & & & & \dTTo<{\lambda_X} & &  & & \\
   & & & & & & & & & \ruTwoar^{\eps_X} & & & \\
   & & & & GX & &
\pile{\rTTo^{G1_X}\\ \Uparrow \scriptstyle \psi_0\\ \rTTo_{1_{GX}}}
& & GX & &\rTTo_{\mu_X} & & FX
\end{diagram}
=
\begin{diagram}[tight,height=1em,width=1.4em]
GX & & \rTTo^{\mu_X} & & FX & &
\pile{\rTTo^{F1_X}\\ \Uparrow\scriptstyle \phi_0 \\ \rTTo_{1_{FX}}}
& & FX & & & & \\
   &\rdTTo(4,4)_{1_{GX}} &       &\ruTwoar_{\eta_X} & & &
& & &\rdTTo(4,4)^{1_{FX}}  & & & \\
   & & & & \dTTo>{\lambda_X} & &= & & \dTTo<{\lambda_X} & &  & & \\
   & & & & & & & & & \ruTwoar^{\eps_X} & & & \\
   & & & & GX & & \rTTo_{1_{GX}} & & GX & &\rTTo_{\mu_X} & & FX
\end{diagram}
\]
and the required equation follows from~\eqref{adj1}.

Equation~\eqref{adj2} and the corresponding property for $\lambda$
imply
\begin{gather*}
\begin{diagram}[tight,height=1em,width=1.5em]
&&&& &&&& GY &&&& &&&&\\
&&&& &&&\ruTTo(4,4)^{Gf}&    &\rdTTo(4,4)^{Gg}&&& &&&&\\
&&&& &&&& \dTwoar>{\psi_2} &&&& &&&&\\
&&&& &&&& &&&& &&&&\\
&&&& GX &&&& \rTTo_{G(f g)} &&&& GZ && \rTTo^{\mu_Z} && FZ\\
&&&\ruTTo(4,4)^{1_{GX}}& &\rdTwoar(8,4)_{\lambda^{-1}_{f g}}&&& &&&&
&\rdTwoar_{\eps_Z} &&\ruTTo(4,4)_{1_{FZ}}&\\
&&&& \uTTo>{\lambda_X} &&&& &&&& \uTTo<{\lambda_Z} &&&&\\
&&&\rdTwoar^{\eta_X}& &&&& &&&& &&&&\\
GX && \rTTo_{\mu_X} && FX &&&& \rTTo_{F(f g)} &&&& FZ &&&&
\end{diagram}\\
=
\begin{diagram}[tight,height=1em,width=1.5em]
&&&&&&&& GY &&&& \rTTo^{1_{GY}} &&&& GY &&&&&&&& \\
&&&&&&&\ruTTo(4,4)^{Gf}&&\rdTwoar(4,4)_{\eps_Y}\rdTTo(8,4)^{\mu_Y}&&&&\rdTwoar(4,4)^{
\eta_Y} &&&
&\rdTwoar(4,8)^{\lambda^{-1}_g}\rdTTo(4,4)^{Gg}&&&&&&& \\
&&&&&&&& \uTTo<{\lambda_Y} &&&& &&&& \uTTo>{\lambda_Y} &&&&&&&&\\
&&&&&&&& &&&& &&&& &&&&&&&&\\
&&&& GX && \rTwoar^{\lambda^{-1}_f} && FY &&&&&&&\rTTo_{1_{FY}}& FY &&&& GZ &&
\rTTo^{\mu_Z} && FZ\\
&&& \ruTTo(4,4)^{1_{GX}}& &&&\ruTTo(4,4)_{Ff}& &&&& \phantom{} &&&&
&\rdTTo(4,4)_{Fg}&&& &\rdTwoar^{\eps_Z} &&
\ruTTo(4,4)_{1_{FZ}}& \\
&&&& \uTTo>{\lambda_X} &&&& &&&& \dTwoar>{\phi_2} &&&&
&&&&\uTTo>{\lambda_Z} &&&&\\
&&&\rdTwoar^{\eta_X}& &&&& &&&& &&&& &&&& &&&&\\
GX &&\rTTo_{\mu_X} && FX &&&& &&&& \rTTo_{F(f g)} &&&& &&&& FZ &&&&
\end{diagram}
\end{gather*}
and the assertion 3(b) for $\mu$ follows.

Let us verify that collections $\eps=(\eps_X)$ and $\eta=(\eta_X)$
determine modifications. Equation~\eqref{adj2} implies that
\[
\begin{diagram}[tight,height=2em,width=1.7em]
FX &&& \rTTo^{Ff} &&& FY &&&& \\
& \rdTTo_{\lambda_X} \rdTTo(6,2)^{\lambda_X} & = &&&&
\dTwoar<{\lambda_f} & \rdTTo(4,2)^{\lambda_Y} &&& \\
\dTTo<{1_{FX}} & \lTwoar^{\eps_X} & GX && \rTTo^{1_{GX}} && GX &&
\rTTo^{Gf} && GY \\
& \ldTTo^{\mu_X} && \ldTwoar(4,2)^{\eta_X} && \ruTTo(6,2)_{\lambda_X} &
\dTwoar<{\lambda_f^{-1}} &&& \ruTTo(4,2)^{\lambda_Y} \ldTwoar>{\eps_Y}
& \dTTo>{\mu_Y} \\
FX &&& \rTTo_{Ff} &&& FY && \rTTo_{1_{FY}} && FY
\end{diagram}
\quad=\qquad
\begin{diagram}[tight,height=2em,width=2.3em]
FX & \rTTo^{Ff} & FY && \\
&&& \rdTTo^{\lambda_Y} & \\
\dTTo<{1_{FX}} & = & \dTTo<{1_{FY}} & \lTwoar^{\eps_Y} & GY \\
&&& \ldTTo_{\mu_Y} & \\
FX & \rTTo_{Ff} & FY &&
\end{diagram}.
\]
This means that \(\eps:\lambda\mu\to1_F:F\to F:\fA\to\fc\) is a
modification.

Also equation~\eqref{adj2} implies that
\[
\begin{diagram}[tight,height=2em,width=2.3em]
&& GX & \rTTo^{Gf} & GY \\
& \ldTTo^{\mu_X} &&& \\
FX & \lTwoar^{\eta_X} & \dTTo>{1_{GX}} & = & \dTTo>{1_{GY}} \\
& \rdTTo_{\lambda_X} &&& \\
&& GX & \rTTo_{Gf} & GY
\end{diagram}
\qquad=\quad
\begin{diagram}[tight,height=2em,width=1.7em]
GX && \rTTo^{1_{GX}} && GX &&& \rTTo^{Gf} &&& GY \\
\dTTo<{\mu_X} & \ldTwoar<{\eta_X} && \ruTTo(4,2)_{\lambda_X} &
\dTwoar>{\lambda_f^{-1}} &&&&& \ruTTo(6,2)^{\lambda_Y}
\ldTwoar(4,2)>{\eps_Y} \ldTTo>{\mu_Y} & \\
FX && \rTTo_{Ff} && FY && \rTTo_{1_{FY}} && FY & \lTwoar^{\eta_Y} &
\dTTo>{1_{GX}} \\
& \rdTTo(4,2)_{\lambda_X} &&& \dTwoar>{\lambda_f} &
\rdTTo(6,2)_{\lambda_Y} &&& = & \rdTTo^{\lambda_Y} & \\
&&&& GX &&& \rTTo_{Gf} &&& GY
\end{diagram}
\]
This means that $\eta:1_G\to\mu\lambda:G\to G:\fA\to\fc$ is a
modification.

Therefore, $\mu$ is a weak 2-transformation quasi-inverse to $\lambda$.
\end{proof}
\fi

\begin{remark}\label{rem-quasi-inverse-2-natural-equivalences}
Clearly, if $\lambda:F\to G:\fA\to\fc$ and $\mu:G\to F:\fA\to\fc$ are
weak 2\n-transformations, quasi-inverse to each other, then both are
2-natural equivalences.
\end{remark}

We shall use the generalization of the classical Yoneda Lemma to
2\n-categories. If we were using strict 2\n-functors and strict
2\n-transformations, we would view 2\n-categories as $cat$-categories,
where $cat$ is the category of categories. This would allow to use one
of the Yoneda structures on 2\n-categories defined by Street and
Walters~\cite[Example~7(1)]{Street:12}, as well as weak Yoneda Lemma
for enriched categories by Eilenberg and
Kelly~\cite[Theorem~I.8.6]{EK},
Kelly~\cite[Section~1.9]{KellyGM:bascec} and strong Yoneda Lemma for
enriched categories by Kelly~\cite[Section~2.4]{KellyGM:bascec}.
However, we need weak 2\n-transformations (and modifications), so we
use the Yoneda Lemma for bicategories obtained by
Street~\cite[(1.9)]{Street:14}. Let us recall the latter statement.

\subsection{\texorpdfstring{The 2-functor $\fA(A,\_)$}
 {The 2-functor A(A,-)}}
Let $\fA$ be a strict 2-category. An arbitrary object $A\in\Ob\fA$
gives rise to a strict 2-functor $\fA(A,\_):\fA\to\Cat$. It is
specified by the following data:
\begin{enumerate}
\item{the function $\Ob\fA\to\Ob\Cat$, $X\mapsto\fA(A,X)$;}
\item{the functor $\fA(A,\_)_{XY}:\fA(X,Y)\to\Cat(\fA(A,X),\fA(A,Y))$
for each pair of objects $X,Y\in\Ob\fA$.}
\end{enumerate}
The functor $\fA(A,\_)_{XY}$ is given as follows. For any 1-morphism
$f:X\to Y$ the functor $\fA(A,f):\fA(A,X)\to\fA(A,Y)$ is given by the
following formulas:
\begin{align*}
(\phi:A\to X) &\mapsto (\phi f:A\to Y), \\
(\pi:\phi\to\psi:A\to X) &\mapsto
(\pi\bull f:\phi f\to \psi f:A\to Y).
\end{align*}
For any 2-morphism $\alpha:f\to g:X\to Y$ the natural transformation
$\fA(A,\alpha):\fA(A,f)\to\fA(A,g):\fA(A,X)\to\fA(A,Y)$ is given
explicitly by its components:
\begin{align*}
(\phi)\fA(A,\alpha) = \phi\bull\alpha:
(\phi)\fA(A,f)=\phi f\to\phi g=(\phi)\fA(A,g),\qquad \phi\in\fA(A,X).
\end{align*}

Let \([\fA,\Cat]\) denote the strict 2\n-category of weak 2\n-functors
\(\fA\to\Cat\), their weak 2\n-transformations and their modifications,
see e.g. \cite[Appendix~A.1.5]{Lyu-sqHopf}.

\begin{lemma}[Yoneda Lemma for bicategories,
 Street~{\cite[(1.9)]{Street:14}}]\label{lem-Yoneda-Lemma-bicategories}
For a homomorphism \(G:\fA\to\Cat\) of bicategories, evaluation at the
identity for each object $A$ of $\fA$ provides the components
\([\fA,\Cat](\fA(A,\_),G)\to GA\) of an equivalence in \([\fA,\Cat]\).
\end{lemma}

We have not found a detailed published proof of the above result in the
existing literature, since it has to be quite lengthy. Curiously, part
of the required statements were formalized and verified by a computer
proof--checker \cite{MR1608340}. On the other hand, in the case of
strict 2-categories one can write down a complete proof in several
pages. For convenience of the reader we decompose it into several
detailed statements, written for a strict 2-category $\fA$, fixed till
the end of this section.
    \ifx\chooseClass2
Proofs of these statements are left to the reader as exercises.
    \fi

A weak 2\n-transformation \(\lambda:\fA(A,\_)\to G:\fA\to\Cat\)
involves, in particular, a functor \(\lambda_A:\fA(A,A)\to GA\).
Evaluating it on the object \(1_A\in\Ob\fA(A,A)\) we get an object
\((1_A)\lambda_A\in\Ob GA\). A modification
\(m:\lambda\to\mu:\fA(A,\_)\to G:\fA\to\Cat\) involves, in particular,
a natural transformation \(m_A:\lambda_A\to\mu_A:\fA(A,A)\to GA\).
Evaluating it on the object \(1_A\in\Ob\fA(A,A)\) we get a morphism
\((1_A)m_A:(1_A)\lambda_A\to(1_A)\mu_A\) of $GA$.

\begin{proposition}\label{pro-ev1A-full-and-faithful}
Let $G:\fA\to\Cat$ be a weak 2\n-functor. Let $A$ be an object of
$\fA$. Then the functor
\begin{align*}
\ev_{1_A}: [\fA,\Cat](\fA(A,\_),G) & \longrightarrow GA, \\
\lambda & \longmapsto (1_A)\lambda_A, \\
m:\lambda \to \mu & \longmapsto (1_A)m_A:(1_A)\lambda_A \to (1_A)\mu_A,
\end{align*}
is full and faithful.
\end{proposition}

\ifx\chooseClass2
    \else
\begin{proof}
Clearly, the assignment \(\ev_{1_A}\) gives a functor. Let us show that
it is faithful. Let two modifications
\(m,n:\lambda\to\mu:\fA(A,\_)\to G:\fA\to\Cat\) be given such that
\((1_A)m_A=(1_A)n_A\). The modification $m$ is a family of natural
transformations \((m_C)_{C\in\Ob\fA}\) satisfying the equation
\begin{equation}
\begin{diagram}
\fA(A,B) & \rTTo^{\fA(A,f)} & \fA(A,C) \\
\dTTo<{\mu_B} \overset{\sss m_B}{\Leftarrow} \dTTo>{\lambda_B} &
\ldTwoar^{\lambda_f} & \dTTo>{\lambda_C} \\
GB & \rTTo^{Gf} & GC
\end{diagram}
\quad = \quad
\begin{diagram}
\fA(A,B) & \rTTo^{\fA(A,f)} & \fA(A,C) \\
\dTTo<{\mu_B} & \ldTwoar^{\mu_f} & \dTTo<{\mu_C}
\overset{\sss m_C}{\Leftarrow} \dTTo>{\lambda_C} \\
GB & \rTTo^{Gf} & GC
\end{diagram}
\label{dia-eq-modif-A(A-)-G}
\end{equation}
for an arbitrary 1\n-morphism \(f:B\to C\) of $\fA$. In particular, it
holds for \(B=A\). Restrict this equation to the object $1_A$ of
\(\fA(A,A)\). Then it gives for an arbitrary 1\n-morphism \(f:A\to C\)
the equation
\begin{multline*}
\bigl[ (f)\lambda_C \rTTo^{(1_A)\lambda_f} ((1_A)\lambda_A)(Gf)
\rTTo^{((1_A)m_A)(Gf)} ((1_A)\mu_A)(Gf) \bigr] \\
= \bigl[ (f)\lambda_C \rTTo^{(f)m_C} (f)\mu_C
\rTTo^{(1_A)\mu_f} ((1_A)\mu_A)(Gf) \bigr].
\end{multline*}
Therefore, the value of $m_C$ on an arbitrary object $f$ of
\(\fA(A,C)\) is completely determined by the morphism \((1_A)m_A\):
\begin{equation*}
(f)m_C = \bigl[(f)\lambda_C \rTTo^{(1_A)\lambda_f} ((1_A)\lambda_A)(Gf)
\rTTo^{((1_A)m_A)(Gf)} ((1_A)\mu_A)(Gf) \rTTo^{(1_A)\mu_f^{-1}}
(f)\mu_C \bigr].
\end{equation*}
Thus, $m=n$ and \(\ev_{1_A}\) is faithful.

Let us prove that \(\ev_{1_A}\) is full. Let
\(\lambda,\mu:\fA(A,\_)\to G:\fA\to\Cat\) be weak 2\n-transformations.
Let \(\phi:(1_A)\lambda_A\to(1_A)\mu_A\) be a morphism of $GA$. We
claim that there is a modification $m:\lambda\to\mu$ such that
\((1_A)m_A=\phi\). The value of $m_C$ on an arbitrary 1\n-morphism
\(f:A\to C\) can be only
\begin{equation}
(f)m_C = \bigl[(f)\lambda_C \rTTo^{(1_A)\lambda_f} ((1_A)\lambda_A)(Gf)
\rTTo^{(\phi)(Gf)} ((1_A)\mu_A)(Gf) \rTTo^{(1_A)\mu_f^{-1}}
(f)\mu_C \bigr],
\label{eq-(f)mC-(phi)Gf}
\end{equation}
as we have seen. Let us verify that, indeed, this formula determines a
modification.

First of all, each $m_C$ is a natural transformation. Indeed, for each
2\n-morphism \(\xi:f\to g:A\to C\) of $\fA$ the following diagram
commutes:
\begin{diagram}[LaTeXeqno]
(f)\lambda_C & \rTTo^{(1_A)\lambda_f} & ((1_A)\lambda_A)(Gf)
& \rTTo^{(\phi)(Gf)} & ((1_A)\mu_A)(Gf) & \rTTo^{(1_A)\mu_f^{-1}}
& (f)\mu_C \\
\dTTo<{(\xi)\lambda_C} && \dTTo<{((1_A)\lambda_A)(G\xi)} &&
\dTTo>{((1_A)\mu_A)(G\xi)} && \dTTo>{(\xi)\mu_C} \\
(g)\lambda_C & \rTTo^{(1_A)\lambda_g} & ((1_A)\lambda_A)(Gg)
& \rTTo^{(\phi)(Gg)} & ((1_A)\mu_A)(Gg) & \rTTo^{(1_A)\mu_g^{-1}}
& (g)\mu_C
\label{dia-mC-nat-trans}
\end{diagram}
The central square commutes because \(G\xi:Gf\to Gg\) is a natural
transformation. Condition 2 for $\lambda$ from
\defref{def-weak-2-trans} implies, in particular, equation
\[
\begin{diagram}[width=4em]
\fA(A,A) & \rTTo^{\fA(A,f)} & \fA(A,C) \\
& \ldTwoar_{\lambda_f} & \\
\dTTo<{\lambda_A} && \dTTo>{\lambda_C} \\
GA &\pile{\rTTo^{Gf}\\ \Downarrow \scriptstyle G\xi\\ \rTTo_{Gg}} & GC
\end{diagram}
=
\begin{diagram}[width=4em]
\fA(A,A) &
\pile{\rTTo^{\fA(A,f)}\\ \Downarrow \scriptstyle
\fA(A,\xi)\\ \rTTo_{\fA(A,g)} }
& \fA(A,C) \\
&& \\
\dTTo<{\lambda_A} & \ldTwoar^{\lambda_g} & \dTTo>{\lambda_C} \\
GA &\rTTo_{Gg} & GC
\end{diagram}.
\]
Restricting this equation to the object $1_A$ of \(\fA(A,A)\) we get an
equation, which expresses precisely commutativity of the left square of
diagram~\eqref{dia-mC-nat-trans}. The right square of
\eqref{dia-mC-nat-trans} commutes by the same reasoning applied to
$\mu$ instead of $\lambda$. Thus, $m_C$ is a natural transformation.

Secondly, we have to prove equation~\eqref{dia-eq-modif-A(A-)-G} for
the family $(m_C)$ and for an arbitrary 1\n-morphism $f:B\to C$ of
$\fA$. On an arbitrary object \(g:A\to B\) of \(\fA(A,B)\) this
equation reads:
\begin{multline*}
\bigl[ (gf)\lambda_C \rTTo^{(g)\lambda_f} ((g)\lambda_B)(Gf)
\rTTo^{((g)m_B)(Gf)} ((g)\mu_B)(Gf) \bigr] \\
= \bigl[ (gf)\lambda_C \rTTo^{(gf)m_C} (gf)\mu_C
\rTTo^{(g)\mu_f} ((g)\mu_B)(Gf) \bigr].
\end{multline*}
Substituting definition~\eqref{eq-(f)mC-(phi)Gf} of $m_C$ we get an
equation, which expresses commutativity of the exterior of the
following diagram:
\begin{equation}
\begin{diagram}[objectstyle=\scriptstyle,nobalance,inline]
((g)\lambda_B)(Gf) & \rTTo^{((1_A)\lambda_g)(Gf)} &
((1_A)\lambda_A)(Gg)(Gf) & \rTTo^{(\phi)(Gg)(Gf)} &
((1_A)\mu_A)(Gg)(Gf) & \rTTo^{((1_A)\mu_g^{-1})(Gf)} &((g)\mu_B)(Gf) \\
\uTTo<{(g)\lambda_f} && \dTTo<{((1_A)\lambda_A)(g,f)\psi_2} &&
\dTTo>{((1_A)\mu_A)(g,f)\psi_2} && \uTTo>{(g)\mu_f} \\
(gf)\lambda_C & \rTTo^{(1_A)\lambda_{gf}} & ((1_A)\lambda_A)(G(gf))
& \rTTo^{(\phi)(G(gf))} & ((1_A)\mu_A)(G(gf))
& \rTTo^{(1_A)\mu_{gf}^{-1}} & (gf)\mu_C
\end{diagram}
\label{dia-3b-4mC-sss}
\end{equation}
The middle square commutes, because $(g,f)\psi_2:(Gg)(Gf)\to G(gf)$ is
a morphism of functors. Property 3(b) of \defref{def-weak-2-trans} for
$\lambda$ implies, in particular, the equation
\begin{equation*}
\begin{diagram}[inline]
&& \fA(A,B) && \\
& \ruTTo^{\fA(A,g)} & \dId>= & \rdTTo^{\fA(A,f)} & \\
\fA(A,A) && \rTTo_{\fA(A,gf)} && \fA(A,C) \\
\dTTo<{\lambda_A} &&& \ldTwoar(4,2)_{\lambda_{gf}} &
\dTTo>{\lambda_C} \\
GA && \rTTo_{G(gf)} && GC
\end{diagram}
\quad = \quad
\begin{diagram}[inline]
&& \fA(A,B) && \\
& \ruTTo^{\fA(A,g)} \ldTwoar(2,4)_{\lambda_g} & \dTTo>{\lambda_B} &
\rdTTo^{\fA(A,f)} & \\
\fA(A,A) && GB & \lTwoar^{\lambda_f} & \fA(A,C) \\
\dTTo<{\lambda_A} & \ruTTo_{Gg} & \dTwoar>{\psi_2} & \rdTTo_{Gf} &
\dTTo>{\lambda_C} \\
GA && \rTTo_{G(gf)} && GC
\end{diagram}
\ .
\end{equation*}
Restricting this equation to the object $1_A$ of $\fA(A,A)$, we will
get precisely the left square of diagram~\eqref{dia-3b-4mC-sss},
therefore, it commutes. The right square of \eqref{dia-3b-4mC-sss}
commutes by the same reasoning applied to $\mu$ instead of $\lambda$.
Therefore, $m$ is a modification, and \(\ev_{1_A}\) is full.
\end{proof}
\fi

\begin{proposition}\label{pro-functor-GA-(ACat)}
For each object $x$ of $GA$ there is a weak 2\n-transformation
\[ \lambda^x = \lambda^{A,x} = \sS{^G}\lambda^{A,x}:
\fA(A,\_) \to (G,\psi_2,\psi_0): \fA \to \Cat,
\]
specified by the family of functors $\lambda^x_C$, $C\in\Ob\fA$:
\begin{align*}
\lambda^x_C: \fA(A,C) & \longrightarrow GC, \\
f: A \to C & \longmapsto (x)(Gf), \\
\xi:f \to g:A \to C & \longmapsto (x)(G\xi):(x)(Gf) \to (x)(Gg),
\end{align*}
and by the family of invertible natural transformations
\[ \lambda^x_f: \fA(A,f)\bull\lambda^x_C \to \lambda^x_B\bull Gf:
\fA(A,B) \to GC,
\]
\(f\in\Ob\fA(B,C)\), which map an object \(g\in\Ob\fA(A,B)\) to the
isomorphism of $GC$:
\begin{equation}
(g)\lambda^x_f \overset{\text{def}}= (x)(g,f)\psi_2^{-1}:
(gf)\lambda^x_C = (x)(G(gf)) \to (x)(Gg)(Gf) = (g)\lambda^x_B(Gf).
\label{eq-(g)-lambda-x-f-def}
\end{equation}
For each morphism \(u:x\to y\) of $GA$ there is a modification
\[ \lambda^u = \lambda^{A,u} = \sS{^G}\lambda^{A,u}:
\lambda^x \to \lambda^y: \fA(A,\_) \to (G,\psi_2,\psi_0): \fA \to \Cat,
\]
specified by the family of natural transformations $\lambda^u_C$,
\(C\in\Ob\fA\):
\begin{gather}
\lambda^u_C: \lambda^x_C \to \lambda^y_C: \fA(A,C) \to GC, \notag \\
\lambda^u_C: (f:A\to C) \longmapsto \bigl((f)\lambda^u_C = (u)(Gf):
(f)\lambda^x_C = (x)(Gf) \to (y)(Gf) = (f)\lambda^y_C\bigr).
\label{eq-lambda-uC}
\end{gather}
The correspondence
\begin{align*}
\Lambda: GA & \longrightarrow [\fA,\Cat](\fA(A,\_),G), \\
x & \longmapsto \lambda^x, \\
u: x \to y & \longmapsto \lambda^u: \lambda^x \to \lambda^y,
\end{align*}
is a functor.
\end{proposition}

\ifx\chooseClass2
    \else
\begin{proof}
As $G:\fA(A,C)\to\Cat(GA,GC)$ is a functor, \(G1_f=1_{Gf}\) for the
unit 2\n-morphism $1_f:f\to f:A\to C$ of $\fA$, and for each pair of
composable 2\n-morphisms \(f \rTTo^\xi g \rTTo^\chi h:A\to C\) of $\fA$
we have
\[ G(\xi\chi) = \bigl(Gf \rTTo^{G\xi} Gg \rTTo^{G\chi} Gh\bigr).
\]
Evaluating these equations on $x$ we get
\((1_f)\lambda^x_C=(x)1_{Gf}=1_{(x)(Gf)}\) and
\((\xi\chi)\lambda^x_C=(\xi\lambda^x_C)(\chi\lambda^x_C)\), thus,
$\lambda^x_C$ is a functor.

We claim that $\lambda^x_f$ given by \eqref{eq-(g)-lambda-x-f-def} is a
natural transformation. Indeed, naturality of $\psi_2$, expressed by
\begin{diagram}
G(gf) & \lTTo^{(g,f)\psi_2} & (Gg)(Gf) \\
\dTTo<{G(\xi\bull f)} & = & \dTTo>{G\xi\bull Gf} \\
G(hf) & \lTTo^{(h,f)\psi_2} & (Gh)(Gf)
\end{diagram}
implies commutativity of
\begin{diagram}
(x)(G(gf)) & \rTTo^{(x)(g,f)\psi_2^{-1}} & (x)(Gg)(Gf) \\
\dTTo<{(x)(G(\xi\bull f))} & = & \dTTo>{(x)(G\xi)(Gf)} \\
(x)(G(hf)) & \rTTo^{(x)(h,f)\psi_2^{-1}} & (x)(Gh)(Gf)
\end{diagram}
which is nothing else, but naturality of $\lambda^x_f$:
\begin{diagram}
(g)\fA(A,f)\lambda^x_C & \rTTo^{(g)\lambda^x_f} & (g)\lambda^x_B(Gf) \\
\dTTo<{(\xi)\fA(A,f)\lambda^x_C} & = & \dTTo>{(\xi)\lambda^x_B(Gf)} \\
(h)\fA(A,f)\lambda^x_C & \rTTo^{(h)\lambda^x_f} & (h)\lambda^x_B(Gf)
\end{diagram}

We claim that
\[ \lambda^x_-: \fA(A,\_)\lambda^x_C \to \lambda^x_B\bull G-:
\fA(A,B) \to GC
\]
is a morphism of functors. That is, for each 2\n-morphism
\(\xi:f\to g:B\to C\) of $\fA$ the following equation holds:
\begin{equation}
\begin{diagram}[width=4em]
\fA(A,B) & \rTTo^{\fA(A,f)} & \fA(A,C) \\
& \ldTwoar_{\lambda^x_f} & \\
\dTTo<{\lambda^x_B} && \dTTo>{\lambda^x_C} \\
GB & \pile{\rTTo^{Gf}\\ \Downarrow \scriptstyle G\xi\\ \rTTo_{Gg}} & GC
\end{diagram}
\quad=\quad
\begin{diagram}[width=4em]
\fA(A,B) &
\pile{\rTTo^{\fA(A,f)}\\ \Downarrow
\scriptstyle \fA(A,\xi) \\ \rTTo_{\fA(A,g)} }
& \fA(A,C) \\
&& \\
\dTTo<{\lambda^x_B} & \ldTwoar^{\lambda^x_g} & \dTTo>{\lambda^x_C} \\
GB & \rTTo_{Gg} & GC
\end{diagram}.
\label{dia-eq-lambda-xf-xg}
\end{equation}
Indeed, for each 1\n-morphism \(h:A\to B\) of $\fA$ we have
\begin{diagram}
(x)(G(hf)) & \rTTo^{(x)(h,f)\psi_2^{-1}} & (x)(Gh)(Gf) \\
\dTTo<{(x)(G(h\bull\xi))} & = & \dTTo>{(x)(Gh)(G\xi)} \\
(x)(G(hg)) & \rTTo^{(x)(h,g)\psi_2^{-1}} & (x)(Gh)(Gg)
\end{diagram}
by naturality of $\psi_2$. Rewriting this equation in the form
\begin{diagram}
(hf)\lambda^x_C & \rTTo^{(h)\lambda^x_f} & (h)\lambda^x_B(Gf) \\
\dTTo<{(h\bull\xi)\lambda^x_C} & = & \dTTo>{((h)\lambda^x_B)(G\xi)} \\
(hg)\lambda^x_C & \rTTo^{(h)\lambda^x_g} & (h)\lambda^x_B(Gg)
\end{diagram}
we deduce that \eqref{dia-eq-lambda-xf-xg} holds on $h$. Therefore,
condition~2 of \defref{def-weak-2-trans} is satisfied.

Let us verify condition~3(a) of \defref{def-weak-2-trans}, that is,
equation
\[
\begin{diagram}[inline,width=4em]
\fA(A,B) & \rTTo^{1_{\fA(A,B)}}_{\fA(A,1_B)} & \fA(A,B) \\
\dTTo<{\lambda^x_B} & \ldTwoar_{\lambda^x_{1_B}} & \dTTo>{\lambda^x_B}
\\
GB & \rTTo_{G1_B} & GB
\end{diagram}
\quad = \quad
\begin{diagram}[inline,width=4em]
\fA(A,B) & \rTTo^{1_{\fA(A,B)}} & \fA(A,B) \\
\dTTo<{\lambda^x_B} & = & \dTTo>{\lambda^x_B} \\
GB & \pile{\rTTo^{1_{GB}}\\ \Downarrow\psi_0\\ \rTTo_{G1_B}} & GB
\end{diagram}.
\]
On an object \(f:A\to B\) of \(\fA(A,B)\) it reads:
\begin{equation}
(f)\lambda^x_{1_B} = ((f)\lambda^x_B)\psi_0:
(f)\lambda^x_B \to (f)\lambda^x_B(G1_B).
\label{eq-(f)-lambdax1B}
\end{equation}
It follows from condition~\eqref{eq-w-2-fun-a1} for $G$,
\[ \bigl[Gf \rTTo^{Gf\bull\psi_0} (Gf)(G1_B) \rTTo^{(f,1_B)\psi_2}
Gf\bigr] = 1_{Gf}: Gf \to Gf: GA \to GB,
\]
which, evaluated on \(x\in\Ob GA\), can be written as
\[ (x)(f,1_B)\psi_2^{-1} = ((x)(Gf))\psi_0: (x)(Gf) \to (x)(Gf)(G1_B).
\]
This is precisely \eqref{eq-(f)-lambdax1B}.

Let us verify condition~3(b) of \defref{def-weak-2-trans}, that is,
equation
\begin{equation}
\begin{diagram}[inline]
&& \fA(A,C) && \\
& \ruTTo^{\fA(A,f)} & \dId>= & \rdTTo^{\fA(A,g)} & \\
\fA(A,B) && \rTTo_{\fA(A,fg)} && \fA(A,D) \\
\dTTo<{\lambda^x_B} &&& \ldTwoar(4,2)_{\lambda^x_{fg}} &
\dTTo>{\lambda_D} \\
GB && \rTTo_{G(fg)} && GD
\end{diagram}
\quad = \quad
\begin{diagram}[inline]
&& \fA(A,C) && \\
& \ruTTo^{\fA(A,f)} \ldTwoar(2,4)_{\lambda^x_f} & \dTTo>{\lambda^x_C} &
\rdTTo^{\fA(A,g)} & \\
\fA(A,B) && GC & \lTwoar^{\lambda^x_g} & \fA(A,D) \\
\dTTo<{\lambda^x_B} & \ruTTo_{Gf} & \dTwoar>{\psi_2} & \rdTTo_{Gg} &
\dTTo>{\lambda^x_D} \\
GB && \rTTo_{G(fg)} && GD
\end{diagram}
\label{dia-eq-3b-lambdax}
\end{equation}
for arbitrary pair of composable 1\n-morphisms \(B\rTTo^f C \rTTo^g D\)
of $\fA$. We have to check this equation on an arbitrary 1\n-morphism
$h:A\to B$. Condition~4(b) of \defref{def-weak-2-func} for $G$ is the
equation
\begin{diagram}
(Gh)(Gf)(Gg) & \rTTo^{Gh\bull(f,g)\psi_2} & (Gh)\bull G(fg) \\
\dTTo<{(h,f)\psi_2\bull Gg} & = & \dTTo>{(h,fg)\psi_2} \\
G(hf).Gg & \rTTo^{(hf,g)\psi_2} & G(hfg)
\end{diagram}
Evaluating it on $x$ we get the equation
\begin{multline*}
(x)(h,fg)\psi_2^{-1} = \bigl[ (x)(G(hfg))
\rTTo^{(x)(hf,g)\psi_2^{-1}} (x)(G(hf))(Gg) \\
\rTTo^{((x)(h,f)\psi_2^{-1})(Gg)} (x)(Gh)(Gf)(Gg)
\rTTo^{((x)(Gh))(f,g)\psi_2} (x)(Gh)(G(fg)) \bigr],
\end{multline*}
which can be rewritten as
\begin{multline*}
(h)\lambda^x_{fg} = \bigl[ (hfg)\lambda^x_D
\rTTo^{(hf)\lambda^x_g} (hf)\lambda^x_C(Gg) \\
\rTTo^{((h)\lambda^x_f)(Gg)} (h)\lambda^x_B(Gf)(Gg)
\rTTo^{((h)\lambda^x_B)(f,g)\psi_2} (h)\lambda^x_B(G(fg)) \bigr].
\end{multline*}
And this is precisely \eqref{dia-eq-3b-lambdax}, evaluated on
\(h:A\to B\).

Therefore, all conditions of \defref{def-weak-2-trans} are satisfied,
and $\lambda^x$ is a weak 2\n-transformation.

Let us show that correspondence~\eqref{eq-lambda-uC} defines a natural
transformation. Indeed, for each 2\n-morphism \(\xi:f\to g:A\to C\) of
$\fA$ the diagram
\begin{diagram}[width=5em]
(f)\lambda^x_C = (x)(Gf) & \rTTo^{(f)\lambda^u_C}_{(u)(Gf)}
& (y)(Gf) = (f)\lambda^y_C \\
\dTTo<{(\xi)\lambda^x_C=}>{(x)(G\xi)} &&
\dTTo<{(y)(G\xi)}>{=(\xi)\lambda^y_C} \\
(g)\lambda^x_C = (x)(Gg) & \rTTo_{(g)\lambda^u_C}^{(u)(Gg)}
& (y)(Gg) = (g)\lambda^y_C
\end{diagram}
commutes due to \(G\xi:Gf\to Gg:GA\to GC\) being a natural
transformation.

We claim that property~2 of \defref{def-modif} holds for \(\lambda^u\).
For an arbitrary 1\n-morphism \(f:B\to C\) of $\fA$ we have to prove
the equation
\[
\begin{diagram}
\fA(A,B) & \rTTo^{\fA(A,f)} & \fA(A,C) \\
\dTTo<{\lambda^y_B} \overset{\sss \lambda^u_B}{\Leftarrow}
\dTTo>{\lambda^x_B} & \ldTwoar^{\lambda^x_f} & \dTTo>{\lambda^x_C} \\
GB & \rTTo^{Gf} & GC
\end{diagram}
\quad = \quad
\begin{diagram}
\fA(A,B) & \rTTo^{\fA(A,f)} & \fA(A,C) \\
\dTTo<{\lambda^y_B} & \ldTwoar^{\lambda^y_f} & \dTTo<{\lambda^y_C}
\overset{\sss \lambda^u_C}{\Leftarrow} \dTTo>{\lambda^x_C} \\
GB & \rTTo^{Gf} & GC
\end{diagram}
\quad.
\]
On the object \(g:A\to B\) of \(\fA(A,B)\) this equation reads
\begin{diagram}[width=7em]
(gf)\lambda^x_C = (x)(G(gf))
& \rTTo^{(g)\lambda^x_f}_{(x)(g,f)\psi_2^{-1}}
& (x)(Gg)(Gf) = (g)\lambda^x_B(Gf) \\
\dTTo<{(gf)\lambda^u_C=}>{(u)(G(gf))} & = &
\dTTo<{(u)(Gg)(Gf)}>{=(g)\lambda^u_B(Gf)} \\
(gf)\lambda^y_C = (y)(G(gf))
& \rTTo_{(g)\lambda^y_f}^{(y)(g,f)\psi_2^{-1}}
& (y)(Gg)(Gf) = (g)\lambda^y_B(Gf)
\end{diagram}
It holds due to \((g,f)\psi_2^{-1}:G(gf)\to(Gg)(Gf):GA\to GC\) being a
natural transformation. Therefore, \(\lambda^u\) is a modification.

The unit morphism \(1_x:x\to x\) of $GA$ goes to the identity
transformation
\[ \lambda^{1_x}_C: (f:A\to C) \longmapsto
\bigl((1_x)(Gf)=1_{(x)(Gf)}: (x)(Gf) \to (x)(Gf)\bigr),
\]
because $Gf$ is a functor. For a pair of composable morphisms
\(x \rTTo^u y \rTTo^v z\) of $GA$ we have
 \(\lambda^{uv}_C=\bigl(\lambda^x_C \rTTo^{\lambda^u_C} \lambda^y_C
 \rTTo^{\lambda^v_C} \lambda^z_C\bigr)\),
since \((uv)(Gf)=(u)(Gf)\cdot(v)(Gf)\) due to $Gf$ being a functor.
Therefore, $\Lambda$ is a functor.
\end{proof}
\fi

The result of Yoneda Lemma~\ref{lem-Yoneda-Lemma-bicategories} for a
strict 2-category $\fA$ can be made more precise as follows.

\begin{proposition}\label{pro-ev1A-Lambda-equivalences}
Functors
\[ \ev_{1_A}:[\fA,\Cat](\fA(A,\_),G) \to GA \quad\text{ and }\quad
\Lambda:GA \to [\fA,\Cat](\fA(A,\_),G)
\]
are equivalences, quasi-inverse to each other.
\end{proposition}

\ifx\chooseClass2
    \else
\begin{proof}
We have
\[ \bigl[GA \rTTo^\Lambda \mbox{} [\fA,\Cat](\fA(A,\_),G)
\rTTo^{\ev_{1_A}} GA\bigr] = G1_A.
\]
Indeed, for any object $x$ of $GA$
\[ (x)\Lambda\ev_{1_A} = (\lambda^x)\ev_{1_A} = (1_A)\lambda^x_A
= (x)(G1_A),
\]
for any morphism $u:x\to y$ of $GA$
\[ (u)\Lambda\ev_{1_A} = (\lambda^u)\ev_{1_A} = (1_A)\lambda^u_A
= (u)(G1_A).
\]
An isomorphism of functors \(\psi_0:1_{GA}\to G1_A\) implies that an
arbitrary object $x$ of $GA$ is isomorphic to
\((x)(G1_A)=((x)\Lambda)\ev_{1_A}\). Thus, \(\ev_{1_A}\) is essentially
surjective on objects. By \propref{pro-ev1A-full-and-faithful}
\(\ev_{1_A}\) is an equivalence. Therefore, $\Lambda$ is isomorphic to
a functor quasi-inverse to \(\ev_{1_A}\). Hence, $\Lambda$ itself is an
equivalence quasi-inverse to \(\ev_{1_A}\).
\end{proof}
\fi

\subsection{\texorpdfstring{Example of strict 2-functor $G=\fA(B,\_)$}
 {Example of strict 2-functor G=A(B,-)}}
Applying \propref{pro-functor-GA-(ACat)} to the strict 2\n-functor
$G=\fA(B,\_):\fA\to\Cat$, we get the following. An arbitrary
1\n-morphism $f:B\to A$ gives rise to  the strict 2-transformation
$f^*=\sS{^{\fA(B,-)}}\lambda^{A,f}:\fA(A,\_)\to\fA(B,\_)$. It is
specified by the family of functors $f^*_C$, $C\in\Ob\fA$:
\begin{align*}
f^*_C = \fA(f,C): \fA(A,C) & \longrightarrow \fA(B,C) \\
(\phi:A\to C) & \longmapsto (f)\fA(B,\phi) = f\phi: B\to C, \\
(\pi:\phi\to\psi:A\to C) & \longmapsto
(f)\fA(B,\pi) = f\bull\pi:f\phi\to f\psi:B\to C.
\end{align*}

An arbitrary 2\n-morphism $\alpha:f\to g:B\to A$ gives rise to the
modification
 $\alpha^*=\sS{^{\fA(B,-)}}\lambda^{A,\alpha}:f^*\to g^*:
 \fA(A,\_)\to\fA(B,\_)$
given by the family of natural transformations
$\alpha^*_C:f^*_C\to g^*_C:\fA(A,C)\to\fA(B,C)$, $C\in\Ob\fA$. The
transformation $\alpha^*_C$ is specified by its components:
\begin{equation*}
\alpha^*_C = \fA(\alpha,C): (\phi:A\to C) \longmapsto
(\alpha)\fA(B,\phi) = \alpha\bull\phi:
(\phi)f^*_C = f\phi\to g\phi = (\phi)g^*_C.
\end{equation*}

By \propref{pro-functor-GA-(ACat)} the correspondence $f\mapsto f^*$,
$\alpha\mapsto\alpha^*$ determines a functor
$Y_{AB}:\fA^{\op}(A,B)=\fA(B,A)\to[\fA,\Cat](\fA(A,\_),\fA(B,\_))$. One
easily verifies that in fact we have a strict 2\n-functor
$Y:\fA^{\op}\to[\fA,\Cat]$.

\begin{corollary}\label{pro-Yab-equivalence}
$Y$ is a local equivalence, i.e., for each pair of objects
$A,B\in\Ob\fA^{\op}$ the functor $Y_{AB}$ is an equivalence.
\end{corollary}

Let us recall also the notion of a birepresentable homomorphism
\(G:\fA\to\Cat\) following Street~\cite[(1.11)]{Street:14}. He
formulates the following statement for an arbitrary bicategory $\fA$,
but we assume that $\fA$ is a strict 2\n-category as usual.

\begin{proposition}\label{pro-representable-2-functor}
Let \(G:\fA\to\Cat\) be a weak 2\n-functor. Then the following
conditions are equivalent:
\begin{enumerate}
\item there exists an object $A$ of $\fA$ and a 2\n-natural equivalence
\(\lambda:\fA(A,\_)\to G\);

\item there exists an object $A$ of $\fA$ and an object $x$ of $GA$
such that the weak 2\n-transformation \(\lambda^x:\fA(A,\_)\to G\) is a
2\n-natural equivalence.
\end{enumerate}
\end{proposition}

\begin{proof}
Clearly, the second property implies the first one. Assume that
condition~1) holds. By \propref{pro-ev1A-Lambda-equivalences} the weak
2\n-transformation $\lambda$ is isomorphic to $\lambda^x$ for some
$x\in\Ob GA$. By \propref{pro-2-natural-quasi-inverse} $\lambda$ is a
quasi-invertible 1\n-morphism of $[\fA,\Cat]$, hence, so is
$\lambda^x$. By \remref{rem-quasi-inverse-2-natural-equivalences}
condition~2) holds.
\end{proof}

\begin{definition}\label{def-representable-by-pair}
A weak 2\n-functor \(G:\fA\to\Cat\) is \emph{representable}
(\emph{birepresentable} in terminology of
Street~\cite[(1.11)]{Street:14}) if it satisfies equivalent conditions
of \propref{pro-representable-2-functor}. A pair $(A,x)$ consisting of
an object $A$ of $\fA$ and an object $x$ of $GA$ is said to
\emph{represent} (\emph{birepresent}) $G$, if
\(\lambda^x=\lambda^{A,x}=\sS{^G}\lambda^{A,x}:\fA(A,\_)\to G\) is a
2\n-natural equivalence.
\end{definition}

\subsection{Uniqueness of the representing pair}
 \label{sec-representing-pair-unique}
It is shown by Street that a representing pair is unique up to an
equivalence in a certain bicategory \cite[(1.10)-(1.11)]{Street:14}.
Let us provide the details in our setting.

Let two pairs $(A,x)$ and $(B,y)$ represent $G$. Then there is a
quasi-inverse to $\lambda^{B,y}:\fA(B,\_)\to G$ weak 2\n-transformation
$\lambda^{B,y-}:G\to\fA(B,\_)$. Define a 2\n-natural equivalence
\(\mu=\lambda^{A,x}\bull\lambda^{B,y-}:\fA(A,\_)\to\fA(B,\_)\). It is
isomorphic to the 2\n-transformation \(\sS{^{\fA(B,-)}}\lambda^{A,f}\)
for some \(f\in\Ob\fA(B,A)\). There is an invertible modification $m$:
\begin{diagram}
\fA(A,\_) && \rTTo^{\sS{^G}\lambda^{A,x}} && G \\
& \rdTTo_{\sS{^{\fA(B,-)}}\lambda^{A,f}} & \dTwoar>m
& \ruTTo_{\sS{^G}\lambda^{B,y}} & \\
&& \fA(B,\_) &&
\end{diagram}
Then \((1_A)m_A:(x)(G1_A)\to(y)(Gf)\) is an isomorphism of $GA$.
Therefore,
\[ x = (x)(1_{GA}) \rTTo^{(x)\psi_0} (x)(G1_A) \rTTo^{(1_A)m_A} (y)(Gf)
\]
is an isomorphism of $GA$. By symmetry we get a 1\n-morphism
\(g:A\to B\) of $\fA$ and an isomorphism \(y \rTTo^\sim (x)(Gg)\) of
$GB$. By construction the strict 2\n-transformations
\begin{align*}
\sS{^{\fA(B,-)}}\lambda^{A,f} &: \fA(A,\_) \to \fA(B,\_), \\
\sS{^{\fA(A,-)}}\lambda^{B,g} &: \fA(B,\_) \to \fA(A,\_)
\end{align*}
are quasi-inverse to each other. In particular,
\begin{multline*}
1_A \simeq
((1_A)\sS{^{\fA(B,-)}}\lambda^{A,f}_A)\sS{^{\fA(A,-)}}\lambda^{B,g}_A
= ((f)\fA(B,1_A))\sS{^{\fA(A,-)}}\lambda^{B,g}_A \\
= (f)\sS{^{\fA(A,-)}}\lambda^{B,g}_A = (g)\fA(A,f) = gf,
\end{multline*}
and by symmetry \(1_B\simeq fg\). Therefore, 1\n-morphisms $f$ and $g$
are quasi-inverse to each other.

Summing up, a pair \((A\in\Ob\fA,x\in\Ob GA)\) representing a weak
2\n-functor \(G:\fA\to\Cat\) is unique up to equivalence $f$ of the
first objects, such that $Gf$ preserves the second object up to an
isomorphism.

\tableofcontents

\end{document}